\documentclass[final,3p,10pt]{elsarticle}

\usepackage{amssymb}
\usepackage{pifont}

\usepackage{amsmath,fixltx2e}
\usepackage{bm}
\usepackage{xcolor} 

\usepackage{placeins} 
\usepackage{tikz}
\usetikzlibrary{patterns}
\usepackage{adjustbox}
\usepackage{subcaption}
\usepackage{varioref}
\usepackage{hyperref}
\usepackage{cleveref} 
\usepackage{booktabs}
\usepackage{fmtcount} 
\usepackage{mathtools}

\usepackage{nomencl}
\makenomenclature

\let\oldcitet\citet
\renewcommand*\citep[1]{({\oldcitet{#1}})}

\newcommand\Nomenclature[4][X]{\nomenclature[#1#2]{#3}{#4}}
\renewcommand\nomgroup[1]{
    \item[
        \hspace{-.7ex}\large\bfseries
        \ifstrequal{#1}{A}{Functions}{%
        \ifstrequal{#1}{B}{Indices}{%
        \ifstrequal{#1}{C}{Groups}{%
        \ifstrequal{#1}{D}{Variables}{%
        \ifstrequal{#1}{E}{Sub and super-indices}{%
        \ifstrequal{#1}{F}{Operators}{%
        }}}}}}
    ]%
    \vspace{2ex}
 }

\Nomenclature[A]{01}{$\delta_{ij}$}{Kronecker delta}
\Nomenclature[A]{02}{$\isotransform$}{Isoparametric transform}
\Nomenclature[B]{01}{$\itype$}{Element type index}
\Nomenclature[B]{02}{$\iele$}{Element index}
\Nomenclature[B]{03}{$\ivar$}{Field variable index}
\Nomenclature[B]{04}{$i, j, k, \ipol, \isigma, \inu$}{Arbitrary summation indices}
\Nomenclature[C]{01}{$\domain$}{Domain}
\Nomenclature[C]{01}{$\stencil$}{Stencil}
\Nomenclature[D]{01}{$\degree$}{Polynomial degree}
\Nomenclature[D]{02}{$\ndim$}{Number of spatial dimensions}
\Nomenclature[D]{03}{$\nvar$}{Number of field variables}
\Nomenclature[D]{04}{$N_{\itype}$}{Number of elements of type $\itype$}
\Nomenclature[D]{05}{$\nodalbasis_{\itype \ipol}$}{Scalar nodal basis polynomial with index $\ipol$ for element type $\itype$}
\Nomenclature[D]{06}{$\modalbasis_{\itype \ipol}$}{Scalar modal basis polynomial with index $\ipol$ for element type $\itype$}
\Nomenclature[D]{07}{$\fluxnodalbasis_{\itype \ipol}$}{Vector nodal basis polynomial with index $\ipol$ for element type $\itype$}
\Nomenclature[D]{08}{$\fluxmodalbasis_{\itype \ipol}$}{Vector modal basis polynomial with index $\ipol$ for element type $\itype$}
\Nomenclature[D]{09}{$\vect{x}$}{Position in the physical space}
\Nomenclature[D]{10}{$\var$}{Conserved variable}
\Nomenclature[D]{11}{$\vargrad$}{Auxiliary gradient variable}
\Nomenclature[D]{12}{$\flux$}{Flux}
\Nomenclature[D]{13}{$\vandermonde$}{Vandermonde matrix}
\Nomenclature[D]{14}{$\jacobian$}{Isoparametric transform Jacobian matrix}
\Nomenclature[D]{15}{$\jacobianrhs$}{Right-hand-side Jacobian matrix}
\Nomenclature[E]{01}{$\atsol{\square}$}{A quantity defined at the solution points}
\Nomenclature[E]{02}{$\atflux{\square}$}{A quantity defined at the flux points}
\Nomenclature[E]{03}{$\atfluxnormal{\square}$}{A normal quantity defined at the flux points}
\Nomenclature[E]{04}{$\atfluxext{\square}$}{A quantity defined at the external flux points}
\Nomenclature[E]{05}{$\atfluxnormalext{\square}$}{A normal quantity defined at the external flux points}
\Nomenclature[E]{06}{$\atfluxint{\square}$}{A quantity defined at the internal flux points}
\Nomenclature[E]{07}{$\atfluxnormalint{\square}$}{A normal quantity defined at the internal flux points}
\Nomenclature[E]{08}{$\atfluxintunique{\square}$}{A quantity defined at the unique internal flux points}
\Nomenclature[E]{09}{$\transformed{\square}$}{A quantity defined in a transformed space}
\Nomenclature[E]{10}{$\unit{\square}$}{A unit magnitude vector}
\Nomenclature[F]{01}{$\commvar_{\ivar}$}{Common value of a conserved variable at the external flux points}
\Nomenclature[F]{02}{$\commflux_{\ivar}$}{Common value of a normal flux at the external flux points}

\newtheorem{remark}{Remark}

\newtheorem{disclaimer}{Disclaimer}

\usepackage[mathscr]{eucal}

\newcommand{\var}{u}
\newcommand{\vargradnovect}{q}
\newcommand{\vargrad}{\vect{\vargradnovect}}
\newcommand{\statevar}{\vect{U}}
\newcommand{\statevargrad}{\vect{Q}}

\newcommand{\vect}[1]{\bm{#1}}
\newcommand{\divergence}{\vect{\nabla} \cdot}

\newcommand{\flux}{\vect{\mathrm{f}}}
\newcommand{\fluxconv}{\vect{\mathrm{f}}^{\text{conv}}}
\newcommand{\fluxvisc}{\vect{\mathrm{f}}^{\text{visc}}}
\newcommand{\fluxnovect}{\mathrm{f}}
\newcommand{\jacobian}{\vect{\mathrm{J}}}
\newcommand{\jacobianrhs}{\mathbb{J}}
\newcommand{\jacobianrhsprime}{\mathbb{G}}
\newcommand{\expjacobianrhs}{\mathbb{E}}

\newcommand{\nodalbasis}{l}
\newcommand{\fluxnodalbasisnovect}{l}
\newcommand{\fluxnodalbasis}{\vect{\fluxnodalbasisnovect}}
\newcommand{\modalbasis}{\psi}
\newcommand{\fluxmodalbasis}{\vect{\psi}}
\newcommand{\correctionnodalbasis}{\atfluxext{\vect{g}}}
\newcommand{\grad}{\vect{\nabla}}

\newcommand{\nvar}{N_{\text{V}}}
\newcommand{\ndim}{N_{\text{D}}}
\newcommand{\determinant}{\text{det  }}
\newcommand{\normal}{\vect{n}}

\newcommand{\unit}[1]{\hat{#1}}

\newcommand{\ncellsinpattern}{N_p}

\newcommand{\stencil}{\mathcal{C}}

\newcommand{\ivar}{\alpha}
\newcommand{\iidim}{i}
\newcommand{\ijdim}{j}
\newcommand{\ipol}{\rho}
\newcommand{\isigma}{\sigma}
\newcommand{\inu}{\nu}
\newcommand{\itype}{e}
\newcommand{\iele}{n}

\newcommand{\isotransform}{\mathcal{M}}

\newcommand{\transformed}[1]{\widetilde{#1}}

\newcommand{\conj}[1]{\overline{#1}}

\newcommand{\atsol}[1]{#1^{(\var)}}
\newcommand{\atflux}[1]{#1^{(f)}}
\newcommand{\atfluxnormal}[1]{#1^{(f\bot)}}
\newcommand{\atfluxint}[1]{#1^{(fi)}}
\newcommand{\atfluxintunique}[1]{#1^{(fiu)}}
\newcommand{\atfluxunique}[1]{#1^{(fu)}}
\newcommand{\atfluxext}[1]{#1^{(fe)}}
\newcommand{\atfluxnormalint}[1]{#1^{(fi\bot)}}
\newcommand{\atfluxnormalext}[1]{#1^{(fe\bot)}}

\newcommand{\vandermonde}{\mathcal{V}}
\newcommand{\commvar}{\mathfrak{C}}
\newcommand{\commonvar}{\mathfrak{C}}
\newcommand{\commflux}{\mathfrak{F}}

\newcommand{\MA}{\mathrm{Ma}}
\newcommand{\RE}{\mathrm{Re}}
\newcommand{\PR}{\mathrm{Pr}}
\newcommand{\PE}{\mathrm{Pe}}

\newcommand{\velocitynovect}{v}
\newcommand{\velocity}{\vect{\velocitynovect}}
\newcommand{\pressure}{p}
\newcommand{\density}{\rho}
\newcommand{\temperature}{\Theta}
\newcommand{\stressTensor}{\tau}

\newcommand{\advecvel}{\vect{\advecvelnovect}}
\newcommand{\advecvelnovect}{c}
\newcommand{\unitary}{\vect{1}}

\newcommand{\degree}{\mathfrak{p}}

\newcommand{\domain}{\Omega}
\newcommand{\norm}[1]{\left\lVert#1\right\rVert}
\newcommand{\numerical}{\delta}
\newcommand{\inner}[1]{\langle#1\rangle}

\newcommand{\penaltyterm}{\eta}
\newcommand{\polynomial}{\mathcal{P}}
\newcommand{\cellsize}{h}
\newcommand{\CFL}{\tau}
\newcommand{\timestep}{\Delta t}
\newcommand{\imag}{\text{I}}
\newcommand{\e}{\text{e}}
\newcommand{\labvec}{\vect{a}}
\newcommand{\patternvar}{U}
\newcommand{\m}{\vect{M}}
\newcommand{\fluxtransformedmatrix}{\transformed{\vect{F}}}
\newcommand{\residualmatrix}{\transformed{\vect{R}}}

\newcommand{\eigenvectormatrix}{\mathbb{W}}

\newcommand{\volume}{\domain}
\newcommand{\ensembleAverage}[1]{\langle#1\rangle}

\newcommand{\xmark}{\ding{55}}

\newcommand{\permutationmatrix}{\vect{P}}
\newcommand{\interpfiu}{\m^{12}}
\newcommand{\gradfi}{\m^{15}}
\newcommand{\gradfe}{\m^{16}}
\newcommand{\interpgradfiu}{\m^{18}}
\newcommand{\divfi}{\m^{13}}
\newcommand{\divfe}{\m^{14}}
\newcommand{\normfi}{\m^{19}}
\newcommand{\tokendim}{\quad \text{with} \quad}

\usepackage{environ}
\NewEnviron{equations}[1]{%
\begin{equation}\begin{split}
  \BODY
\end{split}\label{#1}\end{equation}
}

\NewEnviron{equations*}[1]{%
\begin{equation*}\begin{split}
  \BODY
\end{split}\label{#1}\end{equation*}
}

\newcommand{\Eref}[1]{Eq.~\ref{#1}}
\newcommand{\Fref}[1]{Fig.~\ref{#1}}
\newcommand{\Tref}[1]{Table~\ref{#1}}
\newcommand{\Sref}[1]{Section~\ref{#1}}

\journal{Computer Methods in Applied Mechanics and Engineering}

\begin{document}

\begin{frontmatter}



\title{The Spectral Difference Raviart-Thomas method for two and three-dimensional elements and its connection with the Flux Reconstruction formulation.}

\author[isae]{G. Sáez-Mischlich}

            
\author[imft]{J. Sierra-Ausín}
\ead{javier.sierra@imft.fr}

\author[isae]{J. Gressier}

\address[isae]{Institut Supérieur de l'Aéronautique et de l'Espace (ISAE-SUPAERO), Toulouse 31400, France}
\address[imft]{Institut de Mécanique des Fluides de Toulouse (IMFT), Toulouse 31400, France}

\begin{abstract}
The purpose of this work is to describe in detail the development of the Spectral Difference Raviart-Thomas (SDRT) formulation for two and three-dimensional tensor-product elements and simplexes.
Through the process, the authors establish the equivalence between the SDRT method and the Flux-Reconstruction (FR) approach under the assumption of the linearity of the flux and the mesh uniformity.
Such a connection allows to build a new family of FR schemes for two and three-dimensional simplexes and also to recover the well-known FR-SD method with tensor-product elements.
In addition, a thorough analysis of the numerical dissipation and dispersion of both aforementioned schemes and the nodal Discontinuous Galerkin FR (FR-DG) method with two and three-dimensional elements is proposed through the use of the combined-mode Fourier approach.
SDRT is shown to posses an enhanced temporal linear stability regarding the FR-DG.
On the contrary, SDRT displays larger dissipation and dispersion errors with respect to FR-DG.
Finally, the study is concluded with a set of numerical experiments, the linear advection-diffusion problem, the Isentropic Euler Vortex and the Taylor-Green Vortex (TGV).
The latter test case shows that SDRT schemes present a non-linear unstable behavior with simplex elements and certain polynomial degrees.
For the sake of completeness, the matrix form of the SDRT method is developed and the computational performance of SDRT with respect to FR schemes is evaluated using GPU architectures.
\end{abstract}




\end{frontmatter}


\pagenumbering{gobble}

\printnomenclature

\section{Introduction}

The ever-increasing demand for numerical accuracy to solve turbulent flows has raised interest in the study and application of high-order numerical methods for unstructured grids \citep{Shu2016}.
Such methods have the potential to enhance the accuracy per degree of freedom of numerical simulations in unstructured grids.
Despite several decades of continuous development and investment, the application of high-order numerical schemes to real-world test cases remains limited.
Several reasons may be attributed to this issue such as difficulties to treat solution discontinuities, instabilities due to diminished numerical dissipation, lack of appropriate tools to generate curved mesh elements for under-resolved configurations, etc.

Within high-order methods for unstructured grids, there is no doubt that Spectral-Element Methods (SEM) are among the most promising spatial discretization schemes for the simulation of turbulent flows.
SEM encompass a plethora of different schemes such as the Nodal Discontinuous Galerkin (NDG) \citep{hesthaven2008nodal}, the Flux Reconstruction (FR) \citep{Huynh2007, Vincent2010}, the Spectral Difference (SD) \citep{Kopriva1996, Liu2006} and the Spectral Volume (SV) \citep{Wang2002}.
These schemes share in common several aspects: they generally do not need quadrature evaluation (unless polynomial de-aliasing techniques are put in place), they are compact by nature and they rely on a description of the numerical solution within each mesh element using a nodal polynomial basis.
Such characteristics allow SEM to drastically improve the performance per degree of freedom when compared to other more common second or low-order numerical methods such as the Finite-Volume (FV) method.
Besides, the SEM formulation is suited to exploit the massive computational power of Graphics Processing Units (GPUs) computational architectures due to the data locality and the possibility of expressing all numerical operations as matrix-matrix and/or matrix-vector products.
GPUs have proven to substantially over-perform Central Processing Units (CPUs) architectures with SEM simulations \cite{Witherden2015}.
Despite the apparent numerical and computational beneficial properties of SEM, they present important temporal stability constraints and they lack robustness in under-resolved turbulent flows simulations due to aliasing errors and the lack of sufficient resolution to capture the dissipative scales of the flow.
To alleviate the latter issues, several techniques exist: the Spectral Vanishing Viscosity (SVV) \citep{Manzanero2020}, dealising techniques \citep{Spiegel2015_aliasing}, modal filtering \citep{Glaubitz2017}, application of skew-symmetric formulations \citep{Abe2018}, etc.
Nevertheless, such methods are usually accompanied by an increased computational cost and the need for appropriate tune-in of a certain set of parameters.

The FR method \citep{Huynh2007, Vincent2010} or FRM relies on a solution nodal basis, local to each element, to describe the solution and the discontinuous flux within each element.
Conservation is enforced through the use of common fluxes at element interfaces.
The contribution of such common fluxes to the numerical solution is taken into account by means of the so-called correction functions \cite{VINCENT2011}.
An appropriate choice of such correction functions guarantees the equivalence between the FR formulation and specific NDG schemes \cite{Williams2013tetra, Mengaldo2015, Zwanenburg2016}.
This is particularly important to extend the FR formulation to arbitrary simplex elements, resulting in the so-called FR-DG method.
FR schemes have been applied to numerically solve various systems of conservation laws including the Euler equations \citep{Williams2011}, Navier–Stokes equations \citep{Park2017, Iyer2019}, and their incompressible counterparts \citep{Loppi2018}.
Several implementations of FR are available and have recently demonstrated the possibility to achieve high computational efficiency and scalability on large problems \citep{Witherden2015}.
Within such FR solvers, special focus is placed on the open-source PyFR solver \citep{Witherden2014}, in which the methods described in this work have been implemented.

The foundation for SD schemes was initially introduced by \citet{Kopriva1996} under the name \textit{staggered grid Chebyshev multi-domain methods}.
\citet{Liu2006} adapted the latter method to a more general formulation referred to as the Spectral Difference Method (SDM) which allowed to build stable schemes for both triangular and quadrilateral elements.
The SD method (SDM) is of special interest since the flux function is evaluated with a staggered-grid approach and then projected into a polynomial space one degree higher than that of the solution basis.
To define this basis, an arrangement of internal and external flux points (different from that of the solution points), together with the definition of specific degrees of freedom within these flux points, is proposed for each element.
The SDM has been successfully applied to study non-linear equations with complex physics \cite{Lodato2016, Lodato2019}.
The foundations of the SDM have been conjectured to provide the SDM with additional dealiasing properties compared to other FRM \cite{Cox2021}.
Temporal stability and numerical properties (dissipation and dispersion) of the SDM were initially assessed in \citet{Abeele2008}, showing that the SDM with strictly higher than second-order accuracy were linearly unstable in triangular grids.
A formal stability criterion of the SDM in one-dimensional configurations was discussed in  \citet{Jameson2010}, proving that the SDM is stable in one-dimensional configurations, provided that the internal flux points are located at Gauss-Legendre quadrature points.
Additionally the latter two studies found that, at least for linear cases, the position of the solution points has little to no impact on the accuracy and stability of the method while the definition of the flux points has important implications in those aspects.
The generation of stable high-order SDM for triangular elements was firstly discussed in  \citet{Balan2011}.
Such schemes were built using Raviart-Thomas (RT) basis to build the flux polynomial.
Nevertheless, strictly higher than fourth order schemes were found to be unstable in the former study.
The extension of the SDRT for triangle elements with viscous fluxes was firstly described in \citet{Li2019}.
\citet{Adele2021thesis} proposed the use of certain quadrature points to locate the flux points, proving the existence of stable SDRT schemes for triangular elements for polynomial degrees strictly lower than five, while also providing with a set of flux points which yields stable sixth order SDRT schemes in triangles. 
The latter work also describes the extension of the SDRT for three-dimensional elements, although it does not thoroughly describe the SDRT extension for prismatic elements neither do the results obtained in this work coincide with those presented therein.

This study presents an analysis of the Spectral Difference Raviart-Thomas (SDRT) method for quadrilateral, triangular, hexahedral, tetrahedral and triangular prismatic elements within the open-source PyFR solver \cite{Witherden2014}.
Within this analysis, an equivalence between FR and SDRT will be established and the properties of SDRT schemes will be compared with that of FR-DG.
The study is organized as follows.
\Sref{sec:formulation} introduces the mathematical formulation of SDRT and FR methods, demonstrating the equivalence between the latter two methods under certain conditions.
Next,  \Cref{sec:von_neumann_advection,sec:von_neumann_diffusion}  analyze the dissipation and dispersion properties of the aforementioned schemes with two and three-dimensional elements, using the so-called combined mode approach in the linear advection and linear diffusion equations.
The latter sections also present the linear stability and temporal linear stability criterions for SDRT and FR methods.
To validate the observations in the aforementioned linear analyses, \Sref{sec:numerical_experiments} studies the accuracy of the SDRT and FR-DG schemes through different numerical experiments regarding linear and non-linear test cases.
Finally, the conclusions and future work are illustrated in \Sref{sec:conclusions}.

For the sake of completeness, several appendices have been added to this work to better explain certain aspects of the SDRT method.
\ref{sec:rt_basis} describes the modal flux basis and the location of the flux points of the proposed SDRT method for different element types.
Furthermore, \ref{sec:sdrt_matrix_form} introduces the matrix form of the SDRT method, following the nomenclature of \cite{Witherden2014}.
Next, \ref{sec:perfo_gpu} compares the performance of the SDRT and FR methods with GPU architectures.
At last, \ref{sec:aliasing} discusses the wavenumber aliasing issues of SEM.



\section{Formulation}
\label{sec:formulation}

\begin{disclaimer}
Throughout this work, the SEM notation introduced in \cite{Witherden2014} will be extensively used to describe the numerical methods.
\end{disclaimer}


This work focuses on the solution systems of conservation laws written in the form
\begin{equation}
    \frac{\partial \var_{\ivar}}{\partial t}  + \divergence\flux_{\ivar}  = 0, \vect{x} \in \domain,
    \label{eq:system_conservation_laws}
\end{equation}
where $0 \le \ivar < \nvar$ is the field variable index, $\var_\ivar = \var_\ivar (\vect{x}, t)$ is the correspondent conserved variable with index $\ivar$, $\flux_{\ivar} = \flux_{\ivar} (\var, \grad \var)$ is the flux operator of the considered conserved variable and $\vect{x} = x_i \in \mathbb{R}^{\ndim}$.
The unscripted variables $\var \in \mathbb{R}^{\nvar}$ and $\grad \var \in \mathbb{R}^{\ndim \times \nvar}$ refer to the field variables and their gradients, respectively.
The flux operator is usually split in a convective flux $\fluxconv(\var)$, only dependent on the conserved variables, and a viscous flux $\fluxvisc(\var, \grad \var) = \flux(\var, \grad \var) - \fluxconv(\var)$.
\Eref{eq:system_conservation_laws} may be rewritten as a first-order system as
\begin{equations}{eq:system_conservation_laws_first_order}
\frac{\partial \var_{\ivar}}{\partial t}  + \divergence\flux_{\ivar}  &= 0\\
\vargrad_\ivar &= \grad \var_\ivar,
\end{equations}
where $\vargrad \in \mathbb{R}^{\ndim \times \nvar}$ is the auxiliary gradient variable and its unscripted form $\vargrad$ follows the same convention as that of $\var$ and $\grad \var$.

\begin{figure}
    \centering
    \centering
    \includegraphics[width=1.0\textwidth]{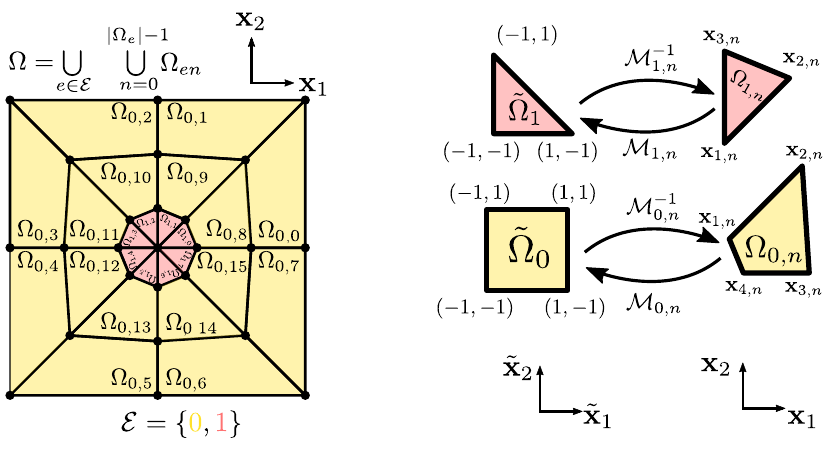}
    \caption{Domain notation}
    \label{fig:Notation_Domain}
\end{figure}

The domain $\domain$ is tessellated with a set of conforming elements, see \Fref{fig:Notation_Domain}.
Each element may be grouped in a set of element types $\mathcal{E}$ such that
\begin{equation}
    \domain = \bigcup_{\itype \in \mathcal{E}} \Omega_\itype \text{,} \quad \domain_\itype = \bigcup_{\iele=0}^{|\domain_\itype| - 1} \domain_{\itype \iele} \text{ and} \quad \bigcap_{\itype \in \mathcal{E}} \bigcap_{\iele=0}^{|\domain_\itype| - 1} \domain_{\itype \iele} = 0.
\end{equation}

Within each element $\domain_{\itype \iele}$, \Eref{eq:system_conservation_laws_first_order} is solved using Spectral-Element methods.
To do so, it is convenient to transform the first-order system to the coordinates of a standard or reference element of each type $\hat{\domain}_\itype$.
This transformation is defined by a mapping function for each element such that
%
\begin{equation}
    \vect{x} = \vect{\isotransform}_{\itype\iele}(\transformed{\vect{x}}) \text{ and} \quad \transformed{\vect{x}} = \vect{\isotransform}_{\itype\iele}^{-1}(\vect{x}) .
    \label{eq:mapping_function}
\end{equation}
In the latter, it is assumed the existence and well-posedness of the inverse of the transformation, i.e., the considered elements are not degenerated.
The associated Jacobian matrix of the mapping function may be written as
\begin{equation}
    \jacobian_{\itype\iele} = J_{\itype\iele\iidim\ijdim} \equiv \frac{\partial \isotransform_{\itype\iele\iidim}}{\partial \transformed{x}_{\ijdim}} \quad \text{and} \quad  J_{\itype\iele} \equiv \determinant \jacobian_{\itype\iele} .
    \label{eq:mapping_jacobian}
\end{equation}
Analogously, the Jacobian of the inverse transformation is introduced as follows
\begin{equation}
    \jacobian_{\itype\iele}^{-1} = J^{-1}_{\itype\iele\iidim\ijdim} \equiv \frac{\partial \isotransform^{-1}_{\itype\iele\iidim}}{\partial x_{\ijdim}} \quad \text{and} \quad  J^{-1}_{\itype\iele} \equiv \determinant \jacobian^{-1}_{\itype\iele} = \frac{1}{J_{\itype\iele}} .
    \label{eq:inverse_mapping_jacobian}
\end{equation}

In what follows, it is supposed that the Jacobian matrices are computed using their analytical values.
Nevertheless, it is worth noting that such a choice might induce free-stream preservation issues in non-uniform grids \cite{Abe2015}.
However, since this work is focused on the analysis of uniform grids, such issues do not occur.
The aforementioned mappings may be used to define the transformed conserved variables, flux operators, and auxiliary variables as \cite{Kopriva2006}
\begin{equations}{eq:transformed_quantities}
\transformed{\var}_{\itype \iele \ivar} &= \transformed{\var}_{\itype \iele \ivar}(\transformed{\vect{x}}, t) = J_{\itype \iele} \var_{\itype \iele \ivar}(\vect{x}, t) ,\\
\transformed{\flux}_{\itype \iele \ivar} &= \transformed{\flux}_{\itype \iele \ivar}(\transformed{\vect{x}}, t) = J_{\itype \iele}(\transformed{\vect{x}}) \jacobian^{-1}_{\itype \iele}(\transformed{\vect{x}}) \flux_{\itype \iele \ivar}(\vect{x}, t) ,\\
\transformed{\vargrad}_{\itype \iele \ivar} &= \transformed{\vargrad}_{\itype \iele \ivar}(\transformed{\vect{x}}, t) = \jacobian_{\itype \iele}^T(\transformed{\vect{x}}) \vargrad(\vect{x}, t) ,
\end{equations}
where the relation between $\vect{x}$ and $\transformed{\vect{x}}$ is given by \Eref{eq:mapping_function}.
Supposing static grids, \Eref{eq:system_conservation_laws_first_order} can be recast as 
\begin{equations}{eq:system_conservation_laws_first_order_transformed}
\frac{\partial \var_{\ivar}}{\partial t}  + J^{-1}_{\itype\iele} \transformed{\vect{\nabla}} \cdot \transformed{\flux}_{\ivar}  &= 0\\
\transformed{\vargrad}_\ivar &= \transformed{\vect{\nabla}} \var_\ivar,
\end{equations}
where $\transformed{\vect{\nabla}} = \partial/\partial \transformed{x}_i$.



%

Spectral-Element Methods use a nodal basis of degree $\degree$, unique to each element type, to describe the numerical solution within each element of a given grid.
For each element type $\itype \in \mathcal{E}$, a set of solution points with position in the reference element $\{ \atsol{\transformed{\vect{x}}_{\itype \ipol}} \}$ such that $0 \le \ipol < \atsol{N}_{\itype}( \degree )$ is defined.
In the latter, $\degree$ is the degree of the nodal basis.
These solution points may be used to build the nodal basis set $\{ \atsol{\nodalbasis_{\itype \ipol}}(\transformed{\vect{x}})   \}$ with the nodal orthogonal property $\atsol{\nodalbasis_{\itype \ipol}}(\atsol{\transformed{\vect{x}}}_{\itype \isigma}) = \delta_{\ipol \isigma}$.

The solution nodal basis may also be expressed through the solution modal basis $\{ \atsol{\modalbasis}_{\itype \ipol} (\transformed{\vect{x}}) \} $.
The relation between the nodal and modal bases is given by the solution Vandermonde matrix, defined as follows
\begin{equation}
    \atsol{\vandermonde}_{\itype \ipol \isigma} = \atsol{\modalbasis}_{\itype \ipol}\left( \atsol{\transformed{\vect{x}}}_{\itype \isigma} \right) .
\end{equation}
Such a matrix links the nodal and modal bases as 
\begin{equation}
    \atsol{\nodalbasis_{\itype \ipol}}(\transformed{\vect{x}}) = \left(\atsol{\vandermonde_{\itype \ipol \isigma}}\right)^{-1} \atsol{\modalbasis}_{\itype \isigma}\left( \transformed{\vect{x}} \right).
\end{equation}
To avoid unwanted numerical errors due to inaccuracies in the evaluation of the inverse Vandermonde matrix, it is recommended to build the solution modal basis using hierarchical orthonormal polynomials.
This allows to drastically reduce the condition number of the solution Vandermonde matrix.
The interested reader is referred to \cite{Karniadakis_2005} for a review of the hierarchical orthornormal bases for tensor-product and simplex elements used in this work.
It is worth mentioning that the nodal solution basis allows to express the conserved variable within a given element as
\begin{equation}
    \var_{\itype \iele \ivar}(\transformed{\vect{x}}) = \atsol{\var}_{\itype \iele \ipol \ivar} \atsol{\nodalbasis}_{\itype \ipol}(\transformed{\vect{x}}) ,
\end{equation}
being $\atsol{\var}_{\itype \iele \ipol \ivar} = \var_{\itype \iele \ivar}(\atsol{\transformed{\vect{x}}}_{\itype \ipol})$.

\subsection{SDRT method}

Together with the solution nodal basis, the SDRT method relies on a vector nodal basis using a set of flux points $\{ \atflux{\transformed{\vect{x}}_{\itype \ipol}} \}$ to compute the transformed flux divergence.
This vector nodal basis is referred to as $\{ \atflux{\fluxnodalbasis}_{\itype \ipol}\left( \transformed{\vect{x}} \right) \}$.
Its correspondent vector modal basis is denoted as $\{ \atflux{\fluxmodalbasis}_{\itype \ipol}\left( \transformed{\vect{x}} \right) \}$.
The flux modal basis lays on the Raviart-Thomas (RT) space of the solution modal basis, i.e.
\begin{equation}
    \transformed{\vect{\nabla}} \cdot \atflux{\fluxmodalbasis} \equiv \atsol{\modalbasis} .
\end{equation}
To build the flux nodal basis, a normal or degree of freedom $\{ \atflux{\unit{\transformed{\normal}}}_{\itype \ipol} \in \mathbb{R}^{\ndim} \}$ is assigned to each flux point.
The interested reader is referred to \ref{sec:rt_basis} for the definition of the RT basis, the location of flux points and the degrees of freedom assigned to each flux point within the different element types utilized in this work.
The flux nodal basis and the degree of freedoms fulfill the nodal condition
\begin{equation}
    \atflux{\fluxnodalbasis_{\itype \ipol}}(\atflux{\transformed{\vect{x}}}_{\itype \isigma}) \cdot \atflux{\unit{\transformed{\normal}}}_{\itype \isigma} = \delta_{\ipol \isigma} .
\end{equation}
Hence, the Vandermonde matrix of the flux basis can be obtained as
\begin{equation}
    \atflux{\vandermonde}_{\itype \ipol \isigma} = \atflux{\fluxmodalbasis}_{\itype \ipol} \left( \atflux{\transformed{\vect{x}}}_{\itype \isigma} \right) \cdot \atflux{\unit{\transformed{\normal}}}_{\itype \isigma} ,
\end{equation}
and its inverse, which links the flux nodal and modal bases as
\begin{equation}
    \atflux{\fluxnodalbasis_{\itype \ipol}}(\transformed{\vect{x}}) = \left(\atflux{\vandermonde_{\itype \ipol \isigma}}\right)^{-1} \atflux{\fluxmodalbasis}_{\itype \isigma}\left( \transformed{\vect{x}} \right) .
\end{equation}
Finding appropriate orthogonal or orthornormal bases to build the flux modal is crucial to reduce the condition number of the flux Vandermonde matrix, thereby avoiding unwanted machine round-off errors.
The interested reader is referred to \ref{sec:rt_basis} for more information on the choice of flux modal bases.
The nodal flux basis allows to express the transformed flux within a given element as
\begin{equation}
    \transformed{\flux}_{\itype \iele \ivar}(\transformed{\vect{x}}) = \atfluxnormal{\transformed{\fluxnovect}}_{\itype \iele \inu \ivar} \atflux{\fluxnodalbasis}_{\itype \inu}(\transformed{\vect{x}}) ,
    \label{eq:flux_nodal_basis_all_points}
\end{equation}
where
\begin{equation}
   \atfluxnormal{\transformed{\fluxnovect}}_{\itype \iele \inu \ivar} = \atflux{\transformed{\flux}}_{\itype \iele \inu \ivar} \cdot \atflux{\unit{\transformed{\normal}}}_{\itype \inu} \text{ and } \atflux{\transformed{\flux}}_{\itype \iele \inu \ivar} = \transformed{\flux}_{\itype \iele \ivar} ( \atflux{\transformed{\vect{x}}_{\itype \inu}} ).
\end{equation}

To ease the implementation and evaluation of the transformed fluxes at element interfaces, the flux points are divided into two sets: the external and internal flux points.
The former set lays on the boundaries of the reference element and its position in the reference element is $\{ \atfluxext{\transformed{\vect{x}}_{\itype \ipol}} \}$ with $0 \le \ipol < \atfluxext{N}_{\itype}$.
The associated flux nodal bases with the external flux points are labeled as $\{ \atfluxext{\fluxnodalbasis}_{\itype \ipol}\left( \transformed{\vect{x}} \right) \}$.
The external flux points are unique and only present a single degree of freedom which coincides with the outward-pointing unit normal vectors of the reference element at the considered external flux points.
Flux points which are internal to the element, referred to as internal flux points and denoted as $\{ \atfluxint{\transformed{\vect{x}}_{\itype \ipol}} \}$ with $0 \le \ipol < \atfluxint{N}_{\itype}$, may be duplicated, i.e. the locations $\{ \atfluxint{\transformed{\vect{x}}_{\itype \ipol}} \}$ might not be unique within a given element.
The associated flux nodal bases to the internal flux points are labeled as $\{ \atfluxint{\fluxnodalbasis}_{\itype \ipol}\left( \transformed{\vect{x}} \right) \}$.
Additionally, the set of unique flux points is referred to as $\{ \atfluxunique{\transformed{\vect{x}}_{\itype \ipol}} \}$ with $0 \le \ipol < \atfluxunique{N}_{\itype}$.
At last, to ease future developments, we define the set of unique flux points $\{ \atfluxintunique{\transformed{\vect{x}}_{\itype \ipol}} \}$ with $0 \le \ipol < \atfluxintunique{N}_{\itype}$.
With such a distinction between internal and external flux points, \Eref{eq:flux_nodal_basis_all_points} may be rewritten as
\begin{equation}
    \transformed{\flux}_{\itype \iele \ivar}(\transformed{\vect{x}}) = \atfluxnormalext{\transformed{\fluxnovect}}_{\itype \iele \inu \ivar} \atfluxext{\fluxnodalbasis}_{\itype \inu}(\transformed{\vect{x}}) + \atfluxnormalint{\transformed{\fluxnovect}}_{\itype \iele \ipol \ivar} \atfluxint{\fluxnodalbasis}_{\itype \ipol}(\transformed{\vect{x}})  .
    \label{eq:flux_nodal_basis_all_points_v2}
\end{equation}
To ensure conservation, the normal flux at external flux points $\atfluxnormalext{\transformed{\fluxnovect}}$ needs to be carefully assessed as it will be shown further in this section.
The SDRT approach differs from the original SDM due to the fact that the interpolatory polynomials are vectors, as opposed to the scalar nodal basis found in the original SDM developed for tensor-product elements.
Nevertheless, it may be demonstrated the SDRT approach for tensor-product cells is equivalent to that of the original SD method \cite{Li2019}.

The first step to compute the flux divergence in the SDRT approach for second-order conservation laws is the computation of the auxiliary gradient variable.
To do so, the conserved variables at solution points are interpolated to the external flux points
\begin{equation}
    \atfluxext{\var}_{\itype \inu \iele \ivar} = \atsol{\var}_{\itype \ipol \iele \ivar} \atsol{\nodalbasis}_{\iele \ipol} \left( \atfluxext{\transformed{\vect{x}}_{\itype \inu}} \right) ,
    \label{eq:interpolate_sol_to_flux}
\end{equation}
and to internal flux points
\begin{equation}
    \atfluxint{\var}_{\itype \inu \iele \ivar} = \atsol{\var}_{\itype \ipol \iele \ivar} \atsol{\nodalbasis}_{\iele \ipol} \left( \atfluxint{\transformed{\vect{x}}_{\itype \inu}} \right) .
    \label{eq:interpolate_sol_to_fluxint}
\end{equation}
Within the latter procedure, computational implementations may take advantage of the fact that only interpolations to the unique internal flux points have to be considered.
The interested reader is referred to \ref{sec:sdrt_matrix_form} for more details on such computational optimizations.

Next, a common value of the conserved variables is selected at external flux points.
\begin{equation}
    \commonvar_{\ivar} \atfluxext{\var}_{\itype \ipol \iele \ivar} = \commonvar_{\ivar} \atfluxext{\var}_{\transformed{\itype \ipol \iele \ivar}} = \commonvar_{\ivar} \left(  \atfluxext{\var}_{\itype \ipol \iele \ivar}, \atfluxext{\var}_{\transformed{\itype \ipol \iele \ivar}}  \right) ,
\end{equation}
where $\commonvar_{\ivar}(\var_L, \var_R)$ is a scalar function which returns a common value from two conserved variables at flux points.
Moreover, the sub-indices $\transformed{\itype \ipol \iele \ivar}$ refer to the element type, flux point and element number adjacent to the flux point $\itype \ipol \iele \ivar$.
Function $\commonvar_{\ivar}$ may be defined through the use of the Local-Discontinuous-Galerkin (LDG) method \cite{Cockburn1998} as
\begin{equation}
    \commonvar_{\ivar}(\var_L, \var_R) = \left( \frac{1}{2} - \beta \right) \var_{L} + \left( \frac{1}{2} + \beta \right) \var_{R} .
    \label{eq:common_var_for_grad}
\end{equation}
Values of $\beta = \pm \frac{1}{2}$ are of interest since it may promote compactness of the schemes in multiple dimensions, although the latter may not be ensured for certain cases on general grids \cite{Peraire2008, Williams2013tri}.
The non-compact Rebay-Bassi 1 (RB1) scheme \cite{Bassi1997} is recovered with $\beta = \frac{1}{2}$.
By using the common values of the conserved variables at flux points, the transformed auxiliary gradient variable is computed as
\begin{equation}
    \atsol{\transformed{\vargrad}}_{\itype \iele \ivar } ( \transformed{\vect{x}}_{\itype \isigma} ) = \atsol{\transformed{\vargrad}}_{\itype \isigma \iele \ivar } =  \atfluxint{\var}_{\itype \inu \iele \ivar}  \atfluxint{\unit{\transformed{\normal}}}_{\itype \inu} \left[\grad \cdot \atfluxint{\fluxnodalbasis}_{\itype \inu }  ( \transformed{\vect{x}}_{\itype \isigma} ) \right] + \commonvar_{\ivar} \atfluxext{\var}_{\itype \ipol \iele \ivar}  \atfluxext{\unit{\transformed{\normal}}}_{\itype \ipol} \left[\grad \cdot \atfluxext{\fluxnodalbasis}_{\itype \ipol }  ( \transformed{\vect{x}}_{\itype \isigma} ) \right] .
\end{equation}
The physical gradients at solution points are related to the reference space gradients as
\begin{equation}
    \atsol{\vargrad}_{\itype \isigma \iele \ivar} = \jacobian^{-T}_{\itype \isigma \iele} \atsol{\transformed{\vargrad}}_{\itype \isigma \iele \ivar } .
\end{equation}
Such gradients have to be interpolated to the flux points using \Eref{eq:interpolate_sol_to_flux} and \Eref{eq:interpolate_sol_to_fluxint} in order to evaluate the transformed fluxes.

At the internal flux points, the transformed flux is defined as
\begin{equation}
    \atfluxnormalint{\transformed{\fluxnovect}}_{\itype \iele \inu \ivar} = \left[\atfluxint{J}_{\itype \inu \iele} (\atfluxint{\jacobian_{\itype \inu \iele}})^{-1} \flux_{\itype \iele \ivar}(\atfluxint{\var}_{\itype \inu \iele}, \atfluxint{\vargrad}_{\itype \inu \iele})\right] \cdot \atfluxint{\unit{\transformed{\normal}}}_{\itype \inu} .
    \label{eq:normal_transformed_flux}
\end{equation}
It is worth mentioning that the implementation may take into account the presence of internal flux points sharing a unique location in the reference element to further optimize the computational performance.
The interested reader is referred to \ref{sec:sdrt_matrix_form} for more information on the implementation aspects and the matrix form of the SDRT method.
On the other hand, a common flux needs to be computed at external flux points to ensure conservation.
This is expressed as
\begin{equation}
    \commflux_{\ivar} \atfluxnormalext{\fluxnovect}_{\itype \isigma \iele \ivar} = - \commflux_{\ivar} \atfluxnormalext{\fluxnovect}_{\transformed{\itype \isigma \iele \ivar}} = \commflux_{\ivar}\left( \atfluxnormalext{\var}_{\itype \isigma \iele}, \atfluxnormalext{\var}_{\transformed{\itype \isigma \iele}}, \atfluxnormalext{\vargrad}_{\itype \isigma \iele}, \atfluxnormalext{\vargrad}_{\transformed{\itype \isigma \iele}},
    \atfluxnormalext{\unit{\vect{n}}}_{\itype \isigma \iele}\right) ,
\end{equation}
where the common flux function $\commflux_{\ivar}$ must be consistent and must present the conservation property, i.e. $\commflux_{\ivar} (\var_L, \var_R, \vargrad_L, \vargrad_R, \unit{\vect{n}}) = - \commflux_{\ivar} (\var_R, \var_L, \vargrad_R, \vargrad_L, -\unit{\vect{n}})$.
Moreover, $\atflux{\unit{\vect{n}}}_{\itype \isigma \iele}$ refers to the physical unit normal at a given external flux point.
The common flux function is usually split into two contributions: the common convective flux and the common viscous flux
\begin{equation}
    \commflux_{\ivar} (\var_L, \var_R, \vargrad_L, \vargrad_R, \unit{\vect{n}}_L) = \commflux^{\text{conv}} (\var_L, \var_R, \unit{\vect{n}}_L)  + \commflux^{\text{visc}} (\var_L, \var_R, \vargrad_L, \vargrad_R, \unit{\vect{n}}_L) .
\end{equation}
The former, which only depends on the conserved variables, may be computed through the Rusanov--Riemann solver
\begin{equation}
    \commflux^{\text{conv}} (\var_L, \var_R, \unit{\vect{n}}_L) = \frac{\unit{\vect{n}}_L}{2}\left( \flux_L^{\text{conv}} + \flux_R^{\text{conv}}  \right) + \frac{|\phi|}{2} \left( \var_L - \var_R \right) .
\end{equation}
where $\phi$ is an estimator of the maximum eigenvalue of the Jacobian of the flux operator.
It is worth mentioning that there exist common convective flux formulations which are specifically tailored to a given system of conservation laws, such as the HLLC \cite{Toro1999} for the Euler/Navier-Stokes equations.
The common viscous flux \cite{Witherden2014} is defined in this work as 
\begin{equation}
    \commflux^{\text{visc}} (\var_L, \var_R, \vargrad_L, \vargrad_R, \unit{\vect{n}}_L) = \unit{\vect{n}}_L \left[ \left( \frac{1}{2} + \beta \right) \flux_{L}^{\text{visc}} + \left( \frac{1}{2} - \beta \right) \flux_{L}^{\text{visc}} \right]  + \penaltyterm \left( \var_L - \var_R \right) .
\end{equation}
where $\penaltyterm$ is a penalty term of the LDG formulation.
The value of $\beta$ utilized in the latter equation is the same as the one used to compute the common value of the conserved variables at external flux points in \Eref{eq:common_var_for_grad}.

Once the value of the common flux at the external points has been assessed, the transformed flux at such points may be computed using
\begin{equation}
    \commflux_{\ivar} \atfluxnormalext{\transformed{\fluxnovect}}_{\itype \isigma \iele \ivar} = \atflux{J}_{\itype \isigma \iele} \norm{\atfluxext{\vect{n}}_{\itype \isigma \iele}} \commflux_{\ivar} \atfluxnormal{\fluxnovect}_{\itype \isigma \iele \ivar} .
\end{equation}

With the values of the transformed flux at internal and external flux points, the divergence at the solution points can be evaluated as
\begin{equation}
    \atsol{\left(\transformed{\vect{\nabla}} \cdot \transformed{\flux}_{\ivar}\right)}_{\itype \ipol \iele \ivar} = \atfluxnormalint{\transformed{\fluxnovect}}_{\itype \iele \inu \ivar} \left[\grad \cdot \atfluxint{\fluxnodalbasis}_{\itype \inu }  ( \vect{x}_{\itype \ipol} ) \right] + \commflux_{\ivar} \atfluxnormalext{\transformed{\fluxnovect}}_{\itype \iele \isigma \ivar} \left[\grad \cdot \atfluxext{\fluxnodalbasis}_{\itype \isigma }  ( \vect{x}_{\itype \ipol} ) \right] .
    \label{eq:divergence_sdrt}
\end{equation}
This allows to rewrite the governing system in a semi-discretized form through the method of lines
\begin{equation}
    \frac{\text{d} \atsol{\var}_{\itype \ipol \iele \ivar}}{\text{d} t} = - J_{\itype \ipol \iele}^{-1 \medskip (\var)} \atsol{\left(\transformed{\vect{\nabla}} \cdot \transformed{\flux}_{\ivar}\right)}_{\itype \ipol \iele \ivar} \quad .
\end{equation}
In this study, the previous equation is solved using explicit time-integration Runge--Kutta (ERK) methods.

\subsection{FR method}

Herein, the FR-DG formulation is described following the notation of \citet{Witherden2014}.
Within this section, it is supposed that the SDRT and FR-DG formulations share the same external flux points distribution.
Nevertheless, in flux-reconstruction methods, no internal flux points are found, i.e. $\atfluxint{N}_{\itype} = 0$.
Therefore, the interpolation procedure to such points is not needed.
To ensure conservation, the flux is split into a discontinuous part, which is represented by the solution nodal basis, and a continuous contribution, which takes into account the difference between the interpolated discontinuous flux at the external flux points and the common flux at those locations.
This may be expressed as
\begin{equation}
    \transformed{\flux}_{\itype \iele \ivar}(\transformed{\vect{x}}) = \atsol{\transformed{\flux}}_{\itype \ipol \iele \ivar} \atsol{\nodalbasis}_{\itype \ipol }(\transformed{\vect{x}}) + \left( \commflux_{\ivar} \atfluxnormalext{\transformed{\fluxnovect}}_{\itype \iele \inu \ivar} - \left[ \atsol{\transformed{\flux}}_{\itype \iele \isigma \ivar} \atsol{\nodalbasis}_{\itype \isigma}(\transformed{\vect{x}}_{\itype \inu}) \right] \cdot \atfluxext{\unit{\transformed{\normal}}}_{\itype \inu}  \right) \correctionnodalbasis_{\itype \inu}(\transformed{\vect{x}}) ,
    \label{eq:flux_fr}
\end{equation}
where $\correctionnodalbasis_{\itype \inu}(\transformed{\vect{x}})$ is a vector correction function whose divergence sits on the polynomial space of the solution nodal basis and which satisfies
\begin{equation}
    \correctionnodalbasis_{\itype \inu}(\transformed{\vect{x}}_{\itype \isigma}) \cdot \unit{\vect{\normal}}_{\itype \isigma} = \delta_{\inu \isigma} .
\end{equation}
With such a polynomial description of the flux, the computation of the divergence (needed to update the solution values at solution points) is straightforward.
The interested reader may refer to \cite{Witherden2014} for the full FR formulation specification.

Certain vector correction function formulations choices allow to recover a nodal DG scheme from the FR formalism \cite{Wang2009, Castonguay2011, Williams2013tetra} for arbitrary simplex elements.
In this work, schemes using the latter correction functions are referred to as FR-DG schemes.
It is worth mentioning that there exist other choices of correction functions, usually referred to as Vincent-Castonguay-Jameson-Huynh (VCJH) schemes, for FR and triangular \cite{Castonguay2011} and tetrahedral \cite{Williams2013tetra} elements which yield energy-stable schemes. In that framework, FR-DG turns out to be the limit case yielding an energy-stable scheme. 

\subsection{Connections between the SDRT and FR formalisms}
\label{sec:sdrt_fr_equivalence}

Both the correction function of the FR formalism and the flux nodal basis SDRT methods lay on the Raviart-Thomas space of the solution basis.
Nevertheless, there exists an important differentiation aspect between FR and SDRT related to the fact that, in FR, the correction functions are not uniquely defined, since there exist fewer external flux points than the cardinal of the RT basis.
Additionally, the FR method approximates the discontinuous flux using the solution nodal basis, while SDRT approximates the latter using the flux nodal basis (which lays on the RT space).
This remark was raised by \cite{Cox2021} as a possible source of aliasing errors of FR  due to non-linearities in regards to SDRT. 

This section aims to propose a link between the FR and SDRT formulations for all the elements types studied in this work, with constant metrics and linear flux.
Such a connection is established using the RT basis of the SDRT method as the correction function, resulting in the FR-SDRT method.
It is worth mentioning that this link had already been established by Huynh \cite{Huynh2007} for tensor-product elements.
Within this section, it is supposed that both the FR and SDRT schemes use the same external flux points and solution points locations.
To demonstrate the analogy between FR-SDRT and SDRT methods in linear test cases with constant metric elements,
let us expand the flux divergence in the FR formulation points at solution points (see \Eref{eq:flux_fr})
\begin{equation}
    \atsol{\left(\transformed{\grad} \cdot \transformed{\flux}\right)}_{\itype \iele \ivar}(\atsol{\transformed{\vect{x}}_{\isigma}})
    =
    \atsol{\transformed{\flux}}_{\itype \ipol \iele \ivar} \cdot  \transformed{\grad}\atsol{\nodalbasis}_{\itype \ipol }(\atsol{\transformed{\vect{x}}_{\isigma}})
    +
    \commflux_{\ivar} \atfluxnormalext{\transformed{\fluxnovect}}_{\itype \iele \inu \ivar} \transformed{\grad} \cdot \correctionnodalbasis_{\itype \inu}(\transformed{\vect{x}}_{\itype \isigma})
    -
    \left[ \atsol{\transformed{\flux}}_{\itype \iele \ipol \ivar} \atsol{\nodalbasis}_{\itype \ipol}(\vect{x}_{\itype \inu}) \right] \cdot \atfluxext{\unit{\transformed{\normal}}}_{\itype \inu} \transformed{\grad} \cdot \correctionnodalbasis_{\itype \inu}(\atsol{\transformed{\vect{x}}_{\isigma}}) .
    \label{eq:divergence_fr}
\end{equation}
If one imposes that the vector correction functions divergence is equal to that of the RT bases associated with the external flux points of the SDRT method then
\begin{equation}
    \transformed{\grad} \cdot  \correctionnodalbasis_{\itype \inu}(\atsol{\transformed{\vect{x}}_{\isigma}}) = \transformed{\grad} \cdot \atfluxext{\fluxnodalbasis}_{\itype \inu}(\atsol{\transformed{\vect{x}}_{\isigma}}) .
\end{equation}
The resulting FR method with such a choice is referred to as FR-SDRT in this work.
The equivalence between FR-SDRT and SDRT schemes may be proven if both Eqs. \ref{eq:divergence_sdrt} and \ref{eq:divergence_fr} are equal.
Hence, if the external flux points location remains the same when using SDRT and FR-SDRT schemes then
\begin{equation}
    \atsol{\transformed{\flux}}_{\itype \ipol \iele \ivar} \cdot \transformed{\grad} \atsol{\nodalbasis}_{\itype \ipol }(\transformed{\vect{x}}_{\itype \isigma})
    -
    \left[ \atsol{\transformed{\flux}}_{\itype \iele \ipol \ivar} \atsol{\nodalbasis}_{\itype \ipol}(\transformed{\vect{x}}_{\itype \inu}) \right] \cdot \atfluxext{\unit{\transformed{\normal}}}_{\itype \inu} \transformed{\grad} \cdot \atfluxext{\fluxnodalbasis}_{\itype \inu}(\transformed{\vect{x}}_{\itype \isigma})
    =
    \atfluxnormalint{\transformed{\fluxnovect}}_{\itype \iele \inu \ivar} \left[\grad \cdot \atfluxint{\fluxnodalbasis}_{\itype \inu }  ( \transformed{\vect{x}}_{\itype \isigma} )\right] .
\end{equation}
The latter equality reduces to
\begin{equation}
    \transformed{\grad} \cdot ( \atsol{\transformed{\flux}}_{\itype \ipol \iele \ivar}  \atsol{\nodalbasis}_{ \itype \ipol })(\atsol{\transformed{\vect{x}}_{\itype \isigma}})
    =
    \transformed{\grad} \cdot (\atfluxnormal{\transformed{\fluxnovect}}_{\itype \inu \iele \ivar}  \atfluxnormal{\fluxnodalbasis}_{\itype \inu})(\atsol{\transformed{\vect{x}}_{\itype \isigma}}) ,
    \label{eq:divergence_equivalence}
\end{equation}
where both nodal expansions represent the polynomial projection of the discontinuous fluxes, which do not take into account the common fluxes, onto the solution and RT bases of a given element.
Such an ansatz may be fulfilled in all elements presenting constant metrics and if the flux is linear.
With these conditions, the flux may be exactly approximated by, at most, a degree $\degree$ modal basis, hence both the RT and solution nodal basis are equivalent, despite the fact that the former uses vector interpolation basis of degree $\degree + 1$ and the latter uses scalar interpolation basis of degree $\degree$.
As the divergence is analytically computed based on the values of the nodal bases, both sides of the ansatz will be equal under the previously described conditions.
\Sref{sec:numerical_experiments_linadvecdiff_equivalence_sdrt_fr} will demonstrate through numerical experiments that FR-SDRT and SDRT schemes are equivalent in linear advection-diffusion cases and two and three-dimensional elements.

\begin{remark}
To further prove the equivalence between the FR-SDRT and SDRT formulations for linear advection-diffusion problems, one must ensure that the gradient computation at the solution points is equivalent in both formulations.
This can be ensured if one uses the same correction functions than that of the fluxes to compute the auxiliary gradient variables.
Such observation may be demonstrated following the same process utilized to determine the flux equivalence in \Eref{eq:divergence_equivalence}, since the calculation of the auxiliary gradient variables follow a similar correction procedure than the fluxes and since the conservative variables are described by polynomials of degree $\degree$ within each element.
This demonstration is left as an exercise to the reader.
\end{remark}

\begin{remark}
Despite the fact that FR-SDRT schemes may be shown to be energy-stable for tensor-product elements \cite{Zwanenburg2016}, it has been observed (through a posteriori analysis) that FR-SDRT schemes do not belong to the family of Vincent–Castonguay–Jameson–Huynh (VCJH) schemes neither for triangles \cite{Castonguay2011} nor for tetrahedrons  \cite{Williams2013tetra}.
Hence, it is thought that the energy stability of FR-SDRT for simplex elements can not be proven using the tools developed in the aforementioned studies.
The a priori demonstration of the energy stability of the FR-SDRT method for simplex elements and its possible connections with filtered DG formalism will be a topic of future research.
\end{remark}

\section{Von-Neumann Analysis}
\label{sec:von_neumann_advection}

This section aims to analyze the dissipation and dispersion errors in linear advection and linear diffusion problems of SEM in two-dimensional and three-dimensional grids, composed of tensor-product or simplex elements.
In particular, a comparison between the SDRT (resp. FR-SDRT) and FR-DG methods will be presented.
The analysis of the dissipation and dispersion of SEM is important to study their associated numerical errors in convective-dominated problems such as those found in turbulent flows.
This analysis can be performed using a modified version of the classical Von-Neumann analysis for FDM and FVM, applied to the numerical solution of linear equations cf.\ \cite{Pereira2020elements,Trojak2020}.

Let us start the analysis by considering the linear advection equation, solved in the periodic domain $\vect{x} \in \domain \equiv [0, L]^3$
\begin{equation}
    \frac{\partial \var(\vect{x}, t)}{\partial t} + \divergence (\advecvel \var(\vect{x}, t)) = 0 ,
    \label{eq:linear_advection}
\end{equation}
where $\advecvel$ is the advection velocity (supposed constant) and defined as 
\begin{equation}
    \advecvel = \advecvelnovect \unit{\unitary}_{\advecvel} ,
\end{equation}
being $\unit{\unitary}_{\advecvel}$ a unitary $\ndim$-dimensional vector computed in three-dimensional configurations as 
\begin{equation}
    \unit{\unitary}_{\advecvel} (\theta_0, \theta_1) = \left( \cos\theta_0 \cos\theta_1, \sin\theta_0, \cos\theta_0 \sin\theta_1 \right) .
\end{equation}
The two-dimensional definition of the latter vector may be obtained by imposing $\theta_1 = 0$.

The exact solution of the linear advection equation is given by the expression
\begin{equation}
    \var(\vect{x}, t) = \var(\vect{x} - \advecvel t, 0) = \var_0(\vect{x} - \advecvel t) ,
\end{equation}
where $\var_0$ is the initial condition.
Supposing that the initial solution is a Fourier mode, i.e. $\var = \e^{\imag \vect{\kappa} \cdot \vect{x}}$, where $\vect{\kappa} = \kappa \unit{\unitary}_{\vect{\kappa}} $ is the spatial wavenumber, then the analytical solution may be expressed as
\begin{equation}
    \var(\vect{x}, t) = \e^{-\imag \omega t} \var_0(\vect{x}) = \e^{-\imag \omega t + \imag \vect{\kappa} \cdot \vect{x}} ,
\end{equation}
with $\omega = \vect{\kappa} \cdot \advecvel$ being the temporal wavenumber.
For future reference, let us define the wavelength as $\lambda = 2\pi/\kappa$.

\begin{figure}
\centering
\begin{subfigure}{0.45\textwidth}
    \centering
    \includegraphics[width=0.95\textwidth]{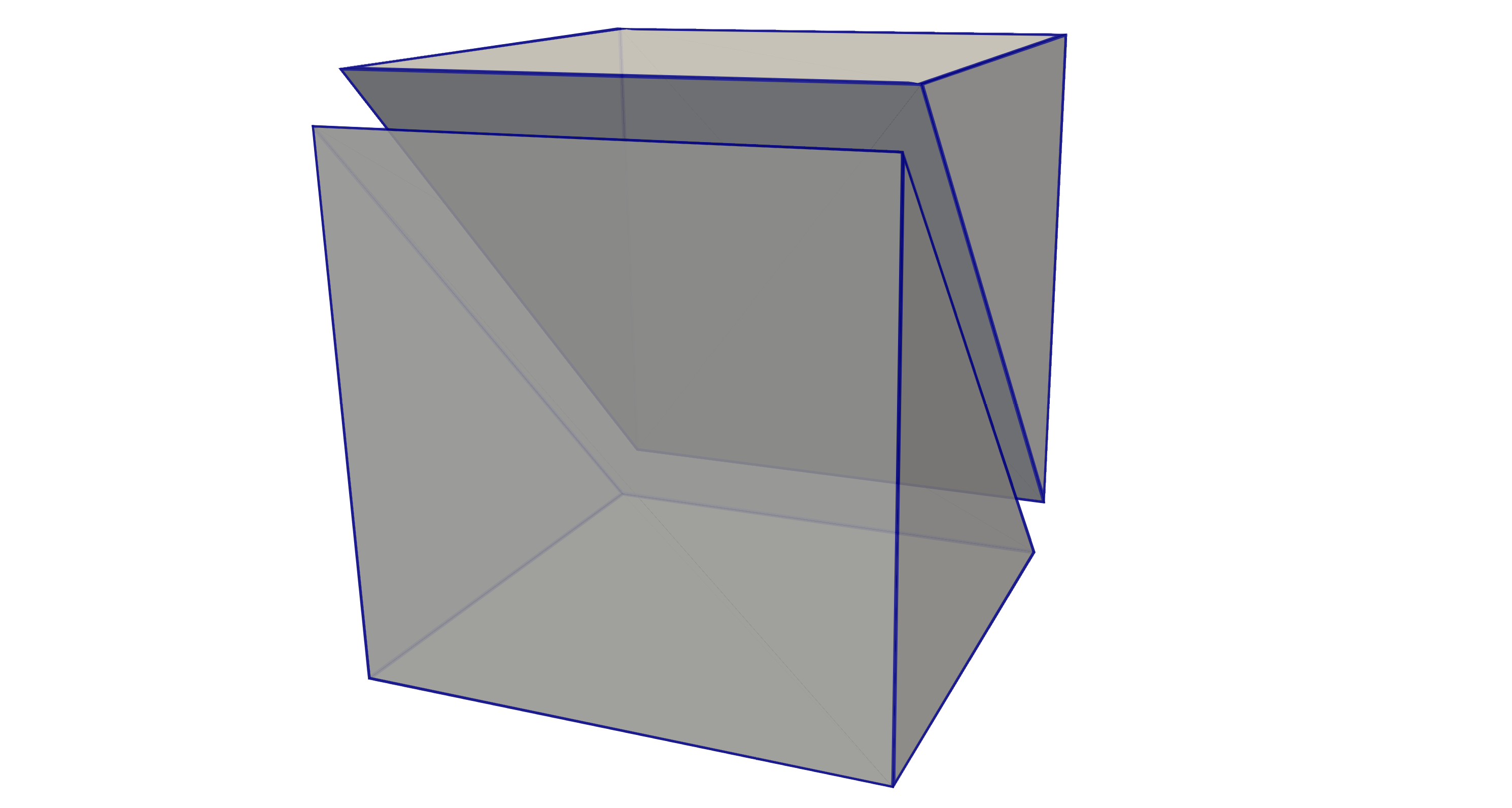}
    \caption{Prisms ($\ncellsinpattern=2$)}
    \label{fig:subdivision_pri}
\end{subfigure}
\begin{subfigure}{0.45\textwidth}
    \centering
    \includegraphics[width=0.95\textwidth]{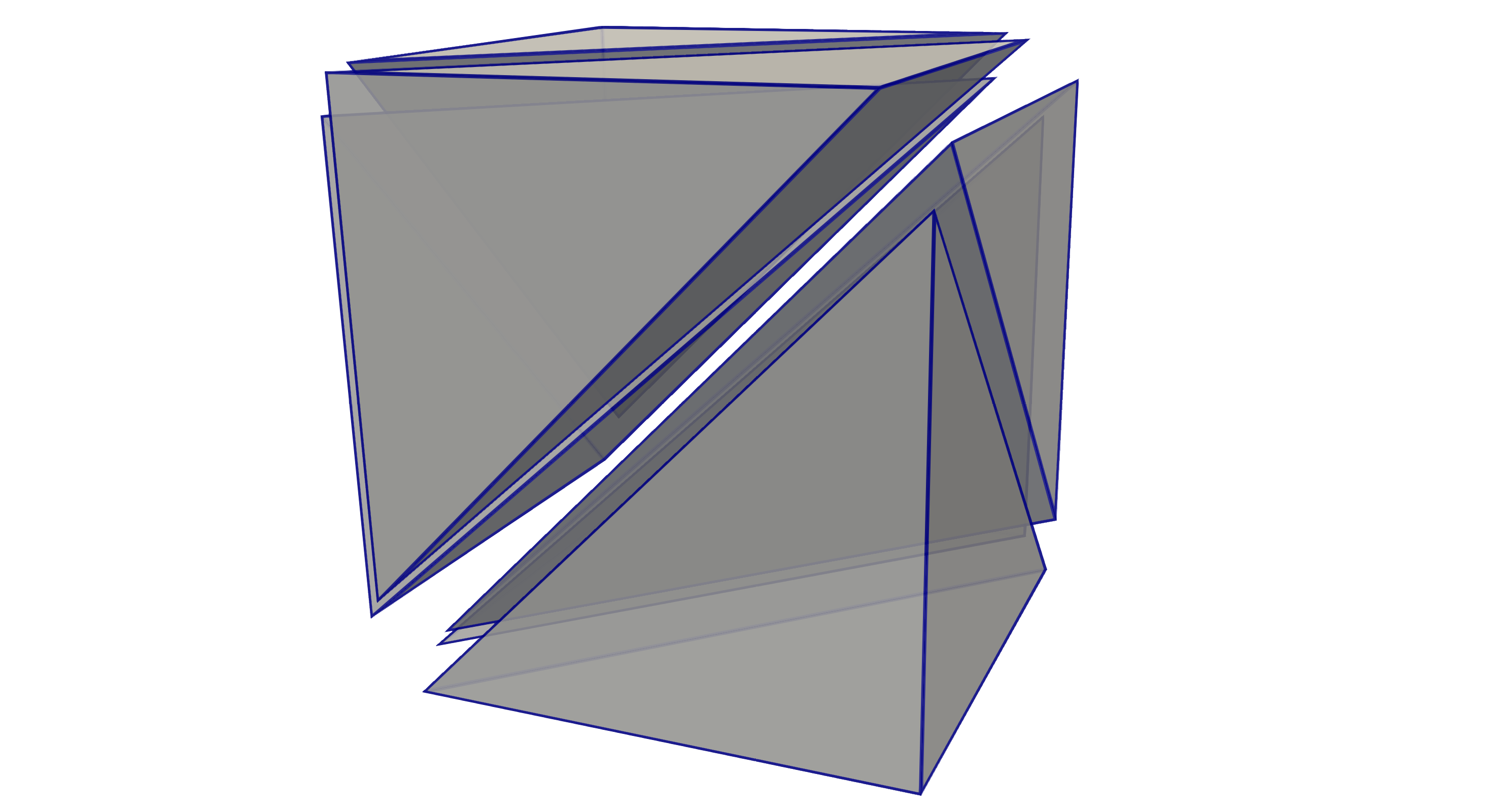}
    \caption{Tetrahedron ($\ncellsinpattern=6$)}
    \label{fig:subdivision_tet}
\end{subfigure}
\caption{Sub-division procedure of a given hexahedron into prismatic (left) and tetrahedron (right) elements.}
\label{fig:subdivision_3D}
\end{figure}

The common flux operator for the linear advection equation is defined in this study using the upwind Lax--Friedrich flux 
\begin{equations}{eq:flux_advection}
    \commflux\left( \var_L, \var_R, \unit{\normal}_L \right) &=
    \frac{\advecvel \cdot \unit{\normal}_L}{2}\left( \var_L + \var_R \right) + \frac{|\advecvel \cdot \unit{\normal}_L|}{2}  \left( \var_L - \var_R \right) \\
    &=\advecvel \cdot \unit{\normal}_L \left[ \frac{(1 + \text{sign}(\advecvel \cdot \unit{\normal}_L))}{2} \var_L + \frac{(1 - \text{sign}(\advecvel \cdot \unit{\normal}_L))}{2} \var_R \right] .
\end{equations}
To solve numerically \Eref{eq:linear_advection}, the domain is triangulated with a uniform mesh of edge size $\cellsize$.
Uniform meshes of simplex elements are generated by sub-dividing the initial mesh made up of tensor-product elements into simplex elements using a similar approach to that shown in \cite{Pereira2020elements}.
The number of sub-divided cells resulting from a unique tensor-product element is referred to as $\ncellsinpattern$ and its value is: $\ncellsinpattern = 1$ for tensor-product elements, $\ncellsinpattern = 2$ for prism and triangular elements and $\ncellsinpattern = 6$ for tetrahedron elements.
See \Fref{fig:subdivision_3D} for a sketch of the subdivision mechanism of hexahedron elements into prismatic and tetrahedron cells.
In what follows, the sub-index $\itype$, referring to the type of element will be dropped as uniform meshes (with a unique element type) are supposed.

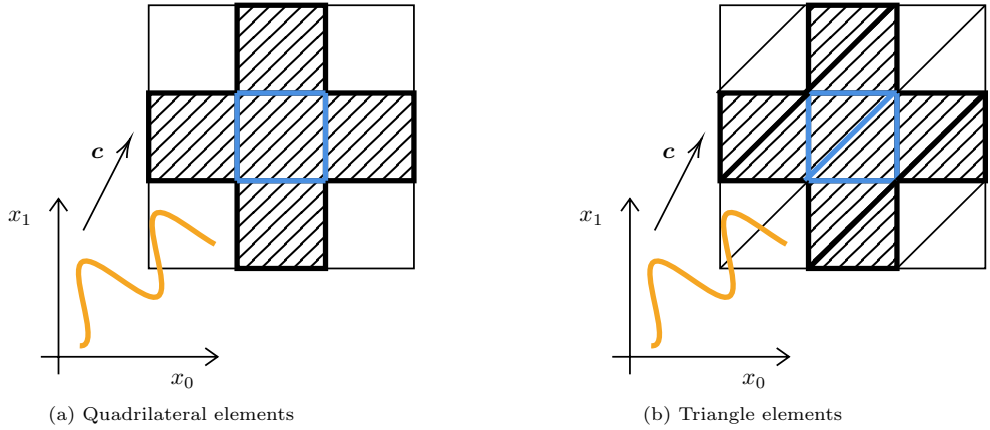
\begin{figure}
    \centering
    \begin{subfigure}{0.45\textwidth}
    \centering
    \adjustbox{width=0.9\textwidth}{
 
\tikzset{
pattern size/.store in=\mcSize, 
pattern size = 5pt,
pattern thickness/.store in=\mcThickness, 
pattern thickness = 0.3pt,
pattern radius/.store in=\mcRadius, 
pattern radius = 1pt}
\makeatletter
\pgfutil@ifundefined{pgf@pattern@name@_nuqpzg2j7}{
\pgfdeclarepatternformonly[\mcThickness,\mcSize]{_nuqpzg2j7}
{\pgfqpoint{0pt}{0pt}}
{\pgfpoint{\mcSize+\mcThickness}{\mcSize+\mcThickness}}
{\pgfpoint{\mcSize}{\mcSize}}
{
\pgfsetcolor{\tikz@pattern@color}
\pgfsetlinewidth{\mcThickness}
\pgfpathmoveto{\pgfqpoint{0pt}{0pt}}
\pgfpathlineto{\pgfpoint{\mcSize+\mcThickness}{\mcSize+\mcThickness}}
\pgfusepath{stroke}
}}
\makeatother

 
\tikzset{
pattern size/.store in=\mcSize, 
pattern size = 5pt,
pattern thickness/.store in=\mcThickness, 
pattern thickness = 0.3pt,
pattern radius/.store in=\mcRadius, 
pattern radius = 1pt}
\makeatletter
\pgfutil@ifundefined{pgf@pattern@name@_3h5gzqq97}{
\pgfdeclarepatternformonly[\mcThickness,\mcSize]{_3h5gzqq97}
{\pgfqpoint{0pt}{0pt}}
{\pgfpoint{\mcSize+\mcThickness}{\mcSize+\mcThickness}}
{\pgfpoint{\mcSize}{\mcSize}}
{
\pgfsetcolor{\tikz@pattern@color}
\pgfsetlinewidth{\mcThickness}
\pgfpathmoveto{\pgfqpoint{0pt}{0pt}}
\pgfpathlineto{\pgfpoint{\mcSize+\mcThickness}{\mcSize+\mcThickness}}
\pgfusepath{stroke}
}}
\makeatother

 
\tikzset{
pattern size/.store in=\mcSize, 
pattern size = 5pt,
pattern thickness/.store in=\mcThickness, 
pattern thickness = 0.3pt,
pattern radius/.store in=\mcRadius, 
pattern radius = 1pt}
\makeatletter
\pgfutil@ifundefined{pgf@pattern@name@_obu523b76}{
\pgfdeclarepatternformonly[\mcThickness,\mcSize]{_obu523b76}
{\pgfqpoint{0pt}{0pt}}
{\pgfpoint{\mcSize+\mcThickness}{\mcSize+\mcThickness}}
{\pgfpoint{\mcSize}{\mcSize}}
{
\pgfsetcolor{\tikz@pattern@color}
\pgfsetlinewidth{\mcThickness}
\pgfpathmoveto{\pgfqpoint{0pt}{0pt}}
\pgfpathlineto{\pgfpoint{\mcSize+\mcThickness}{\mcSize+\mcThickness}}
\pgfusepath{stroke}
}}
\makeatother

 
\tikzset{
pattern size/.store in=\mcSize, 
pattern size = 5pt,
pattern thickness/.store in=\mcThickness, 
pattern thickness = 0.3pt,
pattern radius/.store in=\mcRadius, 
pattern radius = 1pt}
\makeatletter
\pgfutil@ifundefined{pgf@pattern@name@_rsbga484v}{
\pgfdeclarepatternformonly[\mcThickness,\mcSize]{_rsbga484v}
{\pgfqpoint{0pt}{0pt}}
{\pgfpoint{\mcSize+\mcThickness}{\mcSize+\mcThickness}}
{\pgfpoint{\mcSize}{\mcSize}}
{
\pgfsetcolor{\tikz@pattern@color}
\pgfsetlinewidth{\mcThickness}
\pgfpathmoveto{\pgfqpoint{0pt}{0pt}}
\pgfpathlineto{\pgfpoint{\mcSize+\mcThickness}{\mcSize+\mcThickness}}
\pgfusepath{stroke}
}}
\makeatother

 
\tikzset{
pattern size/.store in=\mcSize, 
pattern size = 5pt,
pattern thickness/.store in=\mcThickness, 
pattern thickness = 0.3pt,
pattern radius/.store in=\mcRadius, 
pattern radius = 1pt}
\makeatletter
\pgfutil@ifundefined{pgf@pattern@name@_d9okvvcn8}{
\pgfdeclarepatternformonly[\mcThickness,\mcSize]{_d9okvvcn8}
{\pgfqpoint{0pt}{0pt}}
{\pgfpoint{\mcSize+\mcThickness}{\mcSize+\mcThickness}}
{\pgfpoint{\mcSize}{\mcSize}}
{
\pgfsetcolor{\tikz@pattern@color}
\pgfsetlinewidth{\mcThickness}
\pgfpathmoveto{\pgfqpoint{0pt}{0pt}}
\pgfpathlineto{\pgfpoint{\mcSize+\mcThickness}{\mcSize+\mcThickness}}
\pgfusepath{stroke}
}}
\makeatother
\tikzset{every picture/.style={line width=0.75pt}} 

\begin{tikzpicture}[x=0.75pt,y=0.75pt,yscale=-1,xscale=1]

\draw   (151,99.67) -- (201,99.67) -- (201,149.67) -- (151,149.67) -- cycle ;
\draw   (201,99.67) -- (251,99.67) -- (251,149.67) -- (201,149.67) -- cycle ;
\draw   (151,49.67) -- (201,49.67) -- (201,99.67) -- (151,99.67) -- cycle ;
\draw   (201,49.67) -- (251,49.67) -- (251,99.67) -- (201,99.67) -- cycle ;
\draw   (101,49.67) -- (151,49.67) -- (151,99.67) -- (101,99.67) -- cycle ;
\draw  [pattern=_nuqpzg2j7,pattern size=6pt,pattern thickness=0.75pt,pattern radius=0pt, pattern color={rgb, 255:red, 0; green, 0; blue, 0}][line width=2.25]  (101,99.67) -- (151,99.67) -- (151,149.67) -- (101,149.67) -- cycle ;
\draw  [pattern=_3h5gzqq97,pattern size=6pt,pattern thickness=0.75pt,pattern radius=0pt, pattern color={rgb, 255:red, 0; green, 0; blue, 0}][line width=2.25]  (151,149.67) -- (201,149.67) -- (201,199.67) -- (151,199.67) -- cycle ;
\draw   (101,149.67) -- (151,149.67) -- (151,199.67) -- (101,199.67) -- cycle ;
\draw   (201,149.67) -- (251,149.67) -- (251,199.67) -- (201,199.67) -- cycle ;
\draw  [pattern=_obu523b76,pattern size=6pt,pattern thickness=0.75pt,pattern radius=0pt, pattern color={rgb, 255:red, 0; green, 0; blue, 0}][line width=2.25]  (201,99.67) -- (251,99.67) -- (251,149.67) -- (201,149.67) -- cycle ;
\draw  [pattern=_rsbga484v,pattern size=6pt,pattern thickness=0.75pt,pattern radius=0pt, pattern color={rgb, 255:red, 0; green, 0; blue, 0}][line width=2.25]  (151,49.67) -- (201,49.67) -- (201,99.67) -- (151,99.67) -- cycle ;
\draw  [pattern=_d9okvvcn8,pattern size=6pt,pattern thickness=0.75pt,pattern radius=0pt, pattern color={rgb, 255:red, 0; green, 0; blue, 0}] (151,99.67) -- (201,99.67) -- (201,149.67) -- (151,149.67) -- cycle ;
\draw [color={rgb, 255:red, 74; green, 144; blue, 226 }  ,draw opacity=1 ][line width=2.25]    (151,99.67) -- (201,99.67) ;
\draw [color={rgb, 255:red, 74; green, 144; blue, 226 }  ,draw opacity=1 ][line width=2.25]    (151,149.67) -- (201,149.67) ;
\draw [color={rgb, 255:red, 74; green, 144; blue, 226 }  ,draw opacity=1 ][line width=2.25]    (151,149.67) -- (151,99.67) ;
\draw [color={rgb, 255:red, 74; green, 144; blue, 226 }  ,draw opacity=1 ][line width=2.25]    (201,149.67) -- (201,99.67) ;
\draw  (40,250) -- (140,250)(50,160) -- (50,260) (133,245) -- (140,250) -- (133,255) (45,167) -- (50,160) -- (55,167)  ;
\draw  [color={rgb, 255:red, 245; green, 166; blue, 35 }  ,draw opacity=1 ][line width=2.25]  (62.03,243.73) .. controls (63.06,243.81) and (63.95,243.64) .. (64.66,243.16) .. controls (68.43,240.66) and (66.28,230.38) .. (64.01,219.59) .. controls (61.73,208.8) and (59.59,198.52) .. (63.36,196.01) .. controls (67.12,193.51) and (75.77,199.48) .. (84.83,205.76) .. controls (93.9,212.04) and (102.54,218) .. (106.31,215.5) .. controls (110.08,213) and (107.93,202.72) .. (105.66,191.93) .. controls (103.39,181.14) and (101.24,170.85) .. (105.01,168.35) .. controls (108.77,165.85) and (117.42,171.82) .. (126.48,178.09) .. controls (130.74,181.04) and (134.9,183.92) .. (138.51,185.89) ;
\draw    (64,178) -- (89.76,127.28) ;
\draw [shift={(90.67,125.5)}, rotate = 476.93] [color={rgb, 255:red, 0; green, 0; blue, 0 }  ][line width=0.75]    (10.93,-3.29) .. controls (6.95,-1.4) and (3.31,-0.3) .. (0,0) .. controls (3.31,0.3) and (6.95,1.4) .. (10.93,3.29)   ;

\draw (112,258.4) node [anchor=north west][inner sep=0.75pt]    {$x_{0}$};
\draw (20,165.4) node [anchor=north west][inner sep=0.75pt]    {$x_{1}$};
\draw (67,129.4) node [anchor=north west][inner sep=0.75pt]    {$\advecvel$};

\end{tikzpicture}}
    \caption{Quadrilateral elements}
    \label{fig:quad_patterns}
    \end{subfigure}
    \begin{subfigure}{0.45\textwidth}
    \centering
    \adjustbox{width=0.9\textwidth}{
 
\tikzset{
pattern size/.store in=\mcSize, 
pattern size = 5pt,
pattern thickness/.store in=\mcThickness, 
pattern thickness = 0.3pt,
pattern radius/.store in=\mcRadius, 
pattern radius = 1pt}
\makeatletter
\pgfutil@ifundefined{pgf@pattern@name@_w1yfztebq}{
\pgfdeclarepatternformonly[\mcThickness,\mcSize]{_w1yfztebq}
{\pgfqpoint{0pt}{0pt}}
{\pgfpoint{\mcSize+\mcThickness}{\mcSize+\mcThickness}}
{\pgfpoint{\mcSize}{\mcSize}}
{
\pgfsetcolor{\tikz@pattern@color}
\pgfsetlinewidth{\mcThickness}
\pgfpathmoveto{\pgfqpoint{0pt}{0pt}}
\pgfpathlineto{\pgfpoint{\mcSize+\mcThickness}{\mcSize+\mcThickness}}
\pgfusepath{stroke}
}}
\makeatother

 
\tikzset{
pattern size/.store in=\mcSize, 
pattern size = 5pt,
pattern thickness/.store in=\mcThickness, 
pattern thickness = 0.3pt,
pattern radius/.store in=\mcRadius, 
pattern radius = 1pt}
\makeatletter
\pgfutil@ifundefined{pgf@pattern@name@_whfog6t2l}{
\pgfdeclarepatternformonly[\mcThickness,\mcSize]{_whfog6t2l}
{\pgfqpoint{0pt}{0pt}}
{\pgfpoint{\mcSize+\mcThickness}{\mcSize+\mcThickness}}
{\pgfpoint{\mcSize}{\mcSize}}
{
\pgfsetcolor{\tikz@pattern@color}
\pgfsetlinewidth{\mcThickness}
\pgfpathmoveto{\pgfqpoint{0pt}{0pt}}
\pgfpathlineto{\pgfpoint{\mcSize+\mcThickness}{\mcSize+\mcThickness}}
\pgfusepath{stroke}
}}
\makeatother

 
\tikzset{
pattern size/.store in=\mcSize, 
pattern size = 5pt,
pattern thickness/.store in=\mcThickness, 
pattern thickness = 0.3pt,
pattern radius/.store in=\mcRadius, 
pattern radius = 1pt}
\makeatletter
\pgfutil@ifundefined{pgf@pattern@name@_yajbjt0ue}{
\pgfdeclarepatternformonly[\mcThickness,\mcSize]{_yajbjt0ue}
{\pgfqpoint{0pt}{0pt}}
{\pgfpoint{\mcSize+\mcThickness}{\mcSize+\mcThickness}}
{\pgfpoint{\mcSize}{\mcSize}}
{
\pgfsetcolor{\tikz@pattern@color}
\pgfsetlinewidth{\mcThickness}
\pgfpathmoveto{\pgfqpoint{0pt}{0pt}}
\pgfpathlineto{\pgfpoint{\mcSize+\mcThickness}{\mcSize+\mcThickness}}
\pgfusepath{stroke}
}}
\makeatother

 
\tikzset{
pattern size/.store in=\mcSize, 
pattern size = 5pt,
pattern thickness/.store in=\mcThickness, 
pattern thickness = 0.3pt,
pattern radius/.store in=\mcRadius, 
pattern radius = 1pt}
\makeatletter
\pgfutil@ifundefined{pgf@pattern@name@_wflj3nk4r}{
\pgfdeclarepatternformonly[\mcThickness,\mcSize]{_wflj3nk4r}
{\pgfqpoint{0pt}{0pt}}
{\pgfpoint{\mcSize+\mcThickness}{\mcSize+\mcThickness}}
{\pgfpoint{\mcSize}{\mcSize}}
{
\pgfsetcolor{\tikz@pattern@color}
\pgfsetlinewidth{\mcThickness}
\pgfpathmoveto{\pgfqpoint{0pt}{0pt}}
\pgfpathlineto{\pgfpoint{\mcSize+\mcThickness}{\mcSize+\mcThickness}}
\pgfusepath{stroke}
}}
\makeatother

 
\tikzset{
pattern size/.store in=\mcSize, 
pattern size = 5pt,
pattern thickness/.store in=\mcThickness, 
pattern thickness = 0.3pt,
pattern radius/.store in=\mcRadius, 
pattern radius = 1pt}
\makeatletter
\pgfutil@ifundefined{pgf@pattern@name@_z2m3zekzc}{
\pgfdeclarepatternformonly[\mcThickness,\mcSize]{_z2m3zekzc}
{\pgfqpoint{0pt}{0pt}}
{\pgfpoint{\mcSize+\mcThickness}{\mcSize+\mcThickness}}
{\pgfpoint{\mcSize}{\mcSize}}
{
\pgfsetcolor{\tikz@pattern@color}
\pgfsetlinewidth{\mcThickness}
\pgfpathmoveto{\pgfqpoint{0pt}{0pt}}
\pgfpathlineto{\pgfpoint{\mcSize+\mcThickness}{\mcSize+\mcThickness}}
\pgfusepath{stroke}
}}
\makeatother
\tikzset{every picture/.style={line width=0.75pt}} 

\begin{tikzpicture}[x=0.75pt,y=0.75pt,yscale=-1,xscale=1]

\draw   (172,100) -- (222,100) -- (222,150) -- (172,150) -- cycle ;
\draw   (222,100) -- (272,100) -- (272,150) -- (222,150) -- cycle ;
\draw   (172,50) -- (222,50) -- (222,100) -- (172,100) -- cycle ;
\draw   (222,50) -- (272,50) -- (272,100) -- (222,100) -- cycle ;
\draw   (122,50) -- (172,50) -- (172,100) -- (122,100) -- cycle ;
\draw  [pattern=_w1yfztebq,pattern size=6pt,pattern thickness=0.75pt,pattern radius=0pt, pattern color={rgb, 255:red, 0; green, 0; blue, 0}][line width=2.25]  (122,100) -- (172,100) -- (172,150) -- (122,150) -- cycle ;
\draw  [pattern=_whfog6t2l,pattern size=6pt,pattern thickness=0.75pt,pattern radius=0pt, pattern color={rgb, 255:red, 0; green, 0; blue, 0}][line width=2.25]  (172,150) -- (222,150) -- (222,200) -- (172,200) -- cycle ;
\draw   (122,150) -- (172,150) -- (172,200) -- (122,200) -- cycle ;
\draw   (222,150) -- (272,150) -- (272,200) -- (222,200) -- cycle ;
\draw  [pattern=_yajbjt0ue,pattern size=6pt,pattern thickness=0.75pt,pattern radius=0pt, pattern color={rgb, 255:red, 0; green, 0; blue, 0}][line width=2.25]  (222,100) -- (272,100) -- (272,150) -- (222,150) -- cycle ;
\draw  [pattern=_wflj3nk4r,pattern size=6pt,pattern thickness=0.75pt,pattern radius=0pt, pattern color={rgb, 255:red, 0; green, 0; blue, 0}][line width=2.25]  (172,50) -- (222,50) -- (222,100) -- (172,100) -- cycle ;
\draw  [pattern=_z2m3zekzc,pattern size=6pt,pattern thickness=0.75pt,pattern radius=0pt, pattern color={rgb, 255:red, 0; green, 0; blue, 0}] (172,100) -- (222,100) -- (222,150) -- (172,150) -- cycle ;
\draw [color={rgb, 255:red, 74; green, 144; blue, 226 }  ,draw opacity=1 ][line width=2.25]    (172,100) -- (222,100) ;
\draw [color={rgb, 255:red, 74; green, 144; blue, 226 }  ,draw opacity=1 ][line width=2.25]    (172,150) -- (222,150) ;
\draw [color={rgb, 255:red, 74; green, 144; blue, 226 }  ,draw opacity=1 ][line width=2.25]    (172,150) -- (172,100) ;
\draw [color={rgb, 255:red, 74; green, 144; blue, 226 }  ,draw opacity=1 ][line width=2.25]    (222,150) -- (222,100) ;
\draw  (61,250.33) -- (161,250.33)(71,160.33) -- (71,260.33) (154,245.33) -- (161,250.33) -- (154,255.33) (66,167.33) -- (71,160.33) -- (76,167.33)  ;
\draw    (85,178.33) -- (110.76,127.62) ;
\draw [shift={(111.67,125.83)}, rotate = 476.93] [color={rgb, 255:red, 0; green, 0; blue, 0 }  ][line width=0.75]    (10.93,-3.29) .. controls (6.95,-1.4) and (3.31,-0.3) .. (0,0) .. controls (3.31,0.3) and (6.95,1.4) .. (10.93,3.29)   ;
\draw [line width=2.25]    (172,200) -- (222,150) ;
\draw [line width=2.25]    (222,150) -- (272,100) ;
\draw [line width=2.25]    (122,150) -- (172,100) ;
\draw [line width=2.25]    (170,100) -- (220,50) ;
\draw [color={rgb, 255:red, 74; green, 144; blue, 226 }  ,draw opacity=1 ][line width=2.25]    (170,150) -- (220,100) ;
\draw [line width=0.75]    (222,200) -- (272,150) ;
\draw [line width=0.75]    (222,100) -- (272,50) ;
\draw [line width=0.75]    (122,200) -- (172,150) ;
\draw [line width=0.75]    (120,100) -- (170,50) ;
\draw  [color={rgb, 255:red, 245; green, 166; blue, 35 }  ,draw opacity=1 ][line width=2.25]  (83.03,244.06) .. controls (84.06,244.14) and (84.95,243.97) .. (85.66,243.5) .. controls (89.43,240.99) and (87.28,230.71) .. (85.01,219.92) .. controls (82.73,209.13) and (80.59,198.85) .. (84.36,196.35) .. controls (88.12,193.84) and (96.77,199.81) .. (105.83,206.09) .. controls (114.9,212.37) and (123.54,218.34) .. (127.31,215.84) .. controls (131.08,213.33) and (128.93,203.05) .. (126.66,192.26) .. controls (124.39,181.47) and (122.24,171.19) .. (126.01,168.68) .. controls (129.77,166.18) and (138.42,172.15) .. (147.48,178.43) .. controls (151.74,181.38) and (155.9,184.26) .. (159.51,186.22) ;

\draw (133,258.73) node [anchor=north west][inner sep=0.75pt]    {$x_{0}$};
\draw (41,165.73) node [anchor=north west][inner sep=0.75pt]    {$x_{1}$};
\draw (88,129.73) node [anchor=north west][inner sep=0.75pt]    {$\advecvel$};

\end{tikzpicture}}
    \caption{Triangle elements}
    \label{fig:tri_patterns}
    \end{subfigure}
    \caption{Cell patterns utilized to perform the Von-Neumann analysis in two-dimensional uniform meshes of quadrilateral (left) and triangular (right) elements.
    The cells belonging to the root cell pattern are marked in blue, while the wave is depicted in orange.}
    \label{fig:patterns_2D}
\end{figure}

The resulting semi-discrete system, which discretizes \Eref{eq:linear_advection} reads
\begin{equation}
    \frac{\text{d} \var_{\ipol \iele}}{\text{d} t} = \frac{\advecvelnovect}{h} \sum_{i \in \stencil_{\iele}} \jacobianrhs_{\ipol \iele i \inu} \var_{i \inu} ,
    \label{eq:jacobian_rhs}
\end{equation}
where $\jacobianrhs$ is the right-hand-side (RHS) Jacobian matrix and $\stencil_{\iele}$ is a set that stores the indices of the direct face neighbors of an element $\iele$.
The solution of the latter equation may be simplified by taking advantage of the fact that the RHS Jacobian is a circulant block matrix \cite{Olson2014} in uniform meshes with a unique element type and periodic boundary conditions \cite{Vermeire2017}.
The block size is given by the number of solution points within a given cell pattern.
The interested reader is referred to \ref{sec:rt_basis} for a definition of the number solution points and their location within a given reference element.
In the proposed mesh subdivision procedure, based on the decomposition of tensor-product elements into simplex elements, a cell pattern refers to the resulting sub-divided cells from a unique tensor-product element.
Hence, the number of cells within each cell pattern is $\ncellsinpattern$.
Additionally, the number of circulant blocks in the RHS Jacobian is equal to the number of patterns in the initial non-decomposed mesh, i.e. $2\ndim + 1$.
The interested reader may refer to \Fref{fig:patterns_2D} for a sketch representing the mesh patterns on quadrilateral and triangular meshes.
With such definitions, one may take advantage of the block circulant matrices properties to analyze the dissipation and dispersion properties of SEM.
In the following, the index $\iele$ will refer to a given pattern, while the index $\ipol$ is representative of a given solution point within such pattern.

Let the initial condition be a Fourier mode
\begin{equation}
    \var_{\ipol \iele }(0) = \e^{\imag \vect{\kappa} \cdot \vect{x}_{\ipol \iele }} = \e^{\imag \vect{\kappa} \cdot \vect{x}_{\iele }} \patternvar_{\ipol} ,
    \label{eq:initial_condition_fourier_mode}
\end{equation}
where $\vect{x}_{\iele}$ refers to the center of a given pattern and $U_{\ipol} = \e^{\imag \vect{\kappa} \cdot (\vect{x}_{\iele \ipol} - \vect{x}_{\iele}) } \in \mathbb{R}^{\atsol{N}_e\ncellsinpattern}$ is a projection vector of size equal to the number of solution points per cell multiplied by the number of cells within a cell pattern.
For the considered initial condition, such a projection vector is independent of the considered pattern provided that the solution points are appropriately ordered within each pattern.

The structured arrangement of the mesh allows to represent the pattern centers through a label vector $\labvec \in \mathbb{N}^{\ndim}$ as
\begin{equation}
    \vect{\kappa} \cdot \vect{x}_{\iele} - \vect{\kappa} \cdot \vect{x}_0 =  \cellsize \vect{\kappa} \cdot \labvec ,
    \label{eq:structured_arrangement}
\end{equation}
where $\vect{x}_0$ is an arbitrary root cell pattern.
The RHS Jacobian may be diagonalized using the latter expression and relying on \cite[Collorary 20]{Olson2014}.
This allows to redefine \Eref{eq:jacobian_rhs} as 
\begin{equation}
    \frac{\text{d} \patternvar_{\ipol}}{\text{d} t} = \frac{\advecvelnovect}{h} \sum_{\iele \in \stencil} \jacobianrhs_{\ipol \iele \inu} \e^{\imag \vect{\kappa} \cdot \left( \vect{x}_{\iele} - \vect{x}_0 \right) } \patternvar_{\ipol} = \frac{\advecvelnovect}{h} \sum_{\iele \in \stencil} \jacobianrhsprime_{\ipol\inu} \patternvar_{\ipol} ,
    \label{eq:semi_discrete_pattern}
\end{equation}
where $\stencil$ refers to the neighbor cells of the mesh pattern root cell and $\jacobianrhsprime$ is the reduced right-hand-side Jacobian matrix.
This demonstrates that, with a structured uniform mesh, the computation of the solution may be reduced to the assessment of the values $\patternvar_{\ipol}$ within an arbitrary mesh pattern.

The exact solution of \Eref{eq:semi_discrete_pattern} is given by 
\begin{equation}
    \patternvar_{\ipol} (m\timestep) = \expjacobianrhs^m_{\ipol \inu}  \patternvar_{\inu}(0) ,
    \label{eq:semi_discrete_exponential}
\end{equation}
where $\CFL = \frac{\advecvelnovect \timestep}{\cellsize}$ is the Courant–Friedrichs–Lewy (CFL) number and $\expjacobianrhs = \e^{\CFL \jacobianrhsprime}$ is the exponential matrix of the reduced RHS Jacobian multiplied by the CFL number.
RK time integration methods approximate the exponential matrix, resulting in $\expjacobianrhs^{\numerical}$ (see \cite{Castonguay2011, VANHAREN2017}).
The assessment of \Eref{eq:semi_discrete_exponential} using RK methods allows to take into account the temporal discretization errors within the dissipation and dispersion analysis \cite{Vermeire2017, Pereira2020elements}, as well as to compute the temporal linear stability condition of the resulting discrete system.

After the diagonalization of the linear operator $\expjacobianrhs$, the aforementioned equation reduces to
\begin{equation}
    \patternvar_{\ipol} (m \timestep) = \left( \eigenvectormatrix \Lambda^m \eigenvectormatrix^{-1} \right)_{\ipol \inu}  \patternvar_{\inu}(0),
    \label{eq:patternvar_diagonalized}
\end{equation}
where $\Lambda \in \mathbb{R}^{\atsol{N}_e \times \ncellsinpattern}$ is the diagonal eigenvalue matrix and $\eigenvectormatrix$ is the square eigenvector matrix.
With analytical time and spatial discretizations, the eigenvalues of matrix $\Lambda$ are $\e^{-\imag \omega \timestep}$ cf.\ \cite{Huynh2020}.
Deviations of the numerical eigenvalues from such values indicate dissipation and dispersion errors.
The initial value $U_{\inu}(0)$ is not generally an eigenvector of the exponential matrix.
Hence, the numerical solution is governed by the contribution all each eigenmodes.
The energy stored within each eigenvector may be defined as 
\begin{equation}
    \gamma_{\ipol} = \eigenvectormatrix^{-1}_{\ipol \inu}  \patternvar_{\inu}(0) .
\end{equation}
This parameter allows to describe the so-called physical eigenmode which presents the highest energy contribution, i.e. $\max |\gamma_{\ipol}|$.
Heretofore, most studies only focused on the dissipation and dispersion properties of the physical mode (provided by its associated eigenvalue) to assess the numerical properties of SEM.
Nevertheless, it is worth noting that the dissipation and dispersion of the discrete system are only well characterized by the physical mode for non-aliased wavenumbers \cite{VANHAREN2017, ALHAWWARY2018}.
This is related to the fact that the physical mode only coincides with the spectral radius (the eigenmode whose associated eigenvalue presents the highest absolute value) of the discrete system for non-aliased wavenumbers.
For aliased wavenumbers, such a relation does not exist.
If the physical eigenmode is not the spectral radius of the discrete system, then it will be quickly dissipated during an initial transient regime and, at the asymptotic regime, aliasing appears due to the diagonalization properties matrix $\jacobianrhs$. For example, in one-dimensional configurations the latter matrix is uniquely defined in  $\kappa \cellsize \in [0, \pi]$ cf.\ \cite{VANHAREN2017}.
Therefore, the spectral radius for a given wavenumber presents similar aliasing properties as those found in FDM and FVM.
In two and three-dimensional configurations aliasing occurs if
\begin{equation}
    \kappa h \ge \min_{0\leq i \leq \ndim - 1}\left( \frac{\pi}{\unit{\unitary}_{{\vect{\kappa}}_i}} \right) .
    \label{eq:aliasing_limit}
\end{equation}
Nevertheless, there exist several secondary aliasing conditions, due to the fact that aliasing occurs on a dimension per dimension basis.
Hence, two and three-dimensional aliasing behavior is more complex than the one found in one-dimensional configurations.
The interested reader is referred to \ref{sec:aliasing} for a in-depth description on the aliasing phenomena of SEM and its influence in the dissipation and dispersion errors.
Due to all these facts, the behavior of the discrete system for aliased wavenumbers is highly unsteady.

To characterize the dissipation and dispersion of SEM the combined-mode approach cf.\ \cite{ALHAWWARY2018, VANHAREN2017} will be employed.
The dissipation error at a given iteration $m$ may be associated to imaginary part of a numerical temporal wavenumber $\omega^{\numerical}$ defined as
\begin{equation}
    \text{Im}\left(\omega^{\delta}(m)\right) = -\frac{1}{m\CFL} \ln{\frac{\norm{\patternvar(m \timestep)}}{\norm{\patternvar(0)}}},
    \label{eq:dissipation_measure}
\end{equation}
where the norm of the projection vector is computed as
\begin{equation}
    \norm{\patternvar} = \sqrt{ \frac{1}{|\domain_{\iele}|}\int_{\domain_{\iele}}  \patternvar \conj{\patternvar} \, \text{d} \domain} ,
\end{equation}
in any cell pattern $\domain_{\iele}$ of the considered mesh, with $\conj{\patternvar}$ being the conjugate value of ${\patternvar}$.
This integral needs to be assessed using appropriate numerical quadratures, within the different elements generating the given cell pattern, to avoid introducing quadrature errors in the estimation of the dissipation.
In this work, quadratures of degree ten are used to evaluate the former integral.
On the other hand, the real part of the numerical temporal wavanumber may be computed as
\begin{equation}
    \text{Re}\left(\omega^{\delta}(m)\right) = \text{Re}\left(\omega \right) -\frac{1}{m\CFL} \text{arg}\inner{\patternvar(m \timestep) \e^{\imag \omega m \timestep}, \patternvar(0)},
    \label{eq:dispersion_measure}
\end{equation}
where the inner product is defined as
\begin{equation}
    \inner{a, b} = \frac{1}{|\domain_{\iele}|}\int_{\domain_{\iele}} a \conj{b} \, \text{d} \domain.
\end{equation}

With such dissipation and dispersion measures, the estimator of the error of the numerical temporal wavenumber may be defined as 
\begin{equation}
    \mathcal{D}(m) = \frac{1}{\omega}(| \text{Re}(\omega^{\delta} - \omega) | + \imag |\text{Im}(\omega^{\delta} - \omega)|) .
    \label{eq:diss_disp_error_measure}
\end{equation}
The behavior of error estimator with the wavenumber is representative of the order of accuracy of SEM cf.\ \cite{Pereira2020}.
It is worth mentioning that, for a given wave angle configuration, the dissipation and dispersion measures are a function of the number of iterations.
Such a fact differs substantially from what is observed when performing the Von--Neumann analysis of FD or FVM, since in the latter, the dissipation and dispersion errors are related to a unique eigenmode.

\subsection{Dissipation and dispersion}

This section focuses on the comparison between the dissipation and dispersion maps of FR-DG and SDRT schemes in two uniform meshes with quadrilateral or triangular elements.
To ease such a comparison, the advection velocity vector and the wavenumber vector are supposed to be parallel.
Additionally, only configurations with $\theta_0 = \pi/6$ and $\theta_1 = \pi/4$ (if applicable) are considered.
Exhaustive analysis of the dissipation and dispersion maps for non-aligned waves and other wave angle configurations will be presented in future works.

Since the dissipation and dispersion measures depend on the number of iterations $m$, we will first focus on analyzing cases using $m = m_c = 1 /\CFL$ which implies $t = \cellsize/\advecvelnovect$.
Such a value of the physical time is related to the number of iterations needed for the wave to traverse a given cell.
With this choice, the dissipation and dispersion measures refer to the short-term diffusion and dispersion, which are intimately related to the the dissipation and dispersion of the physical eigenmode.
It is worth noting that, with exponential time integration, the numerical errors are independent of the  $\CFL$ number.
The latter is not the case when using analytical exponential time integration methods.

\begin{figure}[h]
    \begin{subfigure}{0.45\textwidth}
        \centering
        \includegraphics[width=0.9\textwidth]{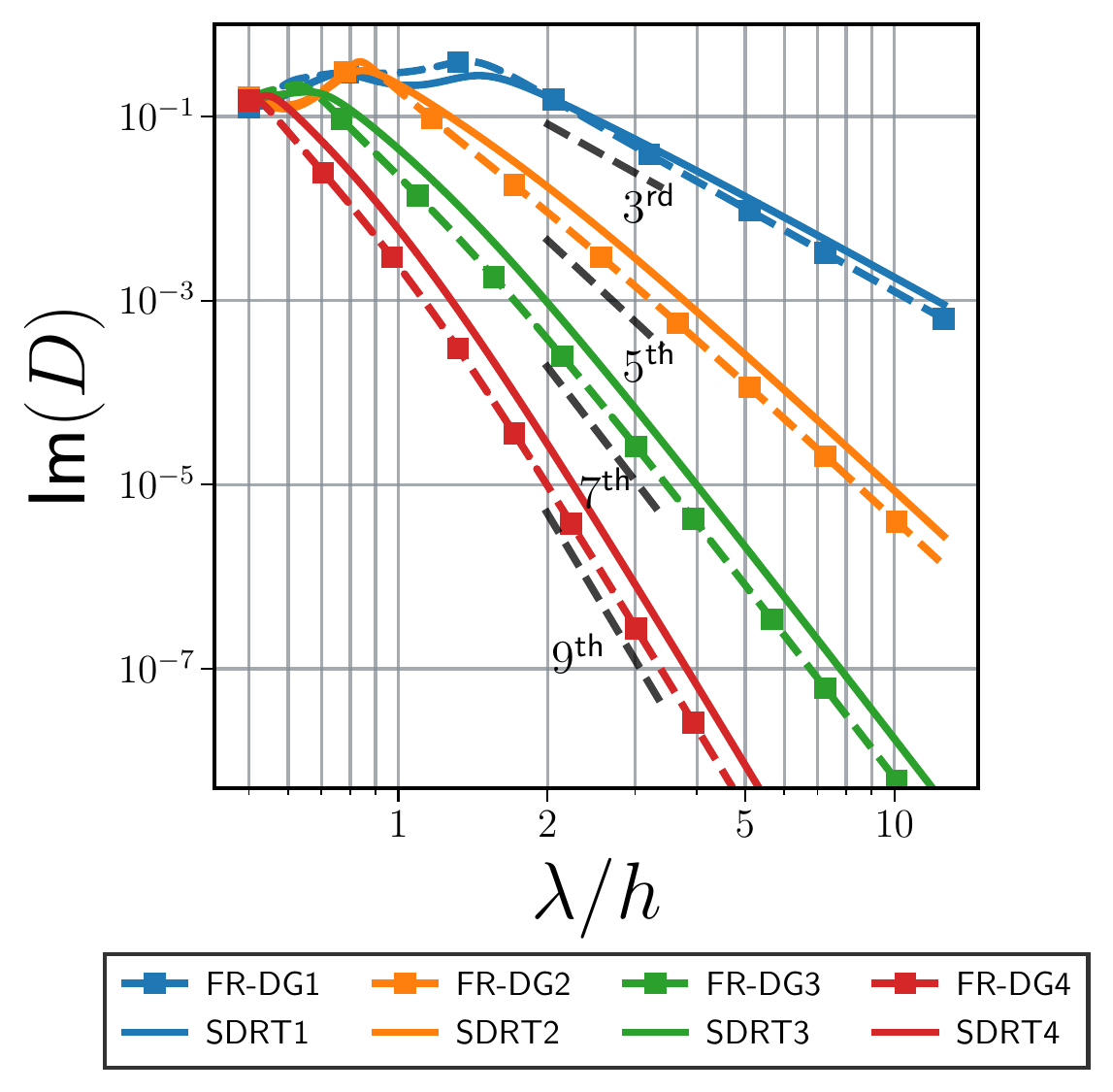}
        \caption{Dissipation error Quadrilateral}
        \label{fig:diss_quad_m1_theta0=pi/6}
    \end{subfigure}
    \begin{subfigure}{0.45\textwidth}
        \centering
        \includegraphics[width=0.9\textwidth]{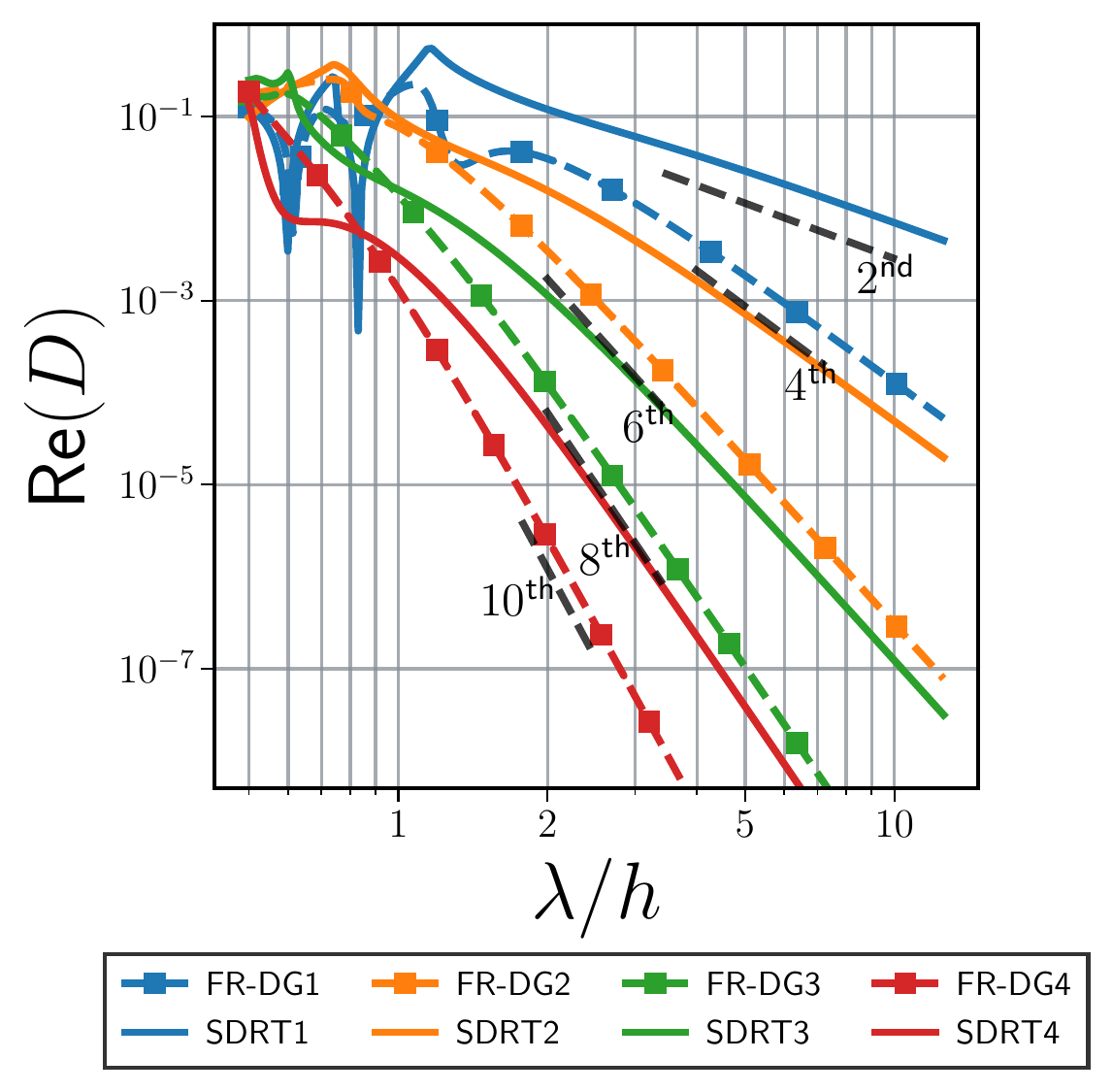}
        \caption{Dispersion error Quadrilateral}
        \label{fig:disp_quad_m1_theta0=pi/6}
    \end{subfigure}\\
    \begin{subfigure}{0.45\textwidth}
        \centering
        \includegraphics[width=0.9\textwidth]{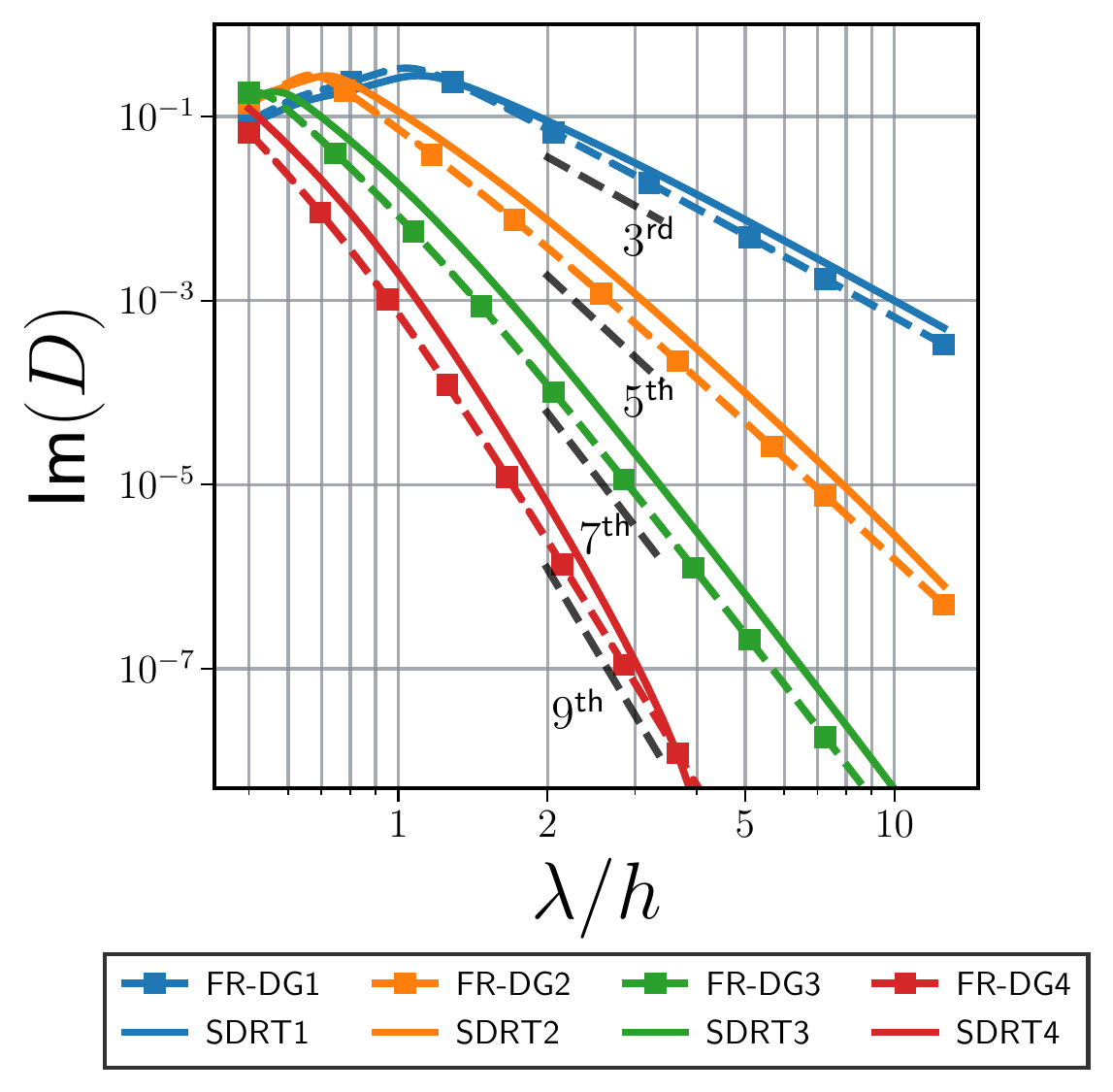}
        \caption{Dissipation error triangles}
        \label{fig:diss_tri_m1_theta0=pi/6}
    \end{subfigure}
    \begin{subfigure}{0.45\textwidth}
        \centering
        \includegraphics[width=0.9\textwidth]{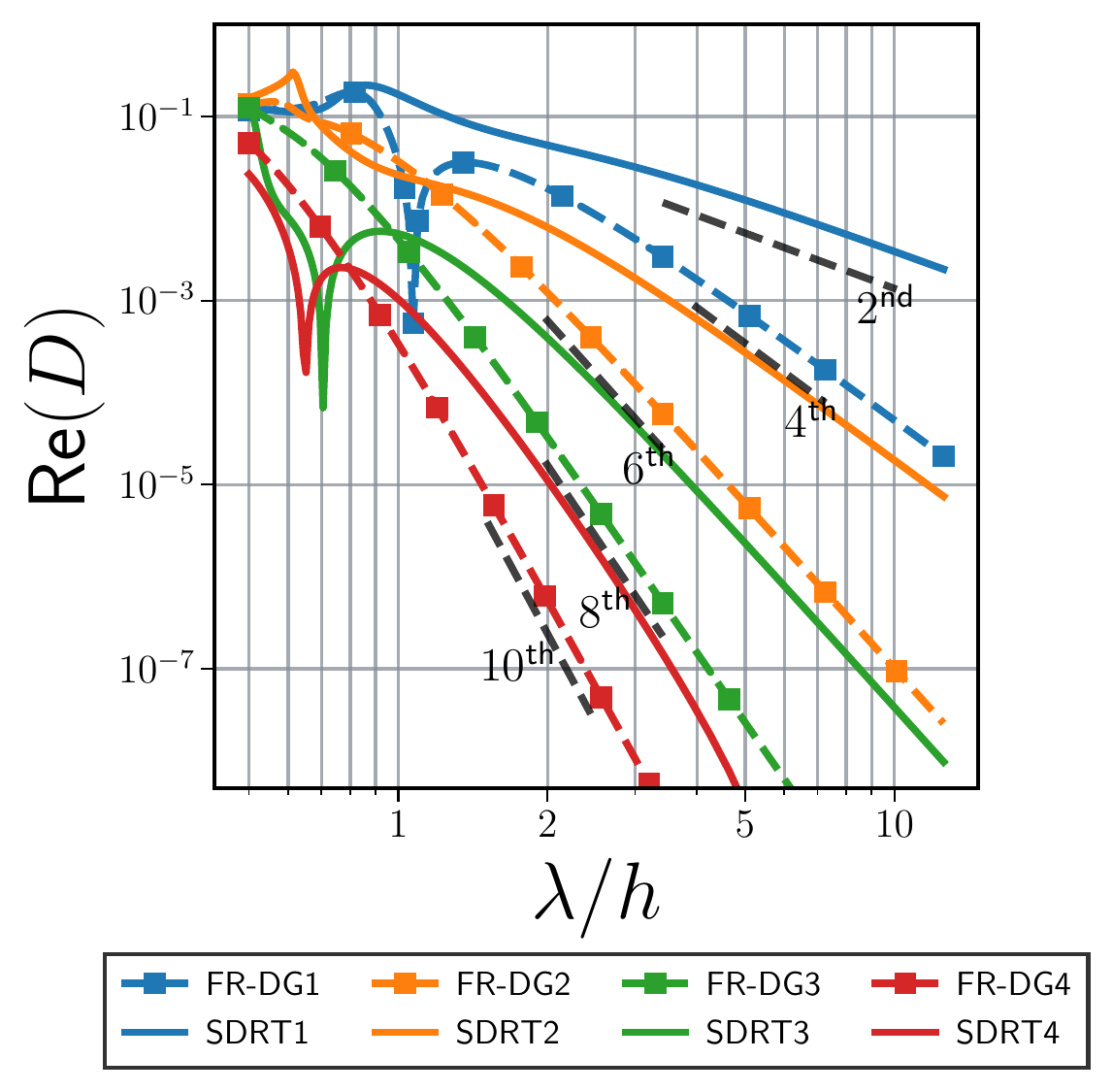}
        \caption{Dispersion error triangles}
        \label{fig:disp_tri_m1_theta0=pi/6}
    \end{subfigure}
    \caption{Dissipation (left) and dispersion (right) errors with exponential time integration obtained with two-dimensional elements, wave angle equal to $\theta_0 = \pi/6$ and $m = m_c$.}
    \label{fig:diss_disp_2D_m1_theta0=pi/6}
\end{figure}

\Fref{fig:diss_disp_2D_m1_theta0=pi/6} represents the short-term dissipation and dispersion error estimator with exponential time integration obtained using quadrilateral elements and a wave angle that is equal to $\theta_0 = \pi/6$ and which is aligned with the advection velocity.
The results are depicted as a function of the cells per wavelength parameter $\lambda / \cellsize$. They show that SDRT method presents increased dissipation and dispersion errors compared to FR-DG for every wavelength value, with the possible exception of the dispersion behavior at aliased wavenumbers.
Moreover, the dissipation errors show $2\degree + 1$ order of accuracy in both the SDRT and FR-DG formulations, while the dispersion converges to $2 \degree$ and $2 \degree + 2$ order in the SDRT and FR-DG methods respectively.
The latter observations concerning the FR-DG formulation had been already found in \cite{Guo_2013, Frean_2019} for two-dimensional meshes of tensor-product elements.
Nevertheless, to the best of the authors' knowledge, the order of accuracy of the dissipation and dispersion maps of FR-DG schemes had not been previously characterized for simplex elements.

\begin{figure}[h]
    \begin{subfigure}{0.45\textwidth}
        \centering
        \includegraphics[width=0.9\textwidth]{figures/figures_advec_diss_disp_v2/diss_quad_advection_m1_theta1=0.5235987755982988_theta2=0.pdf}
        \caption{Dissipation error Quadrilateral $m = m_c$} 
        \label{fig:diss_m1_quad_theta0=pi/6_second_time}
    \end{subfigure}
    \begin{subfigure}{0.45\textwidth}
        \centering
        \includegraphics[width=0.9\textwidth]{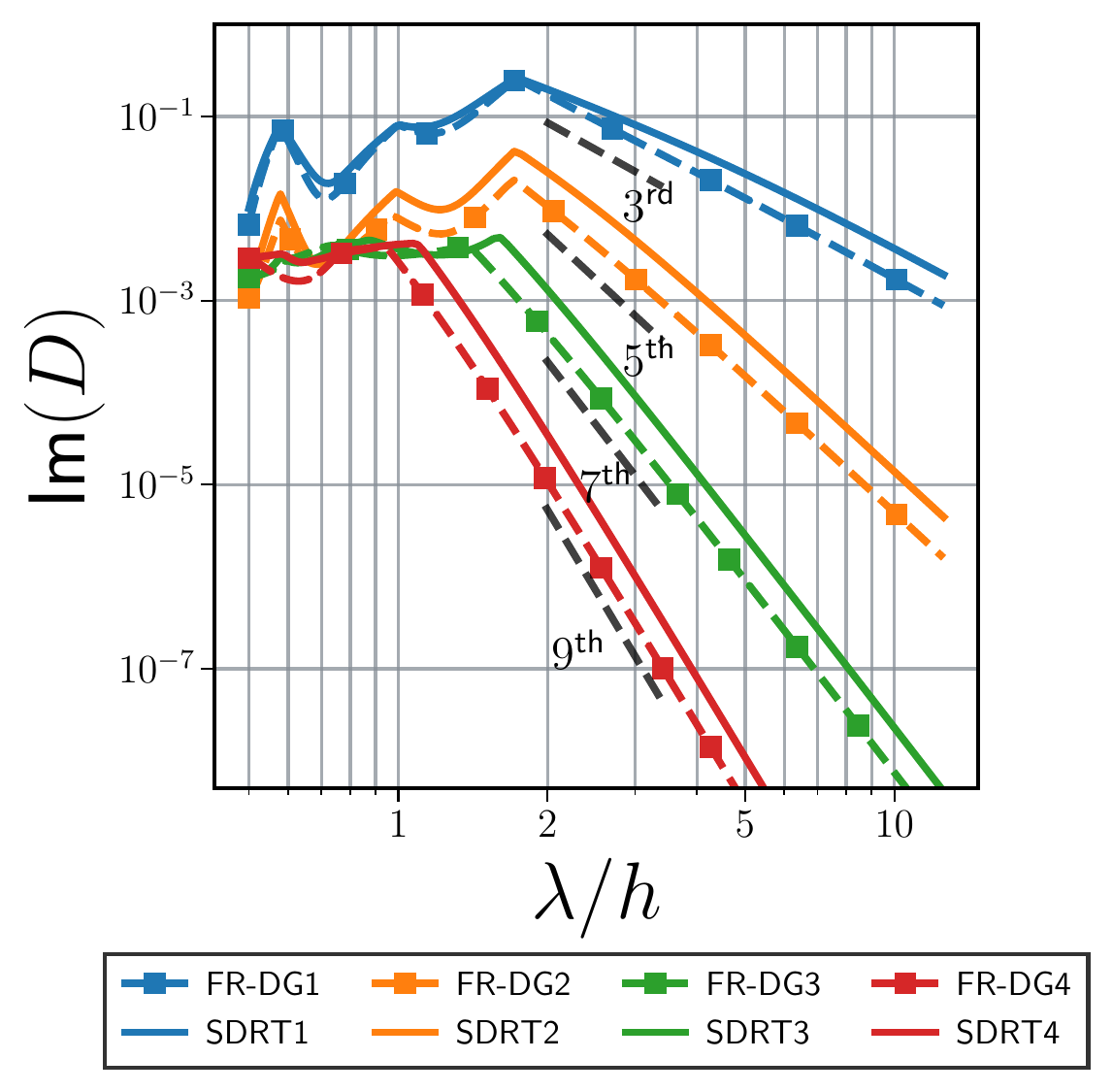}
        \caption{Dissipation error Quadrilateral $m = 400 m_c$}
        \label{fig:diss_m400_quad_theta0=pi/6}
    \end{subfigure}
    \caption{Dissipation errors with exponential time integration obtained with quadrilateral elements, wave angle equal to $\theta_0 = \pi/6$, and $m = m_c$ (left) and $m = 400 m_c$ (right).}
    \label{fig:diss_compare_m_quad_theta0=pi/6}
\end{figure}

\begin{figure}[h]
    \begin{subfigure}{0.45\textwidth}
        \centering
        \includegraphics[width=0.9\textwidth]{figures/figures_advec_diss_disp_v2/diss_tri_advection_m1_theta1=0.5235987755982988_theta2=0.pdf}
        \caption{Dissipation error Triangles $m = m_c$} 
        \label{fig:diss_m1_tri_theta0=pi/6_second_time}
    \end{subfigure}
    \begin{subfigure}{0.45\textwidth}
        \centering
        \includegraphics[width=0.9\textwidth]{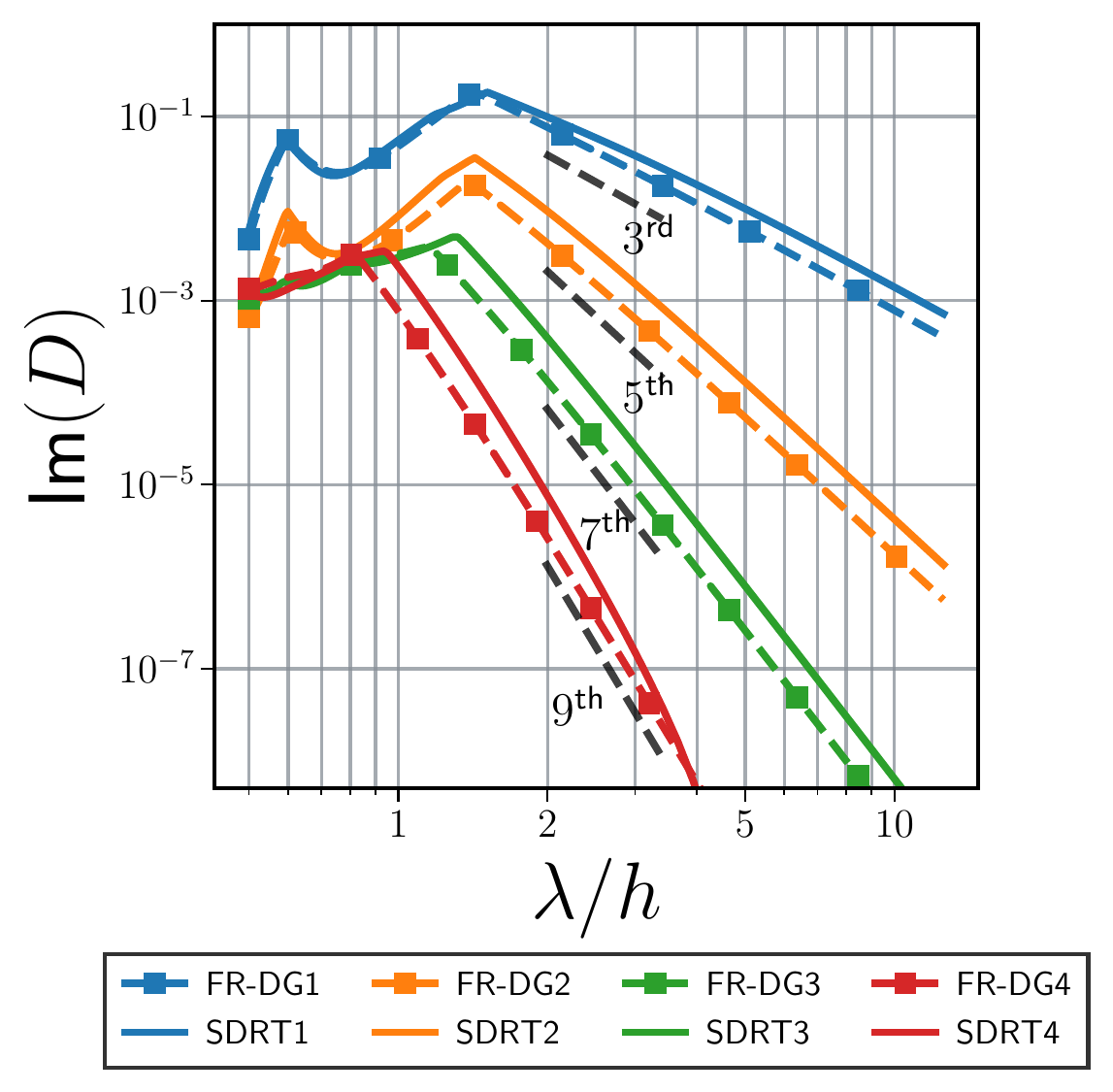}
        \caption{Dissipation error Triangles $m = 400 m_c$}
        \label{fig:diss_m400_tri_theta0=pi/6}
    \end{subfigure}
    \caption{Dissipation errors with exponential time integration obtained with triangular elements, wave angle equal to $\theta_0 = \pi/6$, and $m = m_c$ (left) and $m = 400 m_c$ (right).}
    \label{fig:diss_compare_m_tri_theta0=pi/6}
\end{figure}

To analyze the influence of the choice of the number of iterations $m$ in the dissipation errors, \Fref{fig:diss_compare_m_quad_theta0=pi/6} displays the dissipation errors associated to quadrilateral elements with wave angle that is equal to $\theta_0 = \pi/6$ and for $m = m_c$ and $m = 400 m_c$ coupled with exponential time integration.
Within this figure it is possible to observe that the dissipation errors for $\kappa h < \pi/\cos\theta_0$ or $\lambda / h > 1.73$ (which is the aliasing limit introduced in \Eref{eq:aliasing_limit}) show little to no variations when increasing the amount of iterations performed.
This is related to the fact that the physical eigenmode and the spectral radius coincide in this wavenumber range.
Nevertheless, for aliased wavenumbers $\kappa h > \pi/\cos\theta_0$, the dissipation errors are reduced as the number of iterations increases.
Such an issue can be explained by fact that the physical eigenmode does not coincide with the spectral radius, hence the dissipation errors are heavily influenced by the amount of iterations performed.
For low number of $m$, the dissipation errors are always dominated by the dissipation of the physical eigenmode.
However, for high numbers of $m$, the physical eigenmode is fully dissipated and an aliased behavior is established, which reduces the amount of dissipation errors on a iteration-per-iteration basis.
Nevertheless, the typical aliasing spectrum observed with FD and FVM is not exactly reproduced since this aliasing is established after an important part of initial wave energy (contained in the physical eigenmode) is dissipated.
It is worth noting that as the degree of the schemes increases, the amount of dissipation errors of the physical mode is reduced.
Therefore, more iterations are needed to dissipate the physical eigenmode.
This explains the lack of aliased behavior found within schemes built with $\degree = 4$ at $\kappa h \approx \pi/\cos\theta_0$.
At last, it is worth mentioning that there exist several dissipation local maximums.
Such points are related to secondary aliasing limits when $\kappa h \approx \pi/\sin\theta_0$ (hence $\lambda / h \approx 1$) and/or when $\kappa h \approx  3\pi/\cos\theta_0$ (hence $\lambda / h \approx 0.58$).
For the sake of completeness, \Fref{fig:diss_compare_m_tri_theta0=pi/6} represents the same results with triangular elements.
As the conclusions that can be drawn from this figure are the same than those obtained with quadrilateral elements, no further discussions are added.

\begin{figure}[h]
    \begin{subfigure}{0.45\textwidth}
        \centering
        \includegraphics[width=0.9\textwidth]{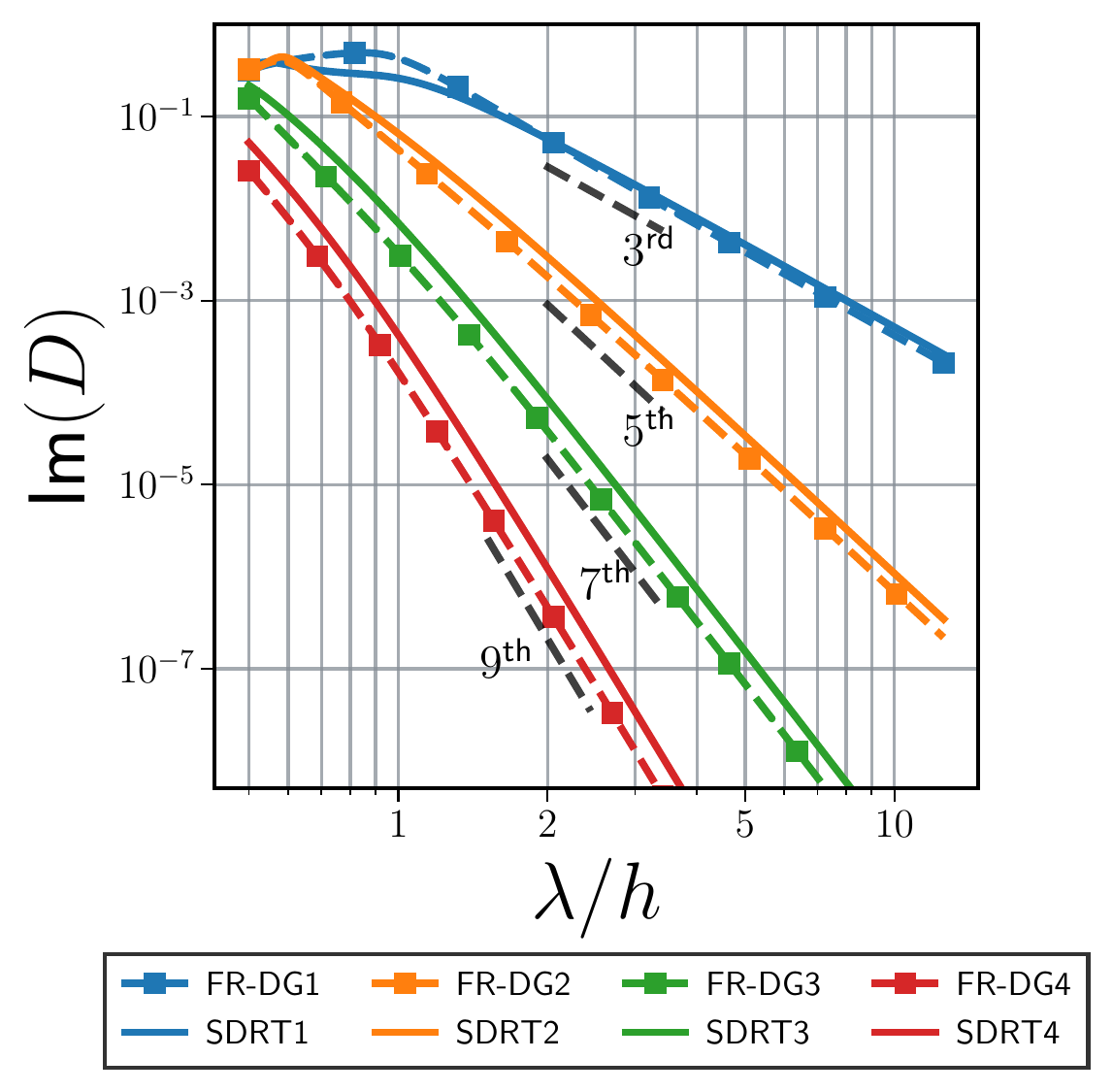}
        \caption{Dissipation error Hexahedra}
        \label{fig:diss_hex_theta0=pi/6_theta1=pi/4}
    \end{subfigure}
    \begin{subfigure}{0.45\textwidth}
        \centering
        \includegraphics[width=0.9\textwidth]{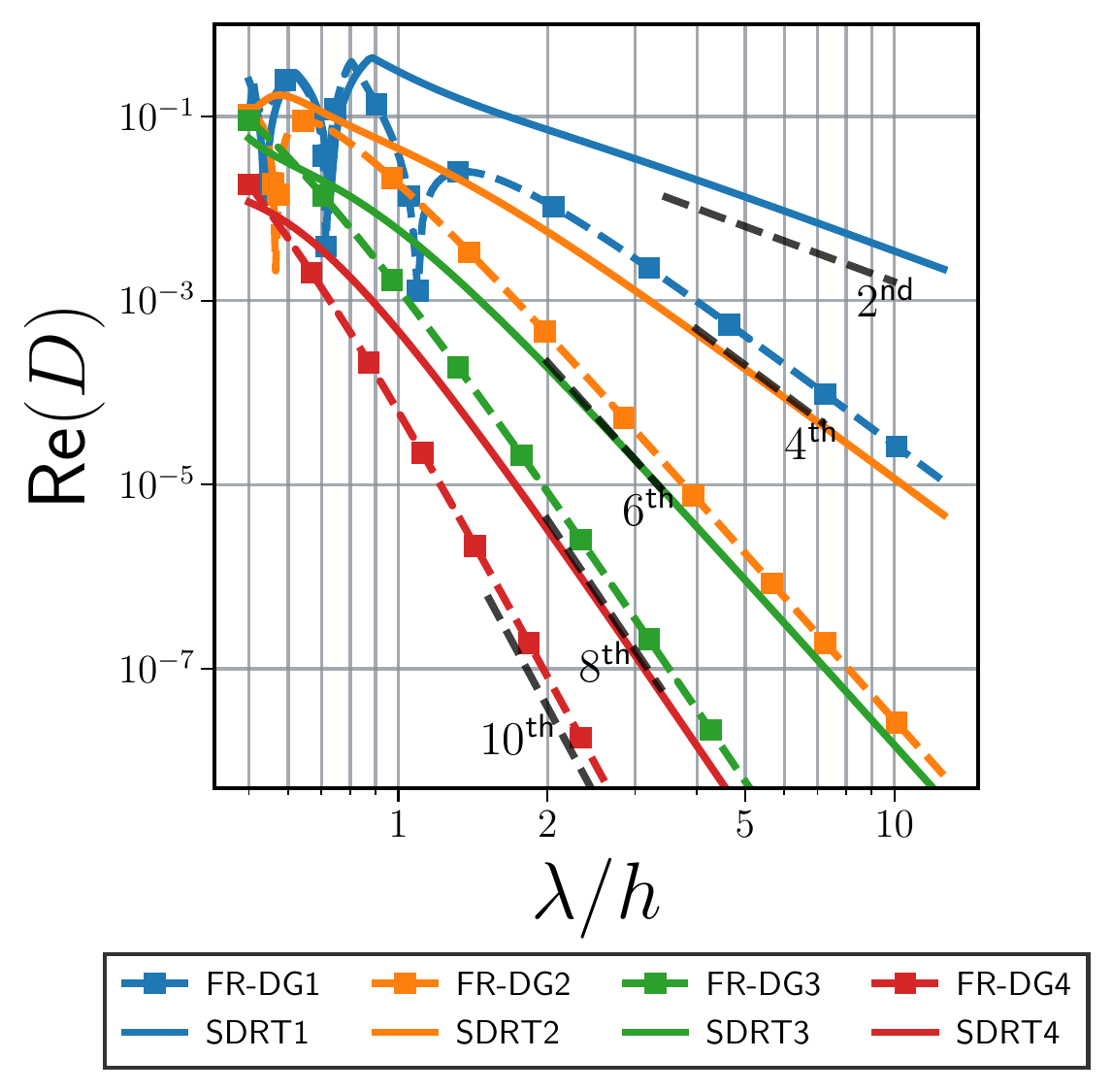}
        \caption{Dispersion error Hexahedra}
        \label{fig:disp_hex_theta0=pi/6_theta1=pi/4}
    \end{subfigure}\\
    \begin{subfigure}{0.45\textwidth}
        \centering
        \includegraphics[width=0.9\textwidth]{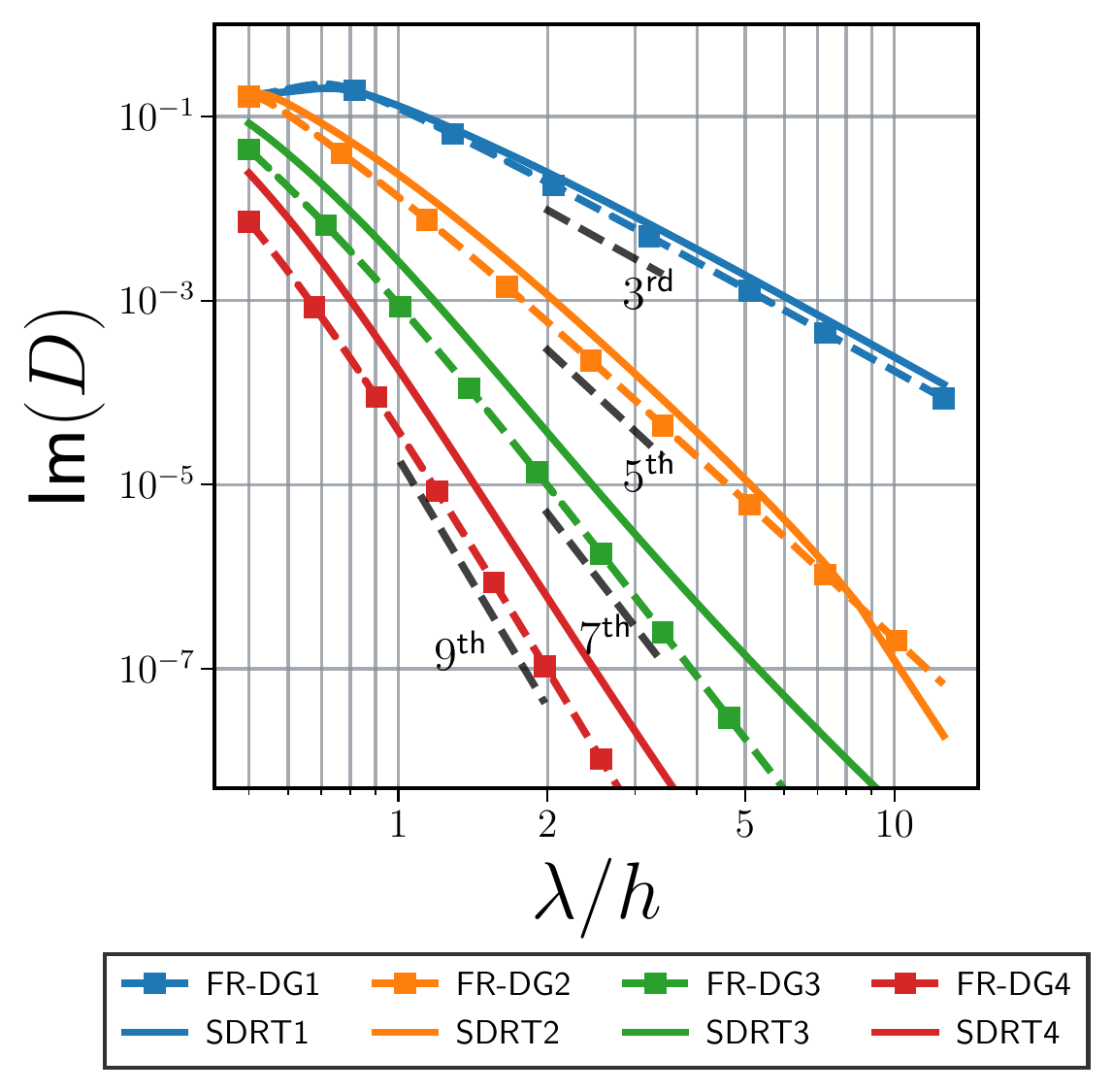}
        \caption{Dissipation error Tetrahedra}
        \label{fig:diss_tet_theta0=pi/6_theta1=pi/4}
    \end{subfigure}
    \begin{subfigure}{0.45\textwidth}
        \centering
        \includegraphics[width=0.9\textwidth]{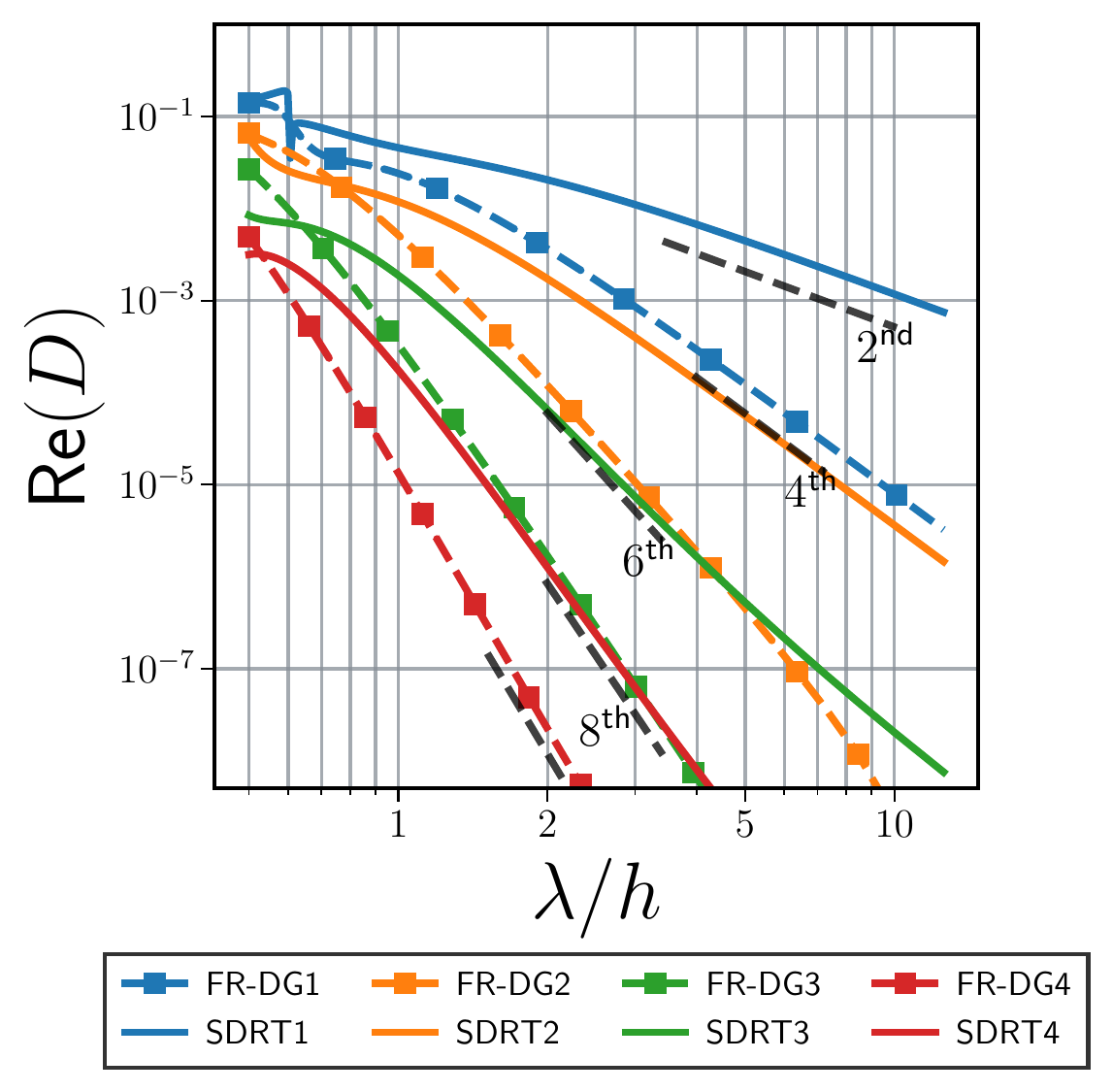}
        \caption{Dispersion error Tetrahedra}
        \label{fig:disp_tet_theta0=pi/6_theta1=pi/4}
    \end{subfigure}\\
    \begin{subfigure}{0.45\textwidth}
        \centering
        \includegraphics[width=0.9\textwidth]{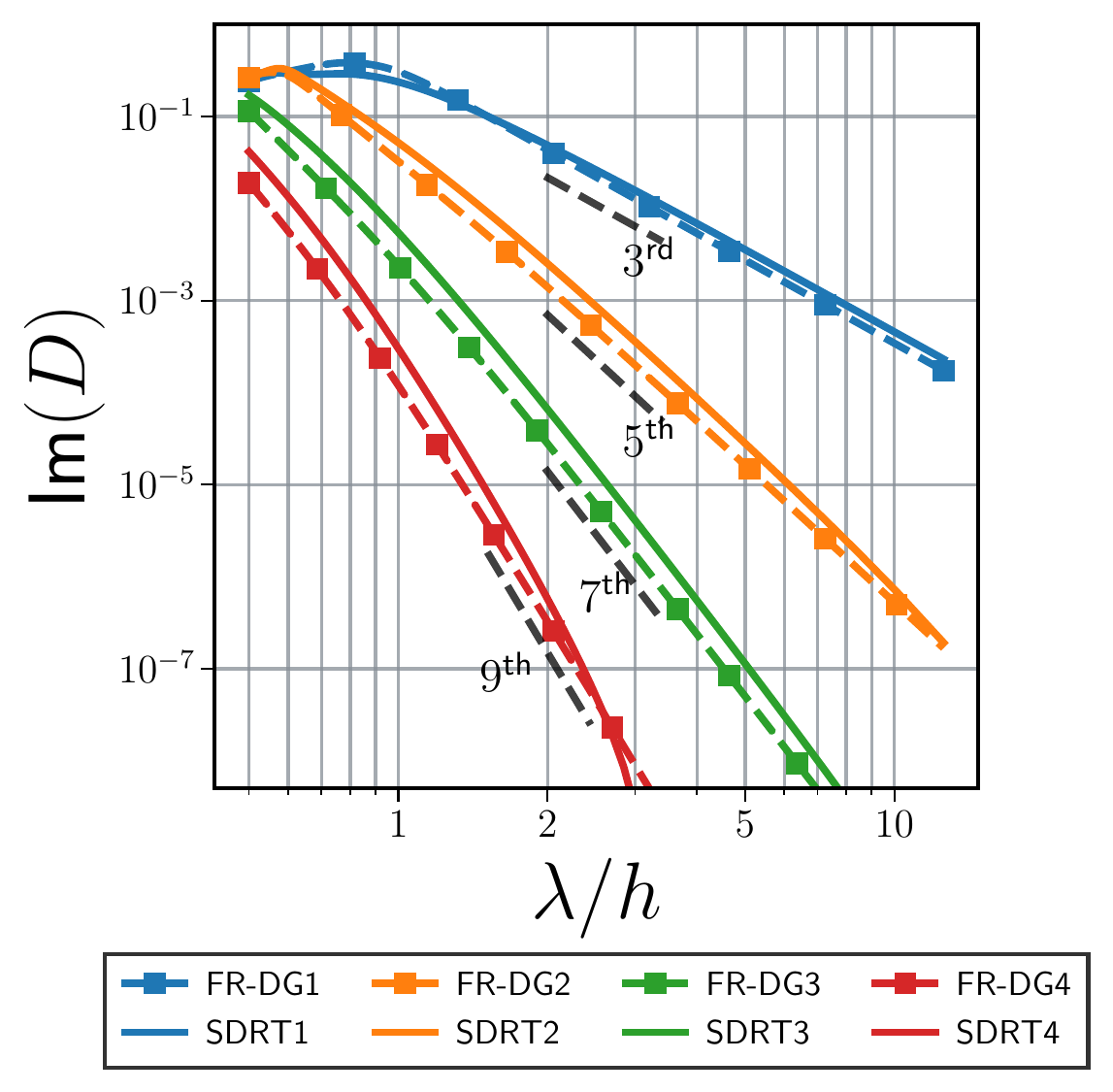}
        \caption{Dissipation error Prisms}
        \label{fig:diss_pri_theta0=pi/6_theta1=pi/4}
    \end{subfigure}
    \begin{subfigure}{0.45\textwidth}
        \centering
        \includegraphics[width=0.9\textwidth]{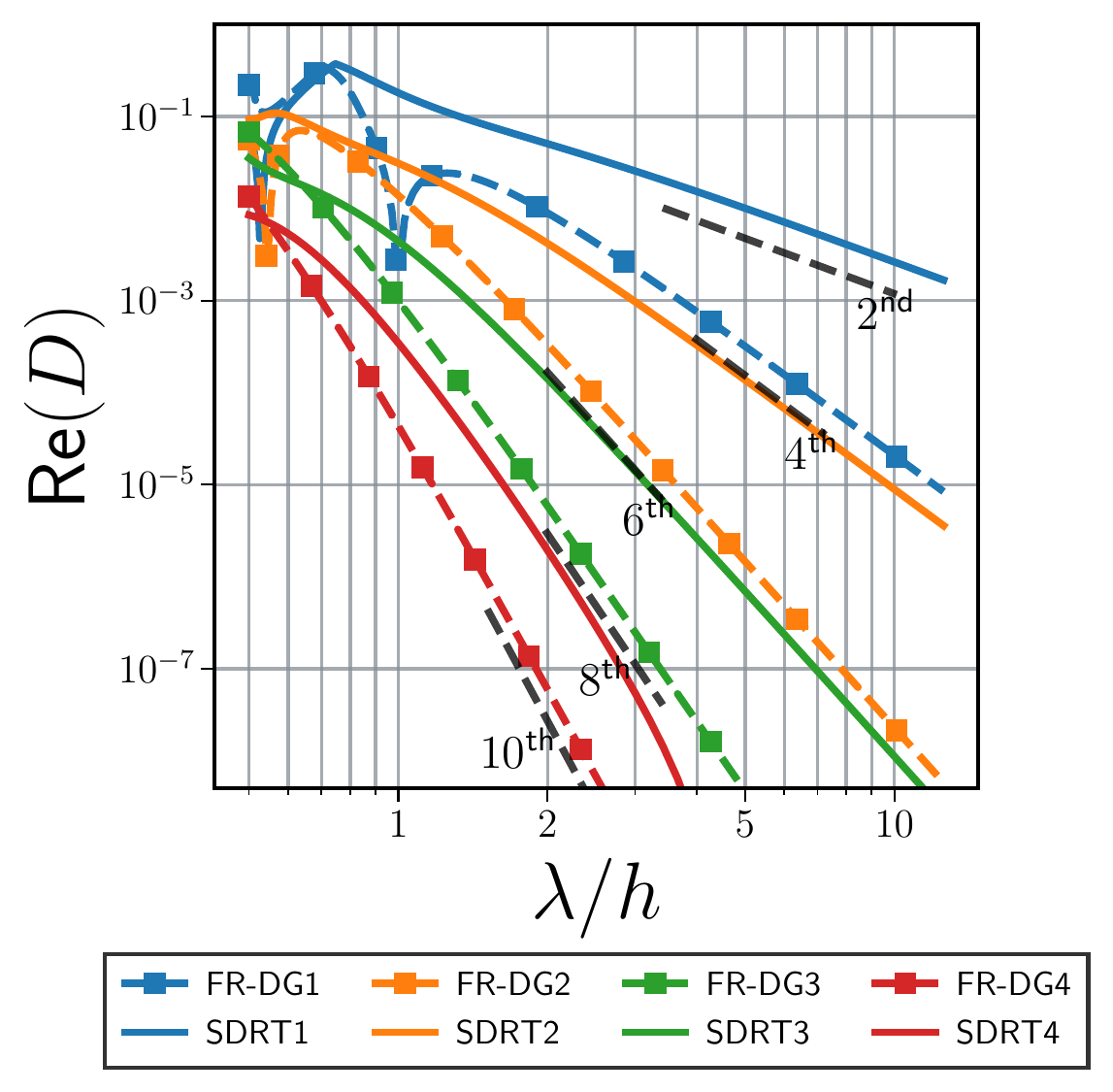}
        \caption{Dispersion error Prisms}
        \label{fig:disp_pri_theta0=pi/6_theta1=pi/4}
    \end{subfigure}
    \caption{Dissipation (left) and dispersion (right) errors with exponential time integration obtained with three-dimensional elements, wave angles equal to $\theta_0 = \pi/6$, $\theta_1 = \pi/4$ and $m = m_c$.}
    \label{fig:diss_disp_3D_m1_theta0=pi/6_theta1=pi/4}
\end{figure}

With regards to three-dimensional elements, \Fref{fig:diss_disp_3D_m1_theta0=pi/6_theta1=pi/4} displays the dissipation and dispersion error estimator with exponential time integration and $m = m_c$ obtained with quadrilateral elements and wave angles $\theta_0 = \pi/6$, $\theta_1 = \pi/4$.
As it was observed in the two-dimensional analysis, SDRT schemes present $2 \degree$ and $2 \degree + 1$ order of accuracy in their dispersion and dissipation maps, respectively, for most elements and polynomial degrees.
On the other hand, FR-DG schemes display $2\degree + 1$ and $2\degree + 2$ order of accuracy in their dispersion and dissipation maps, respectively.
The latter schemes show reduced dispersion and dissipation errors compared to SDRT methods.
It is worth mentioning the reduced order of accuracy of the SDRT3 scheme with tetrahedron elements for well-resolved wavenumbers.
The reason behind this behavior is undetermined, although it is hypothesized that it may be related to the presence of several eigenmodes with similar energy contribution than that of the physical mode and whose associated dissipation and dispersion errors are of lower orders.

\begin{remark}
Although not shown for the sake of brevity, aliasing issues have also been observed in three-dimensional elements.
\end{remark}

Up to this point, the dissipation and dispersion errors analyzed in this work do not take into account the accuracy of the temporal discretization.
In order to consider the temporal discretization defects, the dissipation and dispersion properties will be studied using the classical three stages and third order RK3, four stages fourth order RK4, and the five stages and fourth order RK54 time integrator \cite{Carpenter1994Fourthorder2R}.
Only triangular elements will be considered since, in our studies, all elements showed a similar influence of the temporal discretization errors in the dissipation and dispersion errors.
The aforementioned RK schemes approximate the exponential matrix as
\begin{equation}
    \expjacobianrhs^{\numerical}_{\text{RK3}} = \sum_{i=0}^{3} \frac{\left( \tau \jacobianrhsprime \right)^i}{i!} \text{ , } \expjacobianrhs^{\numerical}_{\text{RK4}} = \sum_{i=0}^{4} \frac{\left( \CFL \jacobianrhsprime \right)^i}{i!} \text{ and } \expjacobianrhs^{\numerical}_{\text{RK54}} = \sum_{i=0}^{4} \frac{\left( \CFL \jacobianrhsprime \right)^i}{i!} + \frac{\left( \CFL \jacobianrhsprime \right)^5}{200}.
    \label{eq:rk_exponential}
\end{equation}

\begin{figure}[h]
    \begin{subfigure}{0.45\textwidth}
        \centering
        \includegraphics[width=0.9\textwidth]{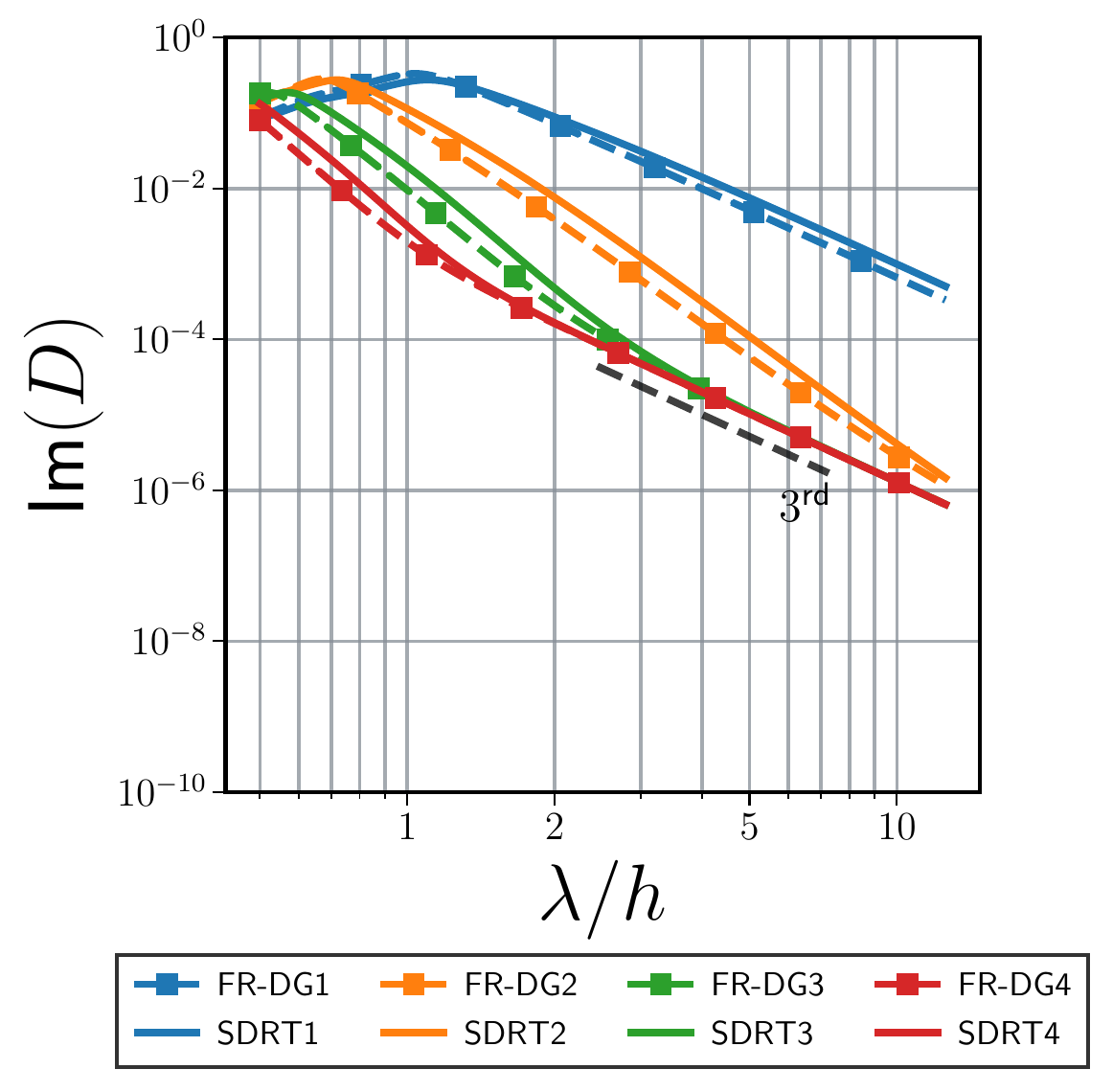}
        \caption{Dissipation error triangles}
        \label{fig:rk3_diss_tri_theta0=pi/6}
    \end{subfigure}
    \begin{subfigure}{0.45\textwidth}
        \centering
        \includegraphics[width=0.9\textwidth]{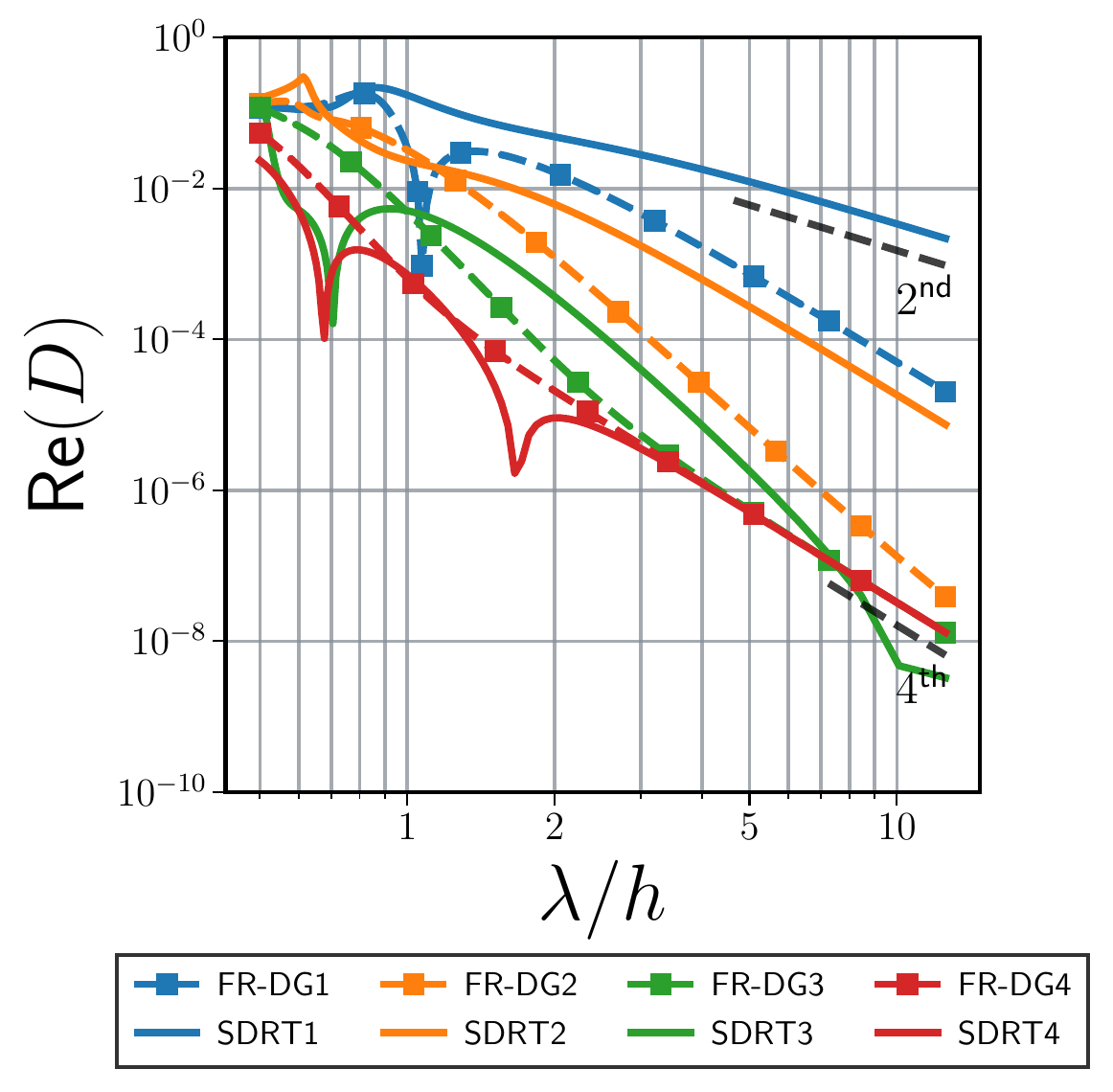}
        \caption{Dispersion error triangles}
        \label{fig:rk3_disp_tri_theta0=pi/6}
    \end{subfigure}
    \caption{Dissipation (left) and dispersion (right) errors with exponential the RK3 time integration scheme obtained with triangular elements, wave angle equal to $\theta_0 = \pi/6$ and $m = m_c$.}
    \label{fig:rk3_diss_disp_2D_theta0=pi/6}
\end{figure}

\Fref{fig:rk3_diss_disp_2D_theta0=pi/6} shows the dissipation and dispersion errors with triangular elements, $\theta_0 = \pi/6$ and the RK3 time integrator with $\CFL = 0.05$ with $m=m_c$.
From these results, it can be clearly observed that the addition of temporal discretization defects to the dissipation and dispersion maps greatly distorts the accuracy of the schemes as the polynomial degree increases.
Furthermore, the order of accuracy of the dissipation and dispersion is reduced to, at most, third and fourth order respectively.
This order of accuracy reduction is consistent with the third order of accuracy of the RK3 time integrator.
Nevertheless, the dissipation and dispersion errors for low $\lambda/h$ remain almost invariant with respect to those obtained with exponential time integration.

\begin{figure}[h]
    \begin{subfigure}{0.45\textwidth}
        \centering
        \includegraphics[width=0.9\textwidth]{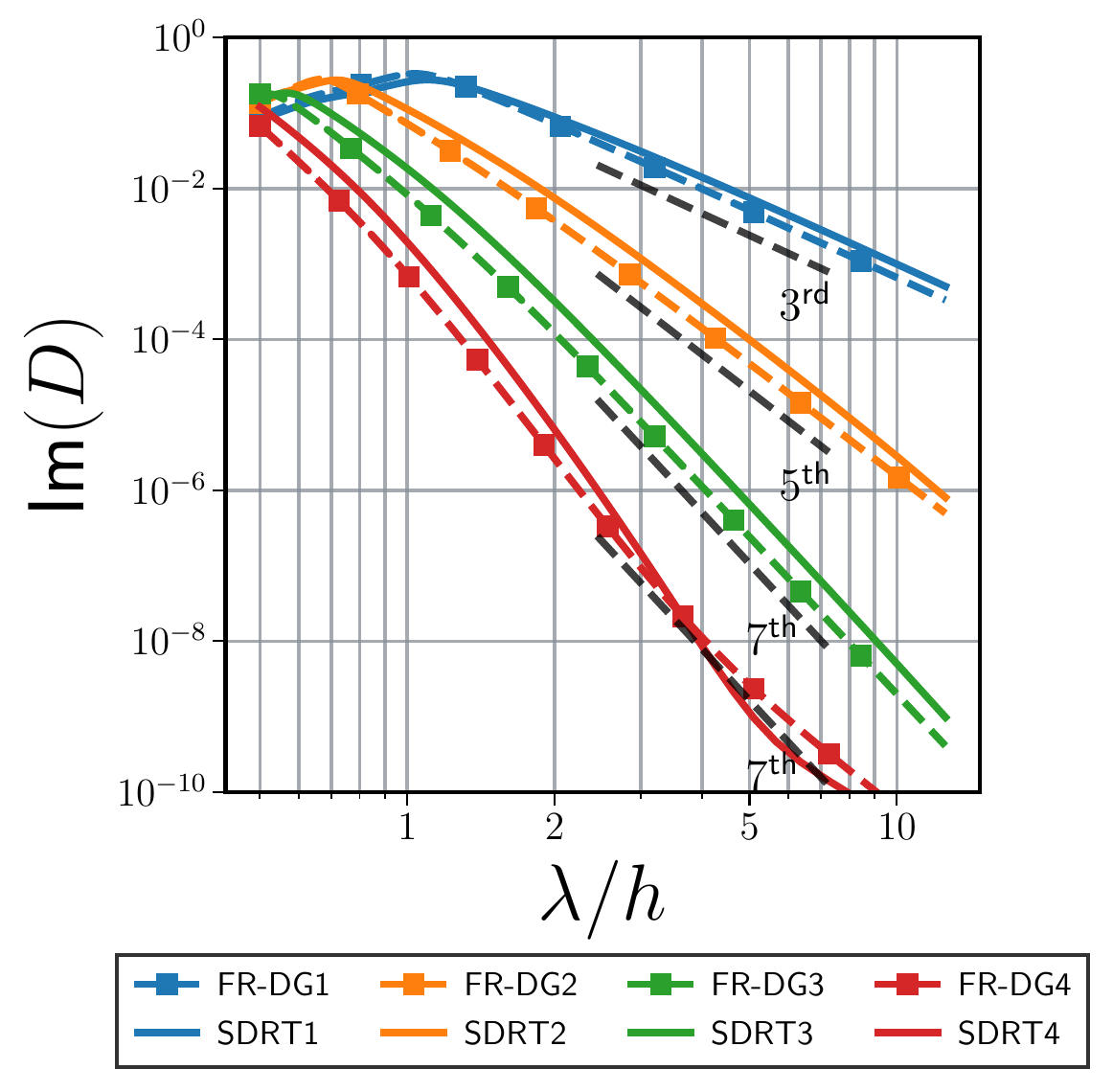}
        \caption{Dissipation error triangles}
        \label{fig:rk54_diss_tri_theta0=pi/6}
    \end{subfigure}
    \begin{subfigure}{0.45\textwidth}
        \centering
        \includegraphics[width=0.9\textwidth]{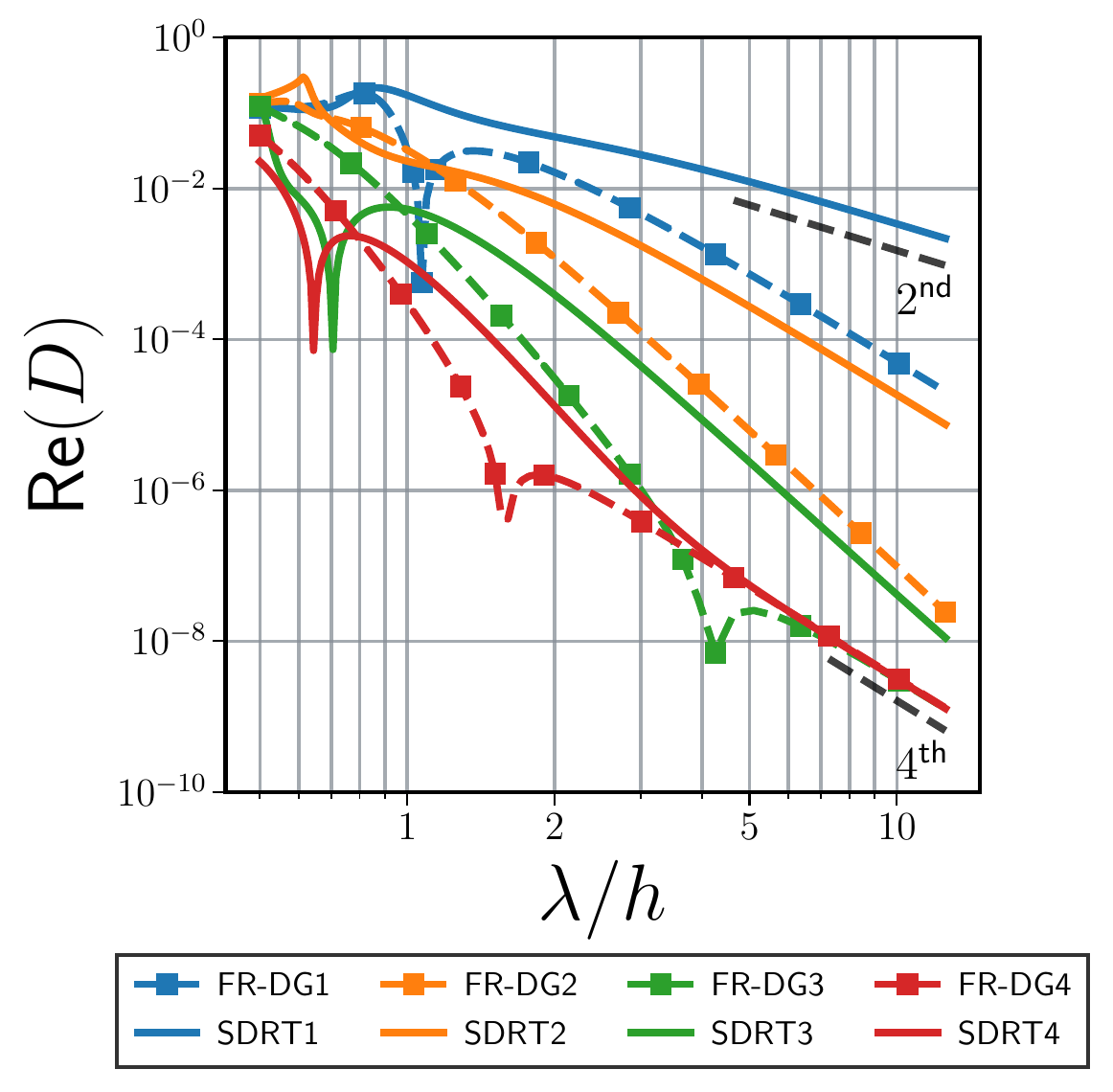}
        \caption{Dispersion error triangles}
        \label{fig:rk54_disp_tri_theta0=pi/6}
    \end{subfigure}
    \caption{Dissipation (left) and dispersion (right) errors with exponential the RK54 time integration scheme obtained with triangular elements, wave angle equal to $\theta_0 = \pi/6$ and $m = m_c$.}
    \label{fig:rk54_diss_disp_2D_theta0=pi/6}
\end{figure}

On the other hand, \Fref{fig:rk54_diss_disp_2D_theta0=pi/6} shows the dissipation and dispersion errors with triangular elements, $\theta_0 = \pi/6$ and the RK54 time integrator with $\CFL = 0.05$ for $m=m_c$.
The results indicates the RK54 time integrator has little influence in the dissipation errors, up to $\lambda / \cellsize \approx 5$ for $\degree = 4$.
Nevertheless, the dispersion errors display fourth order of accuracy for most polynomial degrees, consistent with the formal fourth order of accuracy of the RK54 time integrator.
It is worth mentioning that the influence of the temporal discretization errors in the dissipation and dispersion of SEM can be reduced by decreasing the $\CFL$ number.

\begin{remark}
Although not shown for the sake of brevity, the aliasing issues that appear for high values of the number of iterations are not influenced by the temporal discretization errors of RK schemes.
This can be explained by the fact that the temporal discretization errors only influence the diffusion and dispersion errors of the numerical schemes for well-resolved waves.
\end{remark}

\subsubsection{Temporal linear stability in uniform mesh}
\label{sec:temporal_stability_advection}

Temporal linear stability here refers to the asymptotic stability of the numerical system described in \Eref{eq:semi_discrete_exponential} as a function of the time step.
Boundedness of the solution may be ensured if the absolute value of all eigenvalues of the exponential matrix $\expjacobianrhs^{\numerical}$ is less than one.
Therefore, the discrete linear system defined in \Eref{eq:semi_discrete_exponential} is linearly stable with exponential time integration provided that the RHS Jacobian is stable, i.e. if the spatial discretization is stable.
On the other hand, if the exponential of the RHS Jacobian is approximated through RK methods, there exists a given value of the $\timestep$ or $\CFL$ value for which the spectral radius is strictly greater than one.
Such a $\CFL$ value is referred to as $\CFL_{\text{MAX}}$ and for uniform periodic meshes, it is a function of the reduced wavenumber $\kappa h$, the advection velocity angle and the wavenumber angle.

As the SDRT formulation is equivalent to the original SD method in tensor-product elements, SDRT schemes are not supposed to be linearly unstable with tensor-product elements.
It is worth mentioning that the stability of the SD method in one-dimensional configurations is ensured if the flux points are located at Gauss-Legendre quadrature points \cite{Jameson2010}.
Nevertheless, for triangular, tetrahedron and prismatic elements, it was observed (a posteriori) that SDRT schemes are unstable for $\degree \ge 5$, i.e. the real part of some eigenvalues of $\jacobianrhs$ was found to be substantially higher than machine round-off errors.
Such a fact was also observed in \cite{Adele2021thesis} for triangular elements.
The latter study also proposed a distribution of unique internal flux points which are said to yield stable schemes for $\degree = 5$.
The reason behind such unstable behavior of SDRT schemes in triangles is undetermined and will be further studied in future works.

\Tref{table:cfl_2D_theta=pi/6} shows the $\CFL_{\text{MAX}}$ for quadrilateral and triangular elements, the SDRT and FR-DG schemes and RK3, RK4 \cite{Butcher_2016} and RK54 (five stages, two registers and fourth order) \cite{Carpenter1994Fourthorder2R} scheme. 
The results indicate that SDRT schemes present around 30\% additional higher $\CFL_{\text{MAX}}$ number than that of FR-DG methods.
Moreover, the CFL condition in triangular elements is lower than that of quadrilateral cells.
 
\Tref{table:pri_fr_dg_cfl_theta0=pi/6_theta1=pi/4} shows the $\CFL_{\text{MAX}}$ for hexahedral, tetrahedral and prismatic elements, the SDRT and FR-DG schemes and aforementioned ERK schemes. 
As it was observed for two-dimensional elements, the results indicate that SDRT schemes present around 30\% additional higher $\CFL_{\text{MAX}}$ number than that of FR-DG methods.
Moreover, the CFL condition in tetrahedron and prismatic elements is lower than that of quadrilateral cells.

\begin{table}[!htb]
    \centering
    \begin{subtable}{.45\linewidth}
        \centering
        \begin{tabular}{@{}lccc@{}}
        \toprule
         & \multicolumn{3}{c}{$\CFL_{\text{MAX}}$} \\
         \midrule
         & RK3 & RK4 & RK54 \\
         \midrule
        SDRT1 & 0.459 & 0.509 & 0.695  \\ 
        SDRT2 & 0.235 & 0.281 & 0.392  \\ 
        SDRT3 & 0.149 & 0.165 & 0.247  \\ 
        SDRT4 & 0.102 & 0.119 & 0.173  \\
        \bottomrule
        \end{tabular}
        \caption{SDRT Quadrilateral elements}
        \label{table:quad_sdrt_cfl_theta=pi/6}
    \end{subtable}
    \begin{subtable}{.45\linewidth}
        \centering
        \begin{tabular}{@{}lccc@{}}
        \toprule
         & \multicolumn{3}{c}{$\CFL_{\text{MAX}}$} \\
         \midrule
         & RK3 & RK4 & RK54 \\
         \midrule
        FR-DG1 & 0.306 & 0.339 & 0.507  \\ 
        FR-DG2 & 0.153 & 0.185 & 0.261  \\ 
        FR-DG3 & 0.096 & 0.106 & 0.162  \\ 
        FR-DG4 & 0.065 & 0.078 & 0.113  \\
        \bottomrule
        \end{tabular}
        \caption{FR-DG Quadrilateral elements}
        \label{table:quad_fr_dg_cfl_theta=pi/6}
    \end{subtable}
    \vspace{0.5cm}
    \begin{subtable}{.45\linewidth}
        \centering
        \begin{tabular}{@{}lccc@{}}
        \toprule
         & \multicolumn{3}{c}{$\CFL_{\text{MAX}}$} \\
         \midrule
         & RK3 & RK4 & RK54 \\
         \midrule
        SDRT1 & 0.286 & 0.329 & 0.464  \\ 
        SDRT2 & 0.179 & 0.198 & 0.290  \\ 
        SDRT3 & 0.112 & 0.128 & 0.196  \\ 
        SDRT4 & 0.078 & 0.088 & 0.138  \\ 
        \bottomrule
        \end{tabular}
        \caption{SDRT Triangular elements}
        \label{table:tri_sdrt_cfl_theta=pi/6}
    \end{subtable}
    \begin{subtable}{.45\linewidth}
        \centering
        \begin{tabular}{@{}lccc@{}}
        \toprule
         & \multicolumn{3}{c}{$\CFL_{\text{MAX}}$} \\
         \midrule
         & RK3 & RK4 & RK54 \\
         \midrule
        FR-DG1 & 0.222 & 0.249 & 0.372  \\ 
        FR-DG2 & 0.126 & 0.140 & 0.217  \\ 
        FR-DG3 & 0.086 & 0.096 & 0.146  \\ 
        FR-DG4 & 0.060 & 0.066 & 0.103  \\
        \bottomrule
        \end{tabular}
        \caption{FR-DG Triangular elements}
        \label{table:tri_fr_dg_cfl_theta=pi/6}
    \end{subtable}
    \caption{$\CFL_{\text{MAX}}$ number to ensure linear stability in the linear advection equation with uniform meshes made up of two-dimensional elements with $\theta_0 = \pi/6$.}
    \label{table:cfl_2D_theta=pi/6}
\end{table}

\begin{table}[!htb]
    \centering
    \begin{subtable}{.45\linewidth}
        \centering
        \begin{tabular}{@{}lccc@{}}
        \toprule
         & \multicolumn{3}{c}{$\CFL_{\text{MAX}}$} \\
         \midrule
         & RK3 & RK4 & RK54 \\
         \midrule
        SDRT1 & 0.345 & 0.403 & 0.552  \\ 
        SDRT2 & 0.187 & 0.222 & 0.304  \\ 
        SDRT3 & 0.117 & 0.131 & 0.196  \\ 
        SDRT4 & 0.081 & 0.094 & 0.136  \\
        \bottomrule
        \end{tabular}
        \caption{SDRT Hexahedral}
        \label{table:hex_sdrt_cfl_theta0=pi/6_theta1=pi/4}
    \end{subtable}
    \begin{subtable}{.45\linewidth}
        \centering
        \begin{tabular}{@{}lccc@{}}
        \toprule
         & \multicolumn{3}{c}{$\CFL_{\text{MAX}}$} \\
         \midrule
         & RK3 & RK4 & RK54 \\
         \midrule
        FR-DG1 & 0.237 & 0.269 & 0.394  \\ 
        FR-DG2 & 0.122 & 0.146 & 0.204  \\ 
        FR-DG3 & 0.075 & 0.084 & 0.127  \\ 
        FR-DG4 & 0.052 & 0.061 & 0.088  \\ 
        \bottomrule
        \end{tabular}
        \caption{FR-DG Hexahedral}
        \label{table:hex_fr_dg_cfl_theta0=pi/6_theta1=pi/4}
    \end{subtable}
    \vspace{0.5cm}
    \begin{subtable}{.45\linewidth}
        \centering
        \begin{tabular}{@{}lccc@{}}
        \toprule
         & \multicolumn{3}{c}{$\CFL_{\text{MAX}}$} \\
         \midrule
         & RK3 & RK4 & RK54 \\
         \midrule
        SDRT1 & 0.151 & 0.172 & 0.262  \\ 
        SDRT2 & 0.106 & 0.124 & 0.170  \\ 
        SDRT3 & 0.067 & 0.076 & 0.117  \\ 
        SDRT4 & 0.034 & 0.037 & 0.062  \\
        \bottomrule
        \end{tabular}
        \caption{SDRT Tetrahedra}
        \label{table:tet_sdrt_cfl_theta0=pi/6_theta1=pi/4}
    \end{subtable}
    \begin{subtable}{.45\linewidth}
        \centering
        \begin{tabular}{@{}lccc@{}}
        \toprule
         & \multicolumn{3}{c}{$\CFL_{\text{MAX}}$} \\
         \midrule
         & RK3 & RK4 & RK54 \\
         \midrule
        FR-DG1 & 0.123 & 0.138 & 0.216  \\ 
        FR-DG2 & 0.085 & 0.095 & 0.137  \\ 
        FR-DG3 & 0.055 & 0.061 & 0.097  \\ 
        FR-DG4 & 0.043 & 0.048 & 0.071  \\
        \bottomrule
        \end{tabular}
        \caption{FR-DG Tetrahedra}
        \label{table:tet_fr_dg_cfl_theta0=pi/6_theta1=pi/4}
    \end{subtable}
    \vspace{0.5cm}
    \begin{subtable}{.45\linewidth}
        \centering
        \begin{tabular}{@{}lccc@{}}
        \toprule
         & \multicolumn{3}{c}{$\CFL_{\text{MAX}}$} \\
         \midrule
         & RK3 & RK4 & RK54 \\
         \midrule
        SDRT1 & 0.345 & 0.403 & 0.552  \\ 
        SDRT2 & 0.187 & 0.222 & 0.304  \\ 
        SDRT3 & 0.117 & 0.131 & 0.196  \\ 
        SDRT4 & 0.081 & 0.094 & 0.136  \\ 
        \bottomrule
        \end{tabular}
        \caption{SDRT Prisms}
        \label{table:pri_sdrt_cfl_theta0=pi/6_theta1=pi/4}
    \end{subtable}
    \begin{subtable}{.45\linewidth}
        \centering
        \begin{tabular}{@{}lccc@{}}
        \toprule
         & \multicolumn{3}{c}{$\CFL_{\text{MAX}}$} \\
         \midrule
         & RK3 & RK4 & RK54 \\
         \midrule
        FR-DG1 & 0.237 & 0.269 & 0.394  \\ 
        FR-DG2 & 0.122 & 0.146 & 0.204  \\ 
        FR-DG3 & 0.075 & 0.084 & 0.127  \\ 
        FR-DG4 & 0.052 & 0.061 & 0.088  \\ 
        \bottomrule
        \end{tabular}
        \caption{FR-DG Prisms}
        \label{table:pri_fr_dg_cfl_theta0=pi/6_theta1=pi/4}
    \end{subtable}
    \caption{$\CFL_{\text{MAX}}$ number to ensure linear stability in the linear advection equation with uniform meshes made up of three-dimensional elements with $\theta_0 = \pi/6$ and $\theta_1=\pi/4$.}
\end{table}

Additional tests have been performed utilizing other wave angle configurations and similar conclusions have been obtained.
Therefore, for the sake of brevity, the $\CFL_{\text{MAX}}$ for other wave angle configurations are not depicted in this work.





\FloatBarrier
\section{Von-Neumann analysis with diffusion}
\label{sec:von_neumann_diffusion}


The Von-Neumann analysis depicted in \Sref{sec:von_neumann_advection} for the linear advection equation may also be applied to the linear diffusion equation which reads
\begin{equation}
    \frac{\partial \var(\vect{x}, t)}{\partial t} + \divergence (\mu \grad \var(\vect{x}, t)) = 0 .
    \label{eq:linear_diffusion}
\end{equation}
The dissipation and dispersion errors introduced by SEM when discretizing the latter equation may be predicted using similar tools as those introduced to analyze SEM in the linear advection equation with uniform meshes.
One of the modifications that need to be taken into account is that, for the initial condition considered in \Eref{eq:initial_condition_fourier_mode}, the analytical solution is given by
\begin{equation}
    \var(\vect{x}, t) = \var(\vect{x}, 0)\e^{ -\mu \kappa^2 t} .
\end{equation}
Hence, the temporal wavenumber is $\omega = \imag \mu \kappa^2$.
The discretized version of the diffusion equation reads
\begin{equation}
    \frac{\text{d} \var_{\ipol \iele}}{\text{d} t} = \frac{\mu}{h^2} \sum_{i \in \stencil_{\iele}} \jacobianrhs_{\ipol i \inu} \var_{ i \inu} .
    \label{eq:jacobian_diffusion_rhs}
\end{equation}
For pure-diffusion problems, one may redefine the CFL number as
\begin{equation}
    \CFL = \frac{\mu \timestep}{\cellsize^2} .
\end{equation}

In this study, only the RB1 scheme will be considered to compute the gradients and viscous fluxes, i.e. $\beta = 0$.
Moreover, the penalty term will be set to $\penaltyterm = 0$.
Such choices are motivated due to the fact that an appropriate and symmetric definition of the left and right face connectivity in meshes made up of other than tensor-product elements is cumbersome.
If the latter mesh property is not ensured and the viscous solver is not symmetric (i.e. if $\beta \neq 0$ and/or $\penaltyterm \neq 0$), then the RHS Jacobian matrix will not possess circulant block properties, and therefore, the process described in \Sref{sec:von_neumann_advection} to assess numerical dissipation and dispersion of SEM would not be valid.

It is worth mentioning, that the use of the RB1 formulation has important consequences in the analysis as, due to the non-compactness of the method, the number of neighbors with non-zero contributions in the RHS Jacobian matrix $\jacobianrhs$ is increased with respect to that obtained in pure-advection problems.
In particular, one should consider two face neighbor layers to appropriately carry out the process to compute the dissipation and dispersion errors described in \Sref{sec:von_neumann_advection}.
Hence, in pure-diffusion problems, the set $\stencil_{\iele}$ stores the unique indices of the direct face neighbors of the element $\iele$ and the face neighbors of the aforementioned neighbors.

After computing the eigenvalues from \Eref{eq:patternvar_diagonalized} resulting from the discretization of the diffusion equation, one may compute the numerical temporal wavenumber $\omega^{\numerical}$, using the dissipation and dispersion measures, following Eqs. \ref{eq:dissipation_measure} and \ref{eq:dispersion_measure}.
With such a value of the temporal wavenumber, the error estimator of the numerical temporal wavenumber defined in \Eref{eq:diss_disp_error_measure} may be used to assess the accuracy of the numerical schemes.
As it was found in \cite{Alhawwary_2020}, no dispersion errors were obtained when analyzing the schemes developed in this work.
The reason behind the lack of dispersion errors when dealing with pure-diffusion linear equations in SEM will be studied in future works.

\begin{remark}
Despite the fact that for pure-diffusion problems the RHS Jacobian is independent of the wave angle, both the eigenvalues and eigenvectors shape of \Eref{eq:patternvar_diagonalized} are not.
Therefore, the numerical errors of SEM when discretizing pure-diffusion equations depend on the wave angle.
\end{remark}

\begin{figure}[h]
    \begin{subfigure}{0.45\textwidth}
        \centering
        \includegraphics[width=0.9\textwidth]{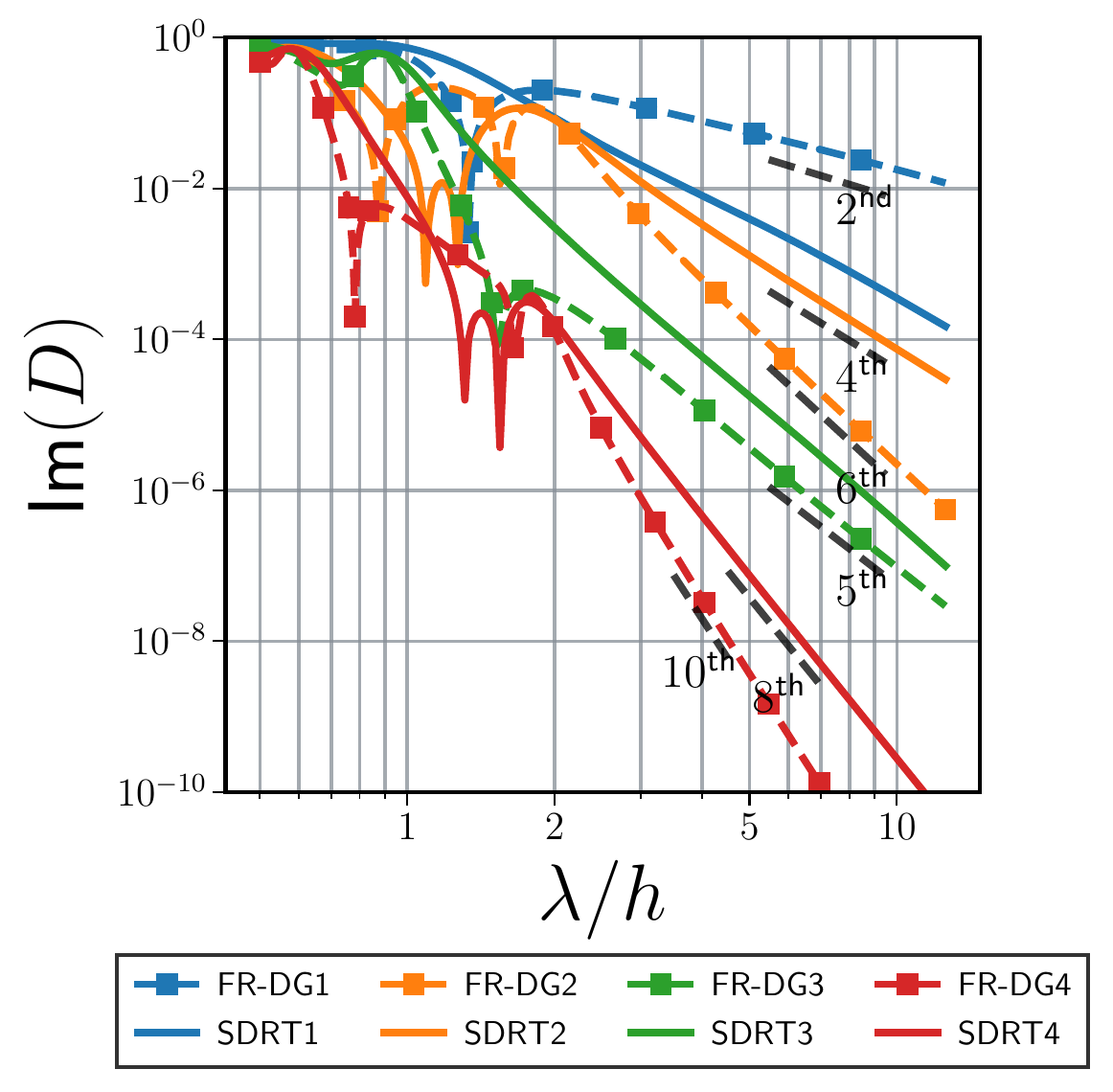}
        \caption{Quadrilaterals}
        \label{fig:diff_diss_quad_theta0=pi/6}
    \end{subfigure}
    \begin{subfigure}{0.45\textwidth}
        \centering
        \includegraphics[width=0.9\textwidth]{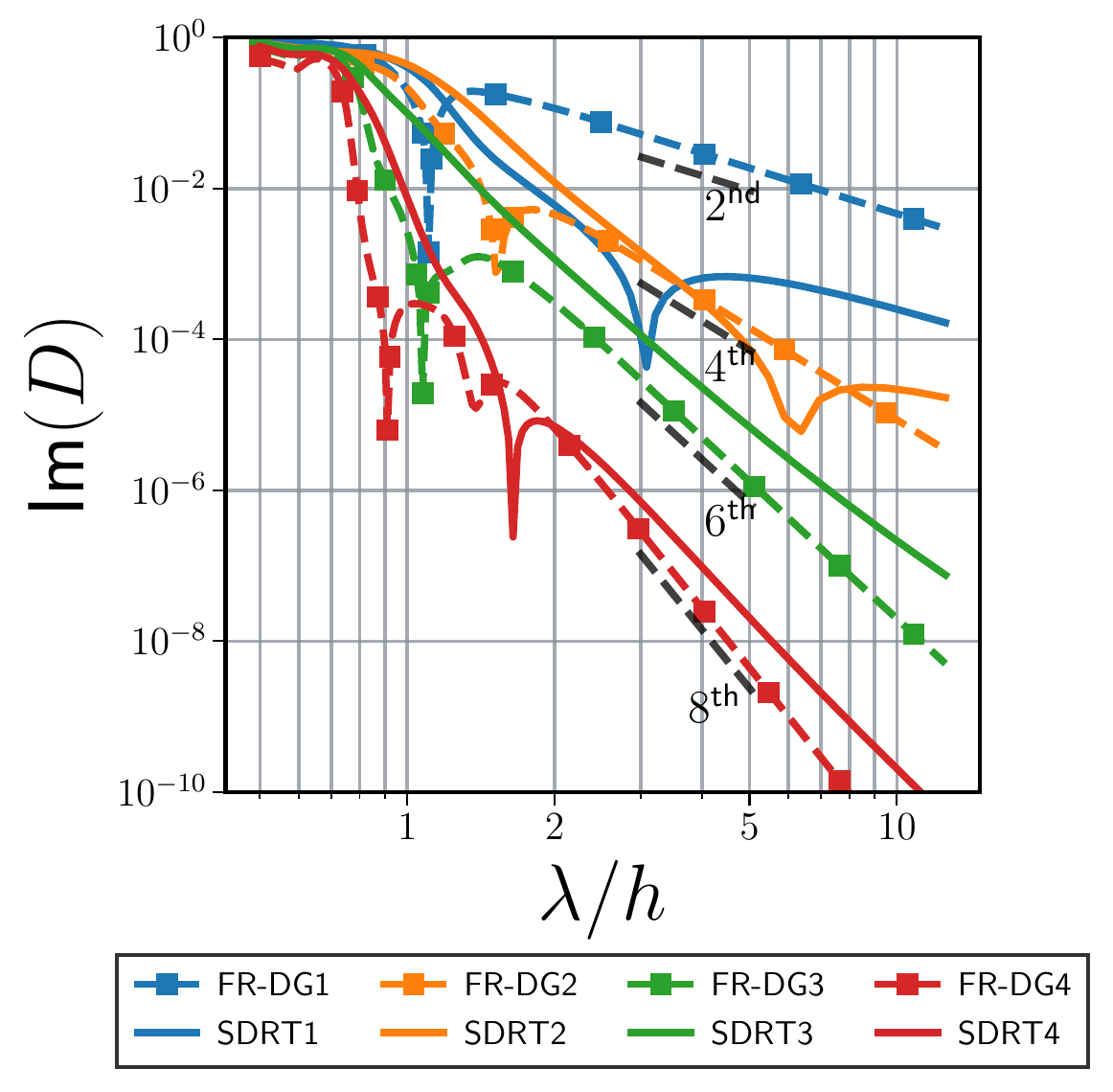}
        \caption{Triangles}
        \label{fig:diff_diss_tri_theta0=pi/6}
    \end{subfigure}
    \caption{Dissipation errors in pure-diffusion problems with exponential time integration obtained with two-dimensional elements, a wave angle equal to $\theta_0 = \pi/6$ and $m=m_c/10$.}
    \label{fig:diff_diss_2D_m/10_theta0=pi/6}
\end{figure}

\Fref{fig:diff_diss_2D_m/10_theta0=pi/6} depicts the dissipation errors obtained using quadrilateral and triangular elements, $\theta_0 = \pi/6$, exponential time integration and $m=m_c/10$.
 When analyzing the order of accuracy of the latter type of elements in pure-diffusion problems it can be observed that both FR-DG and SDRT schemes show $2\degree$ order of accuracy.
 Nevertheless and as it was observed in the pure-advection analysis, the SDRT order of accuracy $\degree = 4$ seems to degrade for well-resolved waves.
 On the other hand, the order of accuracy with quadrilateral elements behaves differently depending on whether $\degree$ is even or odd.
SDRT schemes present $2\degree$ order of accuracy, with the exception of the SDRT1 scheme which shows $2\degree + 2$ order of accuracy.
Additionally, FR-DG methods show $2\degree$ and $2\degree + 2$ order of accuracy for odd and even polynomial degrees respectively.
The reason behind such a behavior of the schemes may be related to error cancellation properties of the RB1 scheme with even polynomial degrees.
 
\begin{figure}[h]
    \begin{subfigure}{0.45\textwidth}
        \centering
        \includegraphics[width=0.9\textwidth]{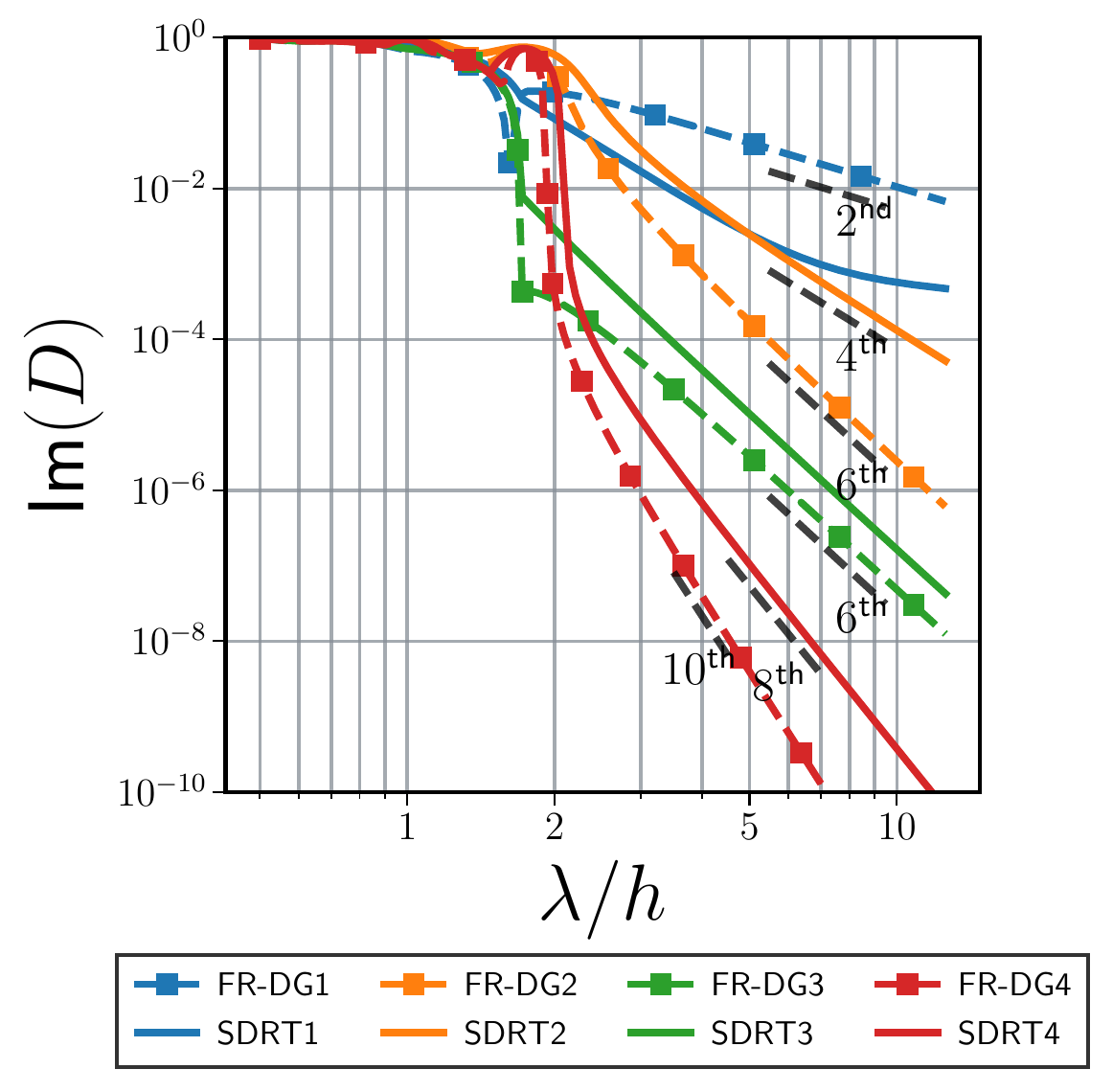}
        \caption{With absolute value.}
        \label{fig:diff_diss_quad_m=10_theta0=pi/6}
    \end{subfigure}
    \begin{subfigure}{0.45\textwidth}
        \centering
        \includegraphics[width=0.9\textwidth]{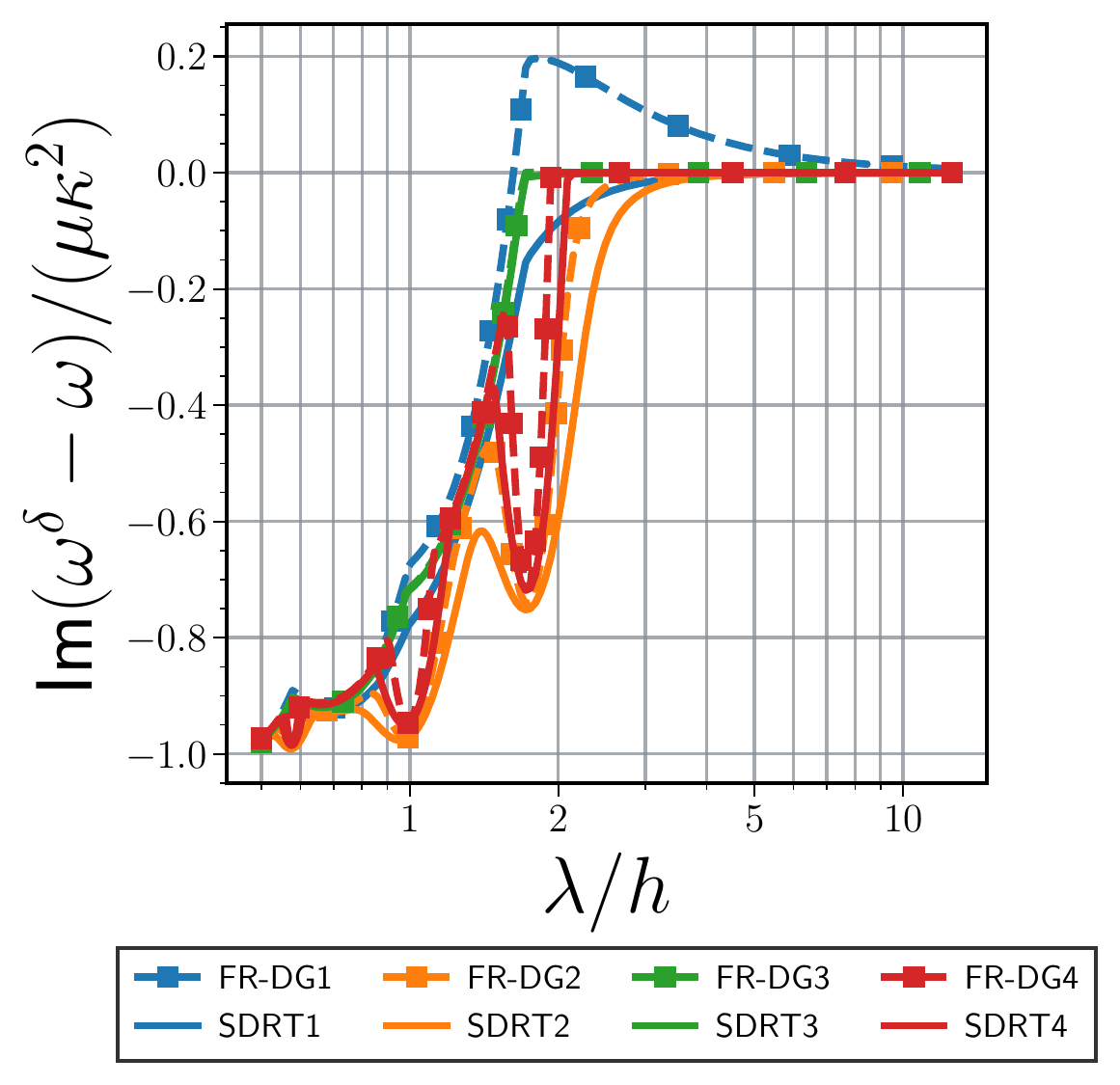}
        \caption{Without absolute value.}
        \label{fig:diff_diss_quad_m=10_noabs_theta0=pi/6}
    \end{subfigure}
    \caption{Dissipation errors in pure-diffusion problems with exponential time integration obtained with quadrilateral elements, a wave angle equal to $\theta_0 = \pi/6$ and $m=10m_c$. The figure on the left was computed using \Eref{eq:diss_disp_error_measure} while the figure on the right was obtained by removing the absolute value operator in \Eref{eq:diss_disp_error_measure}.}
\end{figure}

To assess the influence of the number of iterations in the dissipation maps, \Fref{fig:diff_diss_quad_m=10_theta0=pi/6} illustrates the dissipation errors obtained using quadrilateral  elements, $\theta_0 = \pi/6$, exponential time integration and $m=10 m_c$.
Similarly as it was observed in pure-advection problems, the dissipation errors are not heavily distorted when increasing the number interations until the first aliasing limit $\kappa h < \pi/\cos \theta_0$.
For non-aliased wavenumbers the dissipation errors are very similar to those observed with $m=m_c/10$ in \Fref{fig:diff_diss_quad_theta0=pi/6}, with the exception of the SDRT1 scheme which shows a reduction of its order of accuracy for well-resolved cases.
On the other hand, with aliased wavenumbers the dissipation errors drastically increase for a high number of iterations.
Recently, \cite{Alhawwary_2020} indicated that the dissipation errors for high wavenumbers with pure-diffusion and SEM reduced the dissipation, i.e. the errors are related to the lack of appropriate diffusion provided by the numerical scheme.
To further validate this observation, \Fref{fig:diff_diss_quad_m=10_noabs_theta0=pi/6} represents the dissipation errors from \Eref{eq:diss_disp_error_measure} removing the absolute value from its expression.
The figure illustrates that, for aliased wavenumbers and high number of iterations, the dissipation errors are always negative (except for the FR-DG1 scheme), i.e. the diffusion predicted by the schemes is smaller than expected.
Such an issue is important in conservation laws presenting diffusion terms, since for very high wavenumbers the dissipation of the numerical schemes is close zero, far away from the \textit{physical} dissipation expected for such high wavenumbers.
Additionally, it can also be observed that, the higher the polynomial degree, the lower the numerical diffusion.
Moreover, all SDRT schemes display increased dissipation errors than the FR-DG method, possibly pointing out that SDRT schemes might need of additional numerical dissipation in advection-diffusion problems to avoid the accumulation of energy within the smallest scales.

\begin{figure}
    \begin{subfigure}{0.45\textwidth}
        \centering
        \includegraphics[width=0.9\textwidth]{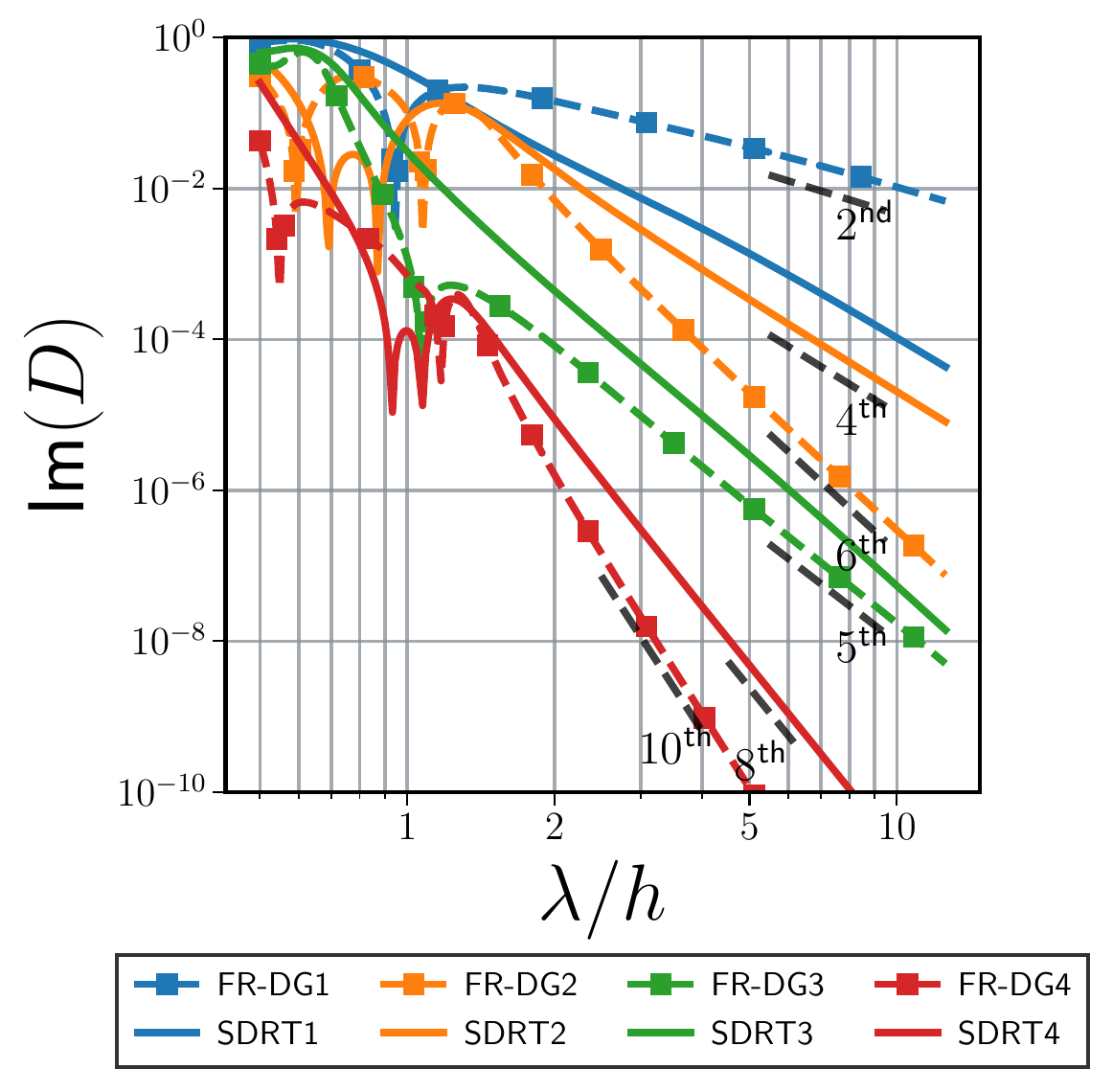}
        \caption{Hexahedral}
        \label{fig:diff_diss_hex_theta0=pi/6_theta1=pi/4}
    \end{subfigure}
    \begin{subfigure}{0.45\textwidth}
        \centering
        \includegraphics[width=0.9\textwidth]{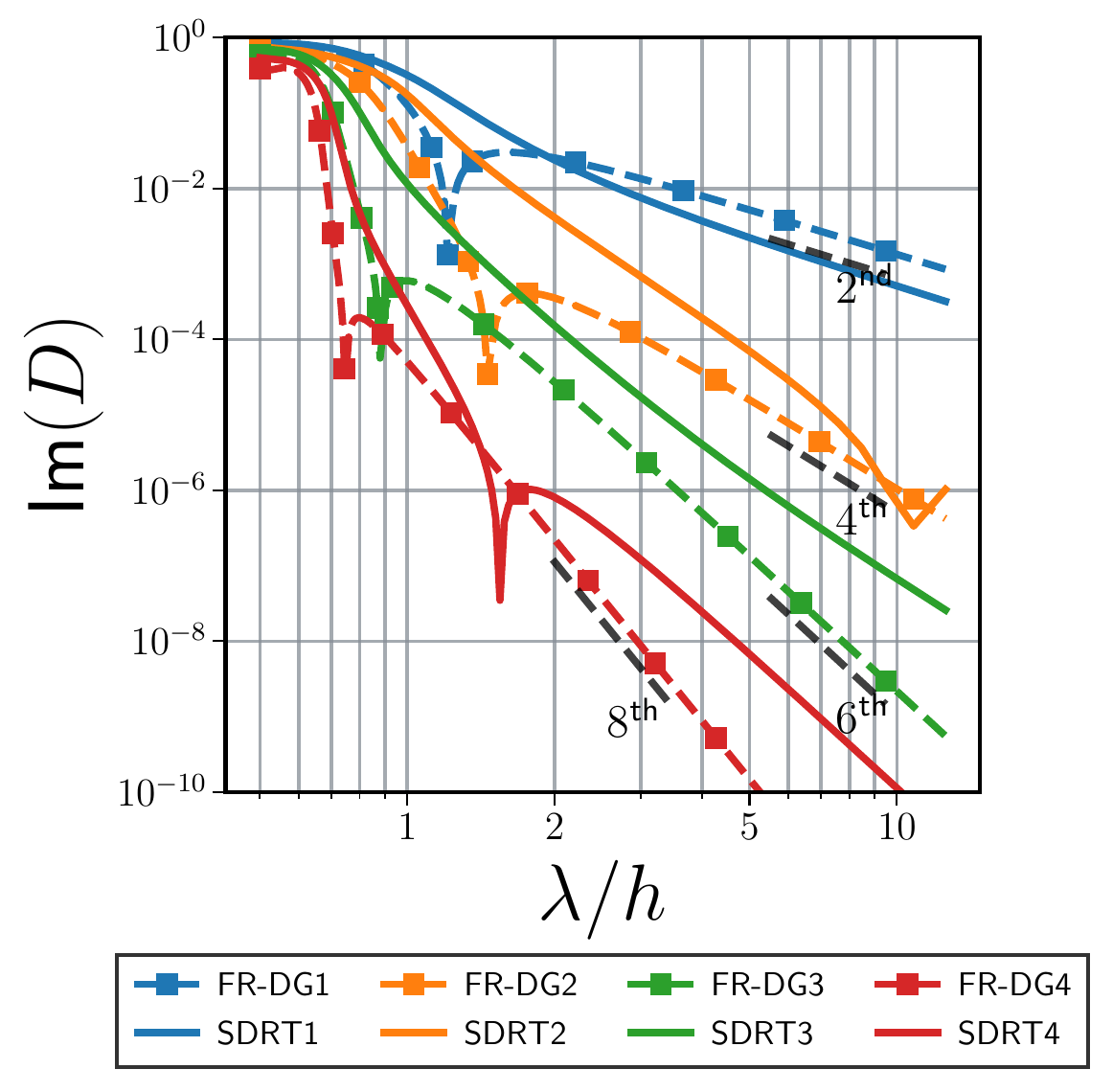}
        \caption{Tetrahedral}
        \label{fig:diff_diss_tet_theta0=pi/6_theta1=pi/4}
    \end{subfigure} \\
    \centering
    \begin{subfigure}{0.45\textwidth}
        \centering
        \includegraphics[width=0.9\textwidth]{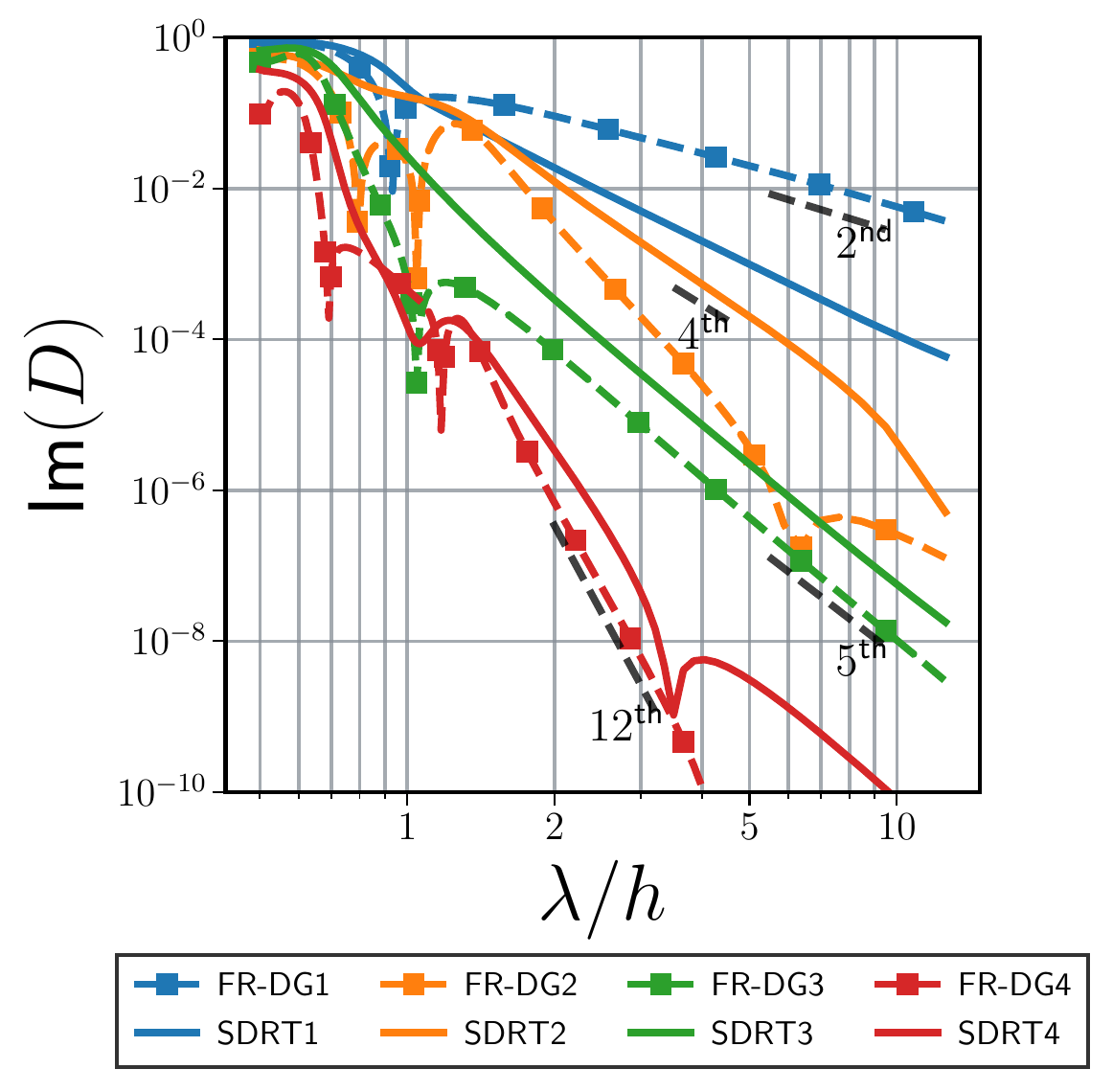}
        \caption{Prisms}
        \label{fig:diff_diss_pri_theta0=pi/6_theta1=pi/4}
    \end{subfigure} \\
    \caption{Dissipation errors in pure-diffusion problems with exponential time integration and $m=m_c/10$ obtained with three-dimensional elements and wave angles equal to $\theta_0 = \pi/6$ and $\theta_1 = \pi/4$.}
    \label{fig:diff_diss_3D_theta0=pi/6_theta1=pi/4}
\end{figure}

Similarly,  \Fref{fig:diff_diss_3D_theta0=pi/6_theta1=pi/4} shows the dissipation errors obtained with three-dimensional elements with $\theta_0 = \pi / 6$ and $\theta_1 = \pi / 4$.
The dissipation error of hexahedral elements is very similar to that observed for quadrilateral elements in the two-dimensional analysis.
FR-DG displays a superconvergent $2\degree + 2$ order of accuracy with even degree polynomials, while its order of accuracy is reduced to $2\degree$ with odd degree polynomials.
SDRT schemes show $2\degree$ order of accuracy.
On the other hand, the behavior of the numerical dissipation with tetrahedral elements and SDRT schemes is slightly different than that observed with triangular elements.
In particular, as the polynomial degree increases the order of accuracy for $\degree \ge 3$ seems to degrade for well-resolved waves.
This issue may be related to the influence of spurious modes and/or the lack of superconvergence properties of the physical/spectral radius eigenmode.
At last, the dissipation errors with prismatic elements display an interesting behavior.
For low values of the cells per wavelength parameter $\lambda / \cellsize < 5$, the dissipation errors resemble those of hexahedron elements, displaying superconvergent behavior for even degree polynomials and FR-DG schemes.
On the other hand, for well-resolved waves, the results are close to those observed with tetrahedron elements.

\begin{remark}
The discussion of discretization errors related to time integration schemes is avoided for the sake of brevity, as similar conclusions than those addressed in \Sref{sec:temporal_stability_advection} have been observed.
\end{remark}

\subsection{Temporal linear stability}
\label{sec:temporal_stability_diffusion}

This section aims to determine the maximum $\CFL_{\text{MAX}}$ number which ensures temporal stability for periodic uniform meshes and pure-diffusion problems.
To compute this parameter, the same process as the one discussed in \Sref{sec:temporal_stability_advection} will be used.
Nevertheless, it is worth noting that as matrix $\jacobianrhs$ is not a function of the wave angle, the $\CFL_{\text{MAX}}$ number is independent of the initial condition.

\Tref{table:diff_cfl_2D_theta=pi/6} displays the $\CFL_{\text{MAX}}$ number for quadrilateral and triangular elements with periodic uniform meshes and pure-diffusion problems.
As it was observed for pure-advection test cases with quadrilateral elements, the $\CFL_{\text{MAX}}$ number of SDRT schemes is substantially higher than that of FR-DG methods.
Such a difference is slightly reduced when considering triangular elements, although it still remains considerably important.
For the sake of brevity, the values of the $\CFL_{\text{MAX}}$ number for three-dimensional elements are not presented.

\begin{table}[!htb]
    \centering
    \begin{subtable}{.45\linewidth}
        \centering
        \begin{tabular}{@{}lccc@{}}
        \toprule
         & \multicolumn{3}{c}{$\CFL_{\text{MAX}}$} \\
         \midrule
         & RK3 & RK4 & RK54 \\
         \midrule
        SDRT1 & 0.2791 & 0.3094 & 0.5174  \\ 
        SDRT2 & 0.0348 & 0.0386 & 0.0646  \\ 
        SDRT3 & 0.0157 & 0.0174 & 0.0291  \\ 
        SDRT4 & 0.0055 & 0.0061 & 0.0103  \\ 
        \bottomrule
        \end{tabular}
        \caption{SDRT Quadrilateral}
        \label{table:diff_quad_sdrt_cfl}
    \end{subtable}
    \begin{subtable}{.45\linewidth}
        \begin{tabular}{@{}lccc@{}}
        \toprule
         & \multicolumn{3}{c}{$\CFL_{\text{MAX}}$} \\
         \midrule
         & RK3 & RK4 & RK54 \\
         \midrule
        FR-DG1 & 0.1570 & 0.1740 & 0.2910  \\ 
        FR-DG2 & 0.0200 & 0.0222 & 0.0371  \\ 
        FR-DG3 & 0.0106 & 0.0117 & 0.0197  \\ 
        FR-DG4 & 0.0032 & 0.0036 & 0.0060  \\ 
        \bottomrule
        \end{tabular}
        \caption{FR-DG Quadrilateral}
        \label{table:diff_quad_fr_dg_cfl}
    \end{subtable}\\
    \vspace{0.5cm}
    \begin{subtable}{.45\linewidth}
        \centering
        \begin{tabular}{@{}lccc@{}}
        \toprule
         & \multicolumn{3}{c}{$\CFL_{\text{MAX}}$} \\
         \midrule
         & RK3 & RK4 & RK54 \\
         \midrule
        SDRT1 & 0.0609 & 0.0675 & 0.1130  \\ 
        SDRT2 & 0.0206 & 0.0228 & 0.0382  \\ 
        SDRT3 & 0.0076 & 0.0085 & 0.0142  \\ 
        SDRT4 & 0.0031 & 0.0035 & 0.0058  \\
        \bottomrule
        \end{tabular}
        \caption{SDRT Triangles}
        \label{table:diff_tri_sdrt_cfl}
    \end{subtable}
    \begin{subtable}{.45\linewidth}
        \begin{tabular}{@{}lccc@{}}
        \toprule
         & \multicolumn{3}{c}{$\CFL_{\text{MAX}}$} \\
         \midrule
         & RK3 & RK4 & RK54 \\
         \midrule
        FR-DG1 & 0.0418 & 0.0464 & 0.0776  \\ 
        FR-DG2 & 0.0144 & 0.0160 & 0.0267  \\ 
        FR-DG3 & 0.0054 & 0.0060 & 0.0101  \\ 
        FR-DG4 & 0.0027 & 0.0030 & 0.0050  \\ 
        \bottomrule
        \end{tabular}
        \caption{FR-DG Triangles}
        \label{table:diff_tri_fr_dg_cfl}
    \end{subtable}
    \caption{$\CFL_{\text{MAX}}$ number to ensure temporal linear stability in the linear diffusion equation with uniform meshes made up of two-dimensional elements.}
    \label{table:diff_cfl_2D_theta=pi/6}
\end{table}

\begin{remark}
As stated in \cite{Watkins2016}, the temporal linear stability criterion for linear advection diffusion conservation laws slightly differs from that of pure-advection or pure-diffusion equations.
\end{remark}

\section{Numerical experiments}
\label{sec:numerical_experiments}

\subsection{Linear Advection Diffusion}
\label{sec:numerical_experiments_linadvecdiff}

To experimentally assess the order-of-accuracy of the SDRT in linear test cases, the linear advection-diffusion equation is solved in two-dimensional and three-dimensional configurations within a domain $\domain \in [0, 2 \pi L]^{\ndim}$ with periodic boundary conditions
\begin{equation}
    \flux(\vect{x}, t) = \advecvel \var(\vect{x}, t)  - \mu \grad \var(\vect{x}, t) ,
\end{equation}
subject to the initial condition
\begin{equation}
    \var(\vect{x}, 0) = \sin{\sum_{i=0}^{\ndim - 1} \frac{ x_i}{L}} .
\end{equation}
The latter initial field is related to the imaginary part of a Fourier mode with wavenumber $\kappa = \sqrt{\ndim}/L$ and wavelength $\lambda = 2 \pi L / \sqrt{\ndim}$.
The analytical solution of the previous problem is given by
\begin{equation}
    \var(\vect{x}, t) = \var(\vect{x} - \advecvel t, 0) \e^{- \mu \ndim t/L^2} .
\end{equation}

The behavior of the linear advection diffusion equation is governed by the Peclet number
\begin{equation}
    \mathrm{Pe} = \frac{\advecvelnovect}{\mu L} \quad .
\end{equation}
Such a non-dimensional parameter plays a similar role to the Reynolds number in the Navier-Stokes equations.
High values of the Peclet number are related to advection-dominated solutions, while low numbers of it imply diffusion-dominated solutions.
To ease the analysis of the numerical errors in future sections, all components of the advection velocity vector are supposed constant.

\subsubsection{Asymptotic order of accuracy with pure advection}

In order to validate the numerical accuracy described in \Sref{sec:von_neumann_advection} and \Sref{sec:von_neumann_diffusion}, it is important to avoid taking into account the initial projection error of the solution on the polynomial basis of the schemes \cite{Guo_2013}.
This is ensured by computing the $L_2$-norm of the solution error using the numerical solution as the reference solution
\begin{equation}
    L_2(\cellsize, m) = \sqrt{\frac{1}{|\volume|} \int_{\volume} \bigg[\var^{\numerical}\left(\vect{x}, \frac{2m\pi L}{\advecvelnovect_0} \right) - \var^{\numerical}\left(\vect{x}, \frac{2(m - 1)\pi L}{\advecvelnovect_0} \right)\bigg]^2 \, \text{d} \volume}
    \label{eq:linadvec_l2_definition}
\end{equation}
where $m$ in the number of periods and $\volume$ is the domain volume.
Such a choice allows to distinguish the numerical errors associated with the physical eigenmode by ensuring the dissipation of all other \textit{spurious} eigenmodes which appear due to the projection of the initial condition on the solution nodal basis.
In this study, the norm defined in \Eref{eq:linadvec_l2_definition} is assessed with a quadrature of sufficient degree, i.e. such that it integrates exactly at least a $2\degree$ polynomial by interpolating the numerical solution to the correspondent quadrature points.

Under certain conditions, the previously described methodology allows to observe the numerical accuracy predicted in \Sref{sec:von_neumann_advection} with linear analysis \cite{Guo_2013}.
If one computes the error norm using the analytical solution as the reference solution, i.e. if $m = 1$, then the schemes will show $\degree + 1$ order of accuracy due to the strong dissipation and dispersion of the eigenmodes other than the physical mode characteristic of the first steps of the simulation \cite{Frean_2019, Guo_2013}.
To minimize the temporal discretization errors RK54 \cite{Carpenter1994Fourthorder2R}  with $\CFL=0.0025$ is employed.
Nevertheless, it is worth mentioning that the use of very higher-order RK methods is recommended to observe the appropriate order of accuracy of spatial discretization schemes cf.\ \cite{Guo_2013}.

\begin{figure}
    \begin{subfigure}{0.45\textwidth}
    \centering
    \includegraphics[width=1.0\textwidth]{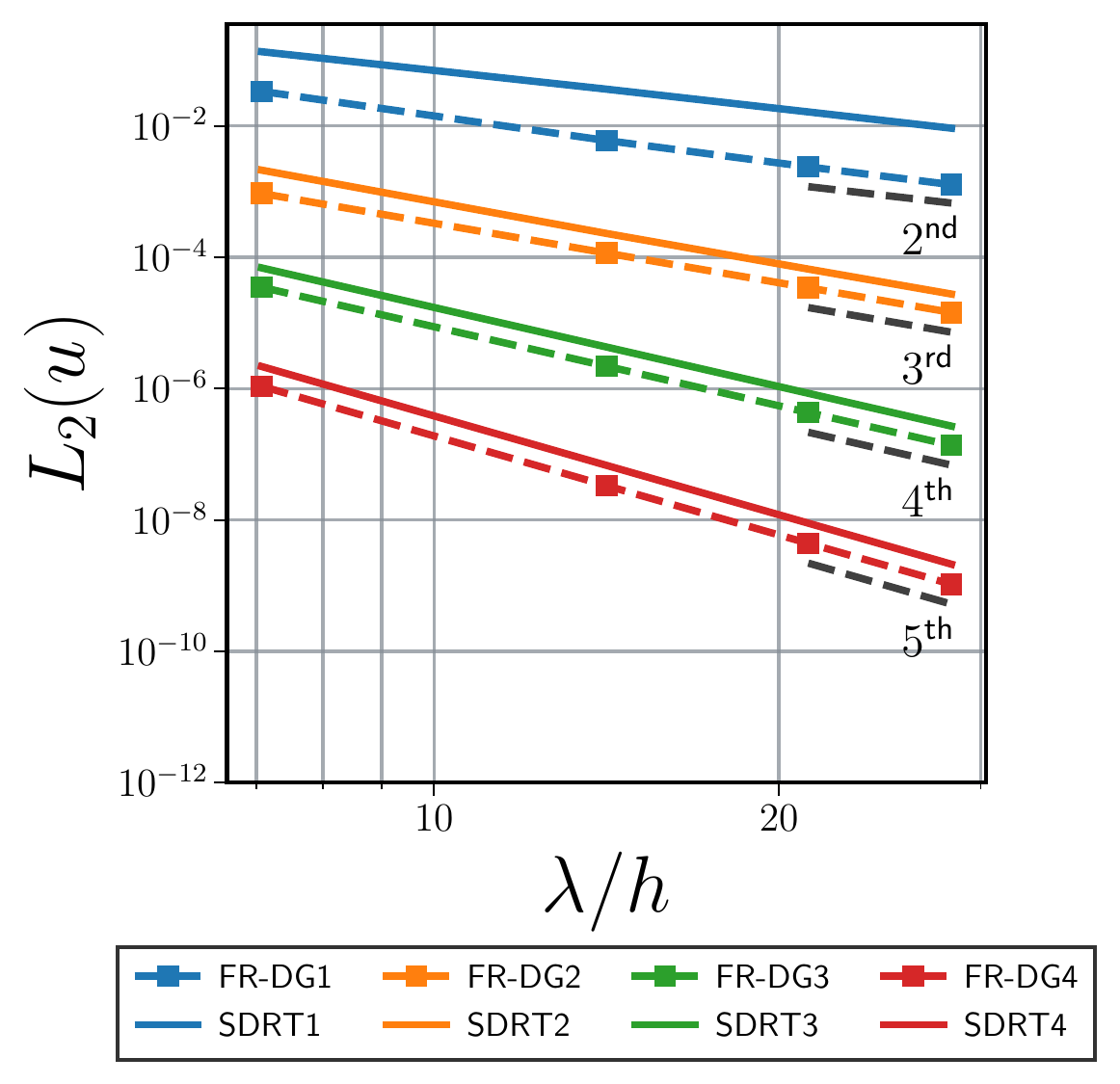}
    \caption{$m = 1$}
    \label{fig:l2_no_ho_quad_sd_fr}    
    \end{subfigure}
    \begin{subfigure}{0.45\textwidth}
    \centering
    \includegraphics[width=1.0\textwidth]{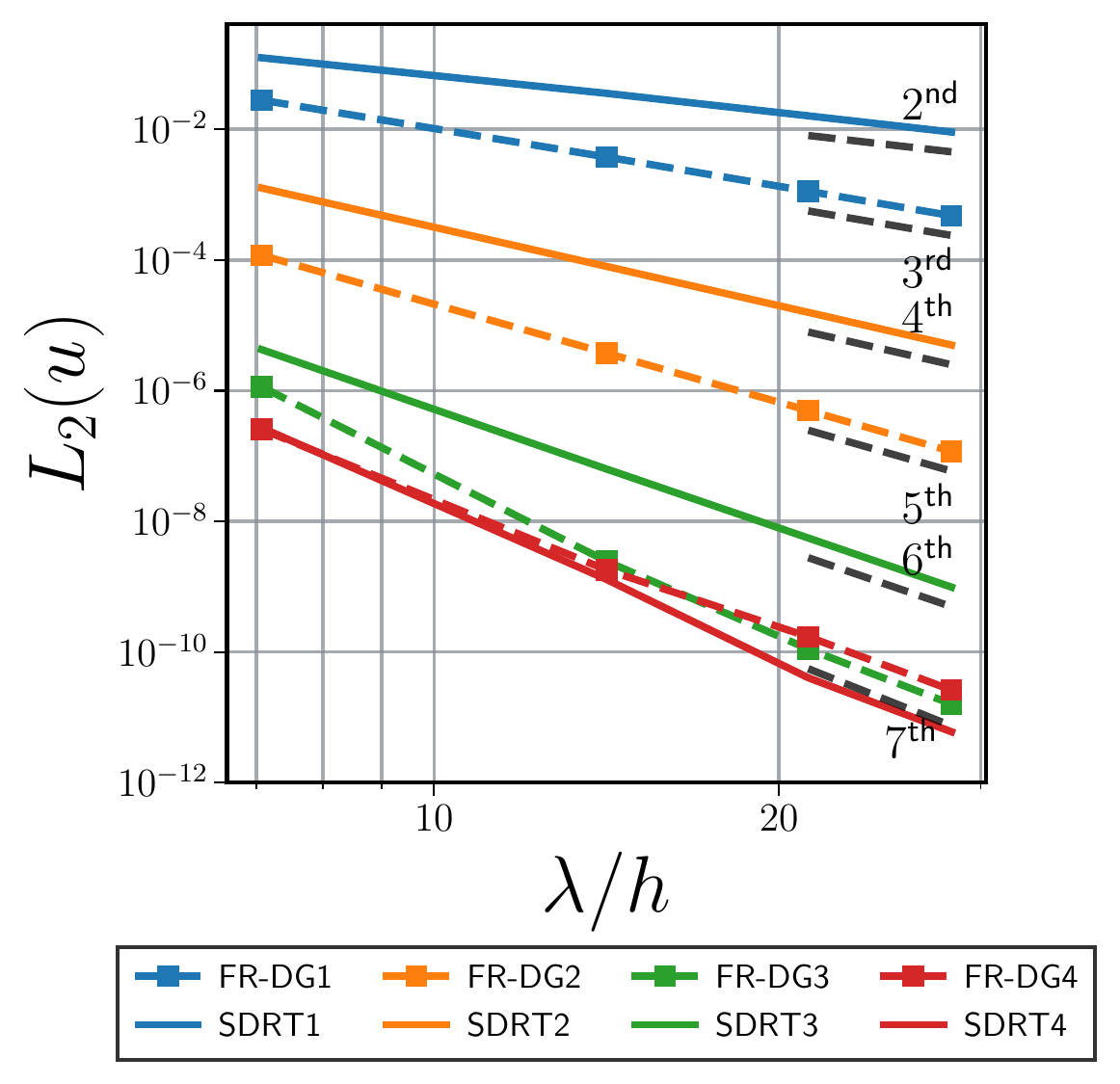}
    \caption{$m = 2$}
    \label{fig:l2_ho_quad_sd_fr} 
    \end{subfigure}
    \caption{$L_2$-norm of the solution error evaluated using \Eref{eq:linadvec_l2_definition} with $m = 1$ (left) and $m = 2$ (right) and obtained using different schemes together with meshes made up of quadrilateral elements.}
    \label{fig:l2_advect_ho_quad}
\end{figure}

\Fref{fig:l2_advect_ho_quad} show the error norm for $\degree= 1, 2, 3$ and $4$ FR-DG and SDRT schemes with quadrilateral elements using $m = 1$ and $m = 2$.
The results illustrate that SDRT and FR-DG schemes display $\degree + 1$ order errors when $m = 1$ as predicted by \citet{Guo_2013}.
Furthermore, when using $m = 2$ the predicted $2\degree$ and $2\degree + 1$ order of accuracy of SDRT and FR-DG schemes in pure advection problems is observed for  $\degree \le 3$.
With $\degree = 4$ the expected order of accuracy is not obtained.
This could be explained due to the lower dissipation rates of the spurious modes and due to the use of a fourth order accurate RK method.
It is worth noting that the data indicate that FR-DG schemes present reduced error values compared to SDRT schemes in quadrilateral elements for both $m = 1$ and $2$.
For the sake of brevity, the computation of the order of accuracy of three-dimensional elements is avoided, as the results in this study were equivalent to those obtained with two-dimensional elements.
The lack of superconvergence when $m = 1$ has been largely discussed in the literature \cite{Guo_2013} and it can be demonstrated that it is related to the time integration errors and the dissipation of spurious eigenmodes resulting from the projection of the initial condition on to the solution nodal basis during the initialization process.

\begin{figure}
    \begin{subfigure}{0.45\textwidth}
    \centering
    \includegraphics[width=1.0\textwidth]{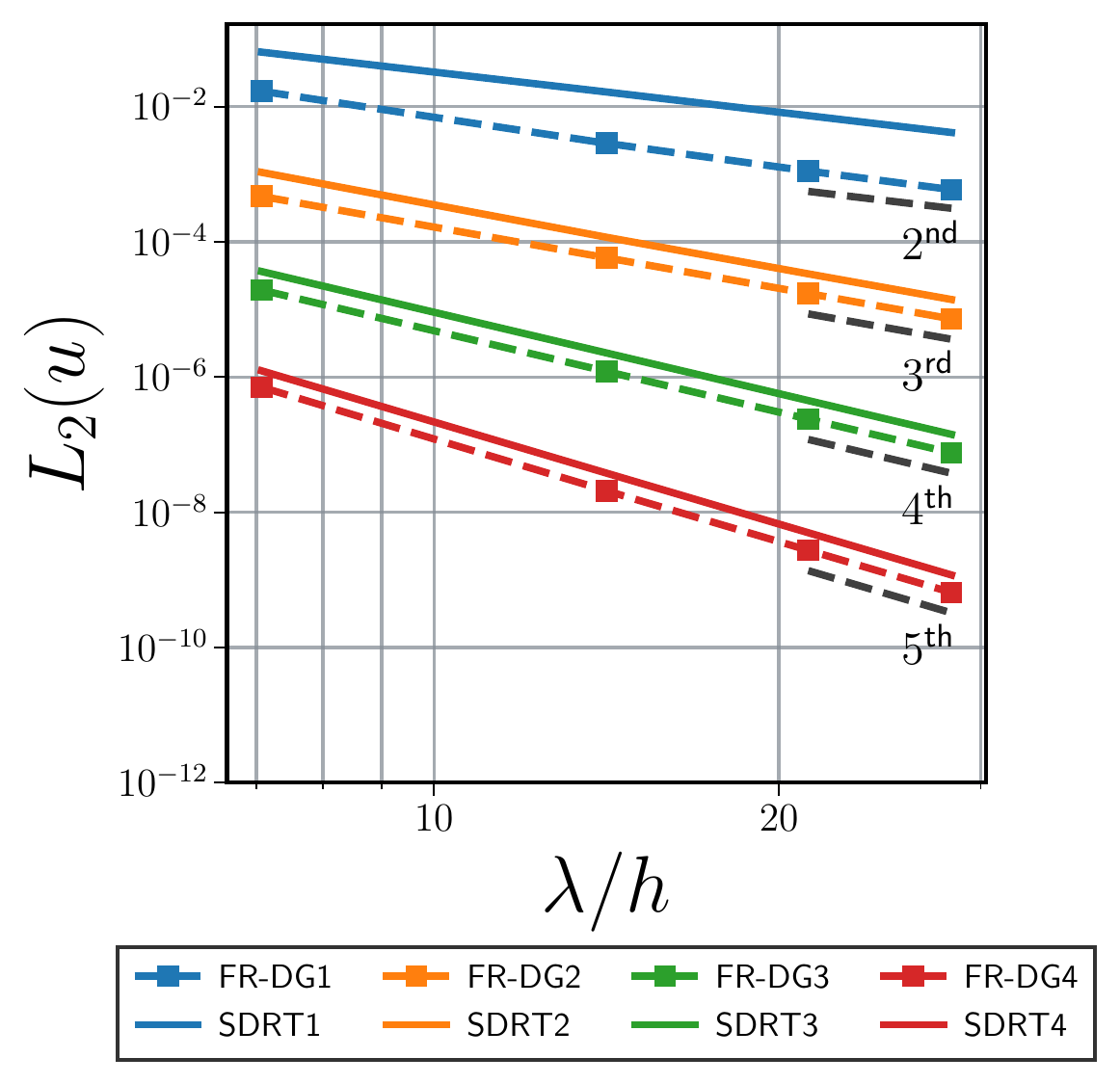}
    \caption{$m = 1$}
    \label{fig:l2_no_ho_tri_sd_fr}    
    \end{subfigure}
    \begin{subfigure}{0.45\textwidth}
    \centering
    \includegraphics[width=1.0\textwidth]{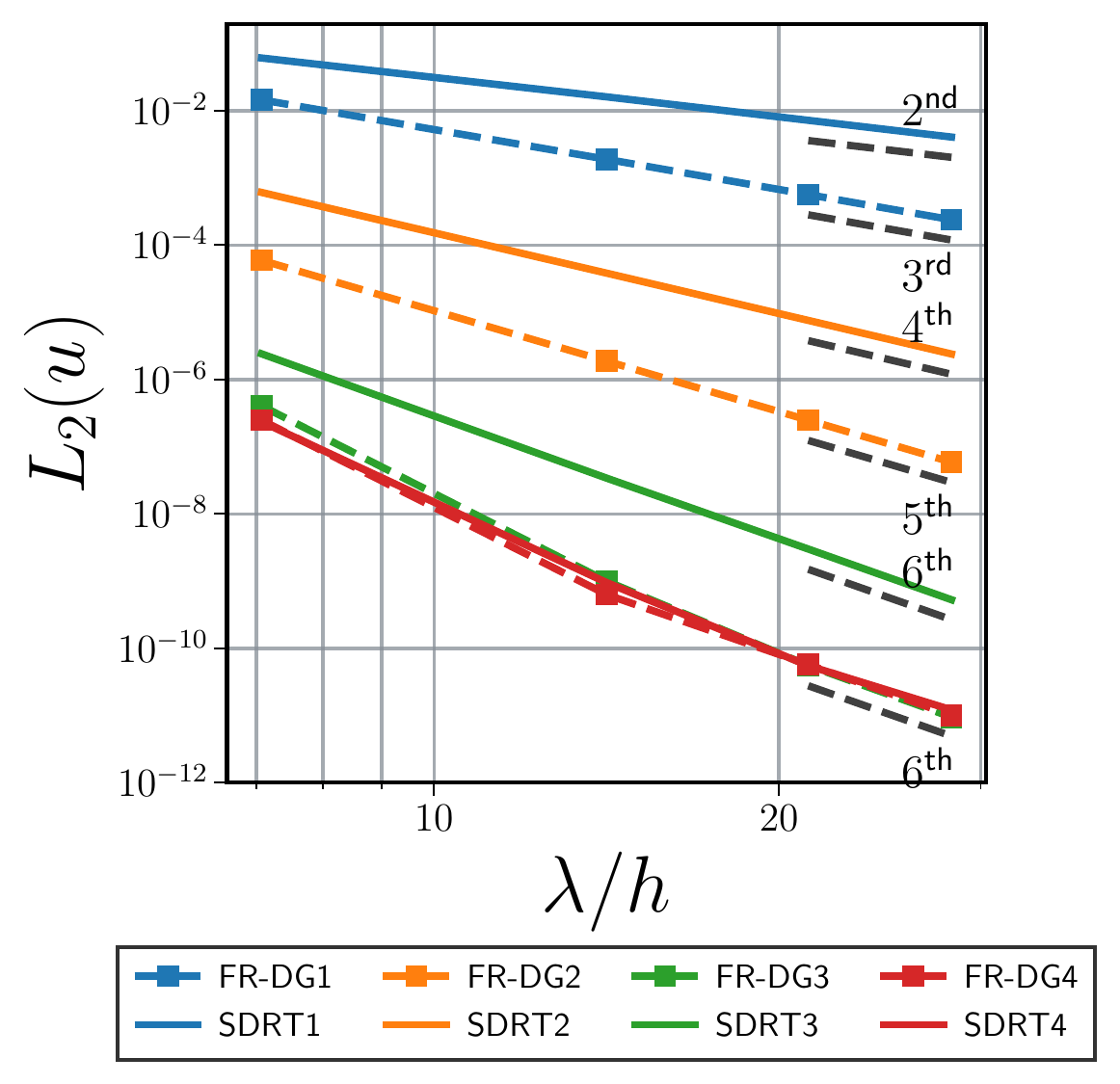}
    \caption{$m = 2$}
    \label{fig:l2_ho_tri_sd_fr} 
    \end{subfigure}
    \caption{$L_2$-norm of the solution error evaluated using \Eref{eq:linadvec_l2_definition} with $m = 1$ (left) and $m = 2$ (right) and obtained using different schemes together with meshes made up of triangular elements.}
    \label{fig:l2_advect_ho_tri}
\end{figure}

\begin{figure}
    \begin{subfigure}{0.45\textwidth}
    \centering
    \includegraphics[width=1.0\textwidth]{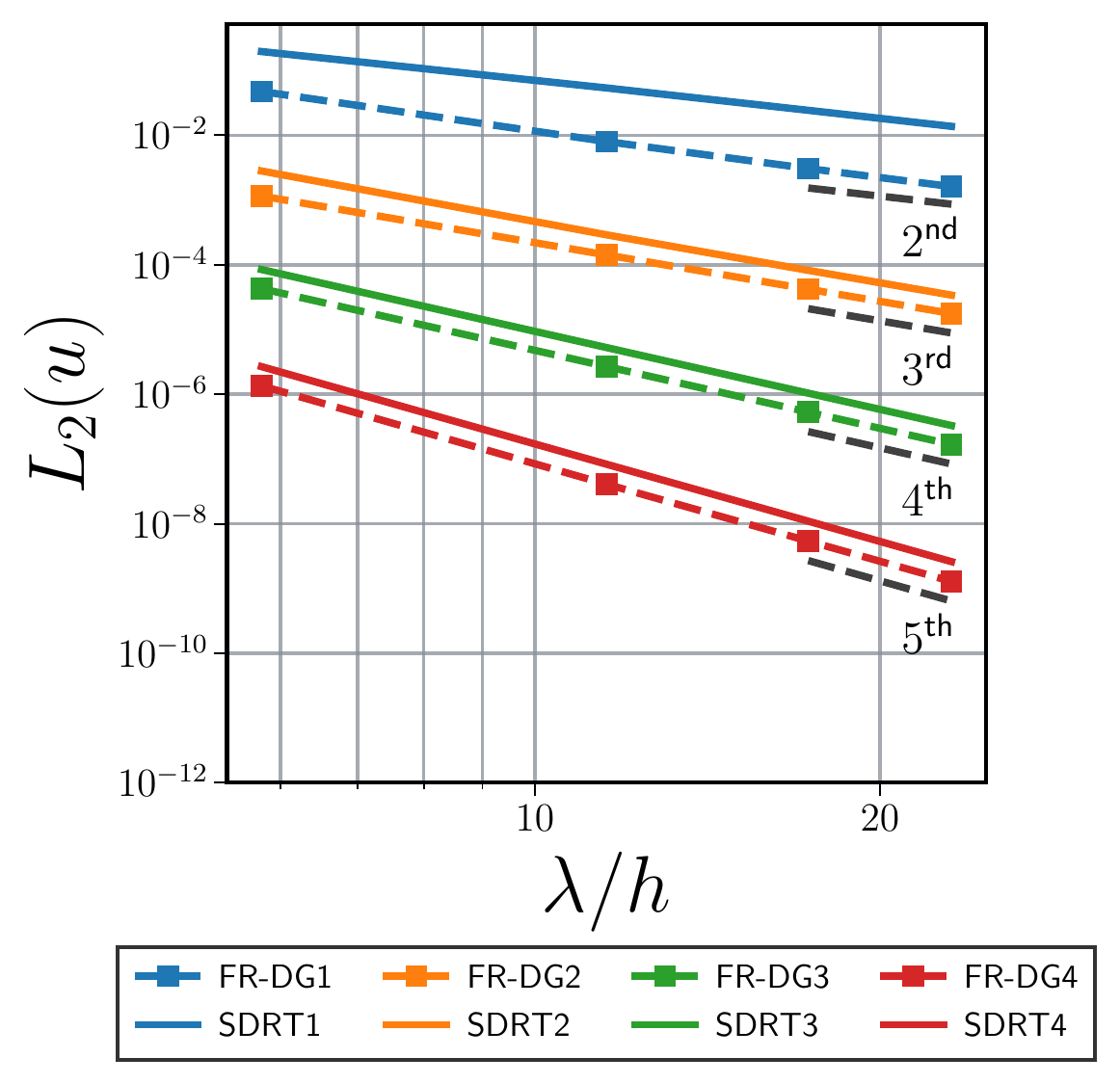}
    \caption{$m = 1$}
    \end{subfigure}
    \begin{subfigure}{0.45\textwidth}
    \centering
    \includegraphics[width=1.0\textwidth]{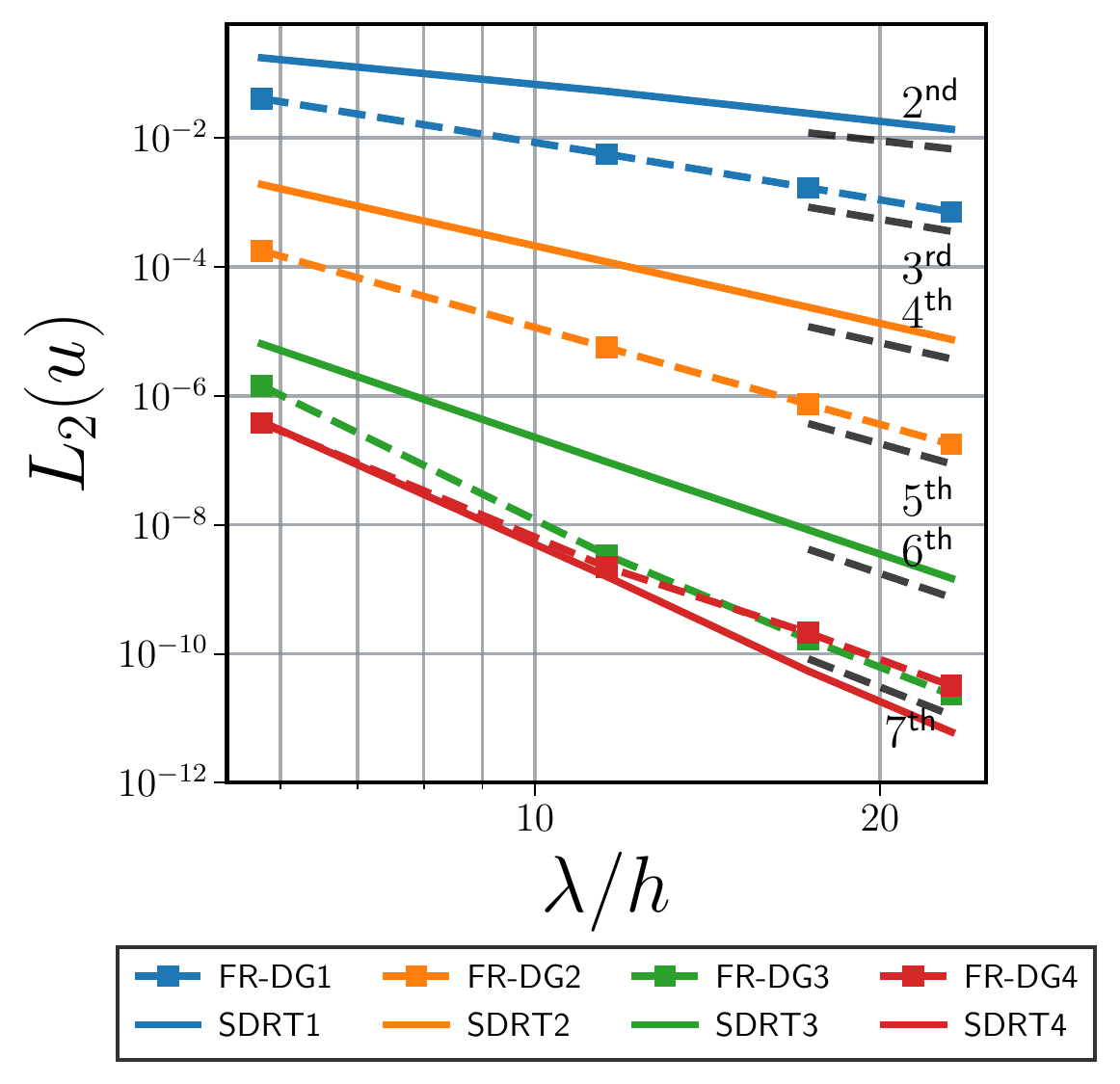}
    \caption{$m = 2$}
    \end{subfigure}
    \caption{$L_2$-norm of the solution error evaluated using \Eref{eq:linadvec_l2_definition} with $m = 1$ (left) and $m = 2$ (right) and obtained using different schemes together with meshes made up of hexahedral elements.}
    \label{fig:l2_advect_ho_hex}
\end{figure}

\begin{figure}
    \begin{subfigure}{0.45\textwidth}
    \centering
    \includegraphics[width=1.0\textwidth]{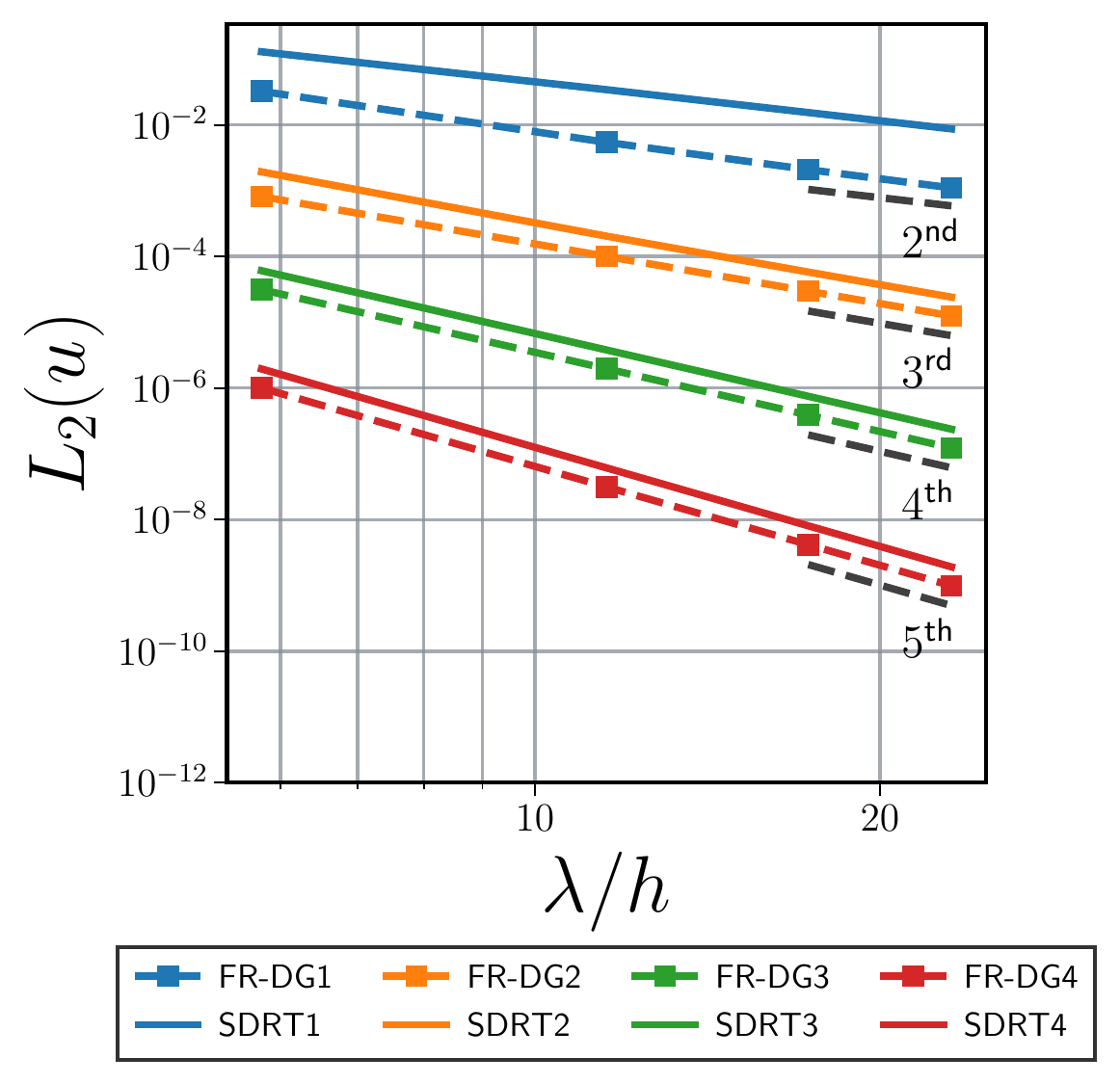}
    \caption{$m = 1$}
    \end{subfigure}
    \begin{subfigure}{0.45\textwidth}
    \centering
    \includegraphics[width=1.0\textwidth]{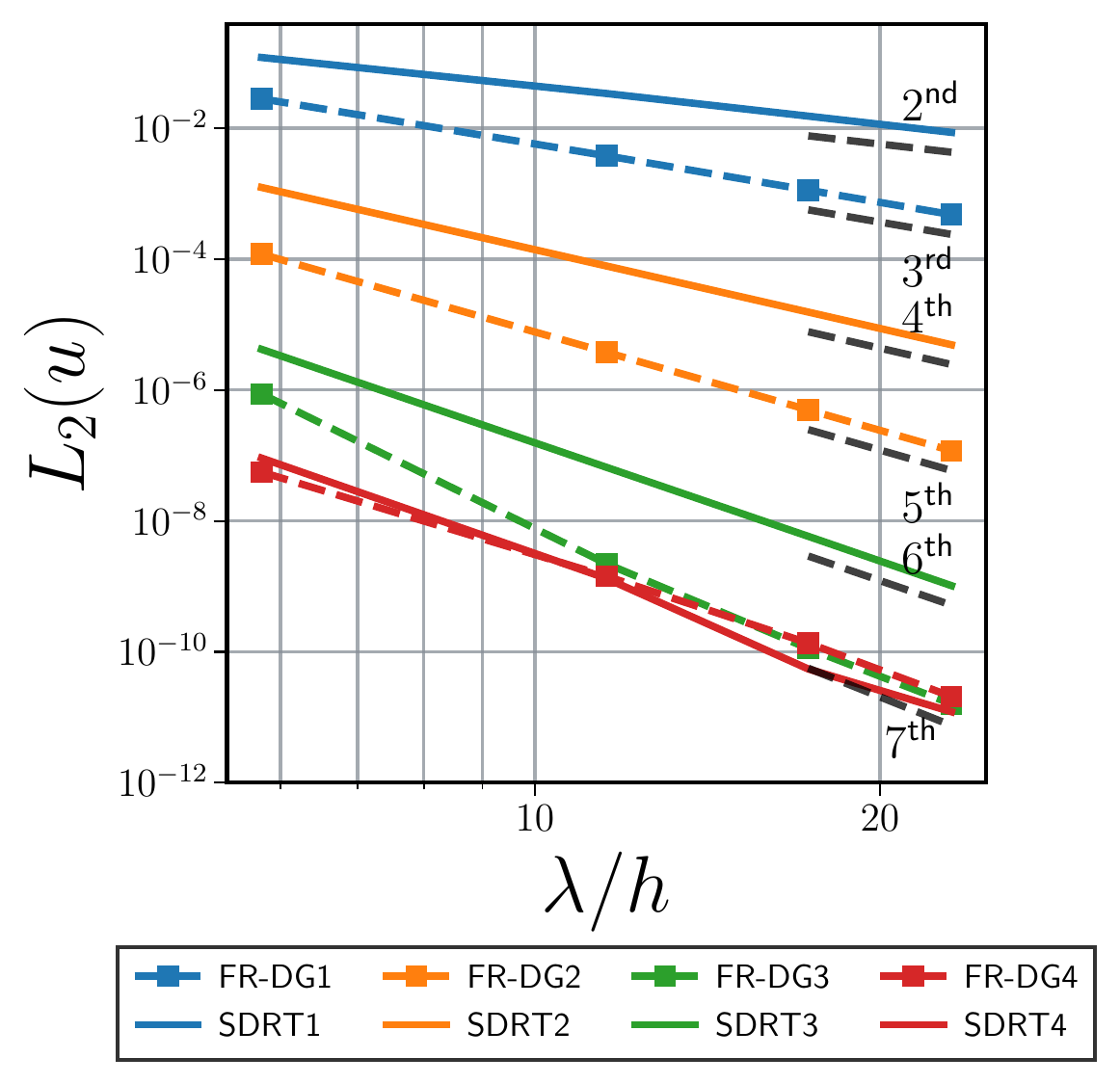}
    \caption{$m = 2$}
    \end{subfigure}
    \caption{$L_2$-norm of the solution error evaluated using \Eref{eq:linadvec_l2_definition} with $m = 1$ (left) and $m = 2$ (right) and obtained using different schemes together with meshes made up of prismatic elements.}
    \label{fig:l2_advect_ho_pri}
\end{figure}

\begin{figure}
    \begin{subfigure}{0.45\textwidth}
    \centering
    \includegraphics[width=1.0\textwidth]{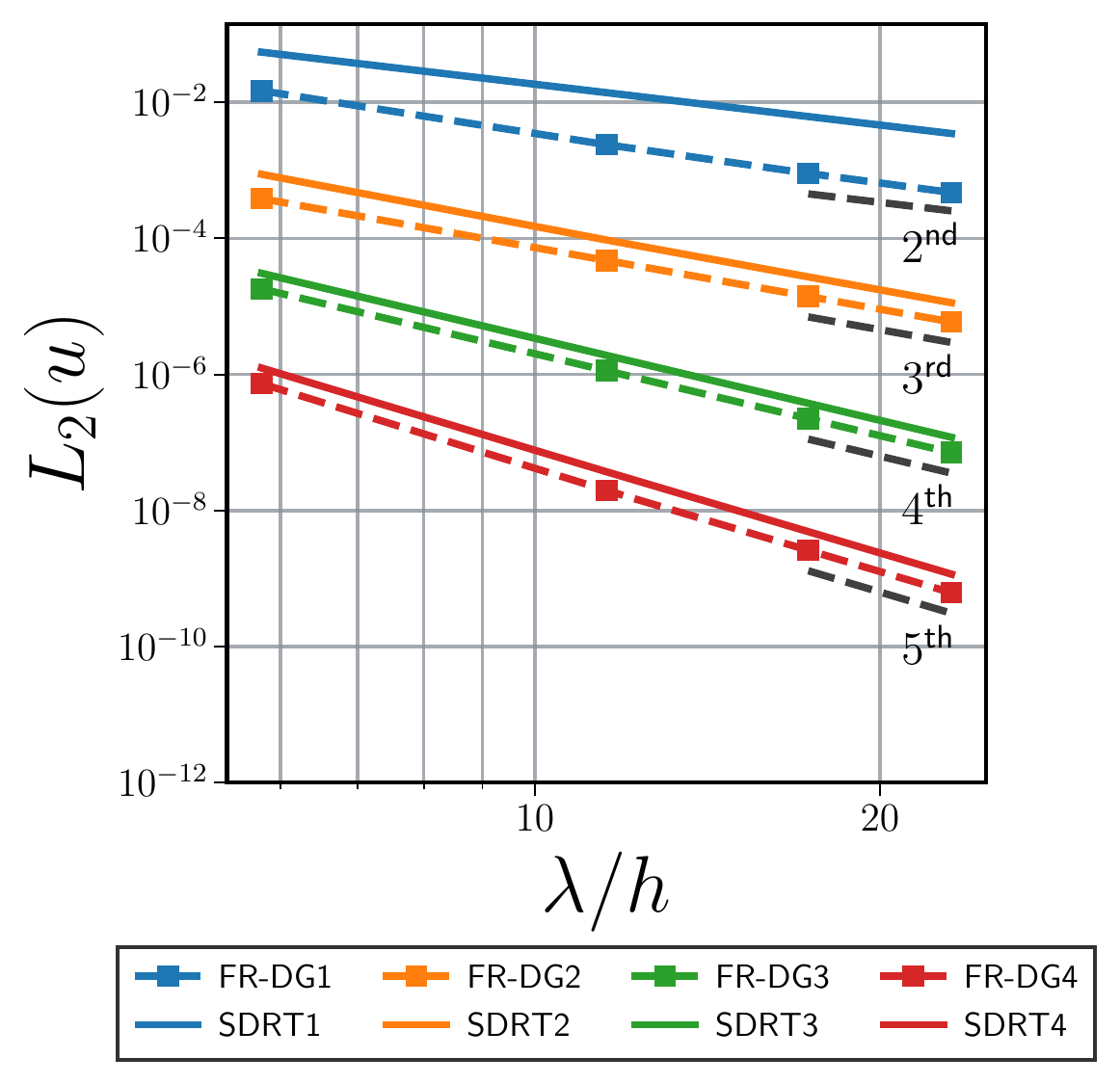}
    \caption{$m = 1$}
    \end{subfigure}
    \begin{subfigure}{0.45\textwidth}
    \centering
    \includegraphics[width=1.0\textwidth]{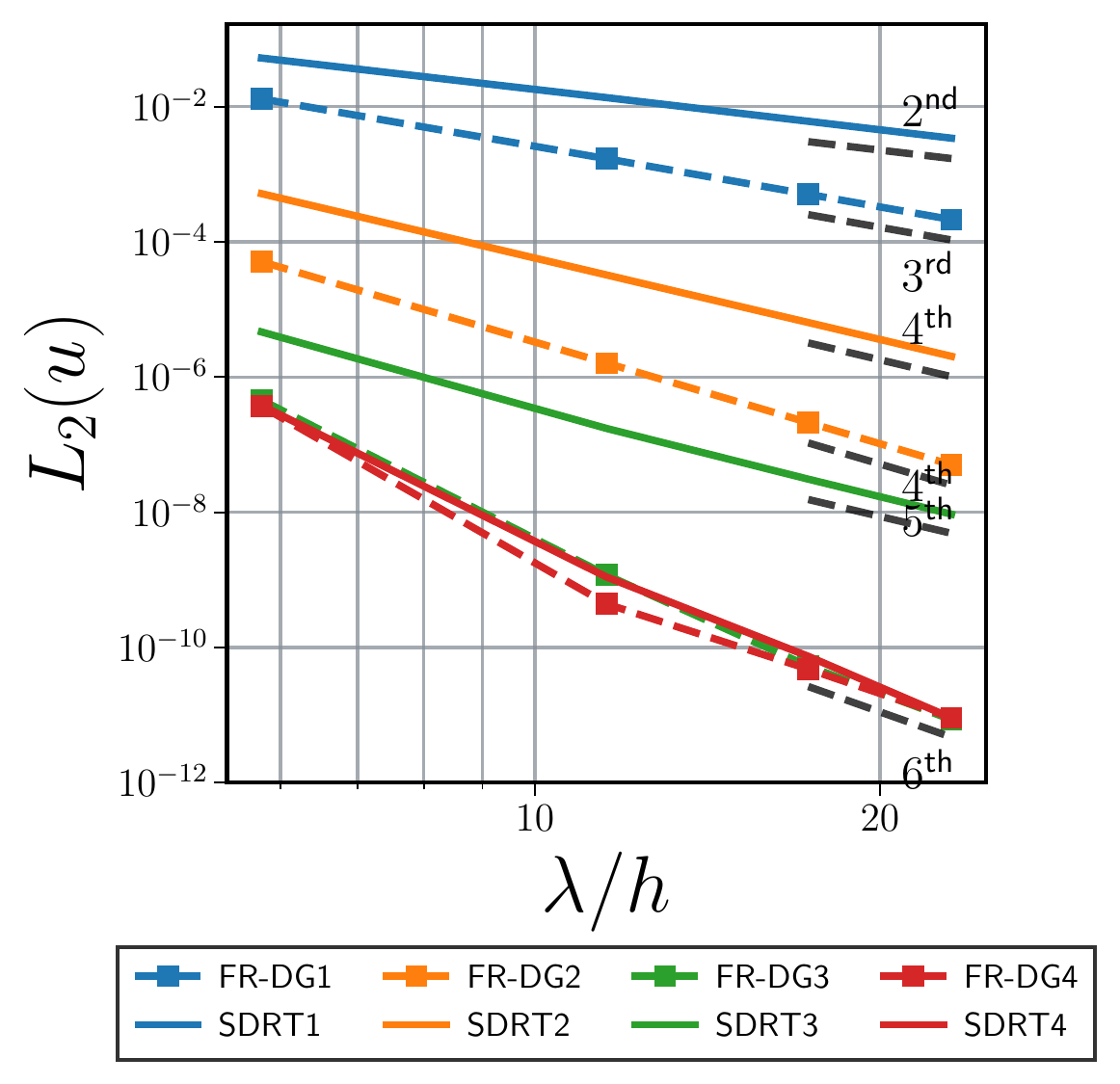}
    \caption{$m = 2$}
    \end{subfigure}
    \caption{$L_2$-norm of the solution error evaluated using \Eref{eq:linadvec_l2_definition} with $m = 1$ (left) and $m = 2$ (right) and obtained using different schemes together with meshes made up of tetrahedral elements.}
    \label{fig:l2_advect_ho_tet}
\end{figure}

For the sake of completeness, Figs. \ref{fig:l2_advect_ho_tri}, \ref{fig:l2_advect_ho_hex}, \ref{fig:l2_advect_ho_pri} and \ref{fig:l2_advect_ho_tet} represent the aforementioned error norm for triangular, hexahedral, prismatic and tetrahedral elements respectively.
The error norm for most elements behaves similarly as that of quadrilateral elements.
Nevertheless, it is worth mentioning the degraded accuracy of error norm computed for $m = 2$ with tetrahedral elements and $\degree = 3$.
Such a behavior was also observed in the analysis of the dispersion errors \Sref{sec:von_neumann_advection} and it may be related to the accuracy and configuration of the different eigenmodes of the SDRT method applied with tetrahedral elements, since this awkward accuracy degradation is not observed for $m = 1$.

\subsubsection{Order of accuracy with $\PE = 10$}

This section aims to validate the order of accuracy of the FR-DG and SDRT schemes utilizing the error norm 
\begin{equation}
    L_2(\cellsize, m) = \sqrt{\frac{1}{|\domain|} \int_{\domain} \bigg[\var^{\numerical}\left(\vect{x}, \frac{2\pi L}{\advecvelnovect_0} \right) - \var\left(\vect{x}, \frac{2\pi L}{\advecvelnovect_0} \right)\bigg]^2 \, \text{d} \domain} ,
    \label{eq:linadvecdiff_l2_definition}
\end{equation}
and Peclet number $\PE = 10$, with $\CFL = 2.5 \cdot 10^{-3}$ (ensuring that both the pure-advection and pure-diffusion $\CFL$ numbers are below the latter threshold).
For the sake of simplicity the results are only presented for $\degree = 3$.
As stated in the previous section, the expected order of accuracy is $\degree + 1 = 4$ due to the initial projection error, despite the fact that, to the best of the authors' knowledge, no published work regarding the asymptotic order of accuracy of pure-diffusion problems discretized with SEM can be found.

\Tref{table:linadvecdiff_l2_sd} displays the error norm obtained with SDRT schemes and all the element types studied in this work, utilizing meshes with a total of $N^{\ndim} \ncellsinpattern$ cells.
As expected, the order of accuracy is $\degree + 1$.
Furthermore, tetrahedral meshes show the lowest error, followed by prismatic and hexahedral meshes.
This can be explained due to the higher number of elements and solution points found in the tetrahedral meshes used in this work.
\Tref{table:linadvecdiff_l2_fr_dg} illustrates the same results obtained with FR-DG methods.
As it was observed for SDRT schemes, the order of accuracy of the schemes is $\degree + 1$, and tetrahedral meshes always present the lowest error norm value.
It is worth mentioning that FR-DG schemes always display lower errors than SDRT schemes for all element types.

\begin{table}[h]
\centering
\begin{tabular}{@{}cccccccccccccccc@{}}
\toprule
  & \multicolumn{2}{c}{quad} & \multicolumn{2}{c}{tri} & \multicolumn{2}{c}{hex} & \multicolumn{2}{c}{pri} & \multicolumn{2}{c}{tet} \\ \cmidrule(lr){2-3} \cmidrule(lr){4-5} \cmidrule(lr){6-7} \cmidrule(lr){8-9} \cmidrule(lr){10-11}  
 & $L_2$ & Order & $L_2$ & Order & $L_2$ & Order & $L_2$ & Order & $L_2$ & Order\\ \midrule 
\multicolumn{1}{l|}{$N = 10$} & 1.32e-05 & - & 5.47e-06 & - & 9.01e-06 & - & 6.00e-06 & - & 2.32e-06 & - \\ 
\multicolumn{1}{l|}{$N = 20$} & 8.96e-07 & 3.88 & 3.04e-07 & 4.17 & 6.14e-07 & 3.87 & 3.92e-07 & 3.93 & 1.33e-07 & 4.12 \\ 
\multicolumn{1}{l|}{$N = 30$} & 1.78e-07 & 3.98 & 5.68e-08 & 4.14 & 1.23e-07 & 3.96 & 7.79e-08 & 3.99 & 2.57e-08 & 4.06 \\ 
\end{tabular}
\caption{$L_2$-norm of the solution error in the linear advection diffusion test case with $\PE = 10$, SDRT schemes and  different element types and meshes. All results were obtained with $\degree = 3$.}
\label{table:linadvecdiff_l2_sd}
\end{table}

\begin{table}[h]
\centering
\begin{tabular}{@{}cccccccccccccccc@{}}
\toprule
  & \multicolumn{2}{c}{quad} & \multicolumn{2}{c}{tri} & \multicolumn{2}{c}{hex} & \multicolumn{2}{c}{pri} & \multicolumn{2}{c}{tet} \\ \cmidrule(lr){2-3} \cmidrule(lr){4-5} \cmidrule(lr){6-7} \cmidrule(lr){8-9} \cmidrule(lr){10-11}  
 & $L_2$ & Order & $L_2$ & Order & $L_2$ & Order & $L_2$ & Order & $L_2$ & Order\\ \midrule 
\multicolumn{1}{l|}{$N = 10$} & 9.76e-06 & - & 3.89e-06 & - & 6.80e-06 & - & 4.50e-06 & - & 1.57e-06 & - \\ 
\multicolumn{1}{l|}{$N = 20$} & 6.80e-07 & 3.84 & 2.42e-07 & 4.01 & 4.78e-07 & 3.83 & 3.08e-07 & 3.87 & 9.80e-08 & 4.00 \\ 
\multicolumn{1}{l|}{$N = 30$} & 1.38e-07 & 3.93 & 4.77e-08 & 4.00 & 9.79e-08 & 3.91 & 6.27e-08 & 3.93 & 1.94e-08 & 4.00 \\ 
\end{tabular}
\caption{$L_2$-norm of the solution error in the linear advection diffusion test case with $\PE = 10$, FR-DG schemes and  different element types and meshes. All results were obtained with $\degree = 3$.}
\label{table:linadvecdiff_l2_fr_dg}
\end{table}

\subsubsection{Equivalence between FR-SDRT and SDRT methods}
\label{sec:numerical_experiments_linadvecdiff_equivalence_sdrt_fr}

This section is devoted to the a posteriori assessment of the equivalence between the FR-SDRT and SDRT methods in two and three-dimensional elements.
This allows to demonstrate the theoretical observations that were drawn in \Sref{sec:sdrt_fr_equivalence}.
To do so, the differences between simulations conducted with $\mathrm{Pe} = 10$, $N = 30$ and FR-SDRT or SDRT schemes will be compared using a $L_{\infty}$ norm evaluated at $t = 2 \pi / \advecvelnovect_0$, i.e. the biggest absolute difference of the numerical solution value at the solution points between FR-SDRT and SDRT simulations.
For the sake of brevity, the results will only be presented for $\degree = 3$ and $\degree = 4$, although the equivalence has been validated a posteriori for all degrees $\degree \in [1, 4]$.
\Tref{table:linadvecdiff_linfty_diff_sdrt_fr_sdrt} and \Tref{table:linadvecdiff_linfty_diff_sdrt_fr_sdrt_p4} illustrate the biggest absolute difference of the solution points value between FR-SDRT and SDRT simulations for different element types and $\degree = 3$ and $\degree = 4$ respectively.
Since this difference is of the same order of magnitude as machine round-off errors for all elements, the equivalence between FR-SDRT and SDRT schemes stated in \Sref{sec:sdrt_fr_equivalence} is also demonstrated through numerical experiments.

\begin{table}[h]
\centering
\begin{tabular}{@{}cccccc@{}}
\toprule
  & \multicolumn{1}{c}{quad} & \multicolumn{1}{c}{tri} & \multicolumn{1}{c}{hex} & \multicolumn{1}{c}{pri} & \multicolumn{1}{c}{tet}\\ \midrule 
$L_\infty$ & 6.11e-16 & 8.26e-16 & 5.93e-16 & 1.25e-15 & 3.30e-16
\end{tabular}
\caption{Absolute difference of the solution points value between FR-SDRT and SDRT simulations for different element types and $\degree = 3$ in the simulation of the linear advection diffusion test case with $\PE = 10$.}
\label{table:linadvecdiff_linfty_diff_sdrt_fr_sdrt}
\end{table}

\begin{table}[h]
\centering
\begin{tabular}{@{}cccccc@{}}
\toprule
  & \multicolumn{1}{c}{quad} & \multicolumn{1}{c}{tri} & \multicolumn{1}{c}{hex} & \multicolumn{1}{c}{pri} & \multicolumn{1}{c}{tet}\\ \midrule 
$L_\infty$ & 6.66e-16 & 1.50e-15 & 3.11e-16 & 2.27e-15 & 5.59e-16
\end{tabular}
\caption{Absolute difference of the solution points value between FR-SDRT and SDRT simulations for different element types and $\degree = 4$ in the simulation of the linear advection diffusion test case with $\PE = 10$.}
\label{table:linadvecdiff_linfty_diff_sdrt_fr_sdrt_p4}
\end{table}

\subsection{Isentropic Euler Vortex}

This section is devoted to the analysis of polynomial aliasing errors which appear in non-linear test cases.
In particular, the Isentropic Euler Vortex \cite{Spiegel2015} will be analyzed following the considerations of \citet{Spiegel2015_aliasing, Cox2021}.
The Isentropic Euler vortex problem is commonly used to test the order of accuracy of numerical methods for conservation laws with non-linear fluxes, in particular the Euler equations.
These equations may be represented as a system of conservation laws that describes the dynamics of inviscid fluids.
The conservative variables of this set of equations are
\begin{equation}
    \var = 
    \begin{pmatrix}
        \rho \\
        \rho \velocity\\
        \rho E
    \end{pmatrix} .
    \label{eq:eulerConservativeVariables}
\end{equation}
where $\density$ is the fluid, $\rho \velocity$ is the fluid momentum, $\velocity$ is the fluid velocity and the total energy is given by $\rho E$ (being $E$ the total energy per unit of mass). 
The flux operator of the Euler equations presents only convective terms and can be written as
\begin{equation}
    \flux(\var) = 
    \begin{pmatrix}
        \rho \velocity \\
        \rho \velocity \otimes \velocity + \pressure \vect{I}\\
        \rho \velocity E + \pressure\velocity 
    \end{pmatrix} .
    \label{eq:fluxOperatorEuler}
\end{equation}
Here, the symbol $\otimes$ represents the dyadic operator, $\pressure$ is the pressure and $\vect{I} \in \mathbb{R}^{\ndim \times \ndim}$ is the identity matrix.
The non-linear system is closed with the equation of state
\begin{equation}
    \rho E = \frac{\pressure}{\gamma - 1} + \frac{1}{2}\rho \velocity \cdot \velocity \quad \text{ and/or } \quad \pressure = \density r \temperature \quad .
\end{equation}
Here, $\gamma$ denotes the adiabatic constant, $r$ is the perfect gas constant and $\temperature$ is the fluid temperature.

In this section, the Isentropic Euler Vortex configuration used in \cite{Witherden2014} is replicated and special emphasis will be placed on the study of aliasing properties of FR-DG, SDRT and FR-SDRT.
The analytical solution of this test case is
\begin{equations}{eq:eulerVortexWitherden}
    \rho(\vect{x}, t) &= \density_{\infty}\left( 1 - \frac{S^2 \MA^2 (\gamma - 1) \e^{2f}}{8 \pi^2} \right)^{1/(\gamma - 1)} \\
    \pressure(\vect{x}, t) &= \pressure_{\infty} \left( \frac{\rho}{\rho_{\infty}} \right)^{\gamma} \\ 
    \velocity(\vect{x}, t) &= \velocitynovect_{\infty} \left(\frac{S \hat{y} \e^{f}}{2 \pi R}, 1 - \frac{S \hat{x} \e^{f}}{2 \pi R}   \right) .
\end{equations}
where $\vect{\hat{x}} = (\hat{x}, \hat{y}) = (x, y - \velocitynovect_{\infty} t)$, $f = (1 - \hat{x}^2 - \hat{y}^2)/(2R^2)$, $S=13.5$ is the vortex strength, $\MA=0.4$ is the free-stream Mach number and $R=1.5$ is the vortex radius.
The free-stream values with $\infty$ sub-index are set as: $\density_{\infty} = 1$, $\velocitynovect_{\infty} = 1$ and $\pressure_{\infty} = \density_{\infty}\velocitynovect_{\infty}^2 / ( \gamma \MA_{\infty}^2)$ to match the conditions imposed in \cite{Witherden2014}.
Simulations are carried out in a domain $\domain \in [-L, L]^2$ with $R/L = 0.075$.
Such a choice of domain size ensures that the velocity and density perturbations are well below the machine round-off errors.
Periodic boundary conditions are imposed at boundaries with constant $y$ coordinate while the limiting values of the initial conditions (i.e. those obtained with $\vect{x} \to \infty$) were used to define the boundary conditions at boundaries with constant $x$ coordinate.
This allows to further reduce numerical errors and instabilities which are introduced when solving the Isentropic Euler Vortex with periodic boundary conditions \cite{Spiegel2015}.
The common fluxes are computed using the Rusanov--Riemann solver.
To measure the order of accuracy the $L_2$ norm of the density error, defined as
\begin{equation}
    L_2(\cellsize, m) = \sqrt{\frac{1}{|\hat{\volume}|} \int_{\hat{\volume}} \bigg[\density^{\numerical}\left(\vect{x}, m t_c \right) - \density\left(\vect{x}, t=0 \right)\bigg]^2 \, \text{d} \volume} ,
    \label{eq:euler_vortex_l2_definition}
\end{equation}
is utilized.
In the latter equation, $t_c = L$ is the characteristic convective time of the vortex, $m$ is an arbitrary positive integer, $\hat{\volume}$ is a the volume of the mesh region comprised within $\vect{x} \in [-L/10, L/10]^2$.
Such a choice of domain used to compute the error norm is considered to avoid taking into account spurious oscillations arising due to the use of periodic boundary conditions that are not representative of the analytical solution \cite{Witherden2014}.
The time step of the simulations is chosen as $\timestep = 1.25 \cdot 10^{-3}$ to ensure that the spatial discretization errors are predominant over temporal discretization errors.

\begin{remark}
As it has been shown in \cite{Spiegel2015, Witherden2014, Cox2021} the order of accuracy obtained from these simulations is heavily dependent on the parameter $m$ and meshes used to evaluate the error norm.
In particular, it was observed that the higher the value of $m$ the higher the order of accuracy, provided that the time step is sufficiently low.
Nevertheless, in the authors' opinion the reason behind this superconvergence behavior of the error norm is not fully determined.
Studies usually explain this observation using the propagation of projection error of the initial condition \cite{Witherden2014}, similar to what it is observed for linear cases \cite{Guo_2013}.
However, it is not clear whether this linear theory can be applied to non-linear problems.
\end{remark}

\begin{figure}
    \begin{subfigure}{0.45\textwidth}
    \centering
    \includegraphics[width=1.0\textwidth]{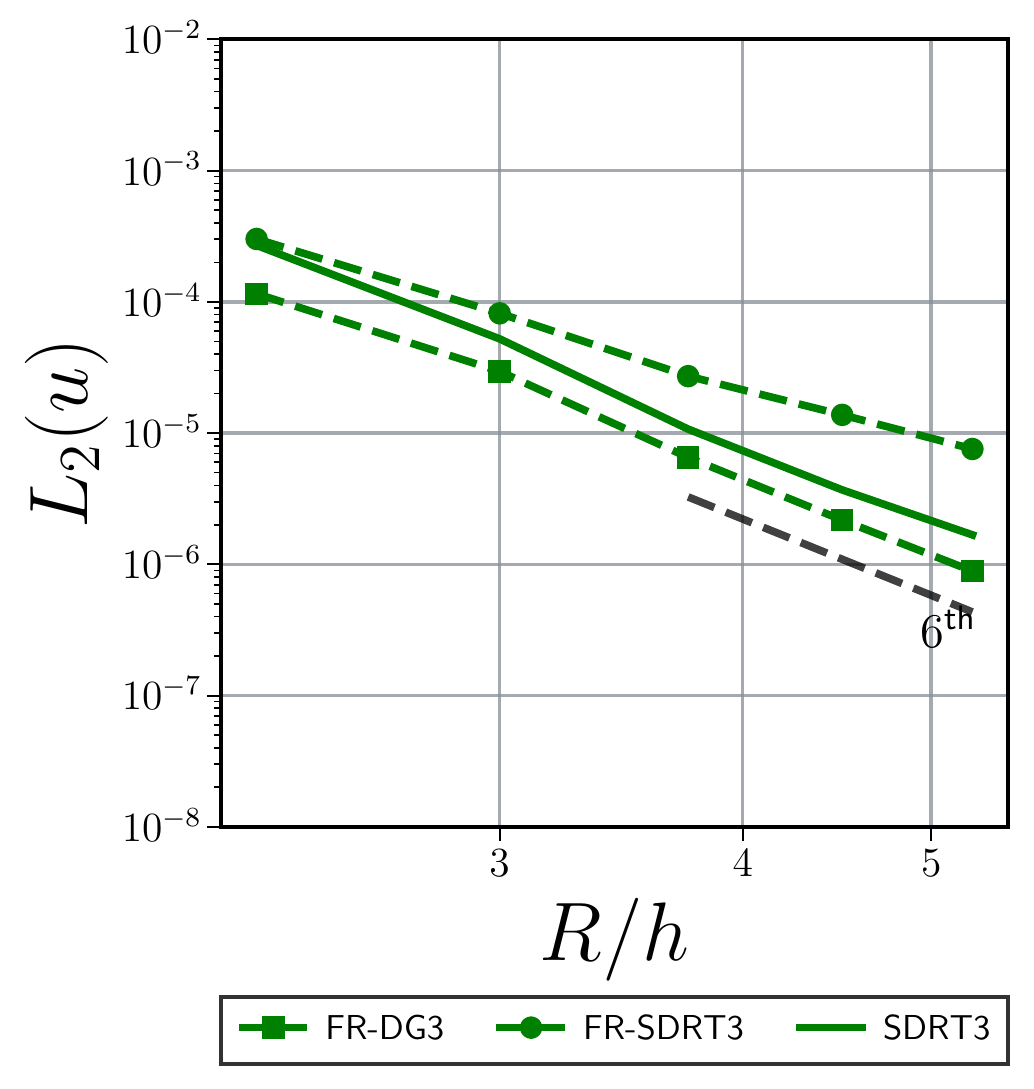}
    \caption{Quadrilateral}
    \end{subfigure}
    \begin{subfigure}{0.45\textwidth}
    \centering
    \includegraphics[width=1.0\textwidth]{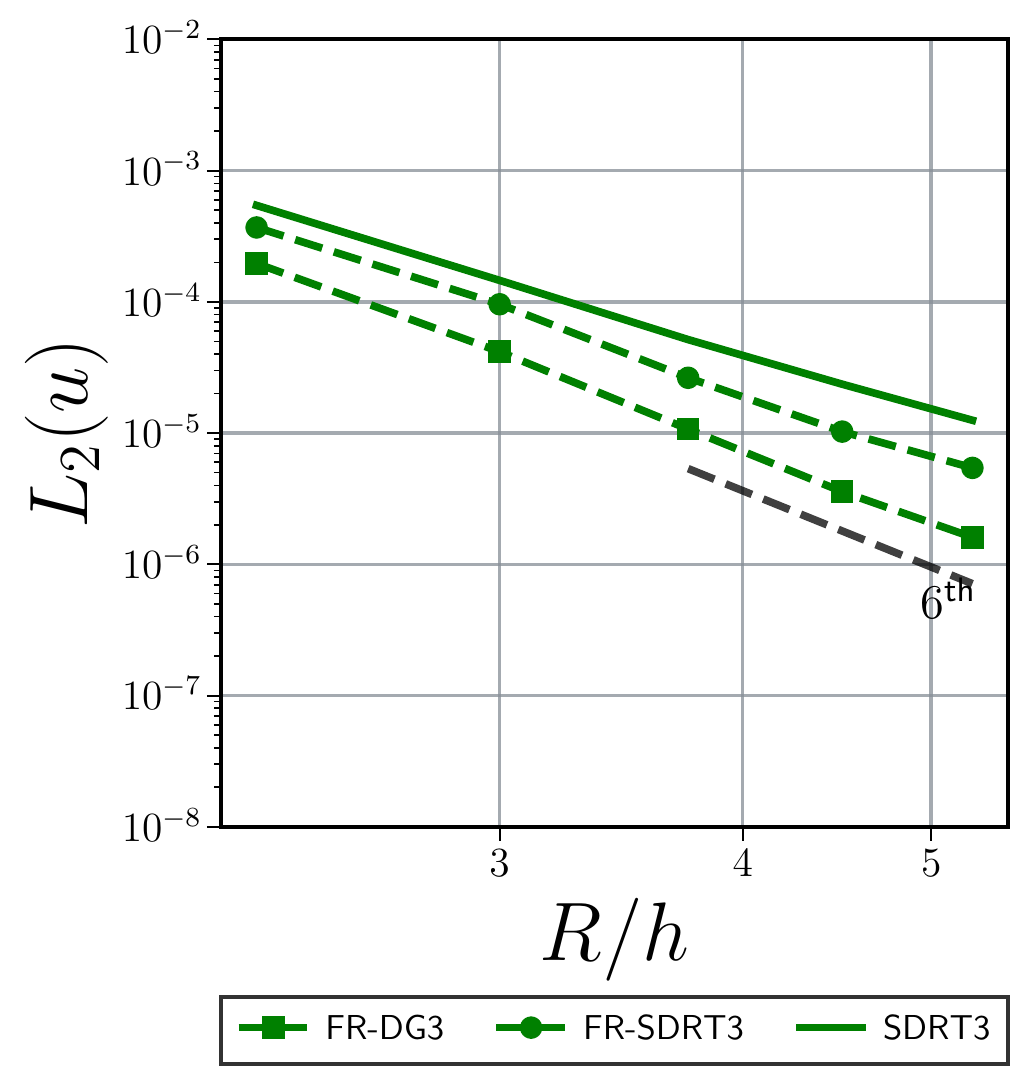}
    \caption{Triangular}
    \end{subfigure}
    \caption{$L_2$-norm of the solution error in the Isentropic Euler Vortex evaluated using \Eref{eq:euler_vortex_l2_definition} with $m = 2$, $\degree = 3$ and obtained using different schemes together with meshes made up of quadrilateral (left) and triangular (right) elements.}
    \label{fig:l2_euler_t80}
\end{figure}

\Fref{fig:l2_euler_t80} represents the $L_2$-norm of the solution error evaluated with $m=2$, $\degree = 3$ and different schemes for quadrilateral and triangular elements.
The results for tensor-product elements show similar errors as those found in \cite{Cox2021}, i.e. FR-DG schemes show the lowest error values followed by SDRT and FR-SDRT schemes.
On the other hand, the results obtained with triangles display that SDRT schemes show the highest error, even higher than FR-SDRT.
The reasoning behind such disparities is not straightforward.
For example, we demonstrated in \Sref{sec:numerical_experiments_linadvecdiff} that FR-SDRT and SDRT schemes are equivalent in linear problems.
However, such equivalence is lost in non-linear problems and hence the differences between simulations performed with FR-SDRT and SDRT may only be explained using non-linear machinery.
In theory, SDRT schemes should present better aliasing properties than FR schemes due to the use of a staggered-grid approach and Raviart-Thomas flux bases to project the fluxes, hence its numerical error should be lower than that of FR-SDRT.
Nevertheless, such an advantage is only observed in this test case when using tensor-product elements and when comparing FR-SDRT and SDRT schemes.
Additionally, FR-DG methods always provide the lowest error values in this test case.
It is even possible that SDRT schemes are non-linearly weakly unstable.
These issues will be studied in future works since, in the authors' opinion, the explanation of such inconsistent conclusions requires tools beyond the linear analysis.

\subsection{Taylor-Green-Vortex}
\label{sec:tgv}

This section analyzes the accuracy of the SDRT and FR-DG methods using the Taylor-Green-Vortex \cite{Wang2013}, which is a standard validation test case used to assess the accuracy of numerical schemes with turbulent-like flows using the Navier-Stokes equation.
The Navier-Stokes equations describe the motion of viscous fluids and share the same convective flux than the Euler equations.
Additionally, they present a viscous flux term which is formulated as
\begin{equation}
    \flux(\var) = 
    \begin{pmatrix}
        \rho \velocity \\
        \rho \velocity \otimes \velocity + \pressure \vect{I}\\
        (\rho E + \pressure)\velocity 
    \end{pmatrix} + \begin{pmatrix}
        0 \\
        - \stressTensor  \\
        - \frac{\mu c_p}{\PR} \nabla \temperature - \stressTensor \cdot \velocity  
    \end{pmatrix} ,
    \label{eq:fluxOperatorNavierStokes}
\end{equation}
where $\mu$ is the kinematic viscosity, $c_p$ is the specific heat of the fluid for constant pressure, $\PR$ is the Prandtl number and $\stressTensor$ is the viscous stress tensor which may be written for Newtonian fluids and under the Stokes hypothesis as
\begin{equation}
    \stressTensor = \mu \left( \grad \velocity + \big( \grad \velocity \big)^T - \frac{2}{3} \grad \cdot \velocity \vect{I} \right) .
\end{equation}

The initial condition of the TGV reads
\begin{equations}{eq:tgvInitCondition}
    \velocitynovect_0 &= \velocitynovect_{\infty} \sin{\frac{x}{L}} \cos{\frac{y}{L}} \cos{\frac{z}{L}} \\
    \velocitynovect_1 &= -\velocitynovect_{\infty} \cos{\frac{x}{L}} \sin{\frac{y}{L}} \cos{\frac{z}{L}} \\
    \velocitynovect_2 &= 0 \\
    \pressure &= \pressure_{\infty} + \frac{\density_{\infty} \velocitynovect_{\infty}^2}{16}\left[\cos{\frac{2x}{L}} + \cos{\frac{2y}{L}}  \right] \left[\cos{\frac{2z}{L}} + 2 \right] \\
    \density &= \frac{\pressure}{r \temperature_{\infty}} .
\end{equations}
In the latter, $\temperature_{\infty}$ is the initial temperature field supposed constant and $\velocitynovect_{\infty}$, $\pressure_{\infty}$ and $\density_{\infty}$ allow to define the non-dimensional parameters Reynolds number $\RE_L$ and Mach number $\MA_{\infty}$ as
\begin{equation}
    \RE_{L} = \frac{\density_{\infty} \velocitynovect_{\infty} L}{\mu} = 1600 \quad , \quad
    \MA_{\infty} = \frac{\velocitynovect_{\infty}}{\sqrt{\gamma \frac{\pressure_{\infty}}{\rho_{\infty}}}} = 10^{-1} \quad \quad .
\end{equation}
Such a value of the Mach number is imposed to obtain results with compressible solvers close to those obtained with incompressible formulations.
Moreover, the flow is supposed to present a constant Prandtl number $\PR = 0.71$.
The simulations are performed in a periodic domain $\Omega \in [-\pi L, \pi L]^3$.

Such a flow configuration experiences a transition to a weakly turbulent state, with the creation of small scales, followed by a decay phase similar to decaying homogeneous turbulence, yet not isotropic according to \cite{Wang2013}.
DNS data from this case is often available in the incompressible limit of the Navier--Stokes equations and for a wide range of Reynolds number (see \cite{vanRees2011} for example).

The comparison of simulations with reference DNS data is usually carried out through the computation of the ensemble average $\langle \square \rangle$ over all the domain of a certain quantity.
This operator is defined as
\begin{equation}
    \langle \square \rangle = \frac{1}{|\Omega|} \int_\Omega \square \text{d}{\Omega} \quad .
\end{equation}
In this work the non-dimensional dissipation of the non-dimensional ensemble average compressible kinetic energy $E^* = \frac{1}{\rho_{\infty} \velocitynovect_{\infty}^2} \rho \velocity \cdot \velocity$ is defined as
\begin{equation}
    \epsilon_1 = -\frac{\text{d} \ensembleAverage{E^*}}{\text{d} t^*} \quad .
\end{equation}
Here, $t^* = t/t_c$ is the non-dimensional time variable and $t_c = L/\velocitynovect_{\infty}$ is the characteristic time.
The analysis of the equation describing the temporal evolution of the kinetic energy allows to estimate numerical errors by the comparison of each of components of the kinetic energy balance equation \cite{Navah2020} which are defined as follows
%
\begin{equations}{eq:dissipationsTaylorGreenVortex4}
    \epsilon_2 &= \frac{2 \mu t_c}{\rho_\infty \velocitynovect_{\infty}^2} \ensembleAverage{\mathbb{S}^d : \mathbb{S}^d} \\
    \epsilon_3 &= -\frac{t_c}{\rho_\infty \velocitynovect_{\infty}^2} \ensembleAverage{\pressure \nabla \cdot \velocity} \\
    \epsilon_1 &= \epsilon_2 + \epsilon_3,
\end{equations}
In the latter, $\mathbb{S}^d = \frac{1}{2}\left(\nabla \velocity + \nabla \velocity^T \right) - \frac{1}{3} \grad \cdot \velocity \vect{I}$ is the deviatoric strain-rate tensor.
Additionally, $\epsilon_2$ is the strain rate dissipation and $\epsilon_3$ is the pressure dissipation.

\begin{remark}
Under the incompressibility hypothesis $\epsilon_3 = 0$, $\epsilon_2$ becomes the classical enstrophy dissipation term \cite{DeBonis2013}.
\end{remark}

Since, all the dissipation terms can be numerically computed from simulation data, deviations of $\epsilon_2 + \epsilon_3$  from $\epsilon_1$ can be linked to numerical dissipation provided that the terms are computed with appropriate numerical quadratures.
In this work, a quadrature of degree 10 is used to assess the ensemble averages.
As in \cite{Navah2020, Cox2021} the error estimator $\epsilon^\numerical$ may be expressed as
\begin{equation}
    \epsilon^\numerical = \epsilon_1 - \epsilon_2 - \epsilon_3 .
\end{equation}

All simulations will be performed with $\beta = 0.5$ and $\tau = 0.1$ which are common choices for simulations of the TGV test case (see \cite{Vermeire2017}) and also add a certain amount of numerical dissipation to the viscous terms which improves the stability of under-resolved simulations.
Furthermore, the common convective fluxes are computed using the Rusanov--Riemann solver.
Additionally, the RK54 time integrator coupled with a PI adaptive controller \cite{Butcher_2016} are utilized to carry out the simulations using adaptive time-stepping.
The latter adaptive time-stepping controller is configured with relative and absolute tolerances equal to $10^{-8}$ and using the $L_{\infty}$ norm as error estimator.

\subsubsection{Validation}

\begin{figure}
    \begin{subfigure}{0.32\textwidth}
    \centering
    \includegraphics[width=0.98\textwidth]{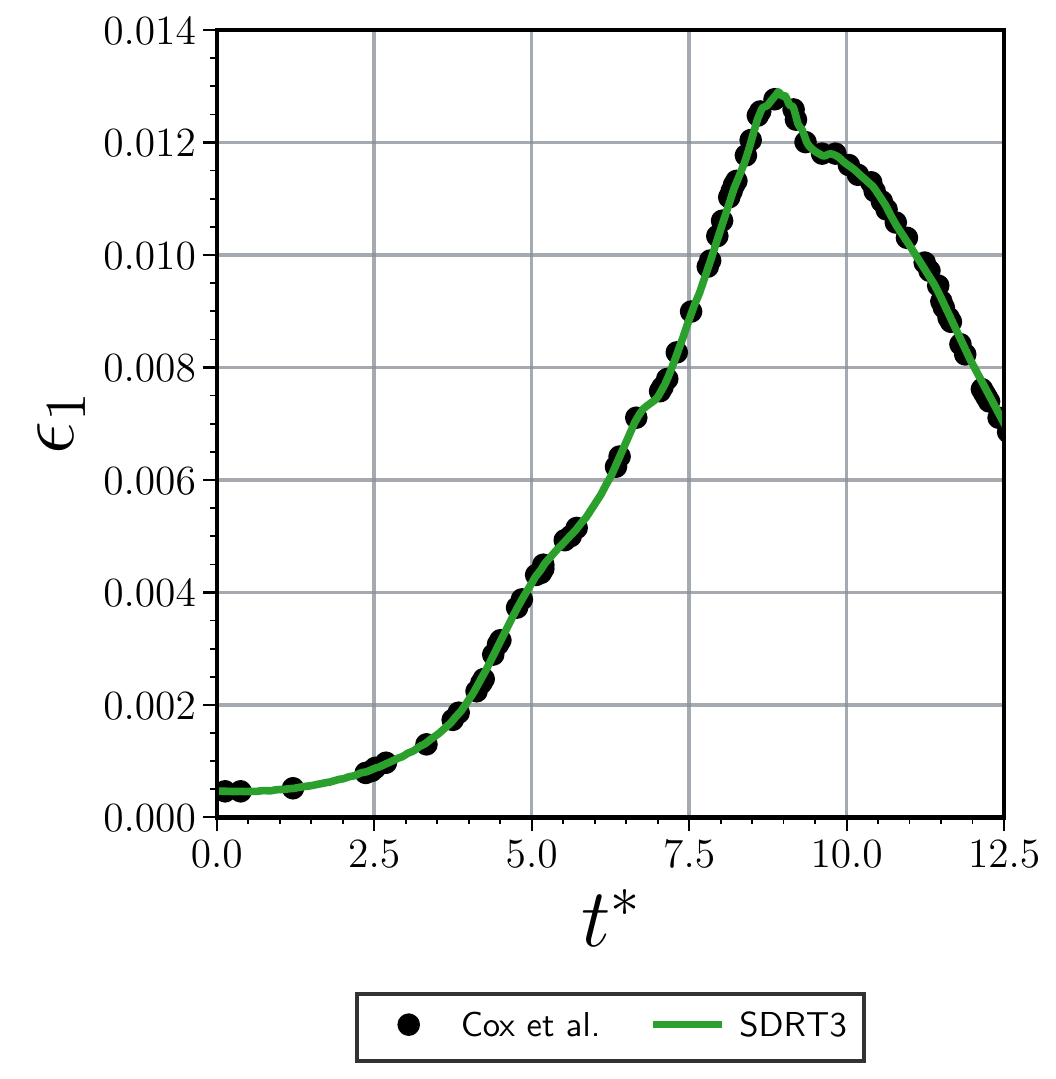}
    \caption{Kinetic energy dissipation}
    \end{subfigure}
    \begin{subfigure}{0.32\textwidth}
    \centering
    \includegraphics[width=0.98\textwidth]{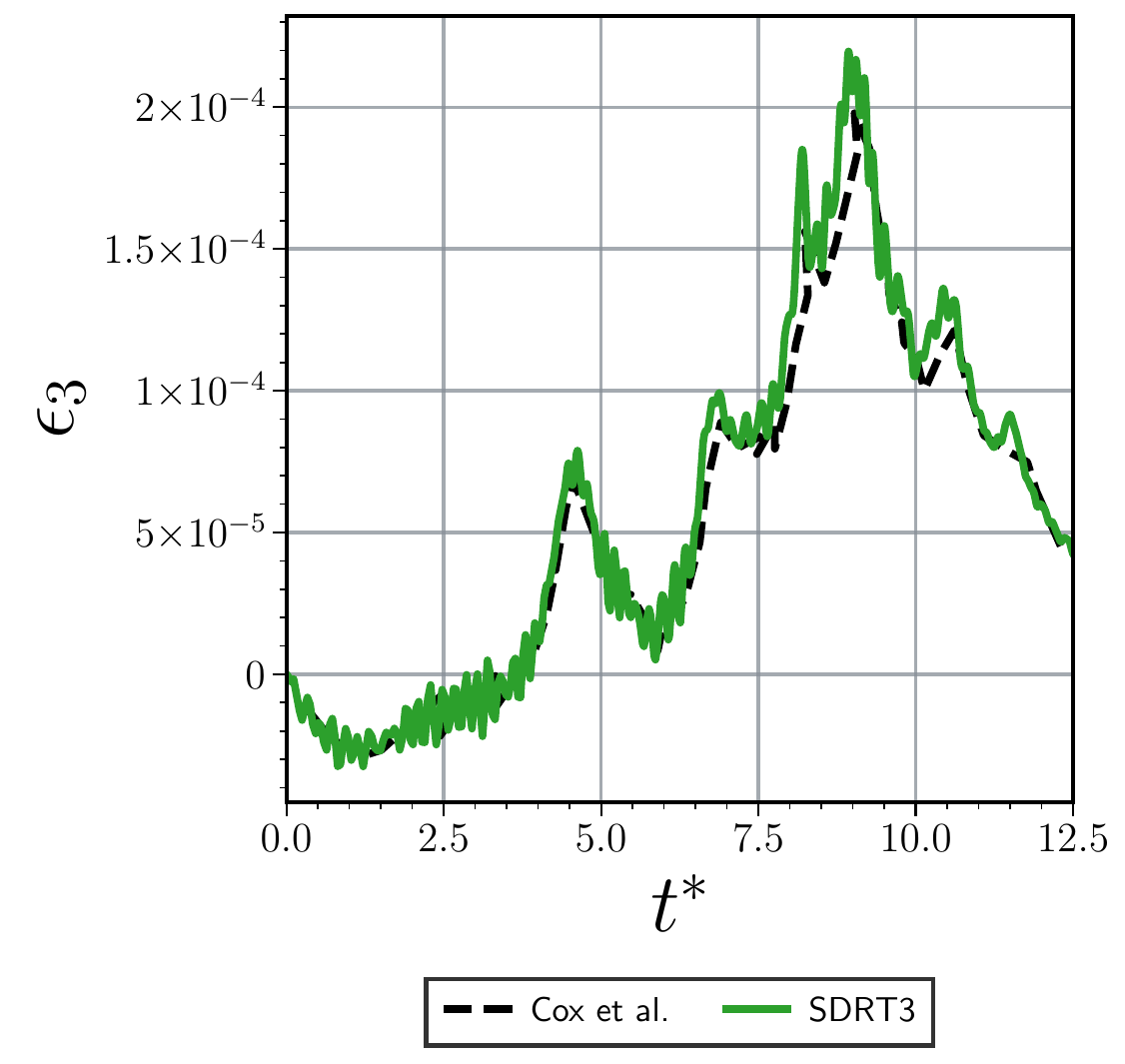}
    \caption{Pressure dissipation}
    \end{subfigure}
    \begin{subfigure}{0.32\textwidth}
    \centering
    \includegraphics[width=0.98\textwidth]{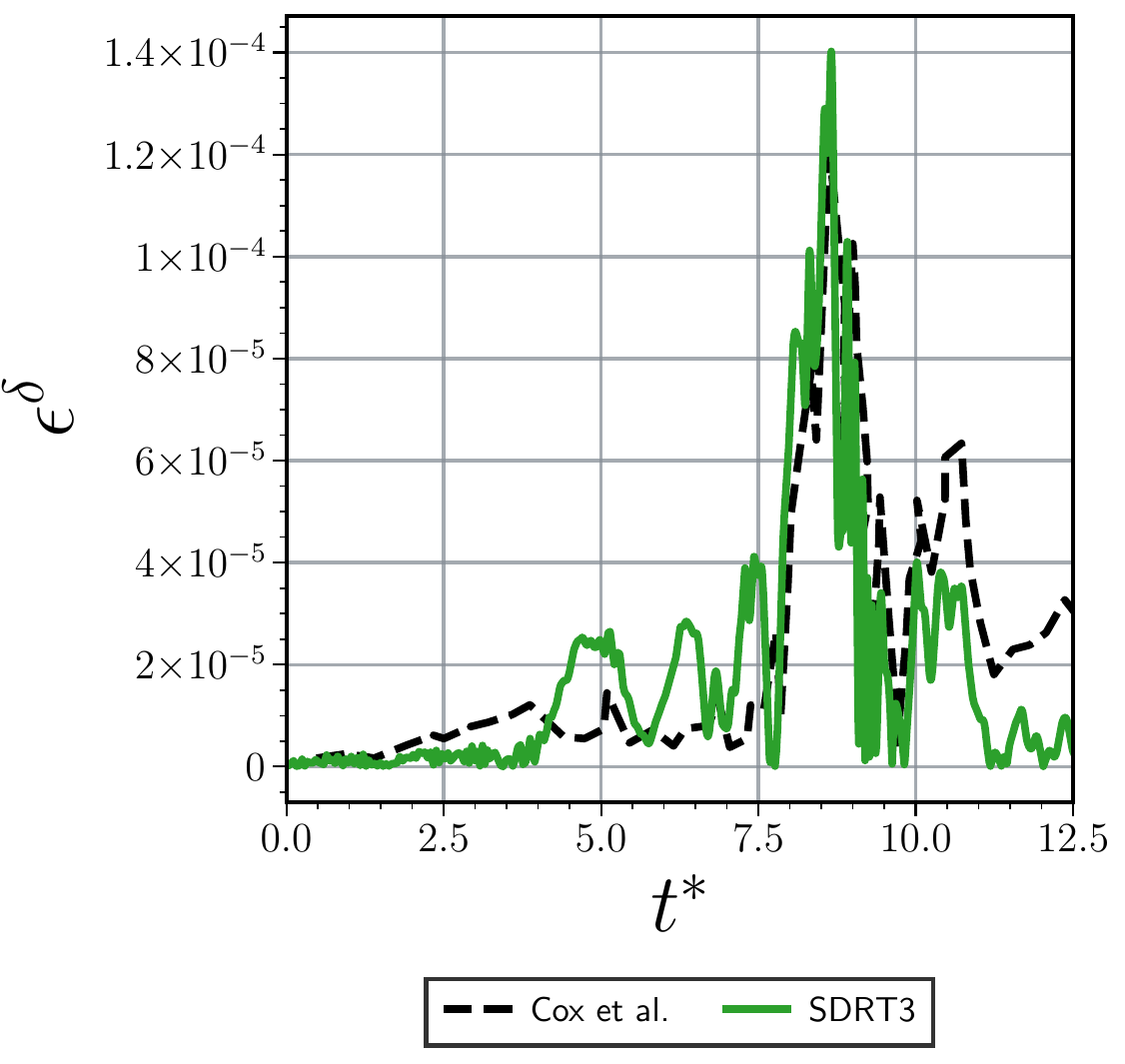}
    \caption{Numerical dissipation}
    \end{subfigure}
    \caption{Temporal evolution of different dissipation terms in the TGV test case obtained with $64^3$ hexahedral elements and the SDRT3 scheme.
    Reference data from \cite{Cox2021} obtained with an in-house SD solver and $\degree = 3$ are used for validation purposes.}
    \label{fig:tgv_sdrt3_256dof_validation}
\end{figure}

To validate the SDRT implementation with hexahedral elements, the TGV test case is simulated using a mesh containing $64^3$ hexahedral elements and $\degree = 3$.
The results are compared with data from \cite{Cox2021} which were obtained with their in-house SD solver implementation.
\Fref{fig:tgv_sdrt3_256dof_validation} represents the different dissipation terms as a function of the non-dimensional time.
The results are very close to each other validating the SDRT implementation in PyFR with hexahedral elements.
Differences found in the numerical dissipation estimator could be related to the use of a different quadrature, the  values of $\beta$ and $\tau$, variations in the adaptive time-stepping method of choice, post processing issues, etc.

\begin{remark}
It is worth mentioning that the computation of the $\epsilon_3$ parameter is highly dependent on the method used to post-process the velocity gradients.
In particular, it was observed that if the parameter $\epsilon_3$ is not assessed using the LDG approach, then its temporal evolution displayed much higher values than those from reference data.
\end{remark}

\subsubsection{Under-resolved configuration}

\begin{figure}
    \begin{subfigure}{0.32\textwidth}
    \centering
    \includegraphics[width=0.95\textwidth]{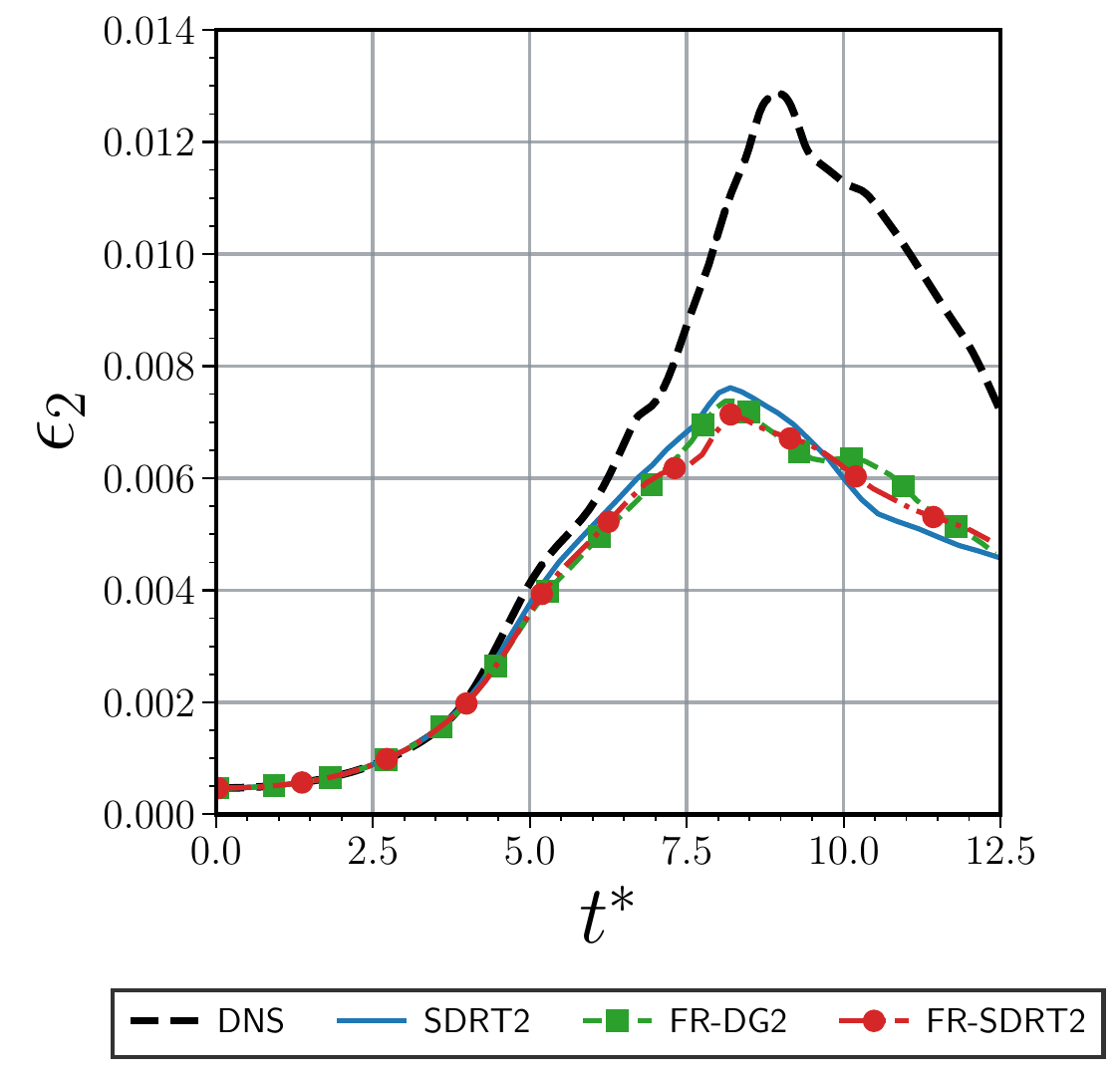}
    \caption{Hexahedral $\degree = 2$}
    \end{subfigure}
    \begin{subfigure}{0.32\textwidth}
    \centering
    \includegraphics[width=0.95\textwidth]{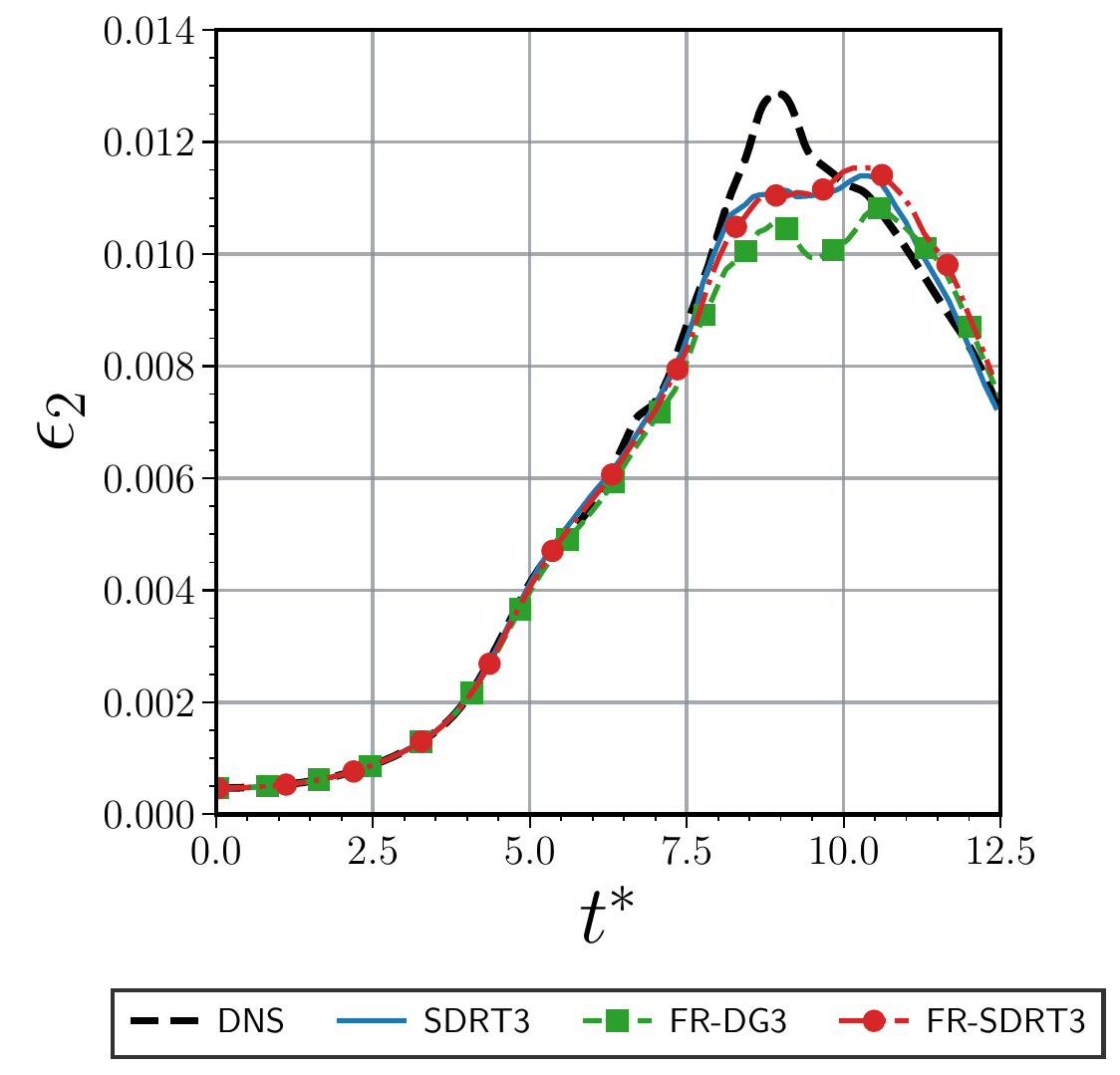}
    \caption{Hexahedral $\degree = 3$}
    \end{subfigure}
    \begin{subfigure}{0.32\textwidth}
    \centering
    \includegraphics[width=0.95\textwidth]{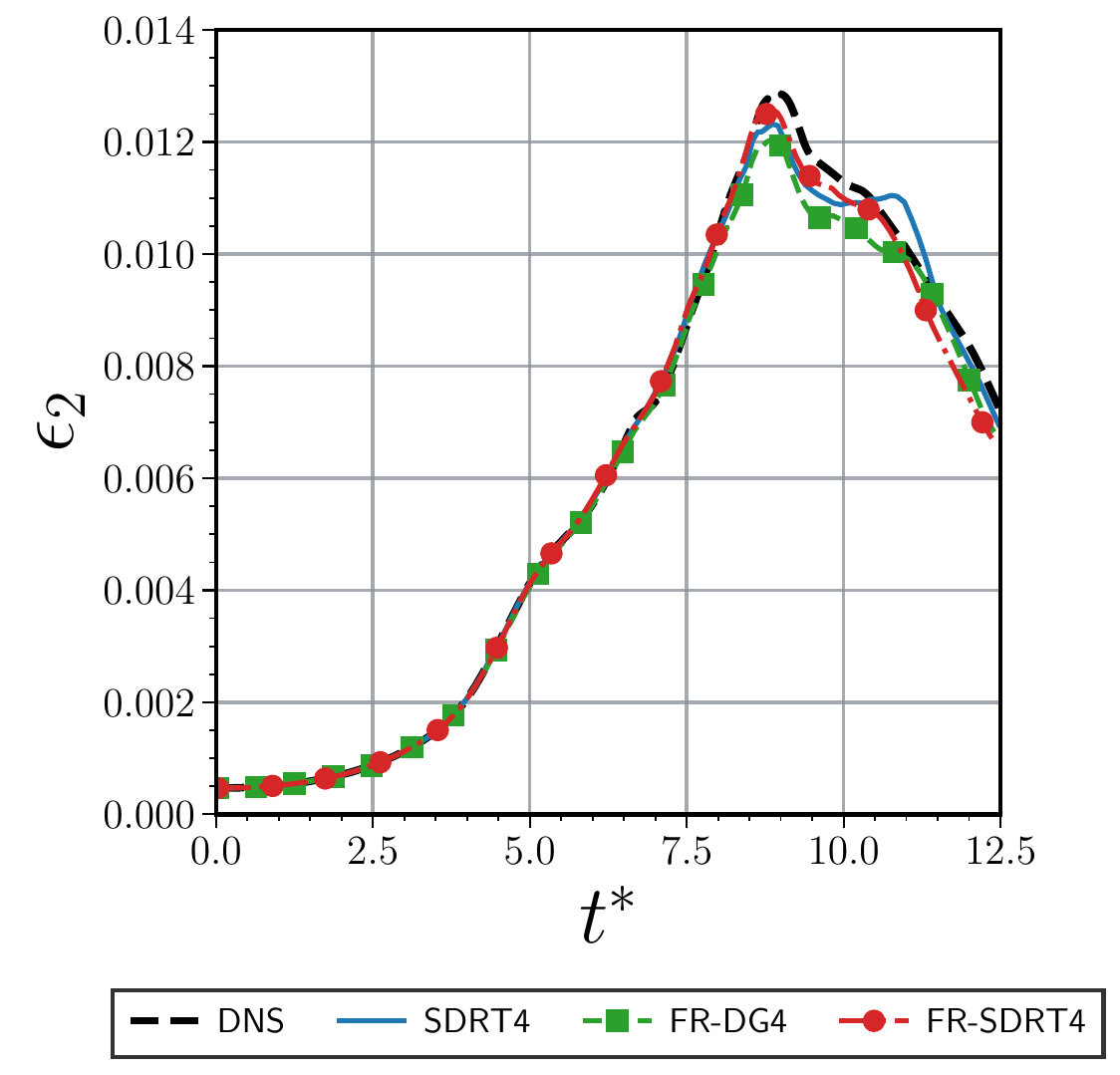}
    \caption{Hexahedral $\degree = 4$}
    \end{subfigure}
    \begin{subfigure}{0.32\textwidth}
    \centering
    \includegraphics[width=0.95\textwidth]{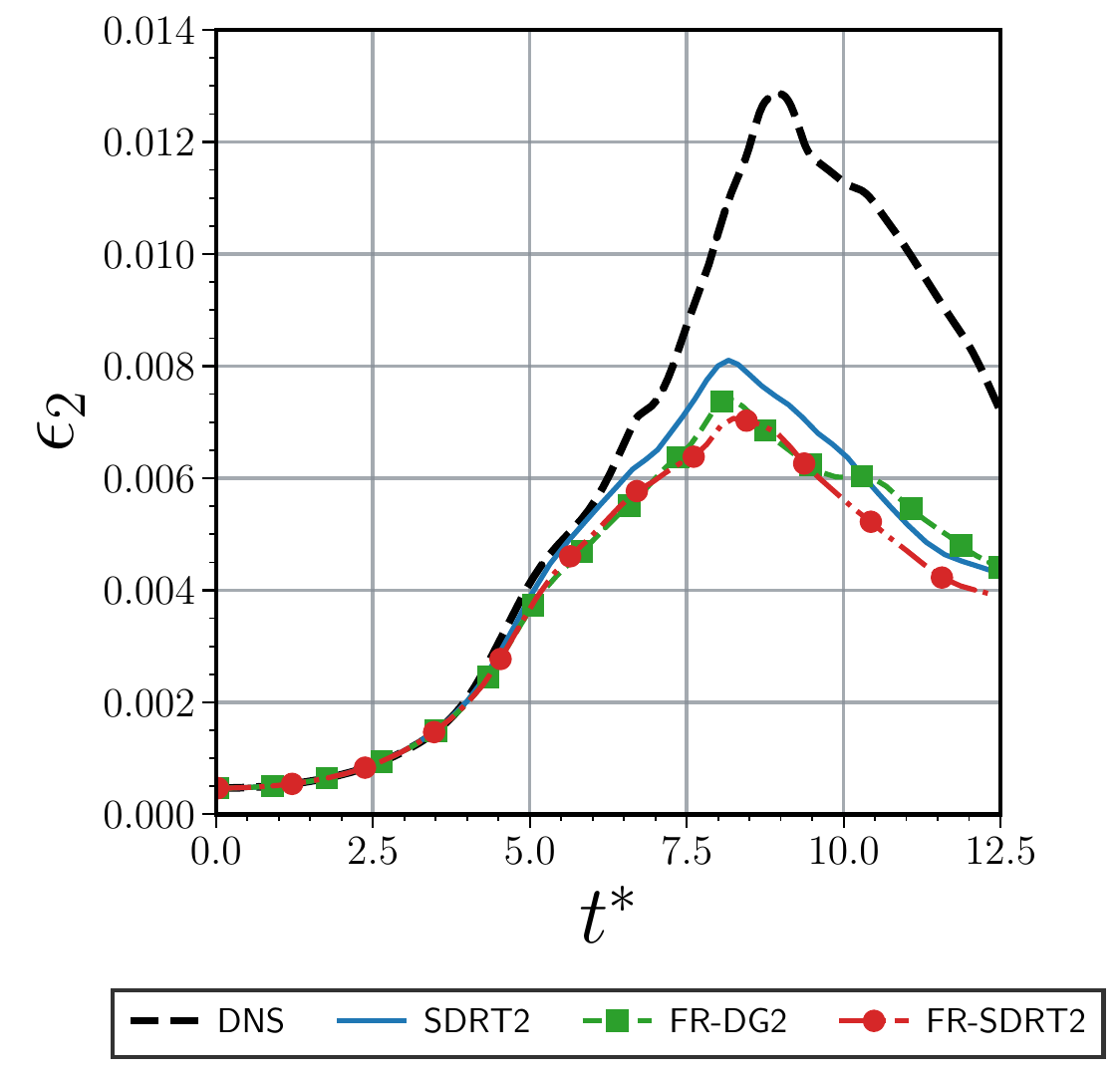}
    \caption{Prismatic $\degree = 2$}
    \end{subfigure}
    \begin{subfigure}{0.32\textwidth}
    \centering
    \includegraphics[width=0.95\textwidth]{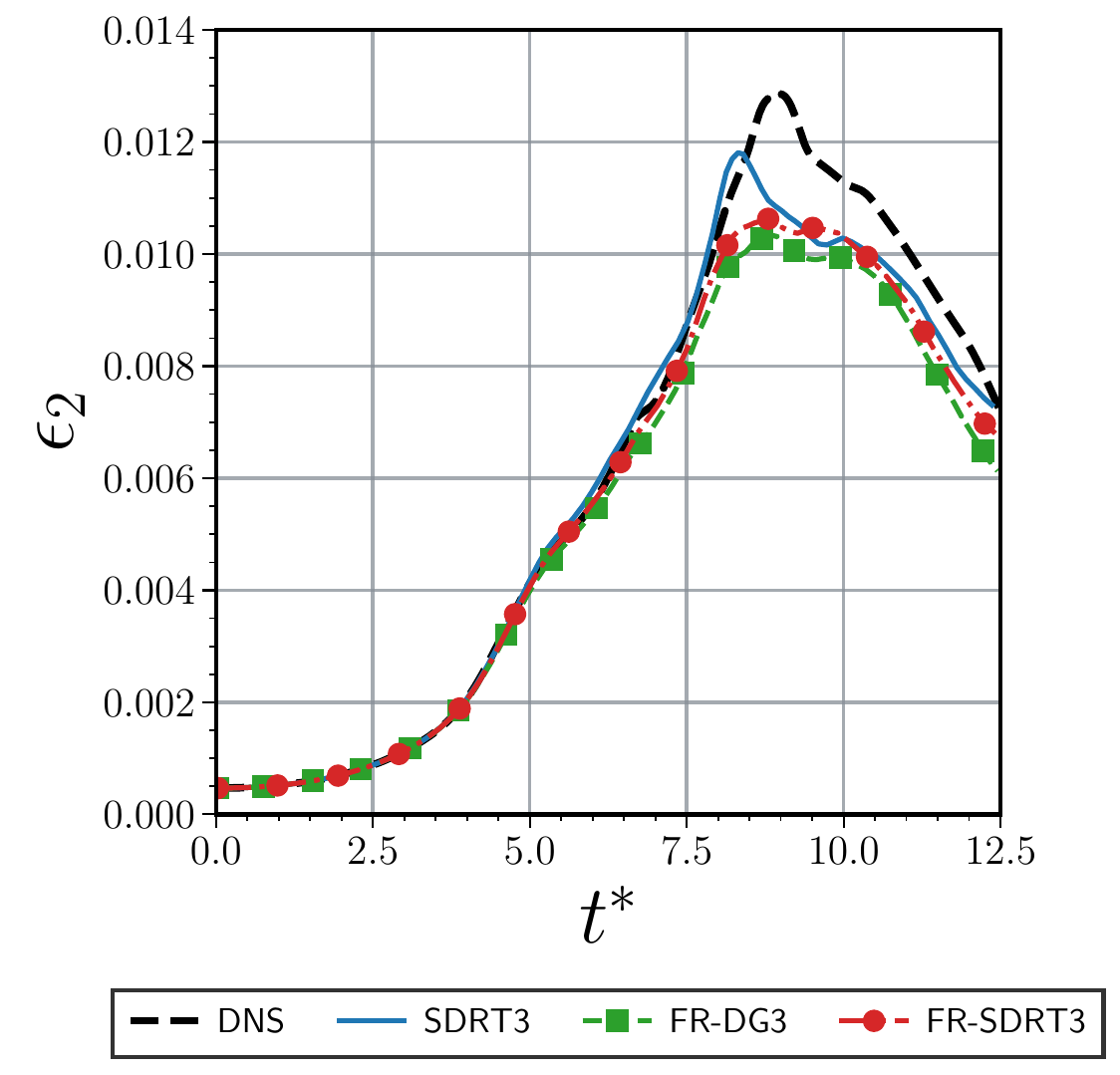}
    \caption{Prismatic $\degree = 3$}
    \end{subfigure}
    \begin{subfigure}{0.32\textwidth}
    \centering
    \includegraphics[width=0.95\textwidth]{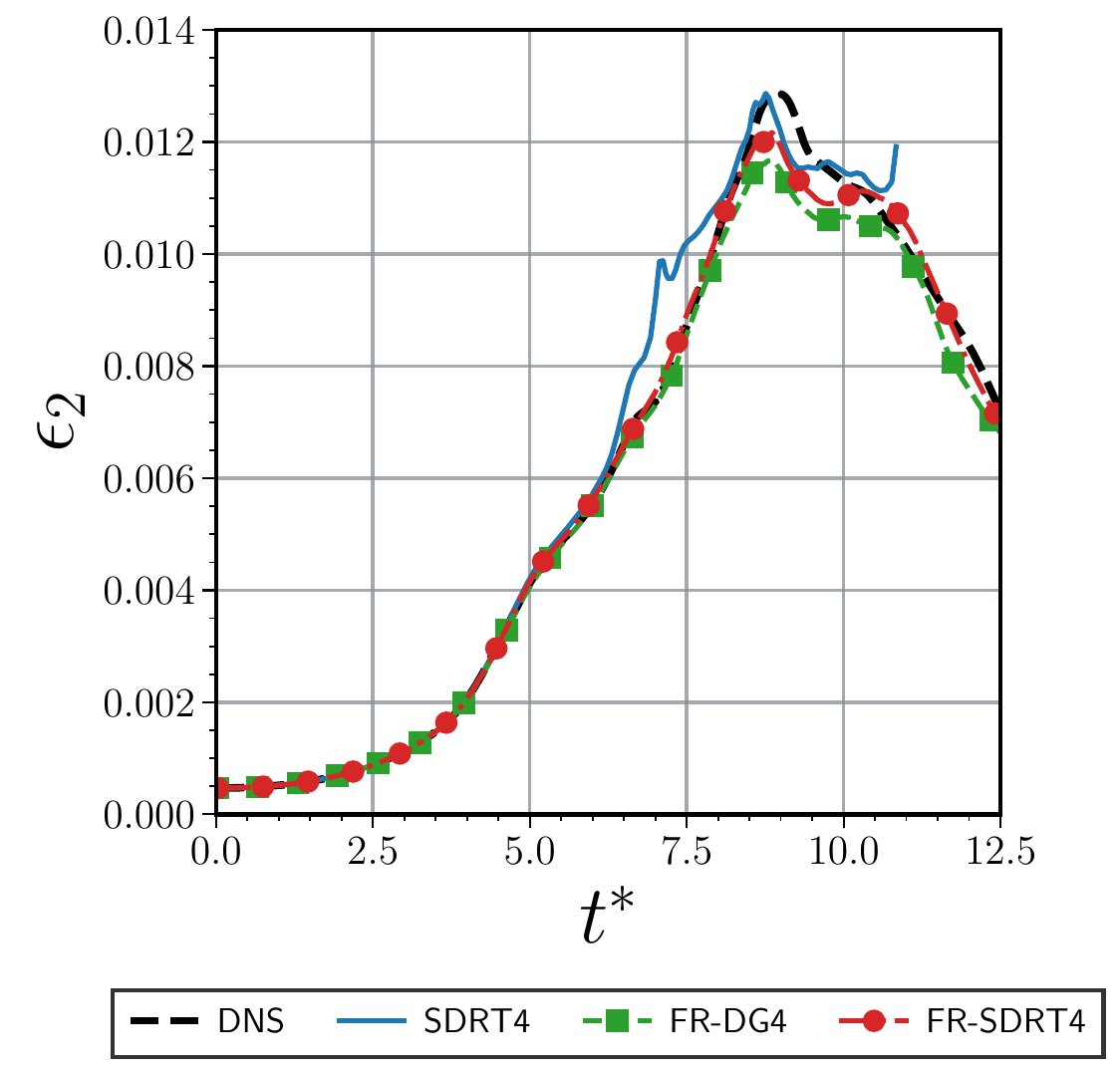}
    \caption{Prismatic $\degree = 4$}
    \end{subfigure}
    \begin{subfigure}{0.32\textwidth}
    \centering
    \includegraphics[width=0.95\textwidth]{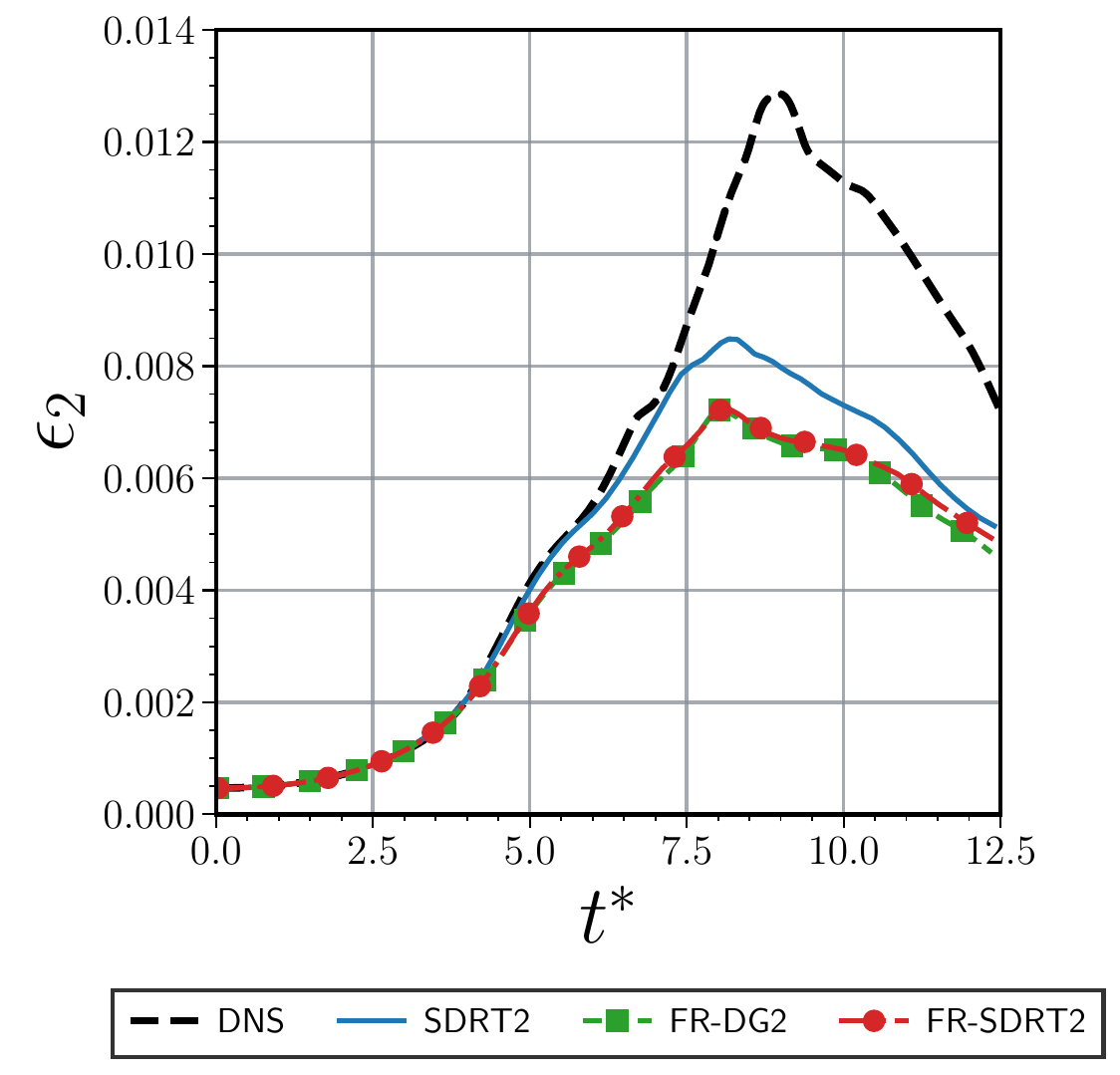}
    \caption{Tetrahedral $\degree = 2$}
    \end{subfigure}
    \begin{subfigure}{0.32\textwidth}
    \centering
    \includegraphics[width=0.95\textwidth]{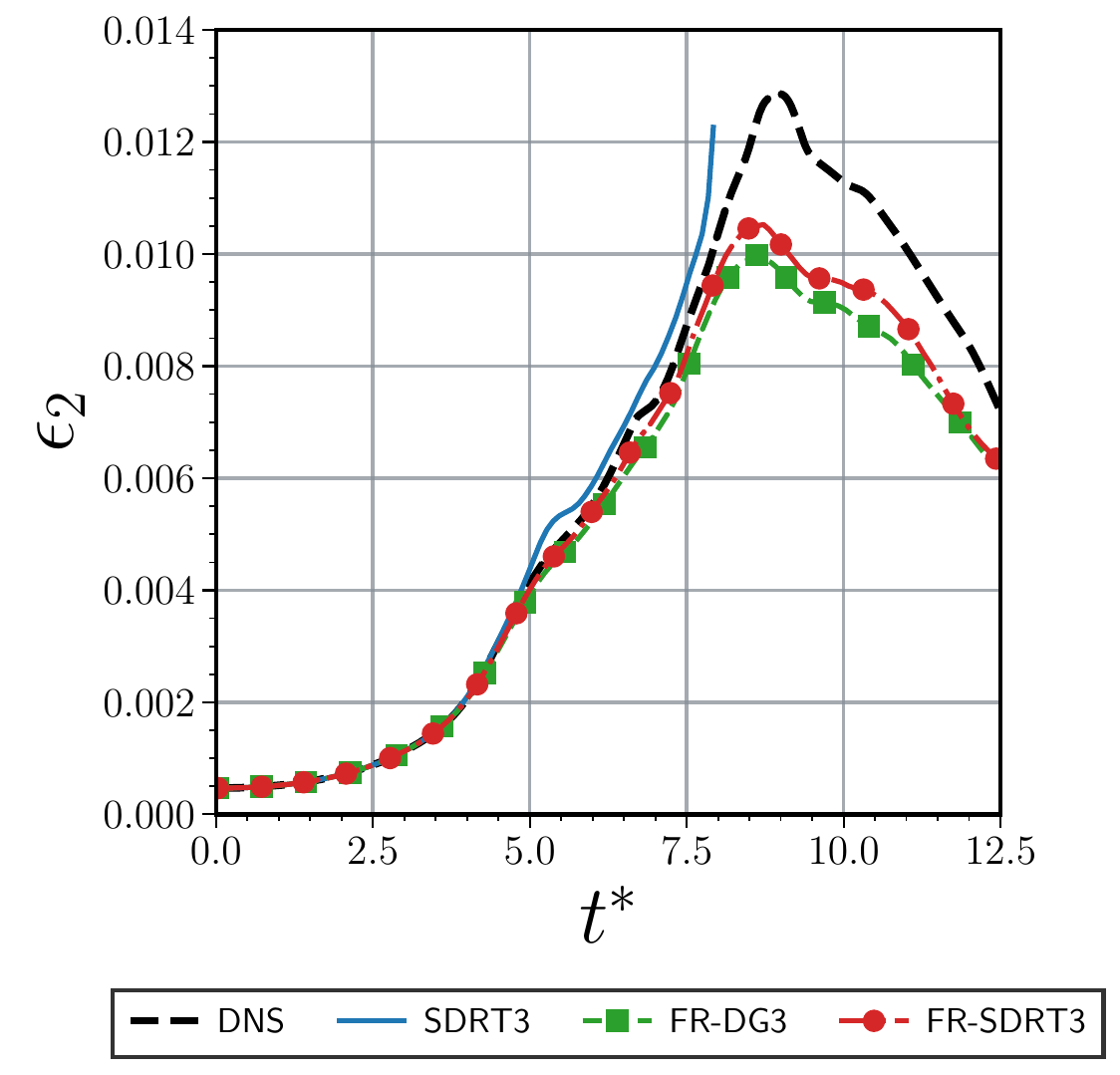}
    \caption{Tetrahedral $\degree = 3$}
    \end{subfigure}
    \begin{subfigure}{0.32\textwidth}
    \centering
    \includegraphics[width=0.95\textwidth]{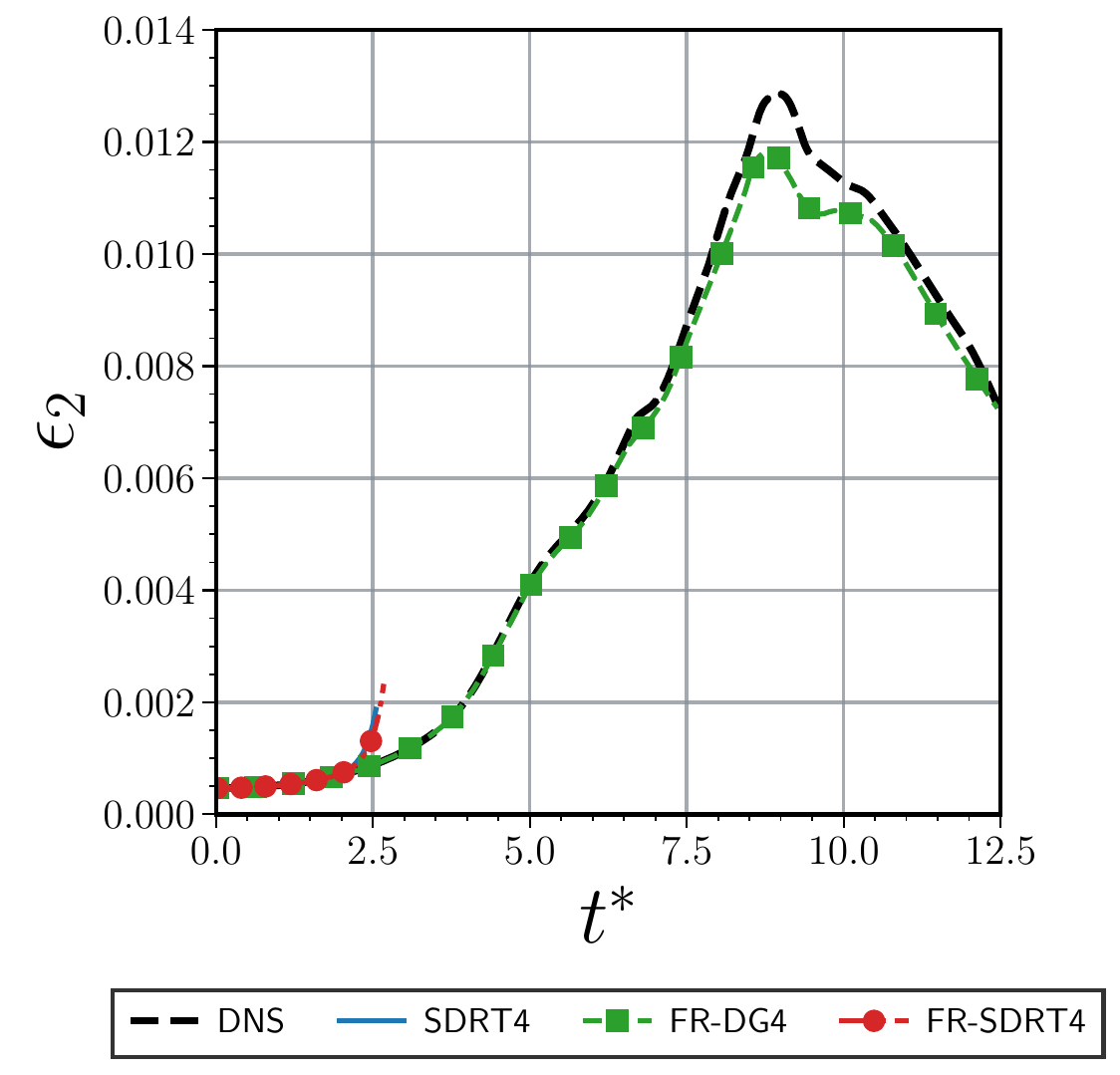}
    \caption{Tetrahedral $\degree = 4$}
    \end{subfigure}
    \caption{Temporal evolution of the $\epsilon_2$ viscous dissipation in the TGV test case with a mesh consisting of $32^3 \ncellsinpattern$, different element types and polynomial degrees.}
    \label{fig:epsilon_2_element_types_orders}
\end{figure}

Herein, results of simulations carried out with $\degree=2, 3$ and $4$, number of cells equal to $32^3 \ncellsinpattern$ and different element types utilizing SDRT, FR-DG and FR-SDRT schemes will be analyzed.
\Fref{fig:epsilon_2_element_types_orders} represents the temporal evolution of the $\epsilon_2$ in the aforementioned configuration for different element types and polynomial degrees.
With hexahedral elements the different schemes display a stable behavior and the FR-SDRT and SDRT schemes display higher values of the the viscous dissipation than FR-DG methods.
Since this viscous dissipation parameter is proportional to the enstrophy under incompressibility assumptions, higher values of the enstrophy values are often related to higher accuracy, provided that the solution does not present local instabilities.
With prismatic elements and $\degree = 2$ the SDRT2 scheme shows the highest value of the $\epsilon_2$ parameter.
Nevertheless for higher polynomial degree the solution displays a certain degree of unstable behavior, showing a very pronounce dissipation maximum at $t^* \approx 8$ with $\degree=3$ and yielding a simulation divergence for $\degree = 4$ at $t^* \approx 11$.
On the other hand, FR-SDRT schemes remain stable and display slightly higher values of the viscous dissipation for $\degree = 3$ and $\degree = 4$.
At last, results with tetrahedral elements are similar to those obtained with $\degree = 2$ and prismatic elements.
However, simulations with SDRT schemes diverge for $\degree > 2$ and those carried out with FR-SDRT methods diverge for $\degree = 4$.

Different tests were performed to assess the root cause of the divergence of simulations performed FR-SDRT and SDRT schemes.
For example, modifying the values of $\beta$, increasing the penalty parameter $\tau$, reducing the tolerances of the PI adaptive time-stepping controller, etc.
Nevertheless, neither of these approaches allowed to stabilize the simulations.
The reason behind the simulations' divergence observed with FR-SDRT and SDRT schemes must be related to non-linear instabilities, since the linear stability of the schemes was demonstrated (for the polynomial degrees studied) in \Sref{sec:von_neumann_advection} and \Sref{sec:von_neumann_diffusion}.
Non-linear stabilities are known to arise due to polynomial aliasing issues or merely due to the non-linear energy stability analysis \cite{Jameson2010}.
Future studies should be carried out to further understand which are the reasons behind the unstable behavior of FR-SDRT and SDRT schemes in non-linear test cases with three-dimensional simplex elements.

\begin{figure}
    \begin{subfigure}{0.32\textwidth}
    \centering
    \includegraphics[width=0.95\textwidth]{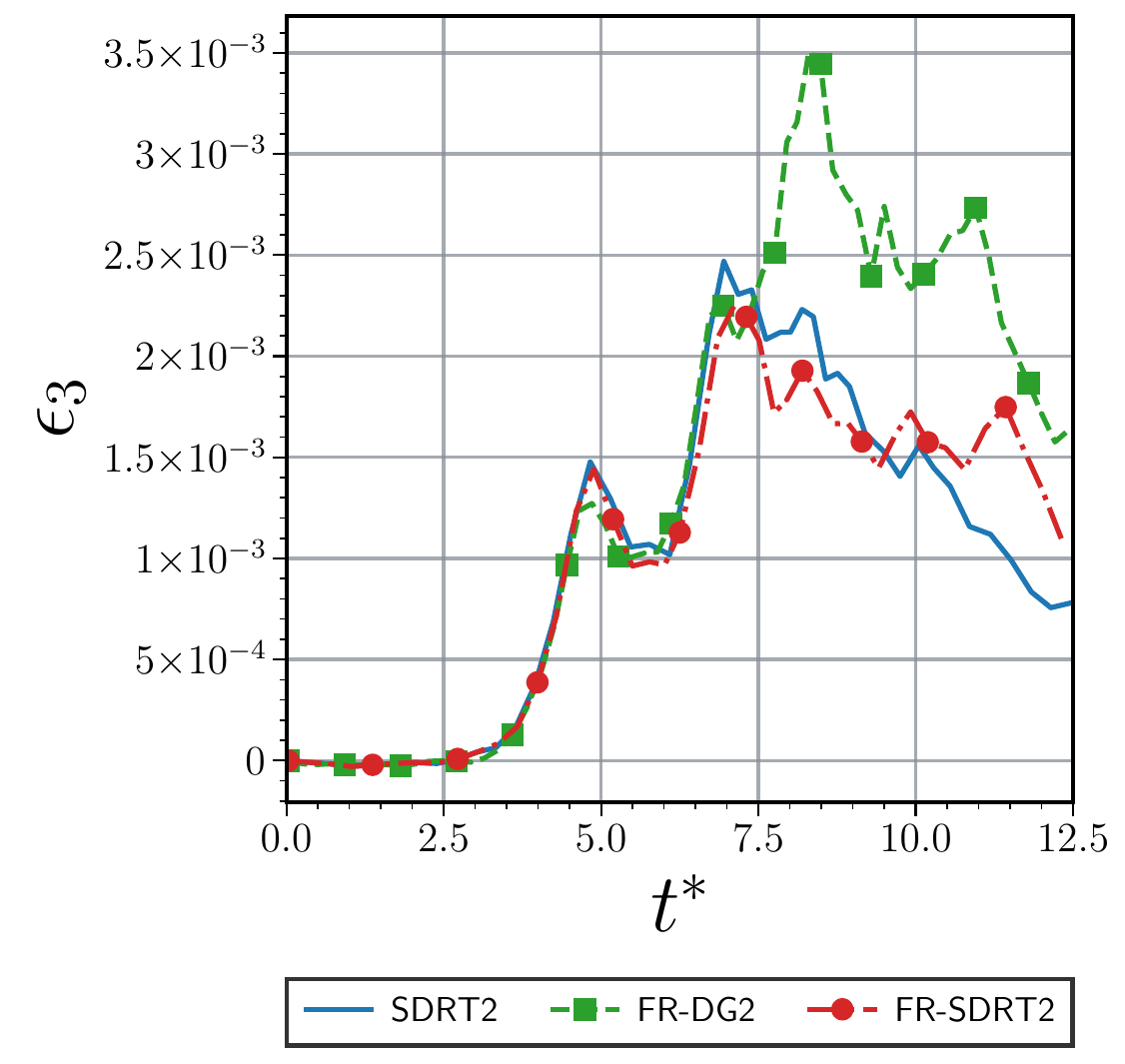}
    \caption{Hexahedral $\degree = 2$}
    \end{subfigure}
    \begin{subfigure}{0.32\textwidth}
    \centering
    \includegraphics[width=0.95\textwidth]{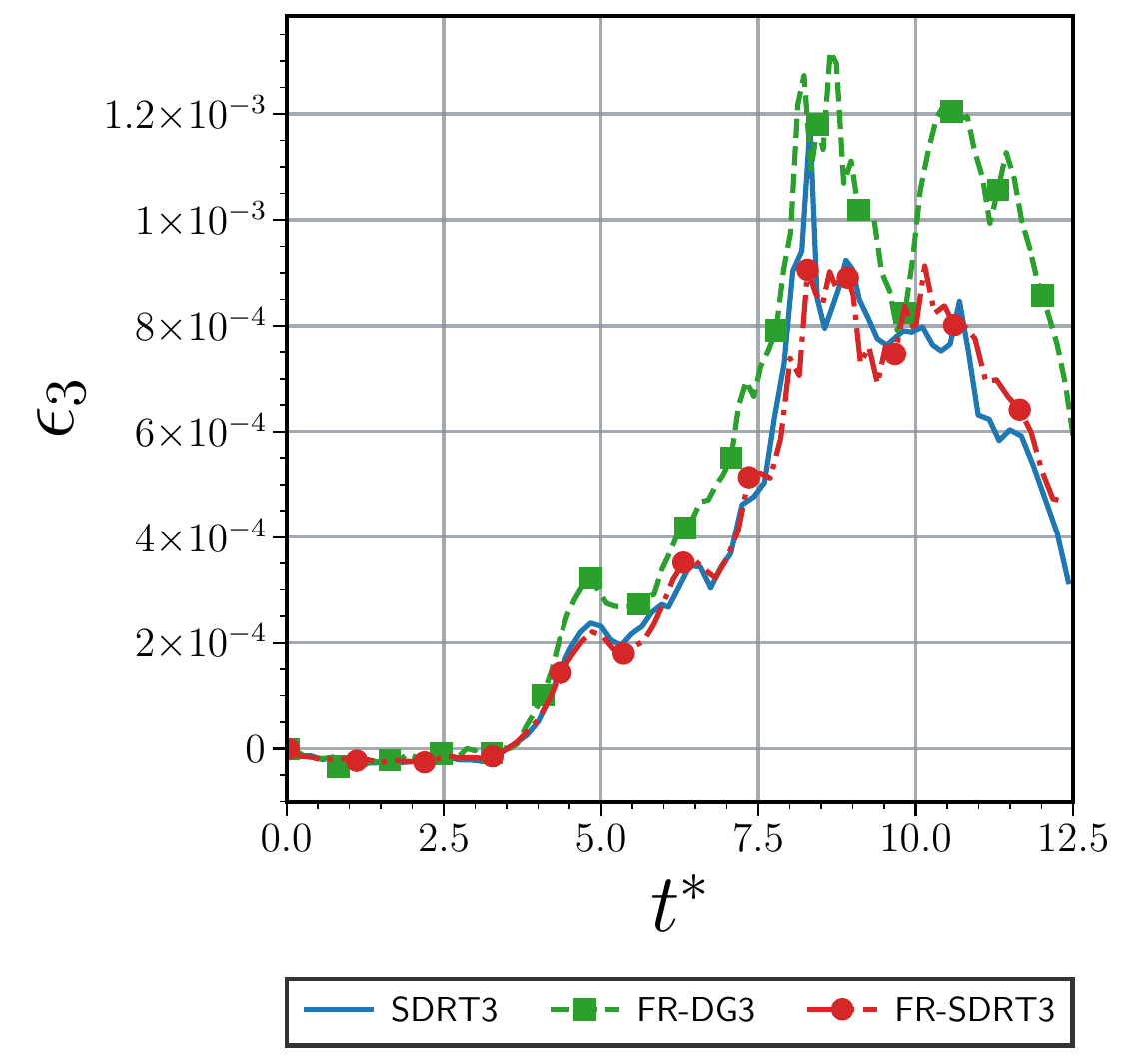}
    \caption{Hexahedral $\degree = 3$}
    \end{subfigure}
    \begin{subfigure}{0.32\textwidth}
    \centering
    \includegraphics[width=0.95\textwidth]{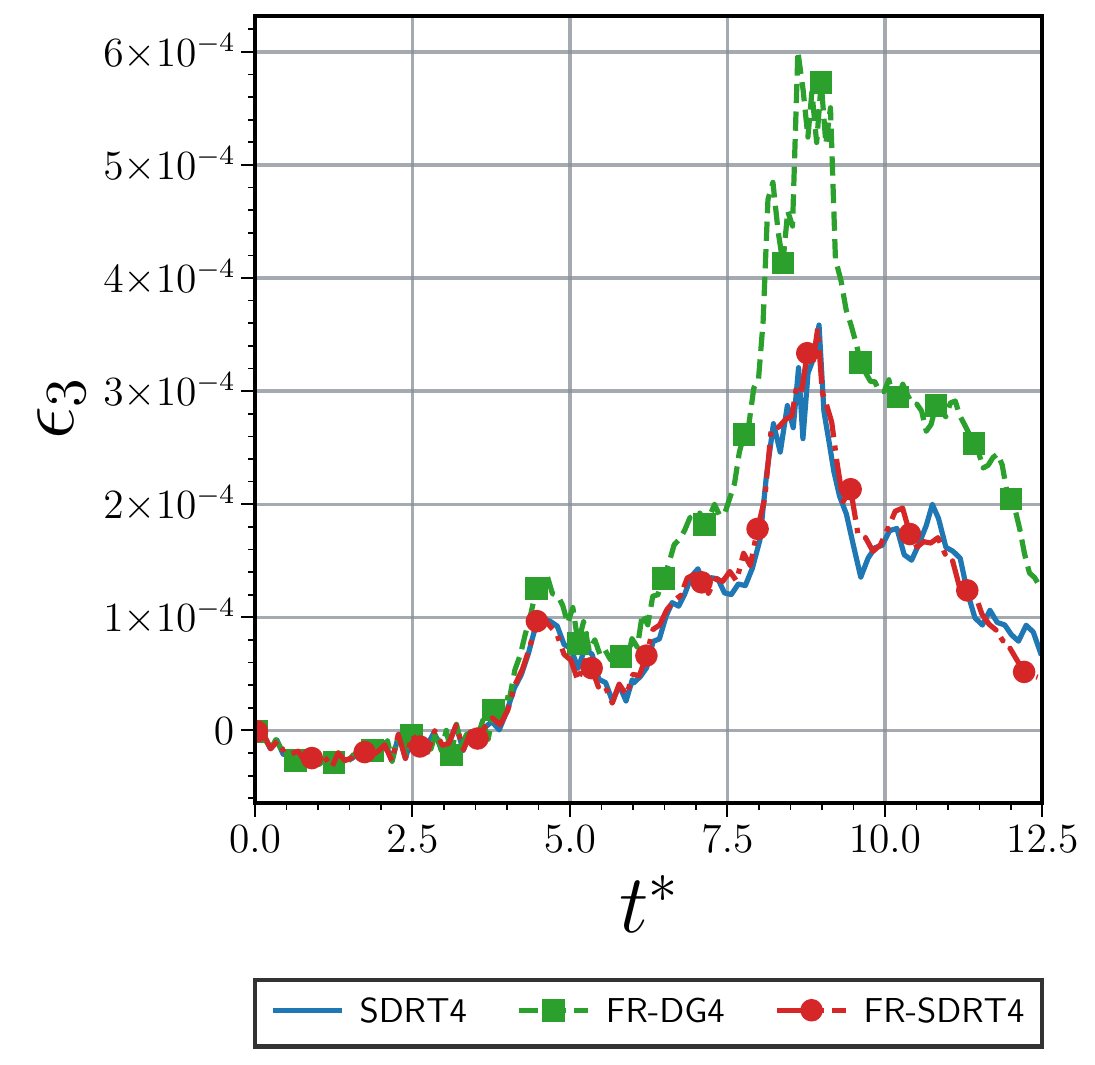}
    \caption{Hexahedral $\degree = 4$}
    \end{subfigure}
    \begin{subfigure}{0.32\textwidth}
    \centering
    \includegraphics[width=0.95\textwidth]{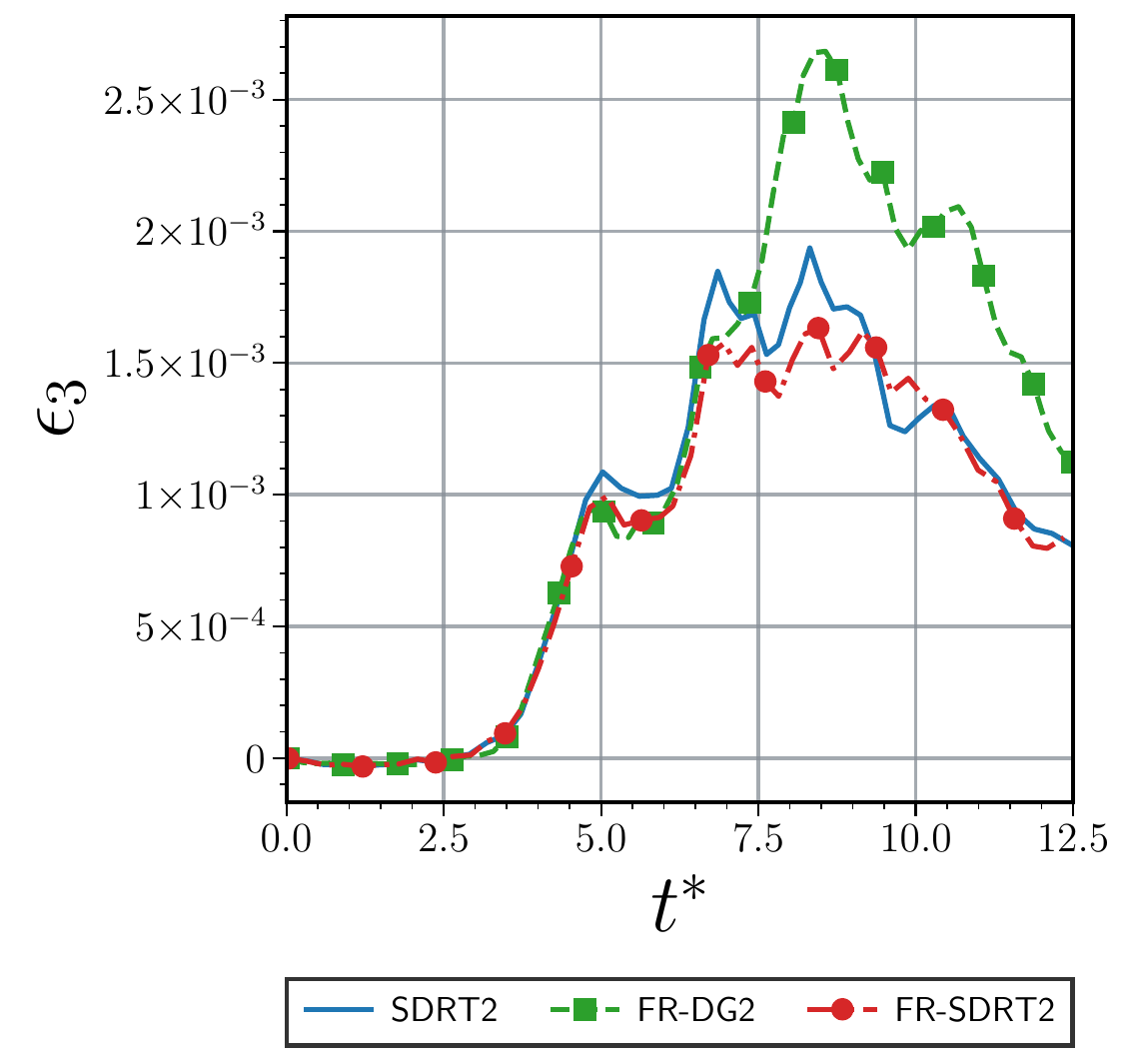}
    \caption{Prismatic $\degree = 2$}
    \end{subfigure}
    \begin{subfigure}{0.32\textwidth}
    \centering
    \includegraphics[width=0.95\textwidth]{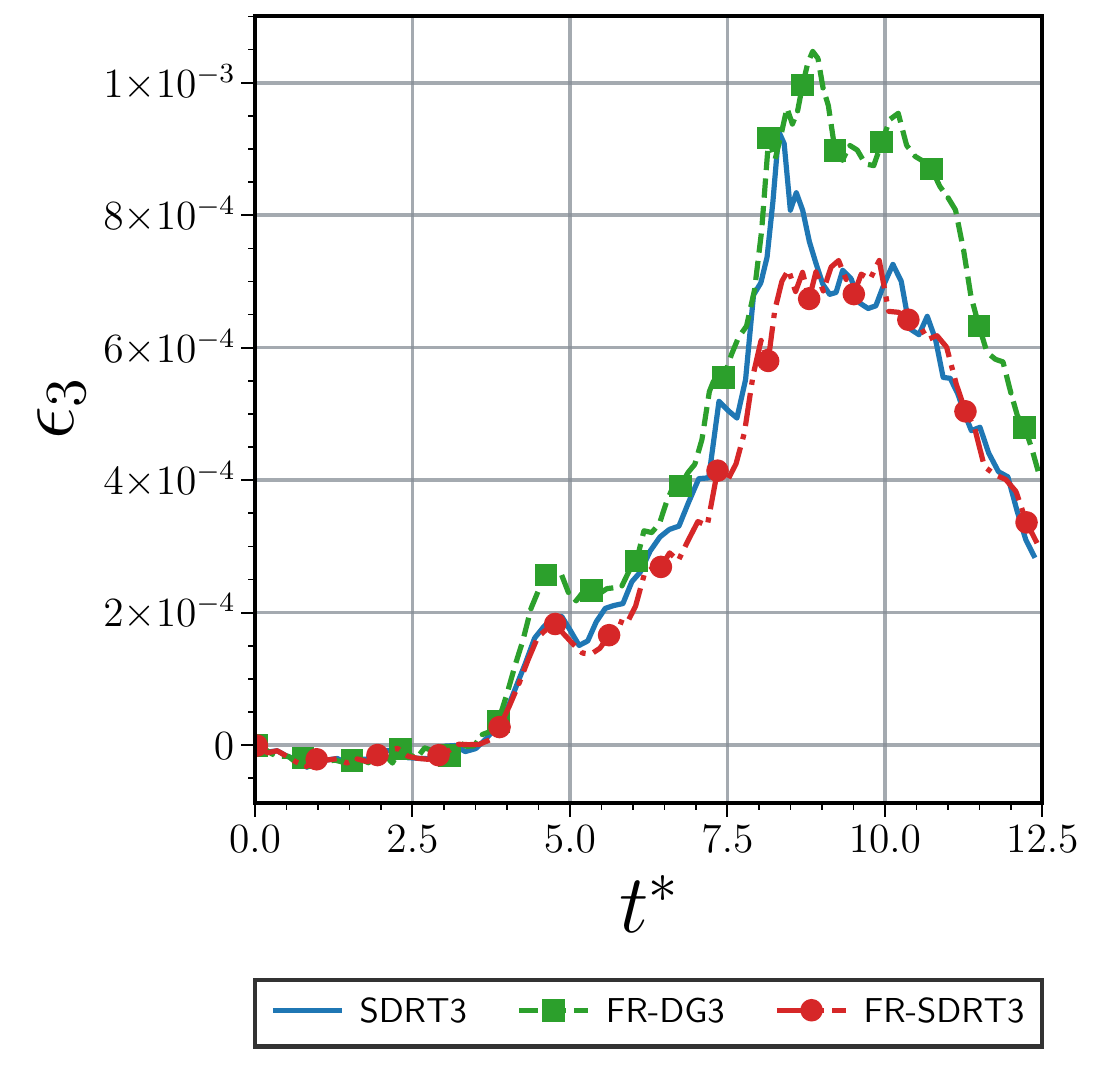}
    \caption{Prismatic $\degree = 3$}
    \end{subfigure}
    \begin{subfigure}{0.32\textwidth}
    \centering
    \includegraphics[width=0.95\textwidth]{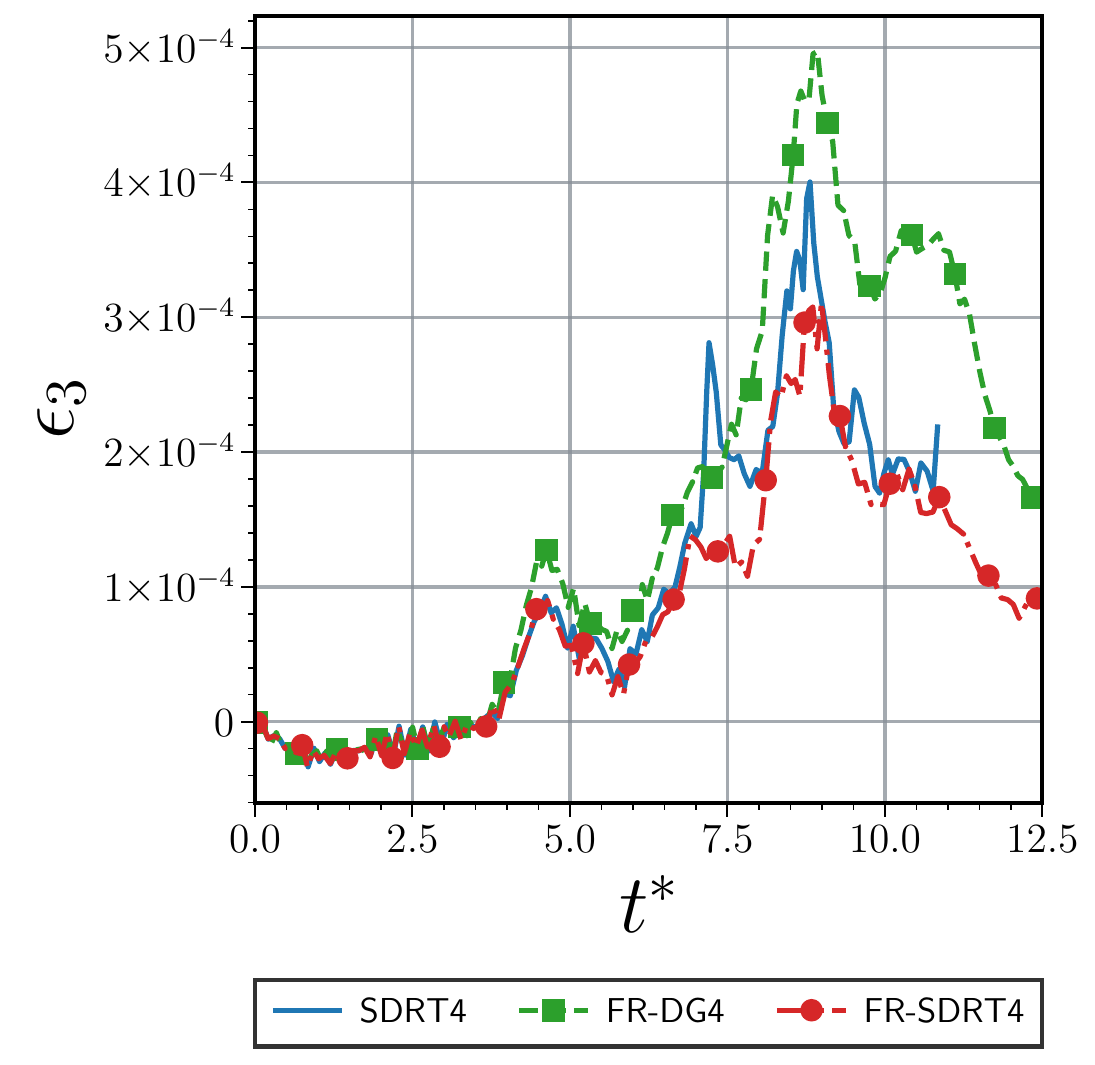}
    \caption{Prismatic $\degree = 4$}
    \end{subfigure}
    \begin{subfigure}{0.32\textwidth}
    \centering
    \includegraphics[width=0.95\textwidth]{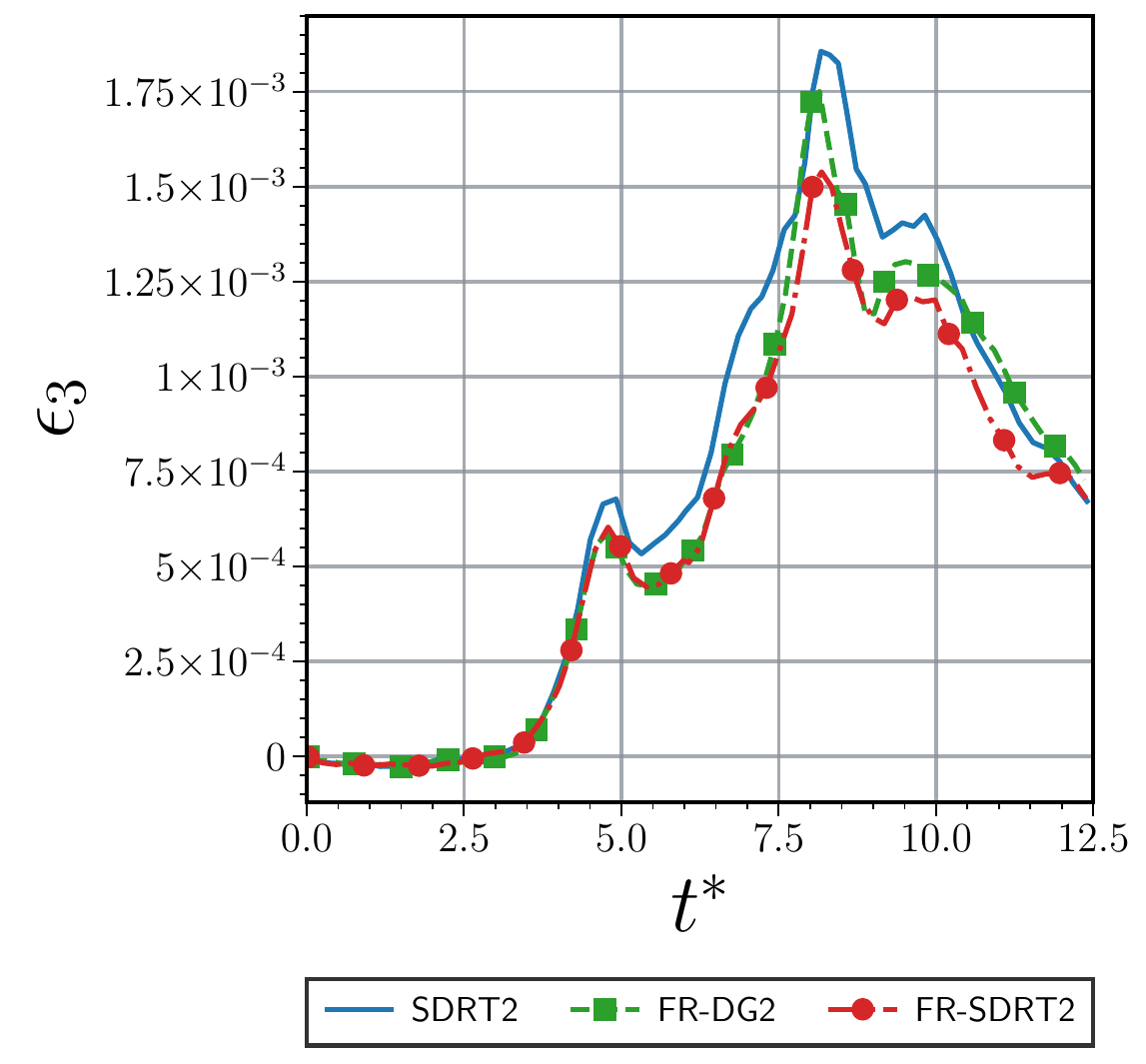}
    \caption{Tetrahedral $\degree = 2$}
    \end{subfigure}
    \begin{subfigure}{0.32\textwidth}
    \centering
    \includegraphics[width=0.95\textwidth]{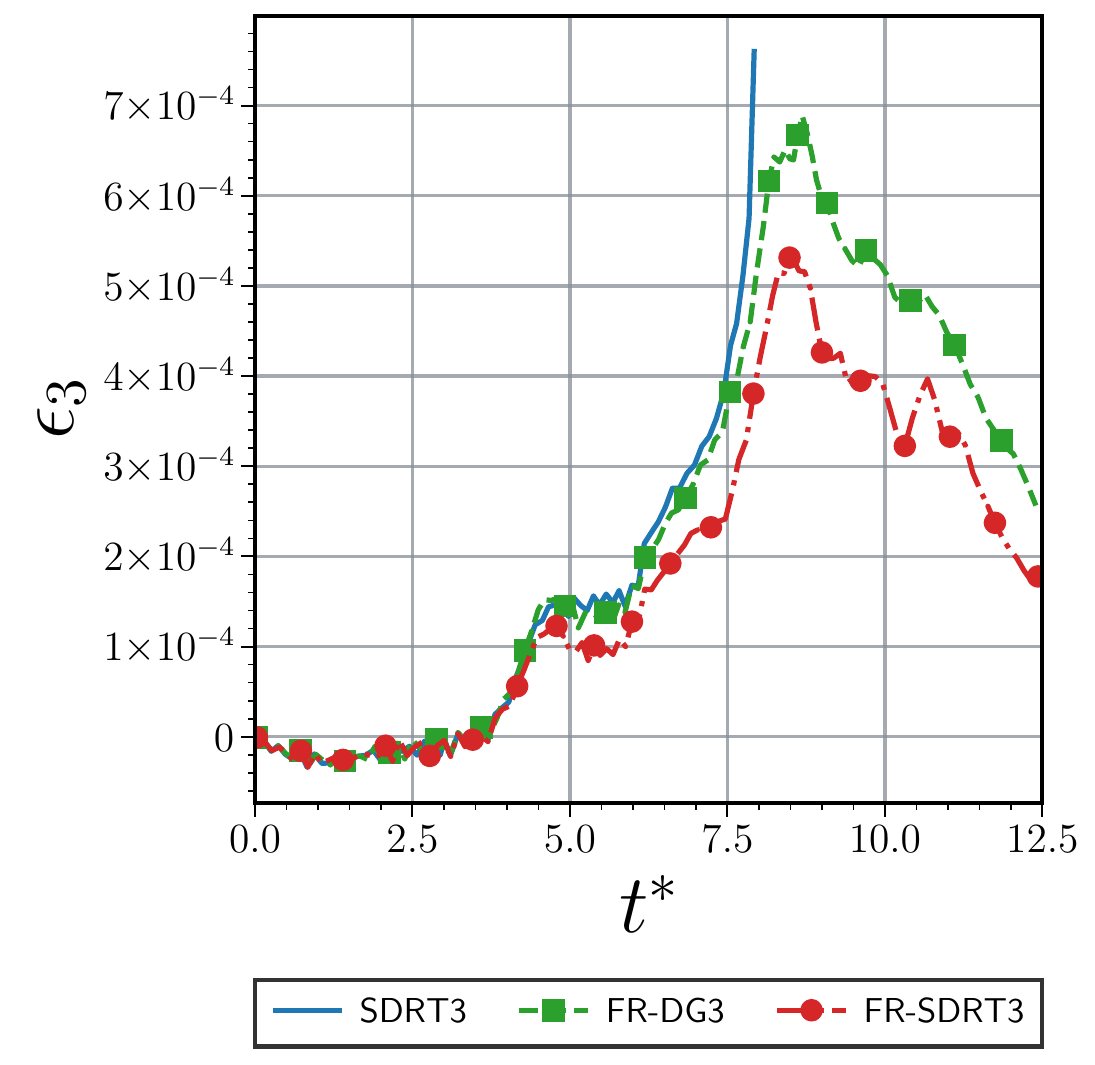}
    \caption{Tetrahedral $\degree = 3$}
    \end{subfigure}
    \begin{subfigure}{0.32\textwidth}
    \centering
    \includegraphics[width=0.95\textwidth]{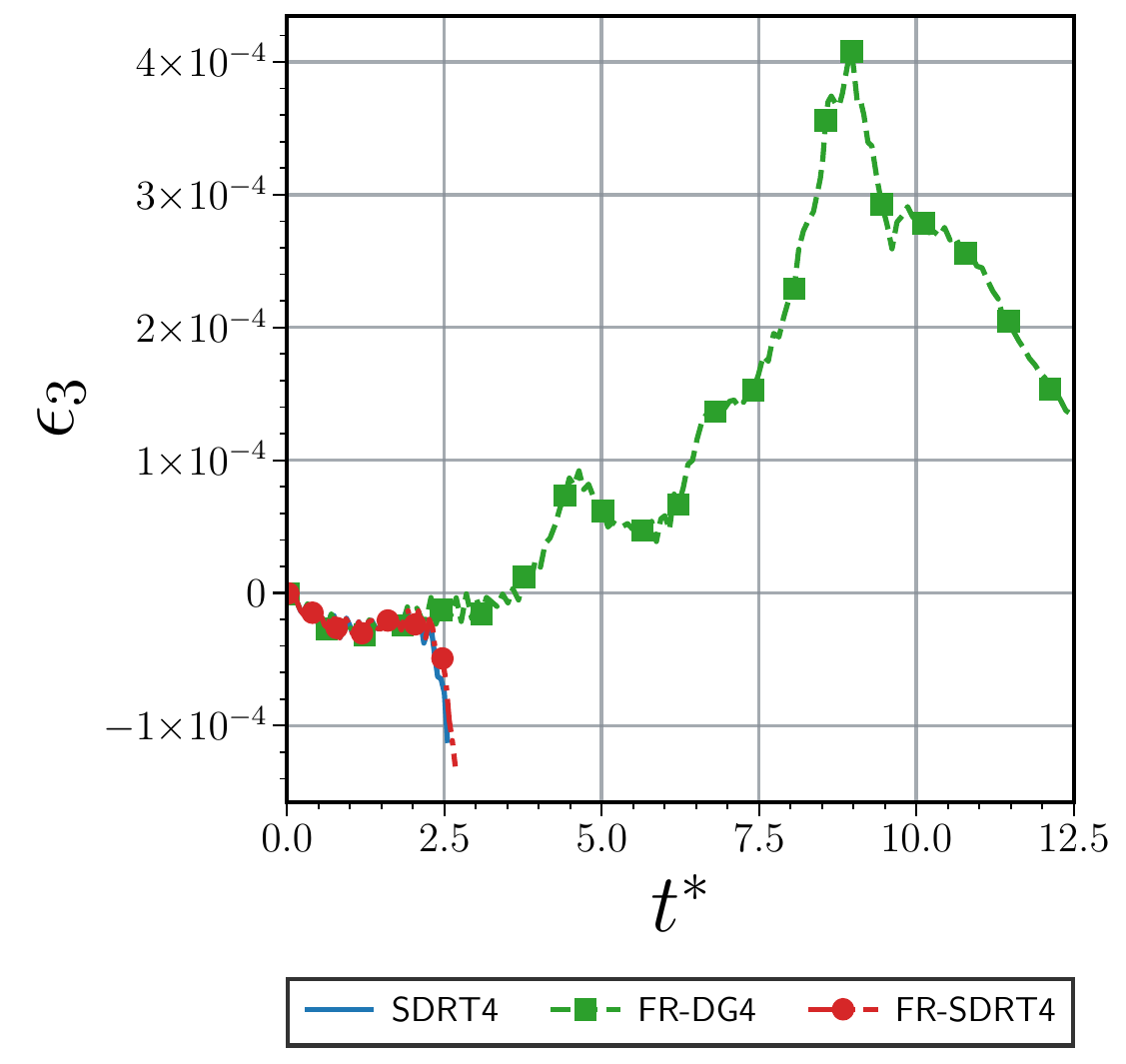}
    \caption{Tetrahedral $\degree = 4$}
    \end{subfigure}
    \caption{Temporal evolution of the $\epsilon_3$ pressure dissipation in the TGV test case with a mesh consisting of $32^3 \ncellsinpattern$, different element types and polynomial degrees.}
    \label{fig:epsilon_3_element_types_orders}
\end{figure}

For the sake of completeness, \Fref{fig:epsilon_3_element_types_orders} represents the temporal evolution of the $\epsilon_3$ in the aforementioned configuration for different element types and polynomial degrees.
The results illustrate that the pressure dissipation term shows lower values when increasing the polynomial degree.
Furthermore, FR-SDRT and SDRT schemes display lower values of the pressure dilation for most configurations, provided that the simulations remained stable.
Such a behavior is often said to imply lower aliasing issues \cite{Cox2021}.

\citet{Park2017} showed that the use of 
 adaptive time-stepping with dealising techniques in FR-DG simulations yielded higher values of the time step.
The reason behind such a fact is not determined.
To assess if this time step increase is also observed when using FR-SDRT and SDRT schemes, \Tref{table:tgv_timestep} represents the average time step of the aforementioned simulations of the TGV test case for different schemes, polynomial degrees and element types.
The results indicate that, as it was illustrated in linear problems, FR-SDRT and SDRT schemes yield higher time step values than FR-DG methods when combined with adaptive time-stepping methods and when using hexahedral and prismatic elements.
With tetrahedral elements, only the FR-SDRT and SDRT schemes with $\degree = 2$ displayed increased time step values.
The increase is substantial, and the ratios are close to those observed in the linear analysis (see \Sref{sec:temporal_stability_advection} and \Sref{sec:temporal_stability_diffusion}).
Such an increased temporal stability is special interesting for unsteady simulations, since it allows to further advance the simulations for a same simulation wall time.
However, as previously discussed, several combinations of element types and schemes resulted in unstable simulations, in particular SDRT3, SDRT4 and FR-SDRT4 showed divergent results with tetrahedral elements, while the SDRT4 scheme was unstable with prismatic elements.

To better compare the performance of the different schemes, it is also important to assess the computational performance of the different methods on an iteration-per-iteration basis.
The interested reader is referred to \ref{sec:perfo_gpu} for an in-depth study on the iteration-per-iteration computational performance of the different schemes with GPUs.

\begin{table}[h]
\centering
\begin{tabular}{@{}cccccccccc@{}}
\toprule
  & \multicolumn{3}{c}{hex} & \multicolumn{3}{c}{pri} & \multicolumn{3}{c}{tet} \\ \cmidrule(lr){2-4} \cmidrule(lr){5-7} \cmidrule(lr){8-10} 
 & SDRT & FR-DG & FR-SDRT & SDRT & FR-DG & FR-SDRT & SDRT & FR-DG & FR-SDRT\\ \midrule 
\multicolumn{1}{l|}{$\mathfrak{p}=2$} & 2.35e-02 & 1.75e-02 & 2.38e-02 & 2.07e-02 & 1.63e-02 & 2.18e-02 & 1.62e-02 & 1.38e-02 & 1.66e-02 \\ 
\multicolumn{1}{l|}{$\mathfrak{p}=3$} & 1.35e-02 & 1.05e-02 & 1.36e-02 & 1.16e-02 & 1.00e-02 & 1.22e-02 & \xmark & 9.21e-03 & 9.15e-03 \\ 
\multicolumn{1}{l|}{$\mathfrak{p}=4$} & 9.02e-03 & 6.62e-03 & 9.14e-03 & \xmark & 6.82e-03 & 7.66e-03 & \xmark & 6.39e-03 & \xmark \\ 
\end{tabular}
\caption{Average time step (in $s$) obtained in the simulation of the TGV with $32^3 \ncellsinpattern$ cells, different element types, polynomial degrees and SDRT, FR-DG and FR-SDRT schemes.
The symbol \xmark \, indicates that the combination of polynomial degree, scheme and element type resulted in unstable simulations.}
\label{table:tgv_timestep}
\end{table}

\section{Conclusions \& perspectives}
\label{sec:conclusions}

In this study, the SDRT formulation has been generalized for two and three-dimensional elements including triangular, tetrahedral and prismatic elements.
Additionally, the equivalence between FR and the SDRT was proven when solving linear equations with uniform mesh and when one utilizes a subset of the RT nodal flux basis to build the vector correction functions of the FR method, resulting in the FR-SDRT formulation.
To the best of the authors' knowledge, aside from the FR-DG method and the VCJH schemes, the FR-SDRT method is the only generalized FR scheme that allows to recover high-order and stable (under certain conditions) schemes with three-dimensional simplex elements.
All these developments were implemented in the open-source PyFR solver.
The dissipation and dispersion of the SDRT method were compared to that of the FR-DG formulation with linear advection and linear diffusion equations with two and three-dimensional tensor-product and simplex elements.
This analyses were performed using a combined-mode method, i.e. by taking into account the information of all eigenmodes in the dissipation and dispersion error measures.
To the best of the authors knowledge, this study is the first to analyze of the dissipation and dispersion properties of SEM with three-dimensional elements using the aforementioned method.
The results showed that the FR-DG maintains $2 \degree + 1$ order of accuracy even with simplex elements, while the SDRT method shows $2 \degree$ order of accuracy.
Nevertheless, the numerical errors of the SDRT method with respect to FR-DG schemes were shown to importantly increase with the polynomial degree of the considered schemes.
Moreover, through temporal linear stability, the SDRT method was shown to provide increased $\CFL_{\text{MAX}}$ values with respect to FR-DG, proving that SDRT schemes might be an appropriate choice to carry out high-order numerical simulations with simplex three-dimensional elements.
Nonetheless, it is worth noting that some SDRT schemes were found to be unstable for triangular,  tetrahedron and prismatic elements for $\degree \ge 5$.
The analytical findings were validated through linear analyses, which demonstrated the predicted order of accuracy of the schemes through numerical experiments.
Additionally, the SDRT, FR-DG and FR-SDRT were tested in the non-linear isentropic Euler vortex test case yielding  non-intuitive results since FR-SDRT schemes were found to be more accurate than SDRT schemes analysis.
At last, the SDRT method with other than tensor-product elements was shown to yield unstable simulations of the Taylor-Green-Vortex for $\degree \ge 3$.
Since these methods were shown to be linearly stable for $\degree < 5$, SDRT schemes are thought to be non-linearly unstable.
Hence, new mathematical tools need to be used to analyze these issues.
It is worth noting that the FR-SDRT method was found to be more stable that SDRT schemes, which is not intuitive.
Future works will be directed towards the extension of the SDRT method for pyramid elements, towards the reduction of the round-off errors introduced in the computation of the inverse Vandermonde matrix in the SDRT method and the theoretical assessment of the linear stability of the SDRT method with simplex elements.

\section*{CRediT authorship contribution statement}

\textbf{G. Sáez-Mischlich}: Conceptualization, Software, Validation, Formal analysis, Writing - original draft.
\textbf{J. Sierra-Ausin}: Conceptualization, Validation, Formal analysis, Writing - original draft.
\textbf{J. Gressier}: Conceptualization, Funding acquisition, Writing - review and editing

\section*{Declaration of Competing Interest}

The authors hereby declare that they have no known competing financial interests or personal relationships that could have appeared to influence the work reported in this paper.

\section*{Acknowledgments}
The authors would like to thank ISAE-SUPAERO for providing with the computational resources and founding to carry out this study.
We would also like to thank Freddie Witherden for his gracious contributions to the development of the open-source PyFR solver and the review of the pull-requests which added certain functionalities needed to perform the studies shown in  this work.





\appendix

\FloatBarrier
\section{Raviart-Thomas (RT) basis and degree of freedom distribution}
\label{sec:rt_basis}

Let $\polynomial$ be a polynomial basis, described in \cite{Bergot2013}, such that
\begin{equations}{eq:polynomial_up_to_degree}
    \polynomial_{\degree}(\transformed{x}_0) &\equiv \text{span}\{ x_0^i \} \text{ with } 0 \le i \le \degree \text{ or }\\
    \polynomial_{\degree}(\transformed{x}_0, \transformed{x}_1) &\equiv \text{span}\{  x_0^i x_1^j \} \text{ with } 0 \le i, j \text{ and } i + j \le \degree  \text{ or }\\
    \polynomial_{\degree}(\transformed{x}_0, \transformed{x}_1, \transformed{x}_2) &\equiv \text{span}\{  x_0^i x_1^j x_2^k \}  \text{ with } 0 \le i, j, k \text{ and } i + j + k \le \degree,
\end{equations}
where each each basis of $\polynomial$ presents $\ndim$ dimensions (depending on the considered element in which these polynomial bases are described).
Within $\ndim$-dimensional bases, only one component is non-zero.
The index of such non-zero component may be deduced from the context and from the examples provided for each elements.
Additionally, let us define $\overline{\polynomial}$ as
\begin{equations}{eq:polynomial_of_unique_degree}
    \overline{\polynomial}_{\degree}(\transformed{x}_0) &\equiv \text{span}\{ x_0^i \} \text{ with } i = \degree \text{ or }\\
    \overline{\polynomial}_{\degree}(\transformed{x}_0, \transformed{x}_1) &\equiv \text{span}\{ x_0^i x_1^j \} \text{ with } i + j = \degree \text{ or }\\
    \overline{\polynomial}_{\degree}(\transformed{x}_0, \transformed{x}_1, \transformed{x}_2) &\equiv \text{span}\{ x_0^i x_1^j x_2^k \} \text{ with } i + j + k = \degree.
\end{equations}
At last, we introduce the polynomial $\mathbb{Q}$
\begin{equations}{eq:polynomial_tensor_product}
    \mathbb{Q}_{n, m}(\transformed{x}_0, \transformed{x}_1) &\equiv \text{span}\{  x_0^i x_1^j \}  \text{ with } 0 \le i, j \le n, m \text{ or }\\
    \mathbb{Q}_{n, m, l}(\transformed{x}_0, \transformed{x}_1, \transformed{x}_2) &\equiv \text{span}\{  x_0^i x_1^j x_2^k \} \text{ with } 0 \le i, j, k \le n, m, l .
\end{equations}
This polynomial also presents $\ndim$-dimensional.

The operator $\times$ indicates the Cartesian product of bases while operator $\oplus$ indicates addition.
The interested reader is referred to the RT basis examples provided within the following sections to better understand the nomenclature.

\subsection{Triangles and Tetrahedrons elements}

The reference triangle and tetrahedron are defined such that 
\begin{equation}
    \transformed{\domain}_{\itype} \in \transformed{x}_i \ge -1; \quad \sum_i \transformed{x}_i \le 0.
\end{equation}


The RT modal bases of a SDRT scheme of degree $\degree$ for triangle and tetrahedron elements are given by
\begin{equation}
    \atflux{\fluxmodalbasis}(\transformed{\vect{x}}) \equiv \polynomial_{\degree}(\transformed{\vect{x}})^{\ndim} \oplus \transformed{\vect{x}} \overline{\polynomial}_{\degree}(\transformed{\vect{x}}) .
\end{equation}
To give an example of such polynomial basis, the RT basis for a triangular element with $\degree = 1$ reads
\begin{equation}
    \fluxmodalbasis(\transformed{\vect{x}}) \equiv
    \begin{bmatrix}
    1 \\ 0
    \end{bmatrix} ,
    \begin{bmatrix}
    \transformed{x} \\ 0
    \end{bmatrix} ,
    \begin{bmatrix}
    \transformed{y} \\ 0
    \end{bmatrix} ,
    \begin{bmatrix}
    0 \\ 1
    \end{bmatrix} ,
    \begin{bmatrix}
    0 \\ \transformed{x}
    \end{bmatrix} ,
    \begin{bmatrix}
    0 \\ \transformed{y}
    \end{bmatrix} ,
    \begin{bmatrix}
    \transformed{x}^2 \\ \transformed{x}\transformed{y}
    \end{bmatrix} ,
    \begin{bmatrix}
    \transformed{x}\transformed{y} \\ \transformed{y}^2
    \end{bmatrix},
\end{equation}
while the RT basis for a tetrahedron element with $\degree = 1$ is given by
\begin{equation}
    \fluxmodalbasis(\transformed{\vect{x}}) \equiv
    \begin{bmatrix}
    1 \\ 0 \\ 0
    \end{bmatrix} ,
    \begin{bmatrix}
    \transformed{x} \\ 0 \\ 0
    \end{bmatrix} ,
    \begin{bmatrix}
    \transformed{y} \\ 0 \\ 0
    \end{bmatrix} ,
    \begin{bmatrix}
    \transformed{z} \\ 0 \\ 0
    \end{bmatrix} ,
    \begin{bmatrix}
    0 \\ 1 \\ 0
    \end{bmatrix} ,
    \begin{bmatrix}
    0 \\ \transformed{x} \\ 0
    \end{bmatrix} ,
    \begin{bmatrix}
    0 \\ \transformed{y} \\ 0
    \end{bmatrix} ,
    \begin{bmatrix}
    0 \\ \transformed{z} \\ 0
    \end{bmatrix} ,
    \begin{bmatrix}
    0 \\ 0 \\ 1
    \end{bmatrix} ,
    \begin{bmatrix}
    0 \\ 0  \\ \transformed{x}
    \end{bmatrix} ,
    \begin{bmatrix}
    0 \\ 0 \\ \transformed{y}
    \end{bmatrix} ,
    \begin{bmatrix}
    0 \\ 0 \\ \transformed{z}
    \end{bmatrix} ,
    \begin{bmatrix}
    \transformed{x}^2 \\ \transformed{x}\transformed{y} \\ \transformed{x}\transformed{z}
    \end{bmatrix} ,
    \begin{bmatrix}
    \transformed{x}\transformed{y} \\ \transformed{y}^2 \\ \transformed{y}\transformed{z}
    \end{bmatrix} ,
    \begin{bmatrix}
    \transformed{x}\transformed{z} \\ \transformed{y}\transformed{z} \\ \transformed{z}^2
    \end{bmatrix}.
\end{equation}

In triangular elements, the orthogonal basis described in \cite{Mozolevski2017} 
are utilized to define the modal basis.
For tetrahedron elements, one may rely on the H(div) hierarchical basis described in \cite{Bergot2013} to reduce the condition number of the flux Vandermonde matrix.
Such choices have been proven crucial to reduce the condition number of the flux Vandermonde matrix in this work.




To distribute the degrees of freedom within triangle and tetrahedron elements, external flux points are related to a single degree of freedom equal to the unitary normal in the reference space at the considered points.
On the other hand, $\ndim$ degrees of freedom oriented along each of the unitary principal axis of the reference element are imposed at unique internal flux points.

Within triangular elements, external flux points are located at Gauss-Legendre quadrature points resulting in $\degree + 1$ points per edge.
On the other hand, internal flux points are located at Williams-Shunn quadrature points \cite{Williams2014} with $\frac{\degree (\degree + 1)}{2}$ unique points.
\Fref{fig:Notation_Triangle_SDRT} represents the distribution of such solution and flux points distribution along the reference triangle element.


For tetrahedron elements, Williams-Shunn quadrature points of degree $\degree$ with $\frac{(\degree + 1)(\degree + 2)}{2}$ points are considered at triangular faces, while internal flux points are set at Williams-Shunn quadrature points of degree $\degree - 1$ with a total of $\frac{\degree(\degree + 1)(\degree + 2)}{6}$ unique points.
\Fref{fig:Notation_Tetra_SDRT} represents the distribution of such solution and flux points distribution along the reference tetrahedron element.


Solution points are also particularized at Williams-Shunn quadrature points such that $\atsol{N_{\itype}} = \frac{(\degree + 1)(\degree + 2)}{2}$ and $\atsol{N_{\itype}} = \frac{(\degree + 1)(\degree + 2)(\degree + 3)}{6}$ for triangular and tetrahedron elements respectively.


\subsection{Tensor-product elements}

The reference tensor-product element is defined as
\begin{equation}
    \transformed{\domain}_{\itype} \in \vect{x} \in [-1, 1]^{\ndim} .
\end{equation}
The RT modal bases of a SDRT scheme of degree $\degree$ for two-dimensional tensor-product elements read
\begin{equation}
    \fluxmodalbasis(\transformed{\vect{x}}) \equiv \mathbb{Q}_{\degree + 1, \degree}(\vect{\transformed{x}}) \times \mathbb{Q}_{\degree, \degree + 1}(\vect{\transformed{x}}),
\end{equation}
and for three-dimensional tensor-product elements
\begin{equation}
    \fluxmodalbasis(\transformed{\vect{x}}) \equiv \mathbb{Q}_{\degree + 1, \degree, \degree}(\vect{\transformed{x}}) \times \mathbb{Q}_{\degree, \degree + 1, \degree}(\vect{\transformed{x}}) \times \mathbb{Q}_{\degree, \degree , \degree + 1}(\vect{\transformed{x}}).
\end{equation}

An example of the RT bases for a two-dimensional tensor-product element with $\degree=1$ is given by
\begin{equation}
    \fluxmodalbasis(\transformed{\vect{x}}) \equiv
    \begin{bmatrix}
    1 \\ 0
    \end{bmatrix} ,
    \begin{bmatrix}
    \transformed{x} \\ 0
    \end{bmatrix} ,
    \begin{bmatrix}
    \transformed{y} \\ 0
    \end{bmatrix} ,
    \begin{bmatrix}
    \transformed{x}\transformed{y} \\ 0
    \end{bmatrix} ,
    \begin{bmatrix}
    \transformed{x}^2 \\ 0
    \end{bmatrix} ,
    \begin{bmatrix}
    \transformed{x}^2\transformed{y} \\ 0
    \end{bmatrix} ,
    \begin{bmatrix}
    0 \\ 1
    \end{bmatrix} ,
    \begin{bmatrix}
    0 \\ \transformed{x}
    \end{bmatrix} ,
    \begin{bmatrix}
    0 \\ \transformed{y}
    \end{bmatrix} ,
    \begin{bmatrix}
    0 \\ \transformed{x}\transformed{y}
    \end{bmatrix} ,
    \begin{bmatrix}
    0 \\ \transformed{y}^2
    \end{bmatrix} ,
    \begin{bmatrix}
    0 \\ \transformed{x}\transformed{y}^2
    \end{bmatrix} .
\end{equation}
On the other hand, the RT basis of hexahedron elements with $\degree = 1$  are
\begin{align*}
    \fluxmodalbasis(\transformed{\vect{x}}) \equiv&
    \begin{bmatrix}
    1 \\ 0 \\ 0
    \end{bmatrix} ,
    \begin{bmatrix}
    \transformed{x} \\ 0 \\ 0
    \end{bmatrix} ,
    \begin{bmatrix}
    \transformed{y} \\ 0 \\ 0
    \end{bmatrix} ,
    \begin{bmatrix}
    \transformed{z} \\ 0 \\ 0
    \end{bmatrix} ,
    \begin{bmatrix}
    \transformed{x}\transformed{y} \\ 0 \\ 0
    \end{bmatrix} ,
    \begin{bmatrix}
    \transformed{x}\transformed{z} \\ 0 \\ 0
    \end{bmatrix} ,
    \begin{bmatrix}
    \transformed{y}\transformed{z} \\ 0 \\ 0
    \end{bmatrix} ,
    \begin{bmatrix}
    \transformed{x}^2 \\ 0 \\ 0
    \end{bmatrix} ,
    \begin{bmatrix}
    \transformed{x}^2\transformed{y} \\ 0 \\ 0
    \end{bmatrix} ,
    \begin{bmatrix}
    \transformed{x}^2\transformed{z} \\ 0 \\ 0
    \end{bmatrix} , \\
    &\begin{bmatrix}
    0 \\ 1 \\ 0
    \end{bmatrix} ,
    \begin{bmatrix}
    0 \\ \transformed{x} \\ 0
    \end{bmatrix} ,
    \begin{bmatrix}
    0 \\ \transformed{y} \\ 0
    \end{bmatrix} ,
    \begin{bmatrix}
    0 \\ \transformed{z} \\ 0
    \end{bmatrix} ,
    \begin{bmatrix}
    0 \\ \transformed{x}\transformed{y} \\ 0
    \end{bmatrix} ,
    \begin{bmatrix}
    0 \\ \transformed{x}\transformed{z} \\ 0
    \end{bmatrix} ,
    \begin{bmatrix}
    0 \\ \transformed{y}\transformed{z} \\ 0
    \end{bmatrix} ,
    \begin{bmatrix}
    0 \\ \transformed{y}^2 \\ 0
    \end{bmatrix} ,
    \begin{bmatrix}
    0 \\ \transformed{x}\transformed{y}^2 \\ 0
    \end{bmatrix} ,
    \begin{bmatrix}
    0 \\ \transformed{z}\transformed{y}^2 \\ 0
    \end{bmatrix} , \\
    &\begin{bmatrix}
    0 \\ 0 \\ 1
    \end{bmatrix} ,
    \begin{bmatrix}
    0 \\ 0 \\ \transformed{x} 
    \end{bmatrix} ,
    \begin{bmatrix}
    0 \\ 0 \\ \transformed{y}
    \end{bmatrix} ,
    \begin{bmatrix}
    0 \\ 0 \\ \transformed{z}
    \end{bmatrix} ,
    \begin{bmatrix}
    0 \\ 0 \\ \transformed{x}\transformed{y}
    \end{bmatrix} ,
    \begin{bmatrix}
    0 \\ 0 \\ \transformed{x}\transformed{z}
    \end{bmatrix} ,
    \begin{bmatrix}
    0 \\ 0 \\ \transformed{y}\transformed{z}
    \end{bmatrix} ,
    \begin{bmatrix}
    0 \\ 0 \\ \transformed{z}^2
    \end{bmatrix} ,
    \begin{bmatrix}
    0 \\ 0 \\ \transformed{x}\transformed{z}^2
    \end{bmatrix} ,
    \begin{bmatrix}
    0 \\ 0 \\ \transformed{y}\transformed{z}^2
    \end{bmatrix}.
\end{align*}

To arrange the location of the flux points, $(\degree + 1)^{\ndim - 1}$ external flux points are located at Gauss-Legendre quadrature points of degree $\degree$ while a total of $\ndim$ sets of $(\degree + 2)(\degree + 1)^{\ndim}$ points are utilized to define the internal flux points.
Such internal flux points are located at Gauss-Legendre quadrature points of degree $\degree - 1$ with a tensor-product arrangement in each direction.
The interested reader is referred to \cite{Liu2006} for further description of such a point arrangement.



\begin{figure}
\centering
\begin{subfigure}{0.45\textwidth}
    \centering
    \includegraphics[width=0.95\textwidth]{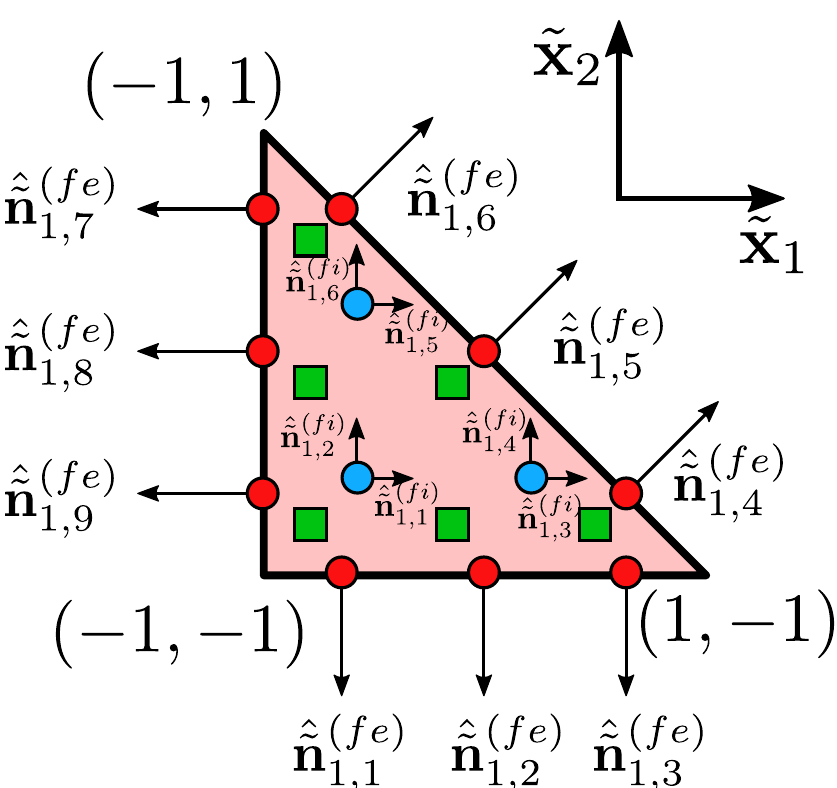}
    \caption{Triangle element ($\degree = 2$)}
    \label{fig:Notation_Triangle_SDRT}
\end{subfigure}
\begin{subfigure}{0.45\textwidth}
    \centering
    \includegraphics[width=0.95\textwidth]{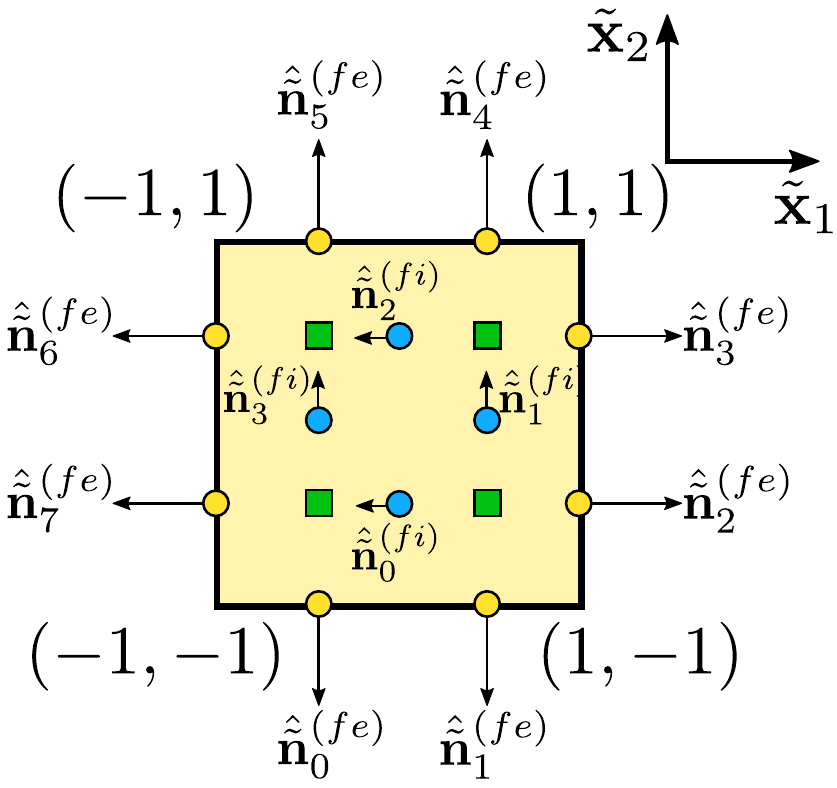}
    \caption{Hexahedra element ($\degree = 1$)}
    \label{fig:Notation_Hexa_SDRT}
\end{subfigure}
\caption{Representation of the different points distribution in the reference domain for the  triangle (a) and quadrilateral (b) element.
Solution points (green squares), external flux points on edges of the triangle (red) and on edges of the square (yellow), internal flux points (blue circles).}
\label{fig:Notation_Element_SDRT_2D}
\end{figure}

For tensor-product elements, there exist $\atflux{N}_{\itype}$ unique physical solution points, each of it has a unique degree of freedom assigned to it.
Solution points are also particularized at Gauss-Legendre quadrature points such that $\atsol{N_{\itype}} = (\degree + 1)^{\ndim}$.
\Fref{fig:Notation_Hexa_SDRT} represents the distribution of such solution and flux points distribution along the reference quadrilateral element.

\begin{remark}
Due to the tensor-product arrangement of the polynomial bases, it is possible to rely on orthonormal or Lagrange polynomials to build the modal basis of tensor-product elements.
This allows to improve the condition number of the Vandermonde matrix associated to the flux nodal basis.
\end{remark}

\subsection{Triangular prismatic elements}

The reference triangular prismatic element is defined as
\begin{equation}
    \transformed{\domain}_{\itype} \in \transformed{x}, \transformed{y} \ge -1 \text{ ; } \transformed{x} + \transformed{y} \le 0 \text{ and } \transformed{z} \in [-1, 1] .
\end{equation}
The RT modal bases of a SDRT scheme of degree $\degree$ for triangular prismatic elements read
\begin{equation}
    \fluxmodalbasis(\transformed{\vect{x}}) \equiv \left[ \polynomial_{\degree}(\transformed{x}_0, \transformed{x}_1)^{2} \oplus (\transformed{x}_0, \transformed{x}_1) \overline{\polynomial}_{\degree}(\transformed{x}_0, \transformed{x}_1) \right] \polynomial_{\degree}(\transformed{x}_2) \times \polynomial_{\degree}(\transformed{x}_0, \transformed{x}_1) \polynomial_{\degree + 1}(\transformed{x}_2) .
\end{equation}
For the sake of clarity, an example of such basis for $\degree = 1$ is provided in the following
\begin{align*}
    \fluxmodalbasis(\transformed{\vect{x}}) \equiv&
    \begin{bmatrix}
    1 \\ 0 \\ 0
    \end{bmatrix} ,
    \begin{bmatrix}
    \transformed{x} \\ 0 \\ 0
    \end{bmatrix} ,
    \begin{bmatrix}
    \transformed{y} \\ 0 \\ 0
    \end{bmatrix} ,
    \begin{bmatrix}
    0 \\ 1 \\ 0
    \end{bmatrix} ,
    \begin{bmatrix}
    0 \\ \transformed{x} \\ 0
    \end{bmatrix} ,
    \begin{bmatrix}
    0 \\ \transformed{y} \\ 0
    \end{bmatrix} ,
    \begin{bmatrix}
    \transformed{x}^2 \\ \transformed{x}\transformed{y} \\ 0
    \end{bmatrix} ,
    \begin{bmatrix}
    \transformed{x}\transformed{y} \\ \transformed{y}^2 \\ 0
    \end{bmatrix} , \\
    &\begin{bmatrix}
    \transformed{z} \\ 0 \\ 0
    \end{bmatrix} ,
    \begin{bmatrix}
    \transformed{x}\transformed{z} \\ 0 \\ 0
    \end{bmatrix} ,
    \begin{bmatrix}
    \transformed{y}\transformed{z} \\ 0 \\ 0
    \end{bmatrix} ,
    \begin{bmatrix}
    0 \\ \transformed{z} \\ 0
    \end{bmatrix} ,
    \begin{bmatrix}
    0 \\ \transformed{x}\transformed{z} \\ 0
    \end{bmatrix} ,
    \begin{bmatrix}
    0 \\ \transformed{y}\transformed{z} \\ 0
    \end{bmatrix} ,
    \begin{bmatrix}
    \transformed{x}^2 \transformed{z} \\ \transformed{x}\transformed{y} \transformed{z} \\ 0
    \end{bmatrix} ,
    \begin{bmatrix}
    \transformed{x}\transformed{y} \transformed{z} \\ \transformed{y}^2 \transformed{z} \\ 0
    \end{bmatrix} ,  \\
    & \begin{bmatrix}
    0 \\ 0 \\ 1
    \end{bmatrix} ,
    \begin{bmatrix}
    0 \\ 0 \\ \transformed{x}
    \end{bmatrix} ,
    \begin{bmatrix}
    0 \\ 0 \\ \transformed{y}
    \end{bmatrix} ,
    \begin{bmatrix}
    0 \\ 0 \\ \transformed{z}
    \end{bmatrix} ,
    \begin{bmatrix}
    0 \\ 0 \\ \transformed{z}\transformed{x}
    \end{bmatrix} ,
    \begin{bmatrix}
    0 \\ 0 \\ \transformed{y}\transformed{z}
    \end{bmatrix} ,
    \begin{bmatrix}
    0 \\ 0 \\ \transformed{z}^2
    \end{bmatrix} ,
    \begin{bmatrix}
    0 \\ 0 \\ \transformed{z}^2\transformed{x}
    \end{bmatrix} ,
    \begin{bmatrix}
    0 \\ 0 \\ \transformed{z}^2\transformed{y}
    \end{bmatrix} .
\end{align*}
To the best of the authors' knowledge no orthonormal RT basis for prismatic elements can be found in the literature, hence the flux Vandermonde matrix presents high condition numbers.


To arrange the external flux points, $(\degree + 1)^{2}$ Gauss-Legendre points are utilized within each of the three quadrilateral faces of the prism, while $\frac{(\degree + 1)(\degree + 2)}{2}$ Williams-Shunn quadrature points are utilized at triangular elements.
Only a single degree of freedom is assigned to external flux points, equal to the normal of the reference-element at the considered flux point

\begin{figure}
\centering
\begin{subfigure}{0.45\textwidth}
    \centering
    \includegraphics[width=0.95\textwidth]{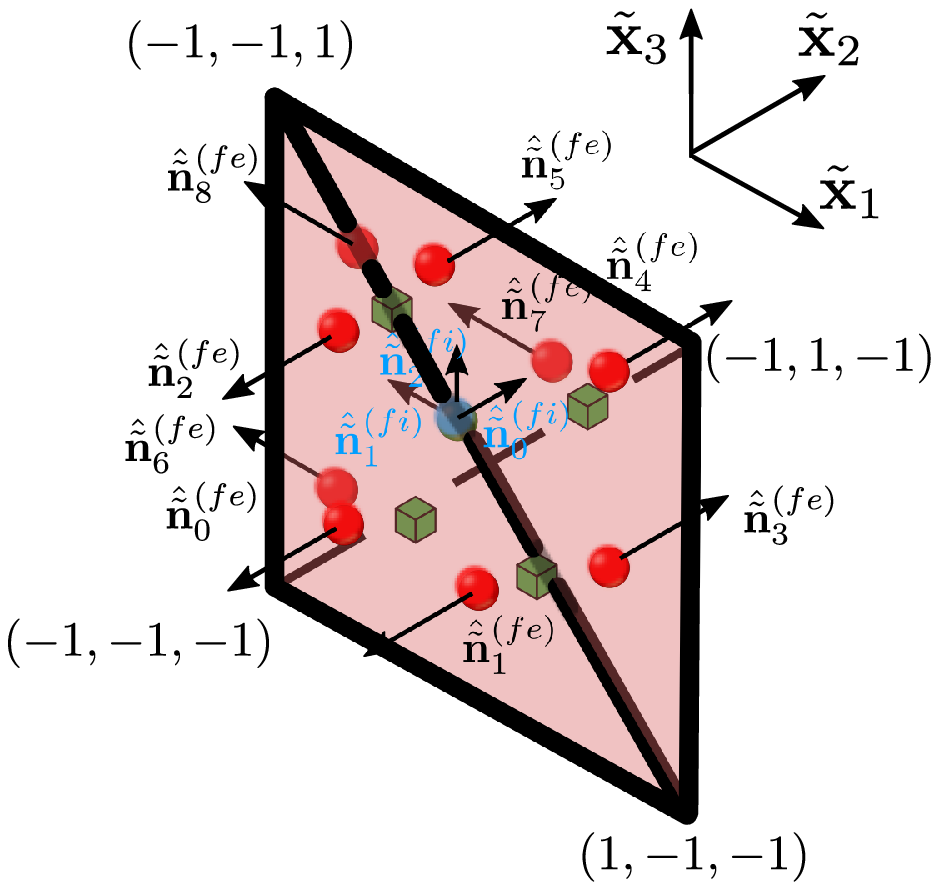}
    \caption{Tetrahedron ($\degree = 1$)}
    \label{fig:Notation_Tetra_SDRT}
\end{subfigure}
\begin{subfigure}{0.45\textwidth}
    \centering
    \includegraphics[width=0.95\textwidth]{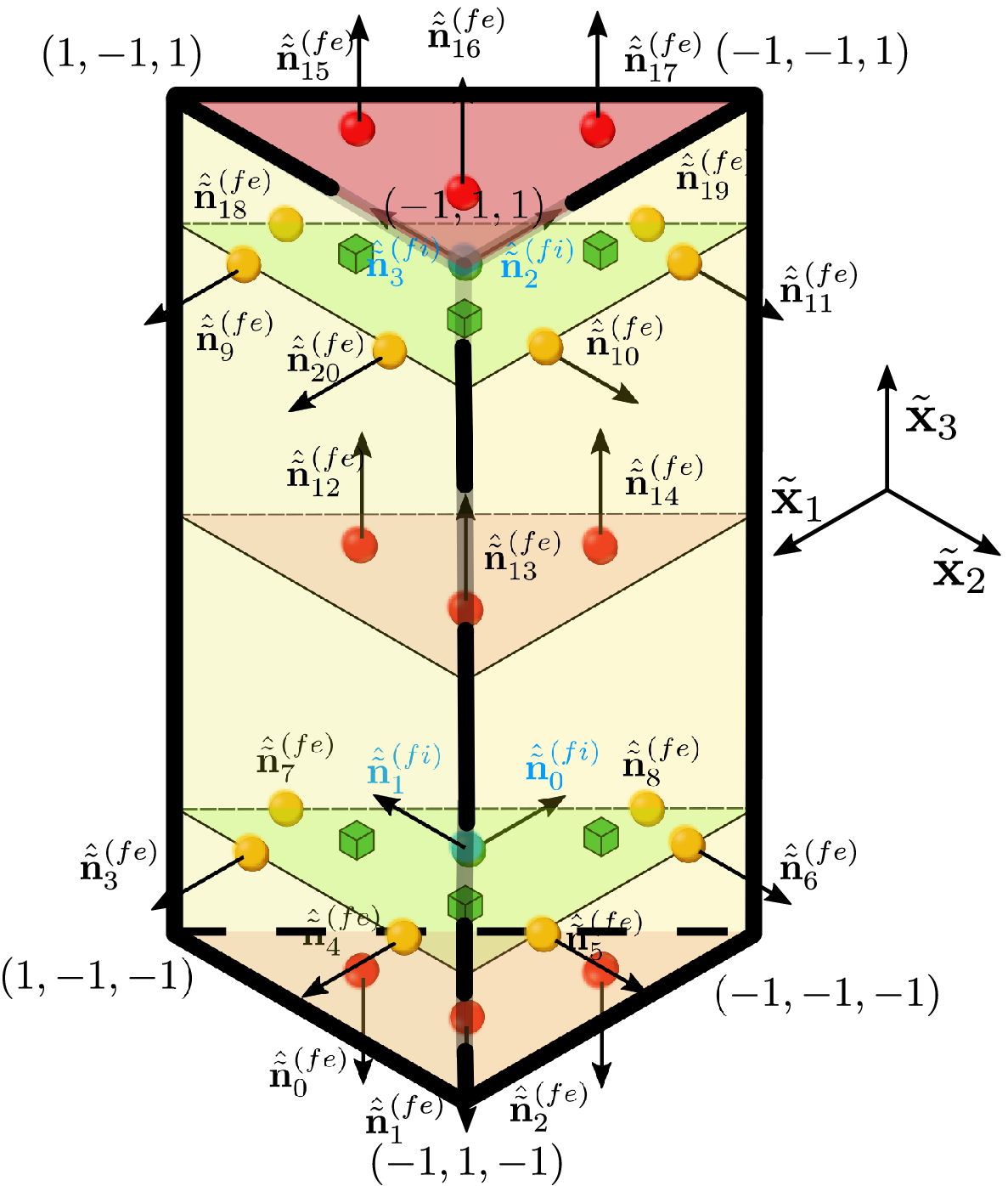}
    \caption{Prism ($\degree = 1$)}
    \label{fig:Notation_Prism_SDRT}
\end{subfigure}
\caption{Representation of the different points distribution in the reference domain for the tetrahedron (a) and prismatic (b) elements. Solution points (green cubes), external flux points (red spheres for triangular faces and yellow spheres for quadrilateral faces), internal flux points (blue spheres).}
\label{fig:Notation_Element_SDRT_3D}
\end{figure}

To arrange internal flux points, two different sets are considered.
The first one is made up of the triangular arrangement of internal flux points duplicated $(p + 1)$ times in the $\transformed{x}_2$ direction using one-dimensional Gauss-Legendre quadrature points in the latter direction.
Two-dimensional degrees of freedom, coinciding with those imposed for internal flux points of triangular elements, are imposed to the aforementioned internal flux points.
The second one consists of $\degree$ duplications of the $\frac{(\degree + 1)(\degree + 2)}{2}$ quadrature points of the triangular faces extruded with $\transformed{x}_2$ coordinate equal to $\degree$ Gauss-Legendre quadrature points.
The latter internal flux points are assigned a unique degree of freedom in the $\transformed{x}_2$ direction.
Such an arrangement may be visualized in \Fref{fig:Notation_Prism_SDRT}.
The first group of internal flux points is represented in blue spheres, while the second one is depicted with red spheres.

Solution points are also particularized at tensor-product configurations of $p + 1$ Gauss-Legendre quadrature points in the extrusion direction and $\frac{(\degree + 1)(\degree + 2)}{2}$ quadrature points of triangles.
This yields a total of $\atsol{N_{\itype}} = \frac{(\degree + 1)^2(\degree + 2)}{2}$ solution points.

\FloatBarrier
\subsection{Summary}
\label{sec:sol_flux_points_tables}

Tables \ref{table:nsolpts_nfpts} and \ref{table:nsolpts_nfpts} represent the number of solution points, external flux points, internal flux points and unique flux points for each element type and as a function of the polynomial degree $\degree$.

\begin{table}[ht]
\centering
\begin{tabular}{@{}lcc@{}}
\toprule
      & $\atsol{N}_{\itype}$ & $\atflux{N}_{\itype}$ \\ \midrule
Triangle & $\frac{(\degree + 1)(\degree + 2)}{2}$ & $(\degree + 1)(\degree + 3)$ \\ \\
Quadrilateral & $(\degree + 1)^2$ & $2(\degree + 1)(\degree + 2)$ \\ \\
Tetrahedron & $\frac{(\degree + 1)(\degree + 2)(\degree + 3)}{6}$ & $\frac{(\degree + 1)(\degree + 2)(\degree + 4)}{2}$ \\ \\
Hexahedron & $(\degree + 1)^3$ & $3(\degree + 1)^2(\degree + 2)$ \\ \\
Prism & $\frac{(\degree + 1)^2(\degree + 2)}{2}$ & $\frac{(\degree + 1)(3 \degree^2 + 12 \degree + 10)}{2}$ \\
         \bottomrule
\end{tabular}
\caption{Number of solution and flux points within each of the elements considered in this work.}
\label{table:nsolpts_nfpts}
\end{table}

\begin{table}[ht]
\centering
\begin{tabular}{@{}lcccc@{}}
\toprule
      & $\atfluxext{N}_{\itype}$ & $\atfluxint{N}_{\itype}$ & $\atfluxintunique{N}_{\itype}$ & $\atfluxint{N}_{\itype}/\atfluxintunique{N}_{\itype}$ \\ \midrule
Triangle & $3(\degree + 1)$ & $\degree(\degree + 1)$ & $\frac{\degree(\degree + 1)}{2}$ & 2 \\ \\
Quadrilateral & $4(\degree + 1)$  & $2\degree(\degree + 1)$ & $2\degree(\degree + 1)$ & 1\\ \\
Tetrahedron & $2(\degree + 1)(\degree + 2)$ & $\frac{\degree(\degree + 1)(\degree + 2)}{2}$ & $\frac{\degree(\degree + 1)(\degree + 2)}{6}$ & 3\\ \\
Hexahedron & $6(\degree + 1)^2$ & $3\degree(\degree + 1)^2$ & $3\degree(\degree + 1)^2$ & 1 \\ \\
Prism & $(\degree + 1)(4 \degree + 5)$ & $\frac{\degree(\degree + 1)(3\degree + 4)}{2}$ & $\frac{\degree (\degree + 1)(2 \degree + 3)}{2}$ & $\frac{3\degree + 4}{2 \degree + 3}$\\
         \bottomrule
\end{tabular}
\caption{Number of external, internal and unique internal flux points within each of the elements considered in this work.}
\label{table:nfpts_extra}
\end{table}







\FloatBarrier
\section{Matrix form of the SDRT method}
\label{sec:sdrt_matrix_form}

As it is the case for the FR formulation \cite{Witherden2014}, the SDRT may be recast in matrix-matrix and matrix-vector operations.
To avoid clashing between the formulation presented in this work and that of the FR formalism in \cite{Witherden2014}, the matrices indices will start at $12$.
To do so, let us start by defining the state matrices of each element type $\iele$ of a given mesh

\begin{equations}{eq:state_matrices}
    \left( \atsol{\statevar}_{\itype} \right)_{\ipol \iele \ivar} &= \atsol{\var}_{\itype \ipol \iele \ivar} \tokendim  \dim \atsol{\statevar}_{\itype} = \atsol{N_{\itype}} \times \nvar |\domain_\itype|, \\
    \left( \atfluxext{\statevar}_{\itype} \right)_{\ipol \iele \ivar} &= \atfluxext{\var}_{\itype \ipol \iele \ivar} \tokendim \dim \atfluxext{\statevar}_{\itype} = \atfluxext{N_{\itype}} \times \nvar |\domain_\itype|, \\
    \left( \atfluxint{\statevar}_{\itype} \right)_{\ipol \iele \ivar} &= \atfluxint{\var}_{\itype \ipol \iele \ivar} \tokendim \dim \atfluxint{\statevar}_{\itype} = \atfluxint{N_{\itype}} \times \nvar |\domain_\itype|, \\
    \left( \atfluxintunique{\statevar}_{\itype} \right)_{\ipol \iele \ivar} &= \atfluxintunique{\var}_{\itype \ipol \iele \ivar} \tokendim \dim \atfluxintunique{\statevar}_{\itype} = \atfluxintunique{N_{\itype}} \times \nvar |\domain_\itype| ,
\end{equations}

The augmented transformation from unique internal flux points to duplicated internal flux points can be carried out using the permutation matrix
\begin{equation}
    \left( \permutationmatrix_{\itype} \right)_{\isigma \ipol } \tokendim \dim \permutationmatrix_{\itype} = \atfluxint{N_{\itype}} \times \atfluxintunique{N_{\itype}} ,
\end{equation}
such that
\begin{equation}
    \atfluxint{\statevar}_{\itype} = \permutationmatrix_{\itype} \atfluxintunique{\statevar}_{\itype} .
\end{equation}
Here, it is worth mentioning that the permutation matrix is the identity matrix for tensor-product elements.

The interpolation procedures may be expressed using the following matrices
\begin{equations}{eq:interpolation_matrices}
    \left( \m^0_{\itype} \right)_{\isigma \ipol} &= \atsol{\nodalbasis}_{\itype \ipol}\left( \atfluxext{\transformed{\vect{x}}}_{\itype \isigma}\right) \tokendim  \dim \m^0_{\itype} = \atfluxext{N_{\itype}} \times \atsol{N_{\itype}}, \\
    \left( \interpfiu_{\itype} \right)_{\isigma \ipol} &= \atsol{\nodalbasis}_{\itype \ipol}\left( \atfluxintunique{\transformed{\vect{x}}}_{\itype \isigma}\right) \tokendim  \dim \interpfiu_{\itype} = \atfluxintunique{N_{\itype}} \times \atsol{N_{\itype}} .
\end{equations}
Hence,
\begin{equations}{eq:interpolation_procedure}
    \atfluxext{\statevar}_{\itype} &= \m^0_{\itype}\atsol{\statevar}_{\itype}, \\
    \atfluxintunique{\statevar}_{\itype} &= \interpfiu_{\itype}\atsol{\statevar}_{\itype}, \\
    \atfluxint{\statevar}_{\itype} &= \interpfiu_{\itype} \permutationmatrix_{\itype} \atsol{\statevar}_{\itype} .
\end{equations}

To compute the transformed gradients at solution points, let
\begin{equations}{eq:gradient_operators}
    \left( \gradfi_{\itype} \right)_{\isigma \ipol} &= \atfluxint{\unit{\transformed{\normal}}}_{\itype \ipol} [\grad \cdot \atfluxint{\fluxnodalbasis}_{\itype \ipol }  ( \atsol{\transformed{\vect{x}}}_{\itype \isigma} )] \tokendim  \dim \gradfi_{\itype} = \ndim \atsol{N_{\itype}} \times \atfluxint{N}_{\itype}, \\
    \left( \gradfe_{\itype} \right)_{\isigma \ipol} &= \atfluxext{\unit{\transformed{\normal}}}_{\itype \ipol} [\grad \cdot \atfluxext{\fluxnodalbasis}_{\itype \ipol }  ( \atsol{\transformed{\vect{x}}}_{\itype \isigma} )] \tokendim  \dim \gradfe_{\itype} = \ndim \atsol{N_{\itype}} \times \atfluxext{N}_{\itype},\\
    \left( \atfluxext{\vect{C}}_{\itype} \right)_{\ipol \iele \ivar} &= \commvar_{\ivar} \atfluxext{\var}_{\itype \ipol \iele \ivar} \tokendim \dim \atfluxext{\statevar}_{\itype} = \atfluxext{N_{\itype}} \times \nvar |\domain_\itype|, \\
    \left( \atsol{\transformed{\statevargrad}}_{\itype} \right)_{\ipol \iele \ivar} &= \atsol{\transformed{\vargrad}}_{\itype \ipol \iele \ivar} \tokendim  \dim \atsol{\transformed{\statevargrad}}_{\itype} = \ndim \atsol{N_{\itype}} \times \nvar |\domain_\itype| .
\end{equations}
Therefore,
\begin{equation}
    \atsol{\transformed{\statevargrad}}_{\itype} = \gradfe_{\itype} \atfluxext{\vect{C}}_{\itype} + \gradfi_{\itype} \atfluxint{\statevar}_{\itype} = \gradfe_{\itype} \atfluxext{\vect{C}}_{\itype} + \gradfi_{\itype} \permutationmatrix_{\itype} \interpfiu_{\itype} \atsol{\statevar}_{\itype} .
\end{equation}

The interpolation of the gradients may be carried out using the following operators
\begin{equations}{eq:interpolate_gradients_matrices}
    \m^{5}_{\itype} &= \text{blockdiag}\left( \m^{0}_{\itype}, \dots, \m^{0}_{\itype} \right) \dim \m^{5}_{\itype} = \ndim \atfluxext{N}_{\itype} \times \ndim \atsol{N}_{\itype} ,\\
    \interpgradfiu_{\itype} &= \text{blockdiag}\left( \interpfiu_{\itype}, \dots, \interpfiu_{\itype} \right) \dim \interpgradfiu_{\itype} = \ndim \atfluxintunique{N}_{\itype} \times \ndim \atsol{N}_{\itype}, \\
    \left( \atsol{\statevargrad}_{\itype} \right)_{\ipol \iele \ivar} &= \atsol{\vargrad}_{\itype \ipol \iele \ivar} \tokendim  \dim \atsol{\statevargrad}_{\itype} = \ndim \atsol{N_{\itype}} \times \nvar |\domain_\itype|, \\
    \left( \atfluxext{\statevargrad}_{\itype} \right)_{\ipol \iele \ivar} &= \atfluxext{\vargrad}_{\itype \ipol \iele \ivar} \tokendim  \dim \atfluxext{\statevargrad}_{\itype} = \ndim \atfluxext{N_{\itype}} \times \nvar |\domain_\itype|, \\
    \left( \atfluxintunique{\statevargrad}_{\itype} \right)_{\ipol \iele \ivar} &= \atfluxintunique{\vargrad}_{\itype \ipol \iele \ivar} \tokendim  \dim \atfluxintunique{\statevargrad}_{\itype} = \ndim \atfluxintunique{N_{\itype}} \times \nvar |\domain_\itype| ,
\end{equations}
hence
\begin{equations}{eq:interpolate_gradients_procedure}
    \atfluxext{\statevargrad}_{\itype} &= \m^{5}_{\itype} \atsol{\statevargrad}_{\itype}, \\
    \atfluxintunique{\statevargrad}_{\itype} &= \interpgradfiu_{\itype} \atsol{\statevargrad}_{\itype} .
\end{equations}

To compute the divergence of the fluxes at the solution points, let us define the following matrices
\begin{equations}{eq:flux_matrices}
    \left( \divfi_{\itype} \right)_{\isigma \ipol} &= 
    \divergence \atfluxint{\fluxnodalbasis}_{\itype \ipol}\left( \atsol{\transformed{\vect{x}}}_{\itype \isigma}\right) \tokendim  \dim \divfi_{\itype} = \atsol{N_{\itype}} \times \atfluxint{N_{\itype}}, \\
    \left( \divfe_{\itype} \right)_{\isigma \ipol} &= \divergence \atfluxext{\fluxnodalbasis}_{\itype \ipol}\left( \atsol{\transformed{\vect{x}}}_{\itype \isigma}\right) \tokendim  \dim \divfe_{\itype} = \atsol{N_{\itype}} \times \atfluxext{N_{\itype}}, \\
    \left( \atfluxnormalext{\transformed{\vect{D}}}_{\itype} \right)_{\ipol \iele \ivar} &= \commflux_{\ivar} \atfluxnormalext{\transformed{\fluxnovect}}_{\itype \ipol \iele \ivar} \tokendim \dim \atfluxnormalext{\transformed{\vect{D}}}_{\itype} = \atfluxext{N_{\itype}} \times \nvar |\domain_\itype|, \\
    \left( \atfluxnormalint{\fluxtransformedmatrix}_{\itype} \right)_{\ipol \iele \ivar} &= \atfluxnormalint{\transformed{\fluxnovect}}_{\itype \ipol \iele \ivar} \tokendim  \atfluxnormalint{\fluxtransformedmatrix}_{\itype} = \atfluxint{N_{\itype}} \times \nvar |\domain_\itype|, \\
    \left( \atfluxintunique{\fluxtransformedmatrix}_{\itype} \right)_{\ipol \iele \ivar} &= \atfluxintunique{\transformed{\flux}}_{\itype \ipol \iele \ivar} \tokendim \dim \atfluxintunique{\fluxtransformedmatrix}_{\itype} = \ndim \atfluxintunique{N_{\itype}} \times \nvar |\domain_\itype|, \\
    \left( \permutationmatrix^{2}_{\itype} \right)_{\isigma \ipol} &= \text{blockdiag}\left( \permutationmatrix_{\itype}, \dots, \permutationmatrix_{\itype} \right) \tokendim \dim \permutationmatrix^{2}_{\itype} = \ndim \atfluxint{N}_{\itype} \times \ndim \atfluxintunique{N}_{\itype},  \\
    \left( \normfi_{\itype} \right)_{\isigma \ipol} &= \atfluxint{\unit{\transformed{\normal}}}_{\itype \inu} \tokendim \dim \normfi_{\itype} = \atfluxint{N}_{\itype} \times \ndim \atfluxintunique{N}_{\itype},  \\
    \left( \atsol{\residualmatrix}_{\itype} \right)_{\ipol \iele \ivar} &= \atsol{(\transformed{\grad} \cdot \transformed{\flux})}_{\itype \ipol \iele \ivar} \tokendim \dim \atsol{\residualmatrix}_{\itype} = \atsol{N_{\itype}} \times \nvar |\domain_\itype|,
\end{equations}
Hence,
\begin{equation}
    \left( \atsol{\residualmatrix}_{\itype} \right)_{\ipol \iele \ivar} = \divfe_{\itype} \atfluxnormalext{\transformed{\vect{D}}}_{\itype} + \divfi_{\itype}   \atfluxnormalint{\fluxtransformedmatrix}_{\itype} = \divfe_{\itype} \atfluxnormalext{\transformed{\vect{D}}}_{\itype} + \divfi_{\itype} \normfi_{\itype} \permutationmatrix^{2}_{\itype}  \atfluxintunique{\fluxtransformedmatrix}_{\itype} .
    \label{eq:residual_sdrt_matrices}
\end{equation}

Numerical solvers may take advantage of the fact that some rows of the product $\normfi_{\itype}  \permutationmatrix^{2}_{\itype}$ may be zero to further optimize the computational performance.
The number of null rows is related to the ratio between the internal flux points the unique internal flux points (see \ref{sec:sol_flux_points_tables} and \ref{sec:rt_basis}).
For example, with tensor-product elements, only one component of the transformed flux at each unique internal flux point is needed to update the residual.
On the other hand, with triangular and tetrahedron elements, all $\ndim$ components of the transformed flux at each unique internal flux point are needed to compute the residual, i.e. no additional optimization is possible.
With prismatic elements, a middle ground between tensor-product and tetrahedron elements is found.
Taking into account such an a priori knowledge of the degrees-of-freedom associated with unique internal flux points allows to further optimize the computational performance by avoiding additional FLOPs related to the multiplication of the metric terms by the fluxes at unique internal flux points, at the cost of additional code complexity.

It is worth mentioning that the FR-SDRT and SDRT methods' equivalence may also be proven using the matrix form of the SDRT method compared to that of the FR formulation described in \cite{Witherden2014}.
This demonstration is left as an exercise for the reader.

\FloatBarrier
\section{Performance Comparison}
\label{sec:perfo_gpu}

In this section a comparison of the perforamance-per-iteration obtained with SDRT, FR-DG and FR-SDRT methods will be presented.
To measure the performance, the parameter PM (measured in $n s$) defined  in \cite{Witherden2014} as
\begin{equation}
    \text{PM} = \frac{T}{N_{\text{iter}} \atsol{N_{\itype}} \nvar |\domain_\itype| N_{\text{RK}}} 10^9 (n s) ,
    \label{eq:performance_per_iter_measure}
\end{equation}
is evaluated in the TGV test case (see \Sref{sec:tgv}).
In the latter equation, $T$ is the wall time needed to perform $N_{\text{iter}}$ time integration steps, $\nvar = 5$ is the number of conserved variables and $N_{\text{RK}} = 5$ is the number of stages of the RK54 time integrator used in the simulations.
The mesh utilized consists of a total of $|\domain_\itype| = 40^3$ tensor-product elements which are subsequently subdivided in $\ncellsinpattern$ cells to obtain meshes made up of prismatic or tetrahedral elements (following the procedure described \Sref{sec:von_neumann_advection}).
This parameter is evaluated four times with $N_{\text{iter}} = 50$ after a the initialization and warmup of the simulation and it is then averaged.
Additionally, the simulations are carried out using constant time step, neglecting possible issues arising from the use of adaptive time stepping.

\Tref{table:perfo_tgv} displays the performance parameter obtained with different type of elements, numerical schemes and polynomial degrees.
The performance was measured using a single NVIDIA V100 GPU.
Since all operations which involve parallel communications in SDRT schemes are equivalent to those found in FR methods, the strong and weak scaling of SDRT and FR schemes should be rather similar.
Hence, the interested reader is referred to the work of \cite{Witherden2015} to further details of the parallel scalability of the SDRT schemes implemented in this work.
The results indicate that FR-DG and FR-SDRT present the same performance, illustrating that FR-SDRT and FR-DG matrices have similar sparsity patterns.
Moreover, the SDRT method shows degraded performance with hexahedral and prismatic elements when compared to FR schemes.
As it can be deduced from \ref{sec:sdrt_matrix_form} this issue is mostly related to additional interpolations that need to be carried out in the SDRT method and also due to the increased size of the divergence correction kernel related to the internal flux points.
Nevertheless, SDRT and FR methods applied tetrahedral elements show very similar performance.
This could be explained due to an improved sparsity pattern of the divergence correction kernel related to the internal flux points with SDRT schemes.
Future studies will analyze the FLOPs associated to each method to further understand the reasoning behind the differences in the performance measure obtained with each element type.

In the study presented herein, the performance is measured on a iteration-per-iteration basis.
Nevertheless, it is worth considering that the time step arising from the adaptive time-stepping method of PyFR yields different values depending on the scheme utilized.
In particular, it was observed (in \Sref{sec:tgv}) that both SDRT and FR-SD schemes present approximately 30\% higher adaptive time step than FR-DG methods in certain configurations.
Such variability of the time step was also observed in \cite{Park2017} when utilizing anti aliasing techniques in FR-DG simulations.
The reason behind such a disparity of adaptive time step values is not yet understood, although it should indicate that unsteady simulations using FR-SDRT and SDRT schemes may be capable of collecting more meaningful data   for a same computational wall time.

\begin{remark}
It is worth recalling that the performance parameter is non-dimensionalized with the number of solution points.
Since, for a given polynomial degree, the element types present a different amount of solution points, the comparison of the data presented in \Tref{table:perfo_tgv} could be slightly misleading.
\end{remark}

\begin{remark}
The performance of the SDRT schemes with prismatic elements could be further optimized by building specific kernels which exploit the number of degrees of freedom associated to each unique internal flux point (as explained \ref{sec:sdrt_matrix_form}).
The results depicted for hexahedral elements already present this optimization.
Since all unique internal flux points contribute to $\ndim$ degrees of freedom in tetrahedral elements, this optimization procedure is not possible for such elements.
\end{remark}

\begin{table}[h]
\centering
\begin{tabular}{@{}cccccccccc@{}}
\toprule
  & \multicolumn{3}{c}{hex} & \multicolumn{3}{c}{pri} & \multicolumn{3}{c}{tet} \\ \cmidrule(lr){2-4} \cmidrule(lr){5-7} \cmidrule(lr){8-10} 
 & SDRT & FR-DG & FR-SDRT & SDRT & FR-DG & FR-SDRT & SDRT & FR-DG & FR-SDRT\\ \midrule 
\multicolumn{1}{l|}{$\mathfrak{p}=1$} & 1.02 & 0.93 & 0.93 & 1.02 & 0.98 & 0.98 & 1.12 & 1.15 & 1.15 \\ 
\multicolumn{1}{l|}{$\mathfrak{p}=2$} & 0.90 & 0.75 & 0.75 & 0.91 & 0.82 & 0.82 & 0.99 & 1.00 & 1.00 \\ 
\multicolumn{1}{l|}{$\mathfrak{p}=3$} & 0.85 & 0.67 & 0.67 & 0.87 & 0.77 & 0.77 & 0.97 & 0.98 & 0.98 \\ 
\multicolumn{1}{l|}{$\mathfrak{p}=4$} & 0.97 & 0.70 & 0.70 & 0.93 & 0.79 & 0.79 & 0.95 & 0.91 & 0.91 \\ 
\end{tabular}
\caption{Performance parameter \Eref{eq:performance_per_iter_measure} values obtained in the TGV test case with a mesh consisting of $40^3 \ncellsinpattern$ cells as a function of the polynomial degree and the element type.
A single NVIDIA V100 GPU was used to evaluate this performance parameter.}
\label{table:perfo_tgv}
\end{table}

\FloatBarrier
\section{Wavenumber Aliasing}
\label{sec:aliasing}


From \Eref{eq:patternvar_diagonalized}, the numerical solution of a linear problem in the asymptotic limit $m \to \infty$ can be written as
\begin{equation}
   \lim_{m \to \infty} \patternvar_{\ipol} (m \timestep) = \eigenvectormatrix_{\ipol SP} \Lambda_{SP}^m \gamma_{SP},
   \label{eq:numerical_solution_m_to_infty}
\end{equation}
provided that the spectral radius is unique.
In the latter equation, the index $SP$ refers to the index of the spectral radius eigenmode.
\Eref{eq:numerical_solution_m_to_infty} indicates that the numerical solution tends to align with the spectral radius eigenvector in the asymptotic limit.
If the spectral radius is not unique, then \Eref{eq:numerical_solution_m_to_infty} would need to be modified by adding the contribution of the eigenmodes whose absolute value of the their associated eigenvalue is equal to the spectral radius.

Using the Matrix-Power-Method theory \cite{VANHAREN2017}, the ratio between the numerical solution norm and the initial condition norm in the asymptotic limit $m \to \infty$ (which measures the dissipation of the numerical solution) is given by
\begin{equation}
    \lim_{m \to \infty} \frac{\norm{\patternvar(m \timestep)}}{\norm{\patternvar(0)}} \approx \frac{ \norm{\eigenvectormatrix_{\ipol SP} \Lambda_{SP}^m \gamma_{SP}}  }{ \norm{\eigenvectormatrix_{\ipol \inu} \Lambda_{\inu}^m \gamma_{\inu}} } .
    \label{eq:dissipation_m_to_infty_vs_0}
\end{equation}
Hence, the dissipation in the asymptotic limit depends on the spectral radius and the initial energy contribution of the spectral radius eigenmode.
If one considers the ratio of dissipation from a given iteration $m$ and a second iteration $m - q$ with $q \ll m$ then
\begin{equation}
    \lim_{m \to \infty} \frac{\norm{\patternvar(m \timestep)}}{\norm{\patternvar((m - q) \timestep)}} \approx \frac{ \norm{\eigenvectormatrix_{\ipol SP} \Lambda_{SP}^m \gamma_{SP}}  }{ \norm{\eigenvectormatrix_{\ipol SP} \Lambda_{SP}^{m - q} \gamma_{SP}} } =  |\Lambda_{SP}|^{q} .
    \label{eq:dissipation_m_to_infty_vs_m_minus_q}
\end{equation}
Therefore, the ratio of dissipation in the asymptotic limit, measured in an iteration-per-iteration basis is just given by the spectral radius.
If the physical eigenmode is equivalent to the spectral radius and its energy contribution is much bigger than that of the other eigenmodes, then Eqs. \ref{eq:dissipation_m_to_infty_vs_0} and \Eref{eq:dissipation_m_to_infty_vs_m_minus_q} are certainly equivalent.
If the physical eigenmode is not coincident with the spectral radius, then the asymptotic dissipation measured from \Eref{eq:dissipation_m_to_infty_vs_0} is heavily perturbed by the fact that the physical eigenmode is not coincident with the spectral radius eigenmode.
Hence, the dissipation will be governed by the initial dissipation of the physical eigenmode until the solution aligns with the spectral radius eigenvector.
Nevertheless, the asymptotic dissipation measured from an iteration-to-iteration perspective will always be given by the spectral radius, independently of its energy contribution to the initial condition.
This might induce wavenumber aliasing issues that had not been previously discussed in the literature.

\begin{remark}
For the sake of brevity, the discussion of the asymptotic behavior of the dispersion is avoided, since similar conclusions to those drawn from the dissipation analysis can be obtained as it was demonstrated in \cite{VANHAREN2017}.
\end{remark}

\subsection{Spectral radius and aliasing}

The numerical solution characteristics of a linear problem in a structured and periodic mesh are governed by the eigenvalues and eigenvectors of matrix $\jacobianrhsprime$ (see \Sref{sec:von_neumann_advection}), which is obtained as
\begin{equation}
    \jacobianrhsprime_{\ipol\inu} = \sum_{\iele \in \stencil} \jacobianrhs_{\ipol \iele \inu} \e^{\imag \vect{\kappa} \cdot \left( \vect{x}_{\iele} - \vect{x}_0 \right) } .
\end{equation}

Let us recall the following wave/wavenumber aliasing identities
\begin{equation}
    \e^{\imag a} = \text{conj}(\e^{-\imag a}) \text{ and } \e^{\imag a} = \e^{\imag( a + 2\pi m)} .
\end{equation}

In one-dimensional problems and when using the structured mesh arrangement illustrated in \Eref{eq:structured_arrangement}, it follows that matrix $\jacobianrhsprime$ is uniquely defined in the range $\kappa h \in [-\pi, \pi]$ due to the aforementioned aliasing identities.
Hence, the spectral radius for any given $\kappa h > \pi$ is related to that obtained with
$ \text{ceil}({\kappa h/\pi}) \pi - \kappa h \text{ if } \text{ceil}({\kappa h/\pi}) \text{ is even }  \text{ and }  \kappa h \mod{ 2 \pi } \text{ if } \text{ceil}({\kappa h/\pi}) \text{ is odd, } $ i.e. the spectral radius is aliased.
Such a fact is also observed in FDM and FVM, although in the latter methods the solution is always characterized by a single eigenmode.
Nevertheless, the energy contribution of each eigenvector in SEM is not uniquely defined in $\kappa h \in [0, \pi]$, i.e. its range is $\kappa h \in [-\infty, \infty)$. 
From a posteriori analyses \cite{Asthana2017}, it can be observed that the physical eigenmode only coincides with the spectral radius for $\kappa h < \pi$.
Therefore, the asymptotic regime dissipation will only be well-characterized by the physical eigenmode for $\kappa h < \pi$.
At $\kappa h = \pi$ the spectral radius is multi-valued, indicating the existence of a branch cut, where the spectral radius ceases to coincide with the physical mode. For other $\kappa h > \pi$, aliasing will appear implying that the numerical dissipation and dispersion are highly unsteady and that the asymptotic dissipation on an iteration-per-iteration basis are reduced since the spectral radius is aliased.
This behavior is similar to that observed in aliased wavenumbers with FDM and FVM.
However, SEM present an additional dissipation mechanism due to the mismatch between the physical and the spectral radius eigenmodes.

\begin{remark}
For the sake of simplicity, the extension of this analysis for two and three-dimensional cases is avoided, since there exist different types of wave aliasing in those configurations.
Nevertheless, as it was demonstrated in \Sref{sec:von_neumann_advection} similar observations regarding the aliasing behavior of the solution may be obtained in two and three-dimensional configurations.
\end{remark}

\bibliographystyle{elsarticle-harv} 
\bibliography{biblio}

\begin{thebibliography}{62}
\expandafter\ifx\csname natexlab\endcsname\relax\def\natexlab#1{#1}\fi
\providecommand{\url}[1]{\texttt{#1}}
\providecommand{\href}[2]{#2}
\providecommand{\path}[1]{#1}
\providecommand{\DOIprefix}{doi:}
\providecommand{\ArXivprefix}{arXiv:}
\providecommand{\URLprefix}{URL: }
\providecommand{\Pubmedprefix}{pmid:}
\providecommand{\doi}[1]{\href{http://dx.doi.org/#1}{\path{#1}}}
\providecommand{\Pubmed}[1]{\href{pmid:#1}{\path{#1}}}
\providecommand{\bibinfo}[2]{#2}
\ifx\xfnm\relax \def\xfnm[#1]{\unskip,\space#1}\fi
\bibitem[{Abe et~al.(2015)Abe, Haga, Nonomura and Fujii}]{Abe2015}
\bibinfo{author}{Abe, Y.}, \bibinfo{author}{Haga, T.},
  \bibinfo{author}{Nonomura, T.}, \bibinfo{author}{Fujii, K.},
  \bibinfo{year}{2015}.
\newblock \bibinfo{title}{On the freestream preservation of high-order
  conservative flux-reconstruction schemes}.
\newblock \bibinfo{journal}{Journal of Computational Physics}
  \bibinfo{volume}{281}, \bibinfo{pages}{28--54}.
\newblock \URLprefix \url{https://doi.org/10.1016/j.jcp.2014.10.011},
  \DOIprefix\doi{10.1016/j.jcp.2014.10.011}.
\bibitem[{Abe et~al.(2018)Abe, Morinaka, Haga, Nonomura, Shibata and
  Miyaji}]{Abe2018}
\bibinfo{author}{Abe, Y.}, \bibinfo{author}{Morinaka, I.},
  \bibinfo{author}{Haga, T.}, \bibinfo{author}{Nonomura, T.},
  \bibinfo{author}{Shibata, H.}, \bibinfo{author}{Miyaji, K.},
  \bibinfo{year}{2018}.
\newblock \bibinfo{title}{Stable, non-dissipative, and conservative
  flux-reconstruction schemes in split forms}.
\newblock \bibinfo{journal}{Journal of Computational Physics}
  \bibinfo{volume}{353}, \bibinfo{pages}{193--227}.
\newblock \URLprefix \url{https://doi.org/10.1016/j.jcp.2017.10.007},
  \DOIprefix\doi{10.1016/j.jcp.2017.10.007}.
\bibitem[{Van~den Abeele et~al.(2008)Van~den Abeele, Lacor and
  Wang}]{Abeele2008}
\bibinfo{author}{Van~den Abeele, K.}, \bibinfo{author}{Lacor, C.},
  \bibinfo{author}{Wang, Z.J.}, \bibinfo{year}{2008}.
\newblock \bibinfo{title}{On the stability and accuracy of the spectral
  difference method}.
\newblock \bibinfo{journal}{Journal of Scientific Computing}
  \bibinfo{volume}{37}, \bibinfo{pages}{162--188}.
\newblock \URLprefix \url{http://dx.doi.org/10.1007/s10915-008-9201-0}.
\bibitem[{Alhawwary and Wang(2018)}]{ALHAWWARY2018}
\bibinfo{author}{Alhawwary, M.}, \bibinfo{author}{Wang, Z.},
  \bibinfo{year}{2018}.
\newblock \bibinfo{title}{Fourier analysis and evaluation of {DG}, {FD} and
  compact difference methods for conservation laws}.
\newblock \bibinfo{journal}{Journal of Computational Physics}
  \bibinfo{volume}{373}, \bibinfo{pages}{835--862}.
\newblock \URLprefix \url{https://doi.org/10.1016/j.jcp.2018.07.018},
  \DOIprefix\doi{10.1016/j.jcp.2018.07.018}.
\bibitem[{Alhawwary and Wang(2020)}]{Alhawwary_2020}
\bibinfo{author}{Alhawwary, M.}, \bibinfo{author}{Wang, Z.},
  \bibinfo{year}{2020}.
\newblock \bibinfo{title}{A combined-mode fourier analysis of {DG} methods for
  linear parabolic problems}.
\newblock \bibinfo{journal}{Journal on Scientific Computing}
  \bibinfo{volume}{42}, \bibinfo{pages}{A3825--A3858}.
\newblock \URLprefix \url{https://doi.org/10.1137\%2F20m1316962},
  \DOIprefix\doi{10.1137/20m1316962}.
\bibitem[{Asthana et~al.(2017)Asthana, Watkins and Jameson}]{Asthana2017}
\bibinfo{author}{Asthana, K.}, \bibinfo{author}{Watkins, J.},
  \bibinfo{author}{Jameson, A.}, \bibinfo{year}{2017}.
\newblock \bibinfo{title}{On consistency and rate of convergence of flux
  reconstruction for time-dependent problems}.
\newblock \bibinfo{journal}{Journal of Computational Physics}
  \bibinfo{volume}{334}, \bibinfo{pages}{367--391}.
\newblock \URLprefix \url{https://doi.org/10.1016/j.jcp.2017.01.008},
  \DOIprefix\doi{10.1016/j.jcp.2017.01.008}.
\bibitem[{Balan et~al.(2011)Balan, May and Sch\"{o}berl}]{Balan2011}
\bibinfo{author}{Balan, A.}, \bibinfo{author}{May, G.},
  \bibinfo{author}{Sch\"{o}berl, J.}, \bibinfo{year}{2011}.
\newblock \bibinfo{title}{A stable spectral difference method for triangles},
  in: \bibinfo{booktitle}{49th {AIAA} Aerospace Sciences Meeting including the
  New Horizons Forum and Aerospace Exposition}, \bibinfo{publisher}{American
  Institute of Aeronautics and Astronautics}.
\newblock \URLprefix \url{https://doi.org/10.2514/6.2011-47},
  \DOIprefix\doi{10.2514/6.2011-47}.
\bibitem[{Bassi and Rebay(1997)}]{Bassi1997}
\bibinfo{author}{Bassi, F.}, \bibinfo{author}{Rebay, S.}, \bibinfo{year}{1997}.
\newblock \bibinfo{title}{High-order accurate discontinuous finite element
  solution of the 2d euler equations}.
\newblock \bibinfo{journal}{Journal of Computational Physics}
  \bibinfo{volume}{138}, \bibinfo{pages}{251--285}.
\newblock \URLprefix \url{https://doi.org/10.1006/jcph.1997.5454},
  \DOIprefix\doi{10.1006/jcph.1997.5454}.
\bibitem[{Bergot and Durufl{\'{e}}(2013)}]{Bergot2013}
\bibinfo{author}{Bergot, M.}, \bibinfo{author}{Durufl{\'{e}}, M.},
  \bibinfo{year}{2013}.
\newblock \bibinfo{title}{Approximation of h(div) with high-order optimal
  finite elements for pyramids, prisms and hexahedra}.
\newblock \bibinfo{journal}{Communications in Computational Physics}
  \bibinfo{volume}{14}, \bibinfo{pages}{1372--1414}.
\newblock \URLprefix \url{https://doi.org/10.4208/cicp.120712.080313a},
  \DOIprefix\doi{10.4208/cicp.120712.080313a}.
\bibitem[{Butcher(2016)}]{Butcher_2016}
\bibinfo{author}{Butcher, J.C.}, \bibinfo{year}{2016}.
\newblock \bibinfo{title}{Numerical Methods for Ordinary Differential
  Equations}.
\newblock \bibinfo{publisher}{John Wiley {\&} Sons, Ltd}.
\newblock \URLprefix \url{https://doi.org/10.1002\%2F9781119121534},
  \DOIprefix\doi{10.1002/9781119121534}.
\bibitem[{Carpenter and Kennedy(1994)}]{Carpenter1994Fourthorder2R}
\bibinfo{author}{Carpenter, M.}, \bibinfo{author}{Kennedy, C.},
  \bibinfo{year}{1994}.
\newblock \bibinfo{title}{Fourth-order 2n-storage runge-kutta schemes}.
\bibitem[{Castonguay et~al.(2011)Castonguay, Vincent and
  Jameson}]{Castonguay2011}
\bibinfo{author}{Castonguay, P.}, \bibinfo{author}{Vincent, P.E.},
  \bibinfo{author}{Jameson, A.}, \bibinfo{year}{2011}.
\newblock \bibinfo{title}{A new class of high-order energy stable flux
  reconstruction schemes for triangular elements}.
\newblock \bibinfo{journal}{Journal of Scientific Computing}
  \bibinfo{volume}{51}, \bibinfo{pages}{224--256}.
\newblock \URLprefix \url{https://doi.org/10.1007/s10915-011-9505-3},
  \DOIprefix\doi{10.1007/s10915-011-9505-3}.
\bibitem[{Cockburn and Shu(1998)}]{Cockburn1998}
\bibinfo{author}{Cockburn, B.}, \bibinfo{author}{Shu, C.W.},
  \bibinfo{year}{1998}.
\newblock \bibinfo{title}{The local discontinuous galerkin method for
  time-dependent convection-diffusion systems}.
\newblock \bibinfo{journal}{{SIAM} Journal on Numerical Analysis}
  \bibinfo{volume}{35}, \bibinfo{pages}{2440--2463}.
\newblock \URLprefix \url{https://doi.org/10.1137/s0036142997316712},
  \DOIprefix\doi{10.1137/s0036142997316712}.
\bibitem[{Cox et~al.(2021)Cox, Trojak, Dzanic, Witherden and Jameson}]{Cox2021}
\bibinfo{author}{Cox, C.}, \bibinfo{author}{Trojak, W.},
  \bibinfo{author}{Dzanic, T.}, \bibinfo{author}{Witherden, F.},
  \bibinfo{author}{Jameson, A.}, \bibinfo{year}{2021}.
\newblock \bibinfo{title}{Accuracy, stability, and performance comparison
  between the spectral difference and flux reconstruction schemes}.
\newblock \bibinfo{journal}{Computers {\&} Fluids} \bibinfo{volume}{221},
  \bibinfo{pages}{104922}.
\newblock \URLprefix \url{https://doi.org/10.1016/j.compfluid.2021.104922},
  \DOIprefix\doi{10.1016/j.compfluid.2021.104922}.
\bibitem[{DeBonis(2013)}]{DeBonis2013}
\bibinfo{author}{DeBonis, J.}, \bibinfo{year}{2013}.
\newblock \bibinfo{title}{Solutions of the taylor-green vortex problem using
  high-resolution explicit finite difference methods}, in:
  \bibinfo{booktitle}{51st {AIAA} Aerospace Sciences Meeting including the New
  Horizons Forum and Aerospace Exposition}, \bibinfo{publisher}{American
  Institute of Aeronautics and Astronautics}.
\newblock \URLprefix \url{https://doi.org/10.2514/6.2013-382},
  \DOIprefix\doi{10.2514/6.2013-382}.
\bibitem[{Frean and Ryan(2019)}]{Frean_2019}
\bibinfo{author}{Frean, D.J.}, \bibinfo{author}{Ryan, J.K.},
  \bibinfo{year}{2019}.
\newblock \bibinfo{title}{Superconvergence and the numerical flux: a study
  using the upwind-biased flux in discontinuous galerkin methods}.
\newblock \bibinfo{journal}{Commun. Appl. Math. Comput.} \bibinfo{volume}{2},
  \bibinfo{pages}{461--486}.
\newblock \URLprefix \url{https://doi.org/10.1007\%2Fs42967-019-00049-2},
  \DOIprefix\doi{10.1007/s42967-019-00049-2}.
\bibitem[{Glaubitz et~al.(2017)Glaubitz, \"{O}ffner and Sonar}]{Glaubitz2017}
\bibinfo{author}{Glaubitz, J.}, \bibinfo{author}{\"{O}ffner, P.},
  \bibinfo{author}{Sonar, T.}, \bibinfo{year}{2017}.
\newblock \bibinfo{title}{Application of modal filtering to a spectral
  difference method}.
\newblock \bibinfo{journal}{Mathematics of Computation} \bibinfo{volume}{87},
  \bibinfo{pages}{175--207}.
\newblock \URLprefix \url{https://doi.org/10.1090/mcom/3257},
  \DOIprefix\doi{10.1090/mcom/3257}.
\bibitem[{Guo et~al.(2013)Guo, Zhong and Qiu}]{Guo_2013}
\bibinfo{author}{Guo, W.}, \bibinfo{author}{Zhong, X.}, \bibinfo{author}{Qiu,
  J.M.}, \bibinfo{year}{2013}.
\newblock \bibinfo{title}{Superconvergence of discontinuous galerkin and local
  discontinuous galerkin methods: Eigen-structure analysis based on fourier
  approach}.
\newblock \bibinfo{journal}{Journal of Computational Physics}
  \bibinfo{volume}{235}, \bibinfo{pages}{458--485}.
\newblock \URLprefix \url{https://doi.org/10.1016\%2Fj.jcp.2012.10.020},
  \DOIprefix\doi{10.1016/j.jcp.2012.10.020}.
\bibitem[{Hesthaven(2008)}]{hesthaven2008nodal}
\bibinfo{author}{Hesthaven, J.}, \bibinfo{year}{2008}.
\newblock \bibinfo{title}{Nodal discontinuous Galerkin methods : algorithms,
  analysis, and applications}.
\newblock \bibinfo{publisher}{Springer}, \bibinfo{address}{New York}.
\bibitem[{Huynh(2007)}]{Huynh2007}
\bibinfo{author}{Huynh, H.T.}, \bibinfo{year}{2007}.
\newblock \bibinfo{title}{A flux reconstruction approach to high-order schemes
  including discontinuous galerkin methods}, in: \bibinfo{booktitle}{18th
  {AIAA} Computational Fluid Dynamics Conference}, \bibinfo{publisher}{American
  Institute of Aeronautics and Astronautics}.
\newblock \URLprefix \url{https://doi.org/10.2514/6.2007-4079},
  \DOIprefix\doi{10.2514/6.2007-4079}.
\bibitem[{Huynh(2020)}]{Huynh2020}
\bibinfo{author}{Huynh, H.T.}, \bibinfo{year}{2020}.
\newblock \bibinfo{title}{Discontinuous galerkin via interpolation: The direct
  flux reconstruction method}.
\newblock \bibinfo{journal}{Journal of Scientific Computing}
  \bibinfo{volume}{82}.
\newblock \URLprefix \url{https://doi.org/10.1007/s10915-020-01175-3},
  \DOIprefix\doi{10.1007/s10915-020-01175-3}.
\bibitem[{Iyer et~al.(2019)Iyer, Witherden, Chernyshenko and
  Vincent}]{Iyer2019}
\bibinfo{author}{Iyer, A.S.}, \bibinfo{author}{Witherden, F.D.},
  \bibinfo{author}{Chernyshenko, S.I.}, \bibinfo{author}{Vincent, P.E.},
  \bibinfo{year}{2019}.
\newblock \bibinfo{title}{Identifying eigenmodes of averaged small-amplitude
  perturbations to turbulent channel flow}.
\newblock \bibinfo{journal}{Journal of Fluid Mechanics} \bibinfo{volume}{875},
  \bibinfo{pages}{758--780}.
\newblock \URLprefix \url{https://doi.org/10.1017/jfm.2019.520},
  \DOIprefix\doi{10.1017/jfm.2019.520}.
\bibitem[{Jameson(2010)}]{Jameson2010}
\bibinfo{author}{Jameson, A.}, \bibinfo{year}{2010}.
\newblock \bibinfo{title}{A proof of the stability of the spectral difference
  method for all orders of accuracy}.
\newblock \bibinfo{journal}{Journal of Scientific Computing}
  \bibinfo{volume}{45}, \bibinfo{pages}{348--358}.
\newblock \URLprefix \url{https://doi.org/10.1007/s10915-009-9339-4},
  \DOIprefix\doi{10.1007/s10915-009-9339-4}.
\bibitem[{Karniadakis and Sherwin(2005)}]{Karniadakis_2005}
\bibinfo{author}{Karniadakis, G.}, \bibinfo{author}{Sherwin, S.},
  \bibinfo{year}{2005}.
\newblock \bibinfo{title}{Spectral/hp Element Methods for Computational Fluid
  Dynamics}.
\newblock \bibinfo{publisher}{Oxford University Press}.
\newblock \URLprefix
  \url{https://doi.org/10.1093\%2Facpro\%3Aoso\%2F9780198528692.001.0001},
  \DOIprefix\doi{10.1093/acprof:oso/9780198528692.001.0001}.
\bibitem[{Kopriva(1996)}]{Kopriva1996}
\bibinfo{author}{Kopriva, D.A.}, \bibinfo{year}{1996}.
\newblock \bibinfo{title}{A conservative staggered-grid chebyshev multidomain
  method for compressible flows. {II}. a semi-structured method}.
\newblock \bibinfo{journal}{Journal of Computational Physics}
  \bibinfo{volume}{128}, \bibinfo{pages}{475--488}.
\newblock \URLprefix \url{https://doi.org/10.1006/jcph.1996.0225},
  \DOIprefix\doi{10.1006/jcph.1996.0225}.
\bibitem[{Kopriva(2006)}]{Kopriva2006}
\bibinfo{author}{Kopriva, D.A.}, \bibinfo{year}{2006}.
\newblock \bibinfo{title}{Metric identities and the discontinuous spectral
  element method on curvilinear meshes}.
\newblock \bibinfo{journal}{Journal of Scientific Computing}
  \bibinfo{volume}{26}, \bibinfo{pages}{301--327}.
\newblock \URLprefix \url{https://doi.org/10.1007/s10915-005-9070-8},
  \DOIprefix\doi{10.1007/s10915-005-9070-8}.
\bibitem[{Li et~al.(2019)Li, Qiu, Liang, Sprague, Xu and Garris}]{Li2019}
\bibinfo{author}{Li, M.}, \bibinfo{author}{Qiu, Z.}, \bibinfo{author}{Liang,
  C.}, \bibinfo{author}{Sprague, M.}, \bibinfo{author}{Xu, M.},
  \bibinfo{author}{Garris, C.A.}, \bibinfo{year}{2019}.
\newblock \bibinfo{title}{A new high-order spectral difference method for
  simulating viscous flows on unstructured grids with mixed-element meshes}.
\newblock \bibinfo{journal}{Computers {\&} Fluids} \bibinfo{volume}{184},
  \bibinfo{pages}{187--198}.
\newblock \URLprefix \url{https://doi.org/10.1016/j.compfluid.2019.03.010},
  \DOIprefix\doi{10.1016/j.compfluid.2019.03.010}.
\bibitem[{Liu et~al.(2006)Liu, Vinokur and Wang}]{Liu2006}
\bibinfo{author}{Liu, Y.}, \bibinfo{author}{Vinokur, M.},
  \bibinfo{author}{Wang, Z.}, \bibinfo{year}{2006}.
\newblock \bibinfo{title}{Spectral difference method for unstructured grids i:
  Basic formulation}.
\newblock \bibinfo{journal}{Journal of Computational Physics}
  \bibinfo{volume}{216}, \bibinfo{pages}{780--801}.
\newblock \URLprefix \url{https://doi.org/10.1016/j.jcp.2006.01.024},
  \DOIprefix\doi{10.1016/j.jcp.2006.01.024}.
\bibitem[{Lodato(2019)}]{Lodato2019}
\bibinfo{author}{Lodato, G.}, \bibinfo{year}{2019}.
\newblock \bibinfo{title}{Characteristic modal shock detection for
  discontinuous finite element methods}.
\newblock \bibinfo{journal}{Computers {\&} Fluids} \bibinfo{volume}{179},
  \bibinfo{pages}{309--333}.
\newblock \URLprefix \url{https://doi.org/10.1016/j.compfluid.2018.11.008},
  \DOIprefix\doi{10.1016/j.compfluid.2018.11.008}.
\bibitem[{Lodato et~al.(2016)Lodato, Vervisch and Clavin}]{Lodato2016}
\bibinfo{author}{Lodato, G.}, \bibinfo{author}{Vervisch, L.},
  \bibinfo{author}{Clavin, P.}, \bibinfo{year}{2016}.
\newblock \bibinfo{title}{Direct numerical simulation of shock wavy-wall
  interaction: analysis of cellular shock structures and flow patterns}.
\newblock \bibinfo{journal}{Journal of Fluid Mechanics} \bibinfo{volume}{789},
  \bibinfo{pages}{221--258}.
\newblock \URLprefix \url{https://doi.org/10.1017/jfm.2015.731},
  \DOIprefix\doi{10.1017/jfm.2015.731}.
\bibitem[{Loppi et~al.(2018)Loppi, Witherden, Jameson and Vincent}]{Loppi2018}
\bibinfo{author}{Loppi, N.}, \bibinfo{author}{Witherden, F.},
  \bibinfo{author}{Jameson, A.}, \bibinfo{author}{Vincent, P.},
  \bibinfo{year}{2018}.
\newblock \bibinfo{title}{A high-order cross-platform incompressible
  navier{\textendash}stokes solver via artificial compressibility with
  application to a turbulent jet}.
\newblock \bibinfo{journal}{Computer Physics Communications}
  \bibinfo{volume}{233}, \bibinfo{pages}{193--205}.
\newblock \URLprefix \url{https://doi.org/10.1016/j.cpc.2018.06.016},
  \DOIprefix\doi{10.1016/j.cpc.2018.06.016}.
\bibitem[{Manzanero et~al.(2020)Manzanero, Ferrer, Rubio and
  Valero}]{Manzanero2020}
\bibinfo{author}{Manzanero, J.}, \bibinfo{author}{Ferrer, E.},
  \bibinfo{author}{Rubio, G.}, \bibinfo{author}{Valero, E.},
  \bibinfo{year}{2020}.
\newblock \bibinfo{title}{Design of a smagorinsky spectral vanishing viscosity
  turbulence model for discontinuous galerkin methods}.
\newblock \bibinfo{journal}{Computers {\&} Fluids} \bibinfo{volume}{200},
  \bibinfo{pages}{104440}.
\newblock \URLprefix \url{https://doi.org/10.1016/j.compfluid.2020.104440},
  \DOIprefix\doi{10.1016/j.compfluid.2020.104440}.
\bibitem[{Mengaldo et~al.(2015)Mengaldo, Grazia, Vincent and
  Sherwin}]{Mengaldo2015}
\bibinfo{author}{Mengaldo, G.}, \bibinfo{author}{Grazia, D.D.},
  \bibinfo{author}{Vincent, P.E.}, \bibinfo{author}{Sherwin, S.J.},
  \bibinfo{year}{2015}.
\newblock \bibinfo{title}{On the connections between discontinuous galerkin and
  flux reconstruction schemes: Extension to curvilinear meshes}.
\newblock \bibinfo{journal}{Journal of Scientific Computing}
  \bibinfo{volume}{67}, \bibinfo{pages}{1272--1292}.
\newblock \URLprefix \url{https://doi.org/10.1007/s10915-015-0119-z},
  \DOIprefix\doi{10.1007/s10915-015-0119-z}.
\bibitem[{Mozolevski and Valmorbida(2017)}]{Mozolevski2017}
\bibinfo{author}{Mozolevski, I.}, \bibinfo{author}{Valmorbida, E.L.},
  \bibinfo{year}{2017}.
\newblock \bibinfo{title}{Efficient equilibrated flux reconstruction in~high
  order raviart-thomas space for discontinuous galerkin methods}, in:
  \bibinfo{booktitle}{Lecture Notes in Computational Science and Engineering}.
  \bibinfo{publisher}{Springer International Publishing}, pp.
  \bibinfo{pages}{467--479}.
\newblock \URLprefix \url{https://doi.org/10.1007/978-3-319-65870-4\_33},
  \DOIprefix\doi{10.1007/978-3-319-65870-4\_33}.
\bibitem[{Navah et~al.(2020)Navah, de~la Llave~Plata and
  Couaillier}]{Navah2020}
\bibinfo{author}{Navah, F.}, \bibinfo{author}{de~la Llave~Plata, M.},
  \bibinfo{author}{Couaillier, V.}, \bibinfo{year}{2020}.
\newblock \bibinfo{title}{A high-order multiscale approach to turbulence for
  compact nodal schemes}.
\newblock \bibinfo{journal}{Computer Methods in Applied Mechanics and
  Engineering} \bibinfo{volume}{363}, \bibinfo{pages}{112885}.
\newblock \URLprefix \url{https://doi.org/10.1016/j.cma.2020.112885},
  \DOIprefix\doi{10.1016/j.cma.2020.112885}.
\bibitem[{Olson et~al.(2014)Olson, Shaw, Shi, Pierre and Parker}]{Olson2014}
\bibinfo{author}{Olson, B.J.}, \bibinfo{author}{Shaw, S.W.},
  \bibinfo{author}{Shi, C.}, \bibinfo{author}{Pierre, C.},
  \bibinfo{author}{Parker, R.G.}, \bibinfo{year}{2014}.
\newblock \bibinfo{title}{Circulant matrices and their application to vibration
  analysis}.
\newblock \bibinfo{journal}{Applied Mechanics Reviews} \bibinfo{volume}{66}.
\newblock \URLprefix \url{https://doi.org/10.1115/1.4027722},
  \DOIprefix\doi{10.1115/1.4027722}.
\bibitem[{Park et~al.(2017)Park, Witherden and Vincent}]{Park2017}
\bibinfo{author}{Park, J.S.}, \bibinfo{author}{Witherden, F.D.},
  \bibinfo{author}{Vincent, P.E.}, \bibinfo{year}{2017}.
\newblock \bibinfo{title}{High-order implicit large-eddy simulations of flow
  over a {NACA}0021 aerofoil}.
\newblock \bibinfo{journal}{{AIAA} Journal} \bibinfo{volume}{55},
  \bibinfo{pages}{2186--2197}.
\newblock \URLprefix \url{https://doi.org/10.2514/1.j055304},
  \DOIprefix\doi{10.2514/1.j055304}.
\bibitem[{Peraire and Persson(2008)}]{Peraire2008}
\bibinfo{author}{Peraire, J.}, \bibinfo{author}{Persson, P.O.},
  \bibinfo{year}{2008}.
\newblock \bibinfo{title}{The compact discontinuous galerkin ({CDG}) method for
  elliptic problems}.
\newblock \bibinfo{journal}{{SIAM} Journal on Scientific Computing}
  \bibinfo{volume}{30}, \bibinfo{pages}{1806--1824}.
\newblock \URLprefix \url{https://doi.org/10.1137/070685518},
  \DOIprefix\doi{10.1137/070685518}.
\bibitem[{Pereira and Vermeire(2020a)}]{Pereira2020}
\bibinfo{author}{Pereira, C.A.}, \bibinfo{author}{Vermeire, B.C.},
  \bibinfo{year}{2020}a.
\newblock \bibinfo{title}{Fully-discrete analysis of high-order spatial
  discretizations with optimal explicit runge{\textendash}kutta methods}.
\newblock \bibinfo{journal}{Journal of Scientific Computing}
  \bibinfo{volume}{83}.
\newblock \URLprefix \url{https://doi.org/10.1007/s10915-020-01243-8},
  \DOIprefix\doi{10.1007/s10915-020-01243-8}.
\bibitem[{Pereira and Vermeire(2020b)}]{Pereira2020elements}
\bibinfo{author}{Pereira, C.A.}, \bibinfo{author}{Vermeire, B.C.},
  \bibinfo{year}{2020}b.
\newblock \bibinfo{title}{Spectral properties of high-order element types for
  implicit large eddy simulation}.
\newblock \bibinfo{journal}{Journal of Scientific Computing}
  \bibinfo{volume}{85}.
\newblock \URLprefix \url{https://doi.org/10.1007/s10915-020-01329-3},
  \DOIprefix\doi{10.1007/s10915-020-01329-3}.
\bibitem[{van Rees et~al.(2011)van Rees, Leonard, Pullin and
  Koumoutsakos}]{vanRees2011}
\bibinfo{author}{van Rees, W.M.}, \bibinfo{author}{Leonard, A.},
  \bibinfo{author}{Pullin, D.}, \bibinfo{author}{Koumoutsakos, P.},
  \bibinfo{year}{2011}.
\newblock \bibinfo{title}{A comparison of vortex and pseudo-spectral methods
  for the simulation of periodic vortical flows at high reynolds numbers}.
\newblock \bibinfo{journal}{Journal of Computational Physics}
  \bibinfo{volume}{230}, \bibinfo{pages}{2794--2805}.
\newblock \URLprefix \url{https://doi.org/10.1016/j.jcp.2010.11.031},
  \DOIprefix\doi{10.1016/j.jcp.2010.11.031}.
\bibitem[{Shu(2016)}]{Shu2016}
\bibinfo{author}{Shu, C.W.}, \bibinfo{year}{2016}.
\newblock \bibinfo{title}{High order {WENO} and {DG} methods for time-dependent
  convection-dominated {PDEs}: A brief survey of several recent developments}.
\newblock \bibinfo{journal}{Journal of Computational Physics}
  \bibinfo{volume}{316}, \bibinfo{pages}{598--613}.
\newblock \URLprefix \url{https://doi.org/10.1016/j.jcp.2016.04.030},
  \DOIprefix\doi{10.1016/j.jcp.2016.04.030}.
\bibitem[{Spiegel et~al.(2015a)Spiegel, Huynh and
  DeBonis}]{Spiegel2015_aliasing}
\bibinfo{author}{Spiegel, S.C.}, \bibinfo{author}{Huynh, H.},
  \bibinfo{author}{DeBonis, J.R.}, \bibinfo{year}{2015}a.
\newblock \bibinfo{title}{De-aliasing through over-integration applied to the
  flux reconstruction and discontinuous galerkin methods}, in:
  \bibinfo{booktitle}{22nd {AIAA} Computational Fluid Dynamics Conference},
  \bibinfo{publisher}{American Institute of Aeronautics and Astronautics}.
\newblock \URLprefix \url{https://doi.org/10.2514/6.2015-2744},
  \DOIprefix\doi{10.2514/6.2015-2744}.
\bibitem[{Spiegel et~al.(2015b)Spiegel, Huynh and DeBonis}]{Spiegel2015}
\bibinfo{author}{Spiegel, S.C.}, \bibinfo{author}{Huynh, H.},
  \bibinfo{author}{DeBonis, J.R.}, \bibinfo{year}{2015}b.
\newblock \bibinfo{title}{A survey of the isentropic euler vortex problem using
  high-order methods}, in: \bibinfo{booktitle}{22nd {AIAA} Computational Fluid
  Dynamics Conference}, \bibinfo{publisher}{American Institute of Aeronautics
  and Astronautics}.
\newblock \URLprefix \url{https://doi.org/10.2514/6.2015-2444},
  \DOIprefix\doi{10.2514/6.2015-2444}.
\bibitem[{Toro(1999)}]{Toro1999}
\bibinfo{author}{Toro, E.F.}, \bibinfo{year}{1999}.
\newblock \bibinfo{title}{The {HLL} and {HLLC} riemann solvers}, in:
  \bibinfo{booktitle}{Riemann Solvers and Numerical Methods for Fluid
  Dynamics}. \bibinfo{publisher}{Springer Berlin Heidelberg}, pp.
  \bibinfo{pages}{315--339}.
\newblock \URLprefix \url{https://doi.org/10.1007/978-3-662-03915-1\_10},
  \DOIprefix\doi{10.1007/978-3-662-03915-1\_10}.
\bibitem[{Trojak et~al.(2020)Trojak, Watson, Scillitoe and Tucker}]{Trojak2020}
\bibinfo{author}{Trojak, W.}, \bibinfo{author}{Watson, R.},
  \bibinfo{author}{Scillitoe, A.}, \bibinfo{author}{Tucker, P.G.},
  \bibinfo{year}{2020}.
\newblock \bibinfo{title}{Effect of mesh quality on flux reconstruction in
  multi-dimensions}.
\newblock \bibinfo{journal}{Journal of Scientific Computing}
  \bibinfo{volume}{82}.
\newblock \URLprefix \url{https://doi.org/10.1007/s10915-020-01184-2},
  \DOIprefix\doi{10.1007/s10915-020-01184-2}.
\bibitem[{Vanharen et~al.(2017)Vanharen, Puigt, Vasseur, Boussuge and
  Sagaut}]{VANHAREN2017}
\bibinfo{author}{Vanharen, J.}, \bibinfo{author}{Puigt, G.},
  \bibinfo{author}{Vasseur, X.}, \bibinfo{author}{Boussuge, J.F.},
  \bibinfo{author}{Sagaut, P.}, \bibinfo{year}{2017}.
\newblock \bibinfo{title}{Revisiting the spectral analysis for high-order
  spectral discontinuous methods}.
\newblock \bibinfo{journal}{Journal of Computational Physics}
  \bibinfo{volume}{337}, \bibinfo{pages}{379--402}.
\newblock \URLprefix \url{https://doi.org/10.1016/j.jcp.2017.02.043},
  \DOIprefix\doi{10.1016/j.jcp.2017.02.043}.
\bibitem[{Veilleux(2021)}]{Adele2021thesis}
\bibinfo{author}{Veilleux, A.}, \bibinfo{year}{2021}.
\newblock \bibinfo{title}{Extension of the Spectral Difference method to
  simplex cells and hybrid grids}.
\newblock Ph.D. thesis. CERFACS.
\bibitem[{Vermeire and Vincent(2017)}]{Vermeire2017}
\bibinfo{author}{Vermeire, B.}, \bibinfo{author}{Vincent, P.},
  \bibinfo{year}{2017}.
\newblock \bibinfo{title}{On the behaviour of fully-discrete flux
  reconstruction schemes}.
\newblock \bibinfo{journal}{Computer Methods in Applied Mechanics and
  Engineering} \bibinfo{volume}{315}, \bibinfo{pages}{1053--1079}.
\newblock \URLprefix \url{https://doi.org/10.1016/j.cma.2016.11.019},
  \DOIprefix\doi{10.1016/j.cma.2016.11.019}.
\bibitem[{Vincent et~al.(2011)Vincent, Castonguay and Jameson}]{VINCENT2011}
\bibinfo{author}{Vincent, P.}, \bibinfo{author}{Castonguay, P.},
  \bibinfo{author}{Jameson, A.}, \bibinfo{year}{2011}.
\newblock \bibinfo{title}{Insights from von neumann analysis of high-order flux
  reconstruction schemes}.
\newblock \bibinfo{journal}{Journal of Computational Physics}
  \bibinfo{volume}{230}, \bibinfo{pages}{8134--8154}.
\newblock \URLprefix \url{https://doi.org/10.1016/j.jcp.2011.07.013},
  \DOIprefix\doi{10.1016/j.jcp.2011.07.013}.
\bibitem[{Vincent et~al.(2010)Vincent, Castonguay and Jameson}]{Vincent2010}
\bibinfo{author}{Vincent, P.E.}, \bibinfo{author}{Castonguay, P.},
  \bibinfo{author}{Jameson, A.}, \bibinfo{year}{2010}.
\newblock \bibinfo{title}{A new class of high-order energy stable flux
  reconstruction schemes}.
\newblock \bibinfo{journal}{Journal of Scientific Computing}
  \bibinfo{volume}{47}, \bibinfo{pages}{50--72}.
\newblock \URLprefix \url{https://doi.org/10.1007/s10915-010-9420-z},
  \DOIprefix\doi{10.1007/s10915-010-9420-z}.
\bibitem[{Wang et~al.(2013)Wang, Fidkowski, Abgrall, Bassi, Caraeni, Cary,
  Deconinck, Hartmann, Hillewaert, Huynh, Kroll, May, Persson, van Leer and
  Visbal}]{Wang2013}
\bibinfo{author}{Wang, Z.}, \bibinfo{author}{Fidkowski, K.},
  \bibinfo{author}{Abgrall, R.}, \bibinfo{author}{Bassi, F.},
  \bibinfo{author}{Caraeni, D.}, \bibinfo{author}{Cary, A.},
  \bibinfo{author}{Deconinck, H.}, \bibinfo{author}{Hartmann, R.},
  \bibinfo{author}{Hillewaert, K.}, \bibinfo{author}{Huynh, H.},
  \bibinfo{author}{Kroll, N.}, \bibinfo{author}{May, G.},
  \bibinfo{author}{Persson, P.O.}, \bibinfo{author}{van Leer, B.},
  \bibinfo{author}{Visbal, M.}, \bibinfo{year}{2013}.
\newblock \bibinfo{title}{High-order {CFD} methods: current status and
  perspective}.
\newblock \bibinfo{journal}{International Journal for Numerical Methods in
  Fluids} \bibinfo{volume}{72}, \bibinfo{pages}{811--845}.
\newblock \URLprefix \url{https://doi.org/10.1002/fld.3767},
  \DOIprefix\doi{10.1002/fld.3767}.
\bibitem[{Wang and Gao(2009)}]{Wang2009}
\bibinfo{author}{Wang, Z.}, \bibinfo{author}{Gao, H.}, \bibinfo{year}{2009}.
\newblock \bibinfo{title}{A unifying lifting collocation penalty formulation
  including the discontinuous galerkin, spectral volume/difference methods for
  conservation laws on mixed grids}.
\newblock \bibinfo{journal}{Journal of Computational Physics}
  \bibinfo{volume}{228}, \bibinfo{pages}{8161--8186}.
\newblock \URLprefix \url{https://doi.org/10.1016/j.jcp.2009.07.036},
  \DOIprefix\doi{10.1016/j.jcp.2009.07.036}.
\bibitem[{Wang and Liu(2002)}]{Wang2002}
\bibinfo{author}{Wang, Z.}, \bibinfo{author}{Liu, Y.}, \bibinfo{year}{2002}.
\newblock \bibinfo{title}{Spectral (finite) volume method for conservation laws
  on unstructured grids}.
\newblock \bibinfo{journal}{Journal of Computational Physics}
  \bibinfo{volume}{179}, \bibinfo{pages}{665--697}.
\newblock \URLprefix \url{https://doi.org/10.1006/jcph.2002.7082},
  \DOIprefix\doi{10.1006/jcph.2002.7082}.
\bibitem[{Watkins et~al.(2016)Watkins, Asthana and Jameson}]{Watkins2016}
\bibinfo{author}{Watkins, J.}, \bibinfo{author}{Asthana, K.},
  \bibinfo{author}{Jameson, A.}, \bibinfo{year}{2016}.
\newblock \bibinfo{title}{A numerical analysis of the nodal discontinuous
  galerkin scheme via flux reconstruction for the advection-diffusion
  equation}.
\newblock \bibinfo{journal}{Computers {\&} Fluids} \bibinfo{volume}{139},
  \bibinfo{pages}{233--247}.
\newblock \URLprefix \url{https://doi.org/10.1016/j.compfluid.2016.09.013},
  \DOIprefix\doi{10.1016/j.compfluid.2016.09.013}.
\bibitem[{Williams et~al.(2011)Williams, Castonguay, Vincent and
  Jameson}]{Williams2011}
\bibinfo{author}{Williams, D.}, \bibinfo{author}{Castonguay, P.},
  \bibinfo{author}{Vincent, P.}, \bibinfo{author}{Jameson, A.},
  \bibinfo{year}{2011}.
\newblock \bibinfo{title}{An extension of energy stable flux reconstruction to
  unsteady, non-linear, viscous problems on mixed grids}, in:
  \bibinfo{booktitle}{20th {AIAA} Computational Fluid Dynamics Conference},
  \bibinfo{publisher}{American Institute of Aeronautics and Astronautics}.
\newblock \URLprefix \url{https://doi.org/10.2514/6.2011-3405},
  \DOIprefix\doi{10.2514/6.2011-3405}.
\bibitem[{Williams et~al.(2013)Williams, Castonguay, Vincent and
  Jameson}]{Williams2013tri}
\bibinfo{author}{Williams, D.}, \bibinfo{author}{Castonguay, P.},
  \bibinfo{author}{Vincent, P.}, \bibinfo{author}{Jameson, A.},
  \bibinfo{year}{2013}.
\newblock \bibinfo{title}{Energy stable flux reconstruction schemes for
  advection{\textendash}diffusion problems on triangles}.
\newblock \bibinfo{journal}{Journal of Computational Physics}
  \bibinfo{volume}{250}, \bibinfo{pages}{53--76}.
\newblock \URLprefix \url{https://doi.org/10.1016/j.jcp.2013.05.007},
  \DOIprefix\doi{10.1016/j.jcp.2013.05.007}.
\bibitem[{Williams et~al.(2014)Williams, Shunn and Jameson}]{Williams2014}
\bibinfo{author}{Williams, D.}, \bibinfo{author}{Shunn, L.},
  \bibinfo{author}{Jameson, A.}, \bibinfo{year}{2014}.
\newblock \bibinfo{title}{Symmetric quadrature rules for simplexes based on
  sphere close packed lattice arrangements}.
\newblock \bibinfo{journal}{Journal of Computational and Applied Mathematics}
  \bibinfo{volume}{266}, \bibinfo{pages}{18--38}.
\newblock \URLprefix \url{https://doi.org/10.1016/j.cam.2014.01.007},
  \DOIprefix\doi{10.1016/j.cam.2014.01.007}.
\bibitem[{Williams and Jameson(2013)}]{Williams2013tetra}
\bibinfo{author}{Williams, D.M.}, \bibinfo{author}{Jameson, A.},
  \bibinfo{year}{2013}.
\newblock \bibinfo{title}{Energy stable flux reconstruction schemes for
  advection{\textendash}diffusion problems on tetrahedra}.
\newblock \bibinfo{journal}{Journal of Scientific Computing}
  \bibinfo{volume}{59}, \bibinfo{pages}{721--759}.
\newblock \URLprefix \url{https://doi.org/10.1007/s10915-013-9780-2},
  \DOIprefix\doi{10.1007/s10915-013-9780-2}.
\bibitem[{Witherden et~al.(2014)Witherden, Farrington and
  Vincent}]{Witherden2014}
\bibinfo{author}{Witherden, F.}, \bibinfo{author}{Farrington, A.},
  \bibinfo{author}{Vincent, P.}, \bibinfo{year}{2014}.
\newblock \bibinfo{title}{{PyFR}: An open source framework for solving
  advection{\textendash}diffusion type problems on streaming architectures
  using the flux reconstruction approach}.
\newblock \bibinfo{journal}{Computer Physics Communications}
  \bibinfo{volume}{185}, \bibinfo{pages}{3028--3040}.
\newblock \URLprefix \url{https://doi.org/10.1016/j.cpc.2014.07.011},
  \DOIprefix\doi{10.1016/j.cpc.2014.07.011}.
\bibitem[{Witherden et~al.(2015)Witherden, Vermeire and
  Vincent}]{Witherden2015}
\bibinfo{author}{Witherden, F.}, \bibinfo{author}{Vermeire, B.},
  \bibinfo{author}{Vincent, P.}, \bibinfo{year}{2015}.
\newblock \bibinfo{title}{Heterogeneous computing on mixed unstructured grids
  with {PyFR}}.
\newblock \bibinfo{journal}{Computers {\&} Fluids} \bibinfo{volume}{120},
  \bibinfo{pages}{173--186}.
\newblock \URLprefix \url{https://doi.org/10.1016/j.compfluid.2015.07.016},
  \DOIprefix\doi{10.1016/j.compfluid.2015.07.016}.
\bibitem[{Zwanenburg and Nadarajah(2016)}]{Zwanenburg2016}
\bibinfo{author}{Zwanenburg, P.}, \bibinfo{author}{Nadarajah, S.},
  \bibinfo{year}{2016}.
\newblock \bibinfo{title}{Equivalence between the energy stable flux
  reconstruction and filtered discontinuous galerkin schemes}.
\newblock \bibinfo{journal}{Journal of Computational Physics}
  \bibinfo{volume}{306}, \bibinfo{pages}{343--369}.
\newblock \URLprefix \url{https://doi.org/10.1016/j.jcp.2015.11.036},
  \DOIprefix\doi{10.1016/j.jcp.2015.11.036}.

\end{thebibliography}



\end{document}